\documentclass[12pt]{article}

\usepackage{amsfonts,amssymb,euscript,latexsym,setspace,%
subeqnarray,eqsection,indent,epic}

\newcommand{\proof}{\par\medskip\noindent{\sc Proof.\ }}
\newcommand{\firstproof}{\par\medskip\noindent{\sc First proof.\ }}
\newcommand{\secondproof}{\par\medskip\noindent{\sc Second proof.\ }}

\newcommand{\fourthproof}{\par\medskip\noindent{\sc Fourth proof.\ }}
\newcommand{\qed}{\quad $\Box$ \medskip \medskip}
\def\proofof#1{\bigskip\noindent{\sc Proof of #1.\ }}

\newtheorem{defin}{Definition}[section]

\newtheorem{PROP}[defin]{Proposition}
\newtheorem{proposition}[defin]{Proposition}

\newtheorem{LMA}[defin]{Lemma}
\newtheorem{lemma}[defin]{Lemma}

\newtheorem{question}[defin]{Question}

\newtheorem{problem}[defin]{Problem}

\newtheorem{THM}[defin]{Theorem}
\newtheorem{theorem}[defin]{Theorem}

\newtheorem{CORO}[defin]{Corollary}
\newtheorem{corollary}[defin]{Corollary}

\newtheorem{conjecture}[defin]{Conjecture}
\newtheorem{noname}[defin]{}

\newcounter{example}[section]
\newenvironment{example}%
{\refstepcounter{example}
 \bigskip\par\noindent{\bf Example \thesection.\arabic{example}.}\quad
}%
{\quad $\Box$}
\def\bexam{\begin{example}}
\def\eexam{\end{example}}

\def\reff#1{(\protect\ref{#1})}
\newcommand{\be}{\begin{equation}}
\newcommand{\ee}{\end{equation}}

\def\spose#1{\hbox to 0pt{#1\hss}}
\def\ltapprox{\mathrel{\spose{\lower 3pt\hbox{$\mathchar"218$}}
 \raise 2.0pt\hbox{$\mathchar"13C$}}}
\def\gtapprox{\mathrel{\spose{\lower 3pt\hbox{$\mathchar"218$}}
 \raise 2.0pt\hbox{$\mathchar"13E$}}}
\font\fourrm  = cmr5 
\def\gtmaybeeq{\mathrel{\spose{\lower 0.75pt\hbox{\kern-1.5pt\fourrm (\quad)}}
  \raise 2.0pt\hbox{$\ge$}}}

\def\half{ {1 \over 2} }

\def\smfrac#1#2{\textstyle{#1\over #2}}
\def\smhalf{ \smfrac{1}{2} }

\newcommand{\tr}{\mathop{\rm tr}\nolimits}
\newcommand{\cotr}{\mathop{\rm cotr}\nolimits}
\newcommand{\diag}{\mathop{\rm diag}\nolimits}
\newcommand{\rank}{\mathop{\rm rank}\nolimits}

\newcommand{\supp}{\mathop{\rm supp}\nolimits}
\newcommand{\per}{\mathop{\rm per}\nolimits}
\newcommand{\ext}{\mathop{\rm ext}\nolimits}
\newcommand{\coext}{\mathop{\rm coext}\nolimits}
\renewcommand{\arg}{{\rm arg}\,}   
\renewcommand{\det}{{\rm det}\,}   
\newcommand{\real}{{\rm Re}\,}     
\newcommand{\re}{{\rm Re}\,}       
\renewcommand{\Re}{{\rm Re}\,}     
\newcommand{\imag}{{\rm Im}\,}     
\renewcommand{\Im}{{\rm Im}\,}     

\newcommand{\restrict}{\upharpoonright}
\newcommand{\implies}{\;\Longrightarrow\;}
\newcommand{\drop}{\setminus}
\renewcommand{\emptyset}{\varnothing}

\def\Z{{\mathbb Z}}
\def\ZZ{{\mathbb Z}}
\def\R{{\mathbb R}}
\def\C{{\mathbb C}}
\def\CC{{\mathbb C}}
\def\N{{\mathbb N}}
\def\NN{{\mathbb N}}
\def\Q{{\mathbb Q}}

\def\scra{\mathcal{A}}
\def\scrb{\mathcal{B}}
\def\B{\mathcal{B}}

\def\scrd{\mathcal{D}}
\def\scre{\mathcal{E}}
\def\scrf{\mathcal{F}}
\def\scrg{\mathcal{G}}
\def\scrh{\mathcal{H}}
\def\scri{\mathcal{I}}
\def\scrj{\mathcal{J}}

\def\scrm{\mathcal{M}}
\def\scro{\mathcal{O}}
\def\scrp{\mathcal{P}}

\def\scrs{\mathcal{S}}

\def\scrw{\mathcal{W}}
\def\scrz{\mathcal{Z}}

\newcommand{\m}{\mathbf{m}}
\newcommand{\n}{\mathbf{n}}
\renewcommand{\k}{\mathbf{k}}
\newcommand{\p}{\mathbf{p}}
\newcommand{\q}{\mathbf{q}}
\renewcommand{\d}{{\mbox{\boldmath $\delta$}}}

\newenvironment{scarray}{
          \textfont0=\scriptfont0
          \scriptfont0=\scriptscriptfont0
          \textfont1=\scriptfont1
          \scriptfont1=\scriptscriptfont1
          \textfont2=\scriptfont2
          \scriptfont2=\scriptscriptfont2
          \textfont3=\scriptfont3
          \scriptfont3=\scriptscriptfont3
        
        \begin{array}{c}}{\end{array}}

\newcommand{\chapquote}[2]{
  \begin{list}{}{
     \setlength{\leftmargin}{0.4\textwidth}
     \setlength{\rightmargin}{0.05\textwidth}
     \setlength{\labelwidth}{0pt}
  }
  \item
  \footnotesize
  #1
  \par\raggedleft
  {\sl --- #2}
  \end{list}
  \normalsize\smallskip
}

\begin{document}

\title{
         \vspace*{-3.7cm}
         Homogeneous Multivariate Polynomials \\
         with the Half-Plane Property
         \vspace*{-3mm}
      }

\author{\small Young-Bin Choe \\[-2mm]
        \small Department of Combinatorics and Optimization \\[-2mm]
        \small University of Waterloo \\[-2mm]
        \small Waterloo, Ontario N2L 3G1, CANADA \\[-2mm]
        \small {\tt ybchoe@math.uwaterloo.ca} \\[2mm]
        \small James G.~Oxley \\[-2mm]
        \small Department of Mathematics \\[-2mm]
        \small Louisiana State University \\[-2mm]
        \small Baton Rouge, LA 70803-4918, USA \\[-2mm]
        \small {\tt oxley@math.lsu.edu} \\[2mm]
        \small Alan D.~Sokal \\[-2mm]
        \small Department of Physics \\[-2mm]
        \small New York University \\[-2mm]
        \small New York, NY 10003, USA \\[-2mm]
        \small {\tt sokal@nyu.edu} \\[2mm]
        \small David G. Wagner \\[-2mm]
        \small Department of Combinatorics and Optimization \\[-2mm]
        \small University of Waterloo \\[-2mm]
        \small Waterloo, Ontario N2L 3G1, CANADA \\[-2mm]
        \small {\tt dgwagner@math.uwaterloo.ca} \\[2mm]
}

\date{February 4, 2002 \\[1mm] revised November 23, 2002}

\maketitle
\thispagestyle{empty}   

\vspace*{-8mm}

\begin{abstract}
A polynomial $P$ in $n$ complex variables is said to have
the ``half-plane property'' (or Hurwitz property)
if it is nonvanishing whenever all the variables
lie in the open right half-plane.
Such polynomials arise in combinatorics, reliability theory,
electrical circuit theory and statistical mechanics.
A particularly important case is when
the polynomial is homogeneous and multiaffine:
then it is the (weighted) generating polynomial of an $r$-uniform set system.
We prove that the support (set of nonzero coefficients)
of a homogeneous multiaffine polynomial with the half-plane property
is necessarily the set of bases of a matroid.
Conversely, we ask:
For which matroids $M$ does the basis generating polynomial $P_{\scrb(M)}$
have the half-plane property?
Not all matroids have the half-plane property,
but we find large classes that do:
all sixth-root-of-unity matroids,
and a subclass of transversal (or cotransversal) matroids
that we call ``nice''.
Furthermore, the class of matroids with the half-plane property
is closed under minors, duality, direct sums, 2-sums,
series and parallel connection, full-rank matroid union,
and some special cases of principal truncation, principal extension,
principal cotruncation and principal coextension.
Our positive results depend on two distinct (and apparently unrelated)
methods for constructing polynomials with the half-plane property:
a determinant construction (exploiting ``energy'' arguments),
and a permanent construction (exploiting the Heilmann--Lieb theorem
on matching polynomials).
We conclude with a list of open questions.
\end{abstract}

\vspace{0.5cm}

\vspace{0.5cm}
\noindent
{\bf 2000 MATHEMATICS SUBJECT CLASSIFICATION:}
Primary 05B35;
Secondary 05C99, 05E99, 15A15, 15A48, 30C15, 32A99, 82B20, 94C05.

\vspace{0.5cm}
\noindent
{\bf RUNNING HEAD:} Polynomials with the Half-Plane Property

\vspace{0.5cm}

\tableofcontents
\clearpage

\section{Introduction}  \label{sec1}

\chapquote{
   It seems that the theory of polynomials, linear in each variable,
   that do not have zeros in a given multidisk or a more general set,
   has a long way to go, and has so far unnoticed connections to
   various other concepts in mathematics.
}{Aimo Hinkkanen \cite[p.~288]{Hinkkanen_97}}

%

Let us consider a connected graph
$G=(V,E)$ as an electrical network,
by assigning to each edge $e \in E$ a complex number $x_e$,
called its {\em conductance}\/ (or {\em admittance}\/).\footnote{
   In this paper all graphs are finite and undirected.
   Loops and multiple edges are allowed unless specified otherwise.
}
The node voltages $\{ \varphi_i \} _{i \in V}$
and current inflows $\{ J_i \} _{i \in V}$
then satisfy the linear system $L(x) \varphi = J$,
where $L(x)$ is the edge-weighted Laplacian matrix for $G$
[that is, $L(x) = B X B^{\rm T}$
 where $B$ is the directed vertex-edge incidence matrix for
 any orientation of $G$, and $X = \diag(\{x_e\})$].
On physical grounds we expect that if $\real x_e > 0$ for all $e$
(i.e.\ every branch is dissipative),
then the network is uniquely solvable once we fix the voltage
at a single reference node $i_0 \in V$,
or in other words that the $i_0$th principal cofactor of $L(x)$ is nonzero.
Now, by the matrix-tree theorem
\cite{Brooks_40,Chaiken_78,Chaiken_82,Zeilberger_85},
each principal cofactor of $L(x)$ equals the spanning-tree sum
\be
 \label{def_TG}
   T_G(x)  \;=\;  \sum_{{\rm trees}\, T \subseteq E}  x^T   \;,
\ee
where we have used the shorthand $x = \{x_e\}_{e \in E}$
and $x^T = \prod_{e \in T} x_e$.
We therefore conjecture:

\begin{theorem}
 \label{thm1.1}
Let $G$ be a connected graph.
Then the spanning-tree polynomial $T_G$ has the ``half-plane property'',
i.e.\ $\real x_e > 0$ for all $e$ implies $T_G(x) \neq 0$.
\end{theorem}

The proof of Theorem~\ref{thm1.1} is not difficult:
Consider any nonzero complex vector $\varphi = \{ \varphi_i \} _{i \in V}$
satisfying $\varphi_{i_0} = 0$.
Because $G$ is connected, we have $B^{\rm T} \varphi \neq 0$.
Therefore, the quantity
\be
   \varphi^* L(x) \varphi  \;=\;
   \varphi^* B X B^{\rm T} \varphi  \;=\;
   \sum_{e \in E} |(B^{\rm T} \varphi)_e|^2 \, x_e
 \label{energy_form}
\ee
has strictly positive real part whenever $\real x_e > 0$ for all $e$;
so in particular $(B X B^{\rm T} \varphi)_i \neq 0$ for some $i \neq i_0$.
It follows that the submatrix of $L(x)$ obtained by suppressing
the $i_0$th row and column is nonsingular,
and so has a nonzero determinant.
Theorem~\ref{thm1.1} then follows from the matrix-tree theorem.\footnote{
   This proof is well known in the circuit-theory literature:
   see e.g.\ \cite[Section 2.7]{Chen_71}
   as well as the related results in
   \cite[pp.~398--401, 430--431 and 850--851]{Desoer_69}
   \cite[pp.~52--53 and 67--69]{Penfield_70}.
   It has, moreover, a natural physical interpretation:
   if $\varphi = \{ \varphi_i \} _{i \in V}$ are the node voltages,
   then the real part of the quadratic form \reff{energy_form}
   is the total power dissipated in the circuit.
   We thank Charles Desoer and Paul Penfield for pointing out references
   \cite{Desoer_69,Penfield_70}.
}

An immediate corollary of Theorem~\ref{thm1.1}
is that the complementary spanning-tree polynomial
\be
 \label{def_TtildeG}
   \widetilde{T}_G(x)  \;=\;  \sum_{{\rm trees}\, T \subseteq E}
       x^{E \setminus T}
\ee
also has the half-plane property, since
\be
   \widetilde{T}_G(x)  \;=\;  x^E \, T_G(1/x)
\ee
and the map $x_e \mapsto 1/x_e$ takes the right half-plane onto itself.

{}From a combinatorial point of view,
the noteworthy fact is that the spanning trees of $G$
constitute the bases of the graphic matroid $M(G)$,
and their complements constitute the bases of the cographic matroid $M^*(G)$.
So $T_G$ and $\widetilde{T}_G$ are the (multivariate)
basis generating polynomials for $M(G)$ and $M^*(G)$, respectively.
This naturally suggests generalizing Theorem~\ref{thm1.1}
to more general matroids and, perhaps, to more general set systems.
Before posing these questions precisely,
we need to fix some notation and terminology.

A {\em set system}\/ (or {\em hypergraph}\/) $\scrs$
on the (finite) ground set $E$
is simply a collection $\scrs$ of subsets of $E$.
Given any set system $\scrs$ on $E$,
we define its {\em (multivariate) generating polynomial}\/ to be
\be
   P_\scrs(x)   \;=\;  \sum_{S \in \scrs} x^S   \;,
\ee
where $x = \{ x_e \} _{e \in E}$ are commuting indeterminates
(which we shall usually take to be complex variables).
The {\em rank}\/ of a set system is the maximum cardinality
of its members (by convention we set rank $= -\infty$ if $\scrs = \emptyset$);
equivalently, it is the degree of the generating polynomial $P_\scrs$.
A set system $\scrs$ is {\em $r$-uniform}\/ if $|S|=r$ for all $S \in \scrs$,
or equivalently if its generating polynomial $P_\scrs$ is
homogeneous of degree $r$.

An {\em abstract simplicial complex}\/ (or {\em complex}\/ for short)
is a set system $\scrs$ satisfying
\begin{quote}
\begin{itemize}
   \item[(I1)]  $\emptyset \in \scrs$.
   \item[(I2)]  If $S \in \scrs$ and $S' \subseteq S$, then $S' \in \scrs$
      (``$\scrs$ is hereditary downwards'').
\end{itemize}
\end{quote}
The members $S$ of a complex are called {\em faces}\/,
and the maximal members (with respect to set-theoretic inclusion)
are called {\em facets}\/.
A complex is called {\em pure}\/ (of rank $r$)
if all its facets have the same cardinality $r$.

A complex $\scrs$ is called a {\em matroid complex}\/
if it satisfies the further condition
\begin{quote}
\begin{itemize}
   \item[(I3)]  If $S_1, S_2 \in \scrs$ and $|S_1| < |S_2|$,
      then there exists an element $e \in S_2 \setminus S_1$ such that
      $S_1 \cup \{e\} \in \scrs$
      (``independence augmentation axiom for matroids'').
\end{itemize}
\end{quote}
The faces of a matroid complex are called the {\em independent sets}\/
of the matroid;
the facets (i.e.\ the maximal independent sets) are called {\em bases}\/.
It is easy to prove that every matroid complex is pure,
i.e.\ all bases have the same cardinality $r$,
called the {\em rank}\/ of the matroid.
Moreover, it is not difficult to show that a collection $\scrb$
of subsets of $E$ is the collection of bases of a matroid on $E$
if and only if it satisfies the following two conditions:
\begin{quote}
\begin{itemize}
   \item[(B1)]  $\scrb$ is nonempty.
   \item[(B2)]  If $B_1, B_2 \in \scrb$ and $x \in B_1 \setminus B_2$,
      then there exists $y \in B_2 \setminus B_1$ such that
      $(B_1 \setminus x) \cup \{y\} \in \scrb$
      (``basis exchange axiom for matroids'').
\end{itemize}
\end{quote}
For more information on matroid theory, see \cite{Oxley_92}.

We shall be particularly interested in the {\em basis generating polynomial}\/
of a matroid $M$,
\be
   P_{\scrb(M)}(x)  \;=\;  \sum_{B \in \scrb(M)}  x^B
   \;.
\ee
We can now pose the following questions concerning
possible extensions of Theorem~\ref{thm1.1}:

\begin{question}
  \label{question1.2}
For which matroids $M$ does the basis generating polynomial $P_{\scrb(M)}$
have the half-plane property?
\end{question}

\noindent
More generally:

\begin{question}
  \label{question1.3}
For which $r$-uniform set systems $\scrs$ does the generating polynomial
$P_{\scrs}$ have the half-plane property?
\end{question}

Our original conjecture was that all matroids
(and no non-matroidal set systems) have the half-plane property.
That would be nice and neat, but it turns out to be false;
and the truth is considerably more interesting and subtle.
Our conjecture is half right:
an $r$-uniform set system with the half-plane property
is necessarily the set of bases of a matroid (Theorem~\ref{thm.matroidal}).
But not every matroid has the half-plane property,
and we do not yet have a complete characterization of those that do.
Nevertheless, we can find large classes of matroids
with the half-plane property:
\begin{itemize}
   \item[(a)]  Every sixth-root-of-unity matroid \cite{Whittle_97}
has the half-plane property
(Theorem~\ref{thm.QA} and Corollary~\ref{cor.determinant}).
This class properly includes the regular matroids,
which in turn properly include the graphic and cographic matroids.
The proof of Theorem~\ref{thm.QA} is, in fact, a direct generalization
of the proof just given for Theorem~\ref{thm1.1}.
   \item[(b)]  Every uniform matroid has the half-plane property,
and indeed has the (stronger) ``Brown--Colbourn property''
(Section~\ref{sec_uniform}).
   \item[(c)]  A significant subclass of transversal matroids
--- those we call ``nice'' ---
have the half-plane property (Section~\ref{sec_transversal}).
Indeed, it may well be true that {\em all}\/ transversal matroids
have the half-plane property, but we have no idea how to prove this.
   \item[(d)]  All matroids of rank or corank at most 2 have the
half-plane property (Corollary~\ref{cor.rank2a}), as do all matroids
on a ground set of at most 6 elements (Proposition~\ref{prop.nle6}).
   \item[(e)]  The class of matroids with the half-plane property
is closed under minors, duality, direct sums, 2-sums,
series and parallel connection, full-rank matroid union,
and some special cases of principal truncation, principal extension,
principal cotruncation and principal coextension
(Section~\ref{sec_constructions}).
\end{itemize}
Moreover, we can show that certain matroids do {\em not}\/ have the
half-plane property:  among these are the Fano matroid $F_7$,
the non-Fano matroid $F_7^-$,
their relaxations $F_7^{--}$, $F_7^{-3}$ and $M(K_4)+e$,
the matroids $P_8$, $P'_8$ and $P''_8$,
the Pappus and non-Pappus matroids,
the free extension (non-Pappus $\drop\,9)+e$,
and all their duals
(Section~\ref{sec_counterexamples}).
The first six of these examples are minor-minimal,
and we conjecture that the others are as well;
but we strongly suspect that this list is incomplete,
and indeed we consider it likely that the set of minor-minimal
non-half-plane-property matroids is infinite.

More generally, we shall consider homogeneous multiaffine polynomials
$P(x) = \sum_{S \subseteq E, |S| = r} a_S x^S$
with arbitrary complex coefficients $a_S$ (not necessarily 0 or 1).
We shall prove two {\em necessary}\/ conditions for $P \not\equiv 0$
to have the half-plane property:
\begin{itemize}
   \item[(a)] $P$ must have the ``same-phase property'',
      i.e.\ all the nonzero coefficients $a_S$ must have the same phase
      (Theorem~\ref{thm.same-phase}).
      So without loss of generality we can assume that all the $a_S$
      are nonnegative.
   \item[(b)] The {\em support}\/
      $\supp(P) = \{S \subseteq E \colon\;  a_S \neq 0 \}$
      must be the collection of bases of a matroid
      (Theorem~\ref{thm.matroidal}).
\end{itemize}
This latter fact is particularly striking:  it shows that matroids
arise {\em naturally}\/ from a consideration of homogeneous multiaffine
polynomials with the half-plane property.
We do not know whether the converse of Theorem~\ref{thm.matroidal}
is true, i.e.\ whether for every matroid $M$ there exists a
homogeneous multiaffine polynomial $P$ with the half-plane property
such that $\supp(P) = \scrb(M)$.
But it is true, at least, for all matroids representable over $\C$
(Corollary~\ref{cor.determinant}).

We shall also also give two {\em sufficient}\/ conditions
for a homogeneous multiaffine polynomial $P$
to have the half-plane property (or be identically zero):
\begin{itemize}
   \item[(a)] {\em Determinant condition}\/ (Theorem~\ref{thm.QA}):
      $a_S = |\det(A \restrict S)|^2$ for some $r \times n$ complex matrix $A$
      [here $n=|E|$, and $A \restrict S$ denotes the square submatrix of $A$
       using the columns indexed by the set $S$].
      This corresponds to $P(x) = \det(A X A^*)$
      where $X = \diag(\{x_e\})$ and ${}^*$ denotes Hermitian conjugate.
   \item[(b)] {\em Permanent condition}\/ (Theorem~\ref{thm.permanent}):
      $a_S = \per(\Lambda \restrict S)$ for some $r \times n$ nonnegative
      matrix $\Lambda$.
      This corresponds to $P(x) = \per(\Lambda X)$.
\end{itemize}
Unfortunately, the relationship between these sufficient conditions
and the half-plane property looks complicated.
Neither family of polynomials contains the other;
their intersection is nonempty;
and their union is a proper subset of the set of
all homogeneous multiaffine polynomials with the half-plane property.

These questions also have a close connection with reliability theory
\cite{Colbourn_87}.
Consider a finite set $E$ of communication channels,
which fail independently with probabilities $\{ q_e \}_{e \in E}$.
Let $\scrs$ be a set system on $E$, whose members we shall interpret as the
sets of failed channels that allow the system as a whole to be operational.
Then the probability that the system is operational
is given by the {\em multivariate reliability polynomial}\/
\be
   {\rm Rel}_\scrs(q)  \;=\;
   \sum_{A \in \scrs} q^A \, (1-q)^{E \setminus A}
 \label{def.relpoly}
\ee
where $q = \{ q_e \}_{e \in E}$.
This is easily related to the multivariate generating polynomial
\be
   P_\scrs(x)   \;=\;   \sum_{A \in \scrs} x^A
\ee
by
\begin{eqnarray}
  {\rm Rel}_\scrs(q)   & = &   (1-q)^E \, P_\scrs\Bigl( {q \over 1-q} \Bigr)
     \\[2mm]
  P_\scrs(x)   & = &   (1+x)^E \, {\rm Rel}_\scrs\Bigl( {x \over 1+x} \Bigr)
\end{eqnarray}
In the reliability context it is natural to assume that
$\scrs$ is a complex,
i.e.\ that $\scrs$ contains $\emptyset$ and is closed under taking subsets.
Indeed, the simplest case arises when $G=(V,E)$
is a connected graph and we declare the system to be operational
if the non-failed edges form a connected spanning subgraph
(this is the ``all-terminal reliability'').
In this case $\scrs$ is the set of complements of
connected spanning subgraphs of $G$,
i.e.\ the family $\scri(M^*(G))$ of independent sets
of the cographic matroid $M^*(G)$.
Our Holy Grail is the following:

\begin{conjecture}[multivariate Brown--Colbourn conjecture
                   \protect\cite{Brown_92,Sokal_01}]
 \label{conj.Brown-Colbourn}
Let $G$ be a connected graph.
If $|q_e| > 1$ for all $e$, then ${\rm Rel}_{\scri(M^*(G))}(q) \neq 0$.
Equivalently, if $G$ is loopless and $\real x_e < -1/2$ for all $e$,
then $P_{\scri(M^*(G))}(x) \neq 0$.
\end{conjecture}

\noindent
(Note that a loop in $G$ corresponds to a coloop in $M^*(G)$.
 Each loop in $G$ has no effect on the reliability polynomial
 but multiplies the independent-set polynomial of $M^*(G)$ by a factor $1+x_e$,
 leading to a root at $x_e = -1$.
 This is why, in discussing the roots of $P_{\scri(M^*(G))}$,
 we need to assume that $G$ is loopless.)

We are at present quite far from a proof of
Conjecture~\ref{conj.Brown-Colbourn}.
One of us has proven Conjecture~\ref{conj.Brown-Colbourn}
for the special case of series-parallel graphs
\cite[Remark 3 in Section 4.1]{Sokal_01};
earlier, another one of us had proven the corresponding univariate result,
i.e.\ when all the $q_e$ take the same value \cite{Wagner_00}.
But series-parallel graphs are a small subset of planar graphs,
and an even smaller subset of all graphs!

Nonetheless one can dream,
and even lacking a proof of Conjecture~\ref{conj.Brown-Colbourn}
it is reasonable to ask whether stronger results might be true.
For example:

\begin{question}
  \label{question1.5}
Which (coloopless) matroids $M$ have the ``Brown--Colbourn property'',
i.e.\ $\real x_e < -1/2$ for all $e$ implies $P_{\scri(M)}(x) \neq 0$?
\end{question}

\noindent
Aside from the cographic (or equivalently, graphic) matroids
of series-parallel graphs, we can prove the Brown--Colbourn property
for uniform matroids $U_{r,n}$ with $0 \le r < n$
(the condition $r < n$ corresponds precisely to forbidding coloops):
this follows immediately from the corresponding univariate result
\cite[Proposition 7.3]{Wagner_00}
combined with the Grace--Walsh--Szeg\"o coincidence theorem
(see Section~\ref{sec_uniform}).
On the other hand, it is easy to show (Corollary~\ref{cor.shiftedHPP})
that the Brown--Colbourn property for $P_{\scri(M)}$
implies the half-plane property for $P_{\scrb(M)}$.
So the results of this paper imply that the Brown--Colbourn property
fails for many matroids (e.g.\ $F_7$, $F_7^-$, \ldots)
and for all pure complexes that are not matroidal.

One purpose of this paper is, therefore, to serve as a ``warm-up''
for an attack on the Brown--Colbourn property,
by studying first a property (the half-plane property)
that is a {\em necessary}\/ condition for the Brown--Colbourn property
and may be easier to characterize.
But the results of our investigation show that
the half-plane property is very interesting in its own right!

Our study of polynomials with the half-plane property
can also be viewed in the wider context of theorems
asserting that some combinatorially interesting class of
multivariate  polynomials are nonvanishing in some large domain
of complex $n$-space.
Theorems of this kind include
the Lee--Yang theorem on ferromagnetic Ising models
and generalizations thereof
\cite{Lee-Yang_52,Simon_74,Lieb-Sokal_81,Hinkkanen_97}
and the Heilmann--Lieb theorem on matching polynomials \cite{Heilmann_72}
(discussed in Section~\ref{sec_Heilmann-Lieb} below).
Useful tools for manipulating such polynomials include
the Grace--Walsh--Szeg\"o coincidence theorem
\cite{Grace_02,Szego_22,Walsh_22,Marden_66}
(see Section~\ref{sec.Grace}),
the Asano contraction lemma \cite{Asano_70,Ruelle_71,Hinkkanen_97}
and the Hinkkanen composition theorem \cite{Hinkkanen_97}
(see Section~\ref{sec.folding}).

The plan of this paper is as follows:
In Section~\ref{sec_polynomials} we discuss the basic properties
of multivariate polynomials with the half-plane property,
including their important connection with
real-part-positive rational functions.
In Section~\ref{sec_local} we discuss a key criterion
that we call the ``local half-plane property''.
In Section~\ref{sec_constructions} we describe a large number
of constructions that preserve the half-plane property;
most of these are motivated by standard operations on matroids.
In Section~\ref{sec_necessary} we provide a necessary and sufficient
condition for a polynomial to have the half-plane property;
this condition will play an important role in finding counterexamples.
We also provide a simple explicit criterion in the rank-2 case.
In Section~\ref{sec_same-phase} we prove that, in a homogeneous polynomial
with the half-plane property, all the nonzero coefficients must have the
same phase.
In Section~\ref{sec_matroidal} we prove that the support of a
homogeneous multiaffine polynomial must be the collection of bases
of a matroid;  we also prove a generalization to the non-multiaffine case.
In Section~\ref{sec_determinant} we give the ``determinant construction'',
and deduce as a corollary that every sixth-root-of-unity matroid
has the half-plane property.
This section also contains some results on $(F,G)$-representability
of matroids that may be of independent interest.
In Section~\ref{sec_uniform} we prove that
every uniform matroid has the half-plane property,
and indeed has the Brown--Colbourn property.
In Section~\ref{sec_transversal} we give the ``permanent construction'',
and deduce as a corollary that a certain subclass of transversal matroids
(those we call ``nice'') have the half-plane property;
we also give numerous examples of nice and non-nice transversal matroids.
In Section~\ref{sec_counterexamples} we give several examples of
matroids that do {\em not}\/ have the half-plane property.
In Section~\ref{sec_numerical} we report our numerical experiments
on matroids for which we have been unable to prove or disprove
the half-plane property.
We conclude, in Section~\ref{sec_open}, with a list of open questions.
In Appendix~\ref{app_matroids} we provide, for the convenience of the reader,
a list of the matroids considered in this paper,
along with brief summaries of their properties.
In Appendix~\ref{app_F1rep} we analyze the $(F,\{1\})$-representability
of matroids.

{\bf Note Added (November 2002)}:
Gordon Royle and one of the authors (A.D.S.)\ have recently discovered
--- to our great surprise --- that the Brown--Colbourn conjecture is false!
The multivariate Brown--Colbourn conjecture is false
already for the simplest non-series-parallel graph,
namely the complete graph $K_4$.
The univariate Brown--Colbourn conjecture is false
for certain simple planar graphs that can be obtained from $K_4$
by parallel and series extension of edges.
This work will be reported separately \cite{Royle-Sokal}.

\section{Polynomials and the half-plane property}   \label{sec_polynomials}

Our main interest is in polynomials that are homogeneous of degree $r$
and are multiaffine (i.e.\ of degree $\le 1$ in each variable separately):
these include the generating polynomials of $r$-uniform set systems
and, more specifically, the basis generating polynomials of
rank-$r$ matroids.
However, many of our results will be valid for homogeneous polynomials
that are not necessarily multiaffine,
or for multiaffine polynomials that are not necessarily homogeneous,
or (sometimes) for general polynomials.
We shall endeavor to state our results in whatever degree of generality
seems most natural, without being pedantic.

\subsection{Basic definitions}

Let $P(x)=\sum_{\bf m} a_{\bf m}x^{\bf m}$ be a polynomial
with complex coefficients in the variables $x=\{x_e\}_{e\in E}$.
We call the finite set $E=E(P)$ the \emph{ground set} of $P$,
and we call its members \emph{elements}.
We shall generally use the letter $n$ to denote the cardinality of $E$;
often we shall simply take $E$ to be $[n] = \{1,\ldots,n\}$.
In $P$, the sum ranges over all \emph{multi-indices} ${\bf m}$
(that is, functions ${\bf m} \colon\, E\rightarrow\NN$),
and only finitely many of the coefficients $a_{\bf m}$ are nonzero.
The \emph{degree} of the monomial $x^{\m}=\prod_{e\in E}x_{e}^{m(e)}$
is $|{\bf m}|=\sum_{e\in E} m(e)$,
and the \emph{degree of $P$} is
$\deg P = \max\{ |{\bf m}| \colon\; a_{\m}\neq 0\}$.
We say that $P$ is \emph{homogeneous of degree $r$}
if $a_{\bf m} = 0$ whenever $|{\bf m}| \neq r$.
For $e \in E$, the \emph{degree of $P$ in $x_e$}
is $\deg_e P = \max\{m(e) \colon\; a_{\m}\neq 0\}$.
We say that $P$ is \emph{affine in $e$} if $\deg_e P \le 1$,
and we say that $P$ is \emph{multiaffine} if it is affine in $e$
for all $e\in E$.
A multiaffine polynomial can be written in the form
$P(x) = \sum_{S \subseteq E} a_S x^S$.

The \emph{support} of a polynomial $P(x)=\sum_{\m} a_{\m}x^{\m}$
is the set of multi-indices with nonzero coefficients:
\begin{equation}
   \supp(P) \;=\; \{\m \colon\; a_{\m}\neq 0\}  \;.
\end{equation}
In the case of a multiaffine polynomial $P(x)=\sum_{S \subseteq E} a_S x^S$,
we identify multi-indices with subsets of $E$, so that
\begin{equation}
   \supp(P) \;=\; \{S \subseteq E \colon\; a_S \neq 0 \}  \;.
\end{equation}

If $D$ is a domain (connected open set) in $\C^n$,
we denote by $\scrf_D$ the set of all functions analytic in $D$
that are either nonvanishing in $D$ or else identically zero.
By Hurwitz's theorem\footnote{
   Hurwitz's theorem states that if $D$ is a domain in $\C^n$
   and $(f_k)$ are nonvanishing analytic functions on $D$
   that converge to $f$ uniformly on compact subsets of $D$,
   then $f$ is either nonvanishing or else identically zero.
   Hurwitz's theorem for $n=1$ is proved in most standard texts
   on the theory of analytic functions of a single complex variable
   (see e.g.\ \cite[p.~176]{Ahlfors_66}).
   Surprisingly, we have been unable to find Hurwitz's theorem
   proven for general $n$ in any standard text on several complex variables
   (but see \cite[p.~306]{Krantz_92} and \cite[p.~337]{Simon_74}).
   So here, for completeness, is the sketch of a proof:
   Suppose that $f(c) = 0$ for some $c = (c_1,\ldots,c_n) \in D$,
   and let $D' \subset D$ be a small polydisc centered at $c$.
   Applying the single-variable Hurwitz theorem,
   we conclude that $f(z_1,c_2,\ldots,c_n) = 0$
   for all $z_1$ such that $(z_1,c_2,\ldots,c_n) \in D'$.
   Applying the same argument repeatedly in the variables $z_2,\ldots,z_n$,
   we conclude that $f$ is identically vanishing on $D'$
   and hence, by analytic continuation, also on $D$.
},
$\scrf_D$ is closed
in the topology of uniform convergence on compact subsets of $D$.
(That is why we included the phrase ``or else identically zero''
 in the definition of $\scrf_D$:  it takes care of a trivial degenerate case
 in a convenient way.)

Let us now suppose that $D = D_1 \times \cdots \times D_n$
where the $D_i$ are domains in $\C$.
If $f \in \scrf_D$ and we fix some of the variables
(say $x_1,\ldots,x_m$) to particular values
$x_1^{(0)} \in D_1$, \ldots, $x_m^{(0)} \in D_m$,
then clearly $f(x_1^{(0)},\ldots,x_m^{(0)}, \,\cdot\,) \in
 \scrf_{D_{m+1} \times \cdots \times D_n}$.
But more is true, at least if $f \in \scrf_D$
is also continuous in the closure
$\overline{D} = \overline{D_1} \times \cdots \times \overline{D_n}$:
then we can fix $x_1^{(0)} \in \overline{D_1}$, \ldots,
$x_m^{(0)} \in \overline{D_m}$ and we still have
$f(x_1^{(0)},\ldots,x_m^{(0)}, \,\cdot\,) \in
 \scrf_{D_{m+1} \times \cdots \times D_n}$
(again by Hurwitz's theorem).
That is, we can fix some variables even on the boundary of the domain,
and as a function of the remaining variables, $f$ must be either
nonvanishing or else identically zero.

We denote by $H$ the open right half-plane
$\{x \in \C \colon\; \real x > 0 \}$,
and by $\overline{H}$ the closed right half-plane
$\{x \in \C \colon\; \real x \ge 0 \}$.
If the polynomial $P$ in $n$ variables belongs to $\scrf_{H^n}$
(i.e.\ is either nonvanishing in the product of open right half-planes
 or else identically zero),
we say that $P$ has the {\em half-plane property}\/.\footnote{
   In the engineering literature, polynomials $P \not\equiv 0$
   with the half-plane property are termed
   {\em widest-sense Hurwitz polynomials}\/:
   see e.g.\ \cite{Fettweis_87}.
}
The discussion of the previous paragraph has the following corollary:

\begin{proposition}
 \label{prop.hurwitz}
Let $P$ have the half-plane property,
let $E'$ be a subset of the ground set $E=E(P)$,
and fix $x' \in \overline{H}^{E'}$.
Then $P(x', \,\cdot\,)$,
considered as a polynomial on the ground set $E \drop E'$,
has the half-plane property.
In particular, this holds when we take $x' = 0$.
\end{proposition}

\subsection{Shifted half-plane property and leading part}
 \label{sec_leadingpart}

Let us denote by $H_{\theta,K}$ the rotated translated open half-plane
$\{x \in \C \colon\; \real (e^{-i\theta} x) > K \}$.
If the polynomial $P$ in $n$ variables belongs to
$\scrf_{H_{\theta,K}^n}$ for some $\theta$ and $K$,
we say that $P$ has the {\em shifted half-plane property}\/.
In particular, if $P \in \scrf_{H_{\pi,1/2}^n}$,
we say that $P$ has the {\em Brown--Colbourn property}\/.

If $P$ is a polynomial of degree $r$,
we denote by $P^\sharp$ the polynomial consisting of those terms in $P$
that have degree $r$,
and we call it the {\em leading part}\/ of $P$.
Clearly $P^\sharp$ is homogeneous of degree $r$.
Moreover, it is easy to see that
\be
   P^\sharp(x)   \;=\;   \lim_{\zeta\to\infty} \zeta^{-r} P(\zeta x)
\ee
as $\zeta$ tends to infinity in $\C$,
uniformly for $x$ in compact subsets of $\C^n$.
Using Hurwitz's theorem, we deduce immediately that:

\begin{proposition}
 \label{prop.shiftedHPP}
If $P$ has the shifted half-plane property,
then $P^\sharp$ has the half-plane property
[and more generally belongs to $\scrf_{H_{\theta',0}^n}$ for all $\theta'$].
\end{proposition}

\begin{corollary}
 \label{cor.shiftedHPP}
Let $M$ be a matroid.
If the independent-set generating polynomial $P_{\scri(M)}$
has the Brown--Colbourn property,
then the basis generating polynomial $P_{\scrb(M)}$
has the half-plane property.
\end{corollary}

\noindent
Henceforth we shall say that a set system $\scrs$ (resp.\ a matroid $M$)
``has the half-plane property'' if its generating polynomial $P_\scrs$
(resp.\ its basis generating polynomial $P_{\scrb(M)}$)
has the half-plane property,
and that a matroid $M$ ``has the Brown--Colbourn property''
if its independent-set generating polynomial $P_{\scri(M)}$ does.


\subsection{Real-part-positive rational functions}

Let $D$ be a domain in $\C^n$,
and let $f$ be a complex-valued function that is analytic on $D$.
We say that $f$ is {\em real-part-positive on $D$}\/
(resp.\ {\em strictly real-part-positive on $D$}\/)
if $\re f(x) \ge 0$ (resp.\ $\re f(x) > 0$) for all $x \in D$.

\begin{lemma}
 \label{lemma.RPP}
Let $f$ be real-part-positive on $D$.
Then either $f$ is strictly real-part-positive on $D$,
or else $f$ is a pure imaginary constant.
(In particular, $f$ is either nonvanishing on $D$ or else identically zero.)
\end{lemma}

\proof
By the open mapping theorem,
either the image $f[D]$ is open in $\CC$
or else $f$ is constant.
Since $f$ is real-part-positive, $f[D]$ is contained in
the closed right half-plane $\overline{H}$.
If $f[D]$ is open, then it is contained in the interior of $\overline{H}$,
i.e.\ in the open right half-plane $H$,
so that $f$ is strictly real-part-positive.
If $f$ is constant, this constant value is either
pure imaginary or else has strictly positive real part.
\qed

\begin{lemma}
 \label{lemma.RPP.composition}
Let $D$ be a domain in $\C^n$, and let $k$ be a positive integer.
Let $f$ be real-part-positive (resp.\ strictly real-part-positive) on $H^k$,
and let $g_1,\ldots,g_k$ be strictly real-part-positive on $D$.
Then $f \circ {\bf g}$ is real-part-positive
(resp.\ strictly real-part-positive) on $D$.
[Here we have used the obvious shorthand ${\bf g} = (g_1,\ldots,g_k)$.]
\end{lemma}

\proof
Trivial.
\qed

\begin{lemma}
 \label{lemma.RPP_implies_HPP}
Let $D$ be a domain in $\C^n$,
and let $g$ and $h$ be analytic functions on $D$.
Suppose that $h$ is nonvanishing on $D$,
and that $g/h$ is real-part-positive on $D$.
Then $g$ is either nonvanishing on $D$ or else identically zero.

In particular, let $P$ and $Q$ be polynomials in $n$ variables,
with $Q \not\equiv 0$.
Suppose that $Q$ has the half-plane property,
and that the rational function $P/Q$ is real-part-positive on $H^n$.
Then $P$ has the half-plane property.
\end{lemma}

\proof
By Lemma~\ref{lemma.RPP},
$g/h$ is either a (pure imaginary) constant function
or else is strictly real-part-positive on $D$.
If $g/h = c$, then $g$ is either identically zero (if $c=0$)
or nonvanishing on $D$ (if $c \neq 0$).
If $g/h$ is strictly real-part-positive on $D$,
then $g$ is manifestly nonvanishing on $D$.
\qed

In this lemma we have {\em assumed}\/ that $h$ is nonvanishing on $D$
in order to guarantee without fuss that $g/h$ is analytic on $D$.
If we drop this assumption, we can still consider $g/h$ as an analytic
function on $D \setminus \scrz(h)$, where
$\scrz(h) = \{ x \in D \colon\, h(x) = 0 \}$
is the zero set of $h$.
[Note that $\scrz(h)$ is a closed set,
 and has empty interior whenever $h \not\equiv 0$.]
It turns out that no generality is gained by this maneuver,
at least when $g$ and $h$ are polynomials
and we exclude the trivial possibility that
$g$ and $h$ contain a common factor:

\begin{lemma}
 \label{lemma.RPP_implies_HPP_2}
Let $P \not\equiv 0$ and $Q \not\equiv 0$ be polynomials in
$n$ complex variables,
with $P$ and $Q$ relatively prime (over $\C$).
Let $D$ be a domain in $\C^n$.
Suppose that the rational function $P/Q$ is real-part-positive on
$D \setminus \scrz(Q)$.
Then in fact $\scrz(Q) \cap D = \emptyset$.
\end{lemma}

\proof
Suppose that for some $z^{(0)} \in D$ we have $Q(z^{(0)}) = 0$.

(a) If $P(z^{(0)}) \neq 0$, then $Q/P$ is analytic in some neighborhood
$U \ni z^{(0)}$ and is nonconstant, so by the open mapping theorem
$(Q/P)[U]$ contains a neighborhood $V \ni 0$.
Therefore, $(Q/P)[U \setminus \scrz(Q)]$ contains $V \setminus \{0\}$,
in violation of the hypothesis that $P/Q$ is real-part-positive on
$D \setminus \scrz(Q)$.

(b) If $P(z^{(0)}) = 0$, then it is known \cite[Theorem 1.3.2]{Rudin_69}
that for every neighborhood $U \ni z^{(0)}$
we have $(P/Q)[U \setminus \scrz(Q)] = \C$,
which again violates the hypothesis that $P/Q$ is real-part-positive on
$D \setminus \scrz(Q)$.
\qed

The most important case for us will be $D=H^n$.
As we shall see, there is a close interplay between
polynomials with the half-plane property
and rational functions that are real-part-positive on $H^n$.\footnote{
   In the engineering literature, rational functions
   that are real-part-positive on $H^n$ are called {\em positive}\/
   (or {\em positive real}\/ if
    their numerator and denominator polynomials have real coefficients):
   see e.g.\ \cite{Brune_31,Ozaki_60,Fettweis_87}.
   We feel, however, that these terms are likely to cause confusion,
   so we prefer the more precise term ``real-part-positive''.
}
One direction of this interplay is given by
Lemmas~\ref{lemma.RPP_implies_HPP} and \ref{lemma.RPP_implies_HPP_2};
the other direction will be given in Proposition~\ref{prop.derivs} below.

\subsection{Derivatives}

If $P$ is a polynomial with ground set $E$
and ${\bf m} = \{ m_e \}_{e \in E}$ is a multi-index,
we define the polynomial $\partial^{\bf m} P$
in the obvious way:
$(\partial^{\bf m} P)(x) =
 \bigl( \prod_{e \in E} (\partial/\partial x_e)^{m_e} \bigr) P(x)$.
We have the following easy but fundamental result:

%

\begin{proposition}
  \label{prop.derivs}
Let $P$ be a polynomial with ground set $E$,
and let $\{ \lambda_e \} _{e \in E}$ be nonnegative real numbers.
Suppose that $P$ has the half-plane property.  Then:
\begin{itemize}
   \item[(a)] The polynomial
       $\sum\limits_{e \in E} \lambda_e \, \partial P/\partial x_e$
       has the half-plane property.
   \item[(b)] Provided that $P \not\equiv 0$, the rational function
       $P^{-1} \sum_{e \in E} \lambda_e \, \partial P/\partial x_e$
       is real-part-positive on $H^E$.
       Indeed, it is strictly real-part-positive on $H^E$ except when
       $\sum_{e \in E} \lambda_e  \,\partial P/\partial x_e \equiv 0$.
\end{itemize}
\end{proposition}

\proof
If $P \equiv 0$ the proposition is trivial, so assume $P \not\equiv 0$.
We shall prove (b);  Lemma~\ref{lemma.RPP_implies_HPP} then implies (a).

Consider first the univariate case $|E|=1$.
If $P$ is a constant function, then $P' \equiv 0$ and the theorem is trivial.
So let $P$ be a univariate polynomial of degree $k \ge 1$
with the half-plane property, i.e.\
$P(x) = C \prod_{i=1}^k (x - \alpha_i)$ with $C \neq 0$
and $\re \alpha_i \le 0$.
Then
\be
   {P'(x) \over P(x)}  \;=\;  \sum_{i=1}^k {1 \over x - \alpha_i}
\ee
is strictly real-part-positive.

Now consider the general multivariate case.
Applying the univariate result to $x_e$ with
$\{x_f\}_{f \neq e}$ held fixed in the open right half-plane,
we conclude that $P^{-1} \partial P/\partial x_e$
is real-part-positive on $H^E$.
The same therefore holds for
$P^{-1} \sum_{e \in E} \lambda_e \, \partial P/\partial x_e$.
Moreover, by Lemma~\ref{lemma.RPP}
this function is either strictly real-part-positive on $H^E$
or else a (pure imaginary) constant.
But it cannot be a nonzero constant,
because $\deg(\sum_{e \in E} \lambda_e \, \partial P/\partial x_e) < \deg P$.
This proves (b) in the general case.
\qed

{\bf Remarks.}
1.  Part (b) of Proposition~\ref{prop.derivs} was proven by
Koga \cite[Theorem 11]{Koga_68}, but it is probably older.

2.  The same proof shows that the differential operator
\be
   \scro  \;=\;
   \alpha \,+\, \sum_{e \in E} \beta_e {\partial \over \partial x_e}
          \,+\, \sum_{e \in E} \gamma_e x_e
 \label{def_scro}
\ee
with $\real \alpha \ge 0$ and $\beta_e,\gamma_e \ge 0$
preserves the half-plane property, since the function
\be
   {\scro P \over P}  \;=\;
   {\alpha P
        \,+\, \sum\limits_{e \in E} \beta_e {\partial P \over \partial x_e}
        \,+\, \sum\limits_{e \in E} \gamma_e x_e P
    \over
    P}
\ee
is real-part-positive.

3.  More generally, let $\scro_1,\ldots,\scro_n$ be differential operators
of the type \reff{def_scro}
with $\real \alpha \ge 0$ and $\beta_e,\gamma_e \ge 0$,
and let $P_1,\ldots,P_n$ have the half-plane property.
For each index $i$, let $Q_i$ be one of the two functions
$P_i$ and $\scro_i P_i$, and let $R_i$ be the other one.
Then
\be
   S  \;=\;  \sum_{i=1}^n \left( \prod_{j \neq i} R_j \right) Q_i
\ee
has the  half-plane property, because
\be
   {S \over \prod\limits_{i=1}^n R_i}  \;=\; \sum_{i=1}^n {Q_i \over R_i}
\ee
is real-part-positive.
[If exactly one of the $R_i$ is identically zero, then $S=Q_i$.
 If two or more of the $R_i$ are identically zero, then $S \equiv 0$.]
This construction was inspired by our reading of the
brief paper of Bose \cite{Bose_72}.


\bigskip

\begin{corollary}
 \label{cor.derivs}
If the polynomial $P$ belongs to $\scrf_{H_{\theta,K}^n}$, then
so does $\partial^{\bf m} P$ for every multi-index ${\bf m}$.
In particular, if $P$ has the half-plane property,
then so does $\partial^{\bf m} P$ for every multi-index ${\bf m}$.
\end{corollary}

Suppose we partition the ground set $E$ into two disjoint subsets
$E'$ and $E''$, and write
\be
   P(x) \;=\;  \sum_{\bf m'} P_{\bf m'}(x'') \, (x')^{\bf m'}
   \;,
 \label{eq.Pcoeffs}
\ee
where of course $x' = \{x_e\}_{e \in E'}$ and $x'' = \{x_e\}_{e \in E''}$.
We then have:

\begin{corollary}
   \label{cor.Pcoeffs}
Let $P$ be written in the form \reff{eq.Pcoeffs}.
If $P$ has the half-plane property,
then so do all of the coefficient functions $P_{\bf m'}$.
\end{corollary}

\proof
We have
$P_{\bf m'}(x'') =
 (1/{\bf m'}!) \, \partial_{x'}^{\bf m'} P(x',x'') \, |_{x'=0}$
where ${\bf m'}! = \prod_{e \in E'} m'_{e}!$.
The result then follows immediately from
Corollary~\ref{cor.derivs} and Proposition~\ref{prop.hurwitz}.
\qed

In fact, a result much stronger than Corollary~\ref{cor.derivs}
was proven two decades ago by Lieb and one of the authors
\cite[Proposition 2.2]{Lieb-Sokal_81}:

\begin{theorem}[Lieb and Sokal \protect\cite{Lieb-Sokal_81}]
  \label{thm.Lieb-Sokal}
Let $\{P_i\}_{i=1}^k$ and $\{Q_i\}_{i=1}^k$
be polynomials in $n$ complex variables,
and define
\begin{eqnarray}
   R(v,w)  & = &   \sum_{i=1}^k P_i(v) Q_i(w)    \\[2mm]
   S(x)    & = &   \sum_{i=1}^k P_i(\partial/\partial x) Q_i(x)
\end{eqnarray}
If $R$ has the half-plane property, then so does $S$.

In particular, we have the case $k=1$:
if $P$ and $Q$ have the half-plane property,
then so does $P(\partial/\partial x) Q(x)$.
\end{theorem}


\subsection{The Grace--Walsh--Szeg\"o coincidence theorem} \label{sec.Grace}

Finally, we shall need a version of the
Grace--Walsh--Szeg\"o coincidence theorem \cite{Grace_02,Szego_22,Walsh_22}.
Let $P(x_1,\ldots,x_n)$ be a multiaffine polynomial that is {\em symmetric}\/
under all permutations of the $x_1,\ldots,x_n$.
Any such polynomial can be written in the form
\be
   P(x_1,\ldots,x_n)  \;=\;
   \sum_{k=0}^n  a_k {n \choose k}^{\! -1} E_k(x_1,\ldots,x_n)
 \label{eq.symmetric_multiaffine}
\ee
where the $E_k$ are the elementary symmetric polynomials,
i.e.\ $E_0 = 1$ and
\be
   E_k(x_1,\ldots,x_n)   \;=\;
   \sum_{1 \le i_1 < i_2 < \ldots < i_k \le n}  x_{i_1} x_{i_2} \cdots x_{i_k}
   \;;
\ee
and conversely, any polynomial of the form \reff{eq.symmetric_multiaffine}
is symmetric and multiaffine.
Now let us define a {\em closed circular region}\/ in $\C$
to be a closed disc (including the degenerate case of a single point),
a closed half-plane, or the closed exterior of a disc.
And let us define an {\em open circular region}\/ in $\C$
to be an open disc, an open half-plane, or the open exterior of a disc.
We then have:

\begin{theorem}[Grace--Walsh--Szeg\"o]
 \label{thm.grace}
Let $P$ be a symmetric multiaffine polynomial in $n$ complex variables,
let $C$ be an open or closed circular region in $\C$,
and let $x_1^{(0)}, \ldots, x_n^{(0)}$ be points in the region $C$.
Suppose, further, that either $\deg P = n$ (i.e.\ $a_n \neq 0$)
or $C$ is convex (or both).
Then there exists at least one point $\xi \in C$ such that
\be
   P(x_1^{(0)},\ldots,x_n^{(0)})  \;=\;  P(\xi,\ldots,\xi)   \;.
 \label{eq.thm.grace}
\ee
\end{theorem}

\proof
Note first that it suffices to prove the theorem for closed circular regions:
for if $C$ is an open circular region, then there exists
a closed circular region $C' \subset C$ that still contains the
points $x_1^{(0)}, \ldots, x_n^{(0)}$.
Now, the standard proof of the Grace--Walsh--Szeg\"o coincidence theorem
(see e.g.\ \cite[Theorem 15.4]{Marden_66})
applies whenever $a_n \neq 0$.
If $a_n = 0$ and $C$ is a closed disc,
then by taking limits from the case $a_n \neq 0$
and using the compactness of $C$,
we can obtain the desired result.\footnote{
   {\sc Proof.}
   Let $P_\epsilon(x_1,\ldots,x_n) = P(x_1,\ldots,x_n) + \epsilon x_1 \cdots x_n$.
   Applying the case $a_n \neq 0$, we conclude that for each $\epsilon \neq 0$
   there exists $\xi_\epsilon \in C$ such that
   $P_\epsilon(x_1^{(0)},\ldots,x_n^{(0)}) =
    P_\epsilon(\xi_\epsilon,\ldots,\xi_\epsilon)$.
   Since $C$ is compact, the family $\{\xi_\epsilon\}$
   must have at least one limit point $\xi \in C$ as $\epsilon \to 0$.
   This limit point satisfies \reff{eq.thm.grace}.
}

Finally, suppose that $a_n = 0$ and $C$ is a closed half-plane.
Without loss of generality we can suppose that $C$ is the
closed right half-plane $\overline{H}$.
Let $\alpha = \max\{|\Im x_i^{(0)}| \colon\, 1 \le i \le n \}$
and $\beta = \max\{\Re x_i^{(0)} \colon\, 1 \le i \le n \}$.
Let $C_R$ be the closed disc of radius $\sqrt{R^2 + \alpha^2}$
centered at $R$.
Note that $C_R$ contains the points $x_1^{(0)}, \ldots, x_n^{(0)}$
whenever $R \ge \beta/2$.
Now apply the theorem to $C_R$:
we obtain a point $\xi_R \in C_R$ satisfying
$P(x_1^{(0)},\ldots,x_n^{(0)}) = P(\xi_R,\ldots,\xi_R)$.
If, for at least one $R$, the point $\xi_R$ lies in $C_R \cap \overline{H}$,
we are done;
and if not, then the points $\xi_R$ all lie in $C_R \setminus \overline{H}$,
and as $R \to \infty$ they must have a limit point $\xi$
in the interval $[-i\alpha,i\alpha]$ of the imaginary axis.
\qed

This theorem is usually applied as follows:
Starting from a univariate polynomial $Q(x) = \sum_k a_k x^k$
and an integer $n \ge \deg Q$,
we construct the unique symmetric multiaffine polynomial
$P(x_1,\ldots,x_n)$ satisfying $P(x,\ldots,x) = Q(x)$,
namely \reff{eq.symmetric_multiaffine}.
In this situation, Theorem~\ref{thm.grace}
implies that if $Q$ is nonvanishing in a circular region $C$,
then $P$ is nonvanishing in $C^n$
(provided that $n = \deg Q$ or $C$ is convex).

More generally, starting from a multivariate polynomial $Q$
in variables $\{x_e\}_{e \in E}$
and integers $n_e \ge \deg_e Q$,
we can introduce variables $\{ x_e^{(j)} \}_{e \in E,\, 1 \le j \le n_e}$
and form the unique multiaffine polynomial $P$ in all these variables
that is symmetric in each family $X_e = \{ x_e^{(j)} \}_{1 \le j \le n_e}$
separately and reduces to $Q$ when $x_e^{(j)} = x_e$ for all $j,e$.
We refer to $P$ as being obtained by the
\emph{$n_e$-fold polarization of $x_e$ in $Q$}, for all $e\in E$.
Then Theorem~\ref{thm.grace}
implies that if $\{C_e\}_{e \in E}$ are circular regions
and $Q(x)$ is nonvanishing when $x_e \in C_e$ for all $e$,
then $P$ is nonvanishing when $x_e^{(j)} \in C_e$ for all $j,e$
(provided that, for each $e$, either $n_e = \deg_e Q$ or $C_e$ is convex).
By this technique, we can often prove theorems for general polynomials $Q$
by reducing them to the special case of multiaffine polynomials $P$
(albeit in a larger number of variables).

\subsection{Algorithms}  \label{sec.algorithms}

It is worth remarking that the half-plane property is algorithmically
decidable, using quantifier-elimination methods
for the theory of real closed fields \cite{Caviness_98}.
Indeed, let $P \not\equiv 0$ be a polynomial in complex variables
$x_1,\ldots,x_n$.  Setting $x_k = a_k + ib_k$ and separating out
real and imaginary parts
\begin{subeqnarray}
   R(a_1,\ldots,a_n,b_1,\ldots,b_n)  & = & \real P(\{a_k + ib_k\})  \\[2mm]
   I(a_1,\ldots,a_n,b_1,\ldots,b_n)  & = & \imag P(\{a_k + ib_k\})
\end{subeqnarray}
the half-plane property for $P$
is immediately seen to be equivalent to the assertion
\be
   \neg (\exists a_1,\ldots,a_n,b_1,\ldots,b_n \in \R) \,
   (a_1 > 0) \wedge \cdots \wedge (a_n > 0)
   \wedge (R(a,b) = 0) \wedge (I(a,b) = 0)
\ee
in the first-order theory of the real numbers,
which is decidable according to a classic result of Tarski \cite{Tarski_51}.
Indeed, by the method of cylindrical algebraic decomposition (CAD)
\cite{Collins_75},
this computation can be performed in time $c_1^{c_2^{n}}$
for suitable constants $c_1$ and $c_2$ (``doubly exponential time'').
Moreover, some more recent algorithms
\cite{Grigorev_88,Heintz_90,Renegar_92,Basu_96}
require only a time $c^n$ (``singly exponential time'').
Unfortunately, this computation seems at present
to be unfeasible in practice even for $n=4$.

%

\section{The local half-plane property}  \label{sec_local}

For any element $e\in E$, the \emph{deletion of $e$ from $P$}
is the polynomial $P^{\drop e}$ on ground set $E \drop e$
that is obtained from $P$ by setting $x_e = 0$.
The \emph{contraction of $e$ from $P$}
is the polynomial $P^{/e}$ on ground set $E$
defined by $P^{/e} = \partial P/\partial x_{e}$.
Note that if $P$ is affine in $e$,
then $P^{/e}$ can also be considered as a polynomial on ground set $E \drop e$
(and we shall usually do so).
We say that an element $e\in E$ is a \emph{loop} of $P$
in case $P^{/e}\equiv 0$,
and is a \emph{coloop} of $P$ in case $P^{\drop e} \equiv 0$.

\begin{PROP}[Deletion/contraction]
 \label{prop3.1}
Let $P$ be a polynomial with the half-plane property.
Then, for every $e\in E(P)$,
both $P^{\drop e}$ and $P^{/e}$ have the half-plane property.
\end{PROP}

\proof
Proposition~\ref{prop.hurwitz}
implies that $P^{\drop e}$ has the half-plane property.
Proposition~\ref{prop.derivs}(a)
implies that $P^{/e}$ has the half-plane property.
\qed

{\bf Remark.}  If $P$ is affine in $e$, then the half-plane property
for the contraction $P^{/e}$ can alternatively be proven using
\be
   P^{/e}(x_{\neq e})  \;=\;  \lim_{x_e \to +\infty} {P(x) \over x_e}
\ee
[see \reff{eq.affine_in_e} below] and Hurwitz's theorem.

\bigskip

Let $P$ be a polynomial and $e \in E(P)$,
and suppose that $P$ is affine in $e$.
Then $P$ can be written in the form
\begin{equation}
   P(x)  \;=\;  P^{\drop e}(x_{\neq e}) + x_e P^{/e}(x_{\neq e})
   \;,
 \label{eq.affine_in_e}
\end{equation}
where we have used the shorthand $x_{\neq e} = \{x_f\}_{f \in E \drop e}$.
In this situation, we say that the pair $(P,e)$ satisfies
the \emph{local half-plane property} if,
whenever $\re x_{f} >0$ for all $f\in E\drop e$,
we have $P^{\drop e}(x_{\neq e})\neq 0$, $P^{/e}(x_{\neq e})\neq 0$ and
\begin{equation}
   \re\,\frac{P^{\drop e}(x_{\neq e})}{P^{/e}(x_{\neq e})} \;\geq\; 0  \;.
 \label{eq.local_HPP}
\end{equation}
Otherwise put, $(P,e)$ has the local half-plane property in case
\begin{itemize}
   \item[(a)] $e$ is neither a loop nor a coloop of $P$,
   \item[(b)] both $P^{\drop e}$ and $P^{/e}$ have the half-plane property, and
   \item[(c)] $P^{\drop e}/P^{/e}$ is real-part-positive on $H^{E\drop e}$.
\end{itemize}
Note also that by Lemma~\ref{lemma.RPP},
there is either strict inequality in \reff{eq.local_HPP}
or else $P^{\drop e}/P^{/e}$ is a nonzero pure imaginary constant
[i.e.\ $P(x) =(1+i\alpha x_{e}) P^{\drop e}(x_{\neq e})$
        for some nonzero real number $\alpha$].

%

We have the following fundamental result:

\begin{THM}
  \label{thm.local_hpp_single_e}
Let $P$ be affine in $e$.  Then the following are equivalent:
\begin{itemize}
   \item[(a)] $P$ has the half-plane property.
   \item[(b)] One of the following four mutually exclusive possibilities holds:
      \begin{itemize}
         \item[(i)]   $P^{/e} \equiv 0$ and $P^{\drop e} \equiv 0$;
         \item[(ii)]  $P^{/e} \equiv 0$, and $P^{\drop e} \not\equiv 0$
             has the half-plane property;
         \item[(iii)] $P^{\drop e} \equiv 0$, and $P^{/e} \not\equiv 0$
             has the half-plane property;
         \item[(iv)]  $(P,e)$ has the local half-plane property,
             i.e.\ $P^{\drop e} \not\equiv 0$ and $P^{/e} \not\equiv 0$
             both have the half-plane property, and
\begin{equation}
   \re\,\frac{P^{\drop e}(x_{\neq e})}{P^{/e}(x_{\neq e})} \;\geq\; 0
 \label{eq.local_HPP2}
\end{equation}
             for all $x_{\neq e} \in H^{E \drop e}$.
       \end{itemize}
   \item[(c)] One of the following three mutually exclusive possibilities holds:
      \begin{itemize}
         \item[(i)]   $P^{/e} \equiv 0$ and $P^{\drop e} \equiv 0$;
         \item[(ii)]  $P^{/e} \equiv 0$, and $P^{\drop e} \not\equiv 0$
             has the half-plane property;
         \item[(iii)] $P^{/e} \not\equiv 0$ has the half-plane property,
             and
\begin{equation}
             \re {P(x) \over P^{/e}(x_{\neq e})} \;>\; 0
\end{equation}
             for all $x \in H^E$.
       \end{itemize}
\end{itemize}
\end{THM}

\proof
(a)$\implies$(b):
By Proposition~\ref{prop3.1}, both $P^{\drop e}$ and $P^{/e}$
have the half-plane property.
If one or both of them is identically zero,
we are in one of the cases (i)--(iii).
If neither is identically zero,
then both $P^{\drop e}(x_{\neq e})$ and $P^{/e}(x_{\neq e})$
are nonzero for all $x_{\neq e} \in H^{E \drop e}$.
Therefore, for each such $x_{\neq e}$,
there is a unique value for $x_e$ such that $P(x) = 0$, namely
\begin{equation}
    x_e  \;=\; - \, \frac{P^{\drop e}(x_{\neq e})}{P^{/e}(x_{\neq e})}
    \;.
\end{equation}
Since $P$ has the half-plane property, we must have
$\re x_{e} \leq 0$, proving (\ref{eq.local_HPP2}).

(b)$\implies$(c):
This is trivial if either (b)(i) or (b)(ii) holds;
so suppose that either (b)(iii) or (b)(iv) holds,
and let $x \in H^E$.
Then $P^{/e}(x_{\neq e})\neq 0$, and
\begin{equation}
   \re\,\frac{P(x)}{P^{/e}(x_{\neq e})}  \;=\;
   \re \left(\frac{P^{\drop e}(x_{\neq e})}{P^{/e}(x_{\neq e})} + x_{e} \right)
   \;>\;  0
\end{equation}
since $\re [P^{\drop e}(x_{\neq e}) / P^{/e}(x_{\neq e})] \ge 0$
and $\re x_e > 0$.

(c)$\implies$(a):
In cases (c)(i) and (c)(ii) this is trivial.
And in case (c)(iii),
we clearly have $P(x) \neq 0$ for all $x \in H^E$,
so it is also trivial.
\qed

\begin{CORO}
   \label{cor.local_hpp}
Let $P$ be a multiaffine polynomial.
Then the following are equivalent:
\begin{itemize}
   \item[(a)] $P$ has the half-plane property.
   \item[(b)] There exists at least one $e \in E(P)$ for which
      condition (b) of Theorem~\ref{thm.local_hpp_single_e} holds.
   \item[(b${}'$)] For all $e \in E(P)$,
      condition (b) of Theorem~\ref{thm.local_hpp_single_e} holds.
   \item[(c)] There exists at least one $e \in E(P)$ for which
      condition (c) of Theorem~\ref{thm.local_hpp_single_e} holds.
   \item[(c${}'$)] For all $e \in E(P)$,
      condition (c) of Theorem~\ref{thm.local_hpp_single_e} holds.
\end{itemize}
\end{CORO}

Now let us drop the hypothesis that $P$ is affine in $e$.
We can still write $P$ in the form
\be
   P(x)  \;=\;  \sum_{k=0}^M P_k(x_{\neq e}) \, x_e^k
 \label{eq.sum_Pk}
\ee
where $M = \deg_e P$.
There is no longer any simple necessary and sufficient condition,
analogous to Theorem~\ref{thm.local_hpp_single_e},
for $P$ to have the half-plane property.
But Fettweis and Basu \cite[Theorem 18]{Fettweis_87}
have proven a very interesting {\em necessary}\/ condition:

\begin{theorem}[Fettweis and Basu \protect\cite{Fettweis_87}]
   \label{thm.fettweis}
Let $P$ be a polynomial with the half-plane property,
written in the form \reff{eq.sum_Pk} for some $e \in E(P)$.
Then each coefficient function $P_k$ has the half-plane property,
and moreover for each $r \ge 0$ we have:
\begin{itemize}
   \item[(a)]  It is impossible to have $P_r \not\equiv 0$,
      $P_{r+1} \equiv P_{r+2} \equiv \ldots \equiv P_{r+\ell} \equiv 0$,
      $P_{r+\ell+1} \not\equiv 0$ with $\ell \ge 2$.
   \item[(b)]  If $P_r \not\equiv 0$, $P_{r+1} \equiv 0$
      and $P_{r+2} \not\equiv 0$, then $P_{r+2}/P_r$ is a
      strictly positive constant.
   \item[(c)]  If $P_r \not\equiv 0$ and $P_{r+1} \not\equiv 0$,
      then $P_{r+1}/P_r$ is real-part-positive on $H^{E \drop e}$.
\end{itemize}
\end{theorem}

To prove Theorem~\ref{thm.fettweis},
we start with a lemma \cite[Lemma 4]{Fettweis_87}:

\begin{lemma}
   \label{lemma.fettweis}
Let $P(x) = \sum_{k=0}^M P_k(x_{\neq e}) \, x_e^k$ be a polynomial.
Let $0 \le r \le s \le M$ and define
\be
   Q(x)  \;=\;  \sum_{k=r}^s {k! \over (k-r)!} \, {(M-k)! \over (s-k)!}
                  \, P_k(x_{\neq e}) \, x_e^k
   \;.
\ee
If $P$ has the half-plane property, then so does $Q$.
\end{lemma}

\proof
Define
\begin{eqnarray}
   P_1(x)  & = &  {\partial^r \over \partial x_e^r} P(x)  \\
   P_2(x)  & = &  x_e^{M-r} P_1(x_{\neq e}, 1/x_e)        \\
   P_3(x)  & = &  {\partial^{M-s} \over \partial x_e^{M-s}} P_2(x)  \\
   Q(x)    & = &  x_e^{s} P_3(x_{\neq e}, 1/x_e)
\end{eqnarray}
By Proposition~\ref{prop.derivs} and
the fact that $x_e \mapsto 1/x_e$ maps the open right half-plane onto itself,
all these operations preserve the half-plane property.
\qed

\proofof{Theorem~\ref{thm.fettweis}}
We have already proven (Corollary~\ref{cor.Pcoeffs})
that all the coefficient functions $P_k$ have the half-plane property,
Now apply Lemma~\ref{lemma.fettweis} with $s=r+\ell+1$
[using $\ell=1$ in (b) and $\ell=0$ in (c)],
yielding a polynomial
\be
   Q(x) \;=\;
   \alpha P_r(x_{\neq e}) \, x_e^r  \,+\,
   \beta P_{r+\ell+1}(x_{\neq e}) \, x_e^{r+\ell+1}
\ee
with $\alpha, \beta > 0$ and $P_r, P_{r+\ell+1} \not\equiv 0$
having the half-plane property.
If $\ell \ge 2$, this polynomial $Q$ never has the half-plane property
[just choose any $x_{\neq e} \in H^{E \drop e}$,
 then there exists a solution to $Q(x) = 0$ with $\re x_e > 0$].
If $\ell = 1$, then $Q$ has the half-plane property
if and only if $P_{r+2}(x_{\neq e})/P_r(x_{\neq e}) > 0$
for all $x_{\neq e} \in H^{E \drop e}$;
but by the open mapping theorem,
this implies that $P_{r+2}/P_r$ is a constant function.
Finally, if $\ell=0$, then $Q$ has the half-plane property
if and only if $\re[P_{r+1}(x_{\neq e})/P_r(x_{\neq e})] \ge 0$
for all $x_{\neq e} \in H^{E \drop e}$.
\qed

\section{Constructions}   \label{sec_constructions}


In the previous section we saw that if a
polynomial $P$ has the half-plane property,
then so do the deletion $P^{\drop e}$ and the contraction $P^{/e}$,
for every element $e\in E$.  In this section we describe some
other constructions that preserve the half-plane property.

Let us recall that the \emph{support} of a
multiaffine polynomial $P(x)=\sum_{S \subseteq E} a_S x^S$
is $\supp(P) = \{S \subseteq E \colon\, a_S \neq 0 \}$.
We shall see later (Theorem~\ref{thm.matroidal})
that if $P \not\equiv 0$ is a homogeneous multiaffine polynomial
with the half-plane property,
then $\supp(P)$ is necessarily the set $\B(M)$ of bases of some matroid $M$
on the ground set $E(P)$.
(Theorem~\ref{thm.matroidal_non-multiaffine}
 generalizes this result to homogeneous polynomials
 that are not necessarily multiaffine.)
Therefore, while describing these constructions on polynomials,
we will also explain, in the homogeneous multiaffine case,
the corresponding operations on matroids.
Indeed, our choice of terminology for these constructions on polynomials
is motivated by analogy with standard terms in matroid theory.

\subsection{Deletion and contraction} \label{sec_constructions.1}

The following proposition is evident from the definitions:

\begin{PROP}
Let $P$ be a homogeneous multiaffine polynomial with $\supp(P) = \B(M)$
for a matroid $M$, and and let $e \in E(P)$.  Then:
\begin{itemize}
   \item[(a)] $e$ is a loop of $P$ if and only if $e$ is a loop of $M$.
   \item[(b)] $e$ is a coloop of $P$ if and only if $e$ is a coloop of $M$.
   \item[(c)] If $e$ is not a loop of $P$, then $\supp(P^{/e}) = \B(M/e)$,
the contraction of $M$ by $e$.
   \item[(d)] If $e$ is not a coloop of $P$, then $\supp(P^{\drop e}) =
\B(M\drop e)$, the deletion of $e$ from $M$.
   \item[(e)] If $e$ is not a loop of $P$, and $P$ is the basis generating
polynomial for $M$, then $P^{/e}$ is the basis generating polynomial for
$M/e$.
   \item[(f)] If $e$ is not a coloop of $P$, and $P$ is the basis
generating polynomial for $M$, then $P^{\drop e}$ is the basis generating
polynomial for $M\drop e$.
\end{itemize}
\end{PROP}

\subsection{Duality}

Given a multiaffine polynomial $P(x)=\sum_{S \subseteq E} a_S x^S$
on the ground set $E$, we define the {\em dual polynomial}\/
\begin{equation}
   P^{*}(x) \;=\; x^{E} P(1/x)  \;=\; \sum_{S \subseteq E} a_S x^{E \drop S}
   \;.
 \label{eq.def_dual}
\end{equation}
We have already encountered (in the Introduction) a special case
of the following principle:

\begin{PROP}[Duality]
 \label{prop.duality}
Let $P$ be a multiaffine polynomial.

If $P$ has the half-plane property,
then the dual polynomial $P^*$ also has the half-plane property.

If $\supp(P)= \B(M)$ for a matroid $M$, then $\supp(P^{*}) =\B(M^{*})$,
where $M^{*}$ is the matroid dual to $M$.

If $P$ is the basis generating polynomial for a matroid $M$,
then $P^*$ is the basis generating polynomial for $M^*$.
\end{PROP}

\proof
The half-plane property for $P^*$ follows immediately from
the fact that $x_e \mapsto 1/x_e$ maps the open right half-plane onto itself.
The remainder follows immediately from the definition of
the dual matroid $M^*$.
\qed

Let us also record, for future reference, some easy identities
relating duality to deletion and contraction:

\begin{LMA}
 \label{lemma.duality}
Let $P$ be a multiaffine polynomial on the ground set $E$,
and let $e \in E$.  Then
\begin{eqnarray}
    (P^{\drop e})^*   & = &   (P^*)^{/e}    \\
    (P^{/e})^*        & = &   (P^*)^{ \drop e}
\end{eqnarray}
where $P^{\drop e}$, $P^{/e}$, $(P^*)^{/e}$ and $(P^*)^{ \drop e}$
are all considered as having ground set $E \drop e$.
\end{LMA}

\subsection{Direct sum}

The next proposition is trivial, and has already been used implicitly:

\begin{PROP}
 \label{prop.product}
Let $P$ and $Q$ be polynomials with the half-plane property.
Then the pointwise product $PQ$
[considered on the ground set $E(P) \cup E(Q)$]
also has the half-plane property.
\end{PROP}

The special case where $P$ and $Q$ have disjoint ground sets
has a simple matroidal interpretation:

\begin{PROP}[Direct sum]
Let $P$ and $Q$ be homogeneous multiaffine polynomials
with $E(P)\cap E(Q)=\emptyset$.

If $\supp(P)=\B(M)$ and $\supp(Q)=\B(N)$ for matroids $M$ and $N$,
then $\supp(PQ)=\B(M\oplus N)$, the direct sum of $M$ and $N$.

If $P$ and $Q$ are the basis generating polynomials for matroids $M$ and $N$,
then $PQ$ is the basis generating polynomial for $M\oplus N$.
\end{PROP}

Clearly, if a matroid $M$ is the direct sum of two matroids $M_1,M_2$
of nonzero rank, then its basis generating polynomial $P_{\scrb(M)}$
is reducible.  A very strong version of the converse is true as well:

\begin{proposition} 
 \label{prop.irreducible}
Let $M$ be a connected matroid with ground set $E$,
and let $P$ be a multiaffine polynomial in the indeterminates
$\{x_e\}_{e \in E}$ with coefficients in an integral domain $R$.
If $\supp(P) = \B(M)$, then $P$ is irreducible over $R$.
\end{proposition}

\noindent
We begin with an easy but crucial lemma:

\begin{lemma}
  \label{lemma.irreducible}
Let $P_1$ and $P_2$ be nonzero polynomials in the indeterminates
$\{x_e\}_{e \in E}$ with coefficients in an integral domain $R$.
Suppose that $P_1 P_2$ is multiaffine.  Then:
\begin{itemize}
   \item[(a)]  There exist {\em disjoint}\/ subsets $E_1,E_2 \subseteq E$
      such that $P_i$ uses only the variables $\{x_e\}_{e \in E_i}$
      ($i=1,2$).
   \item[(b)]  $P_1$ and $P_2$ are both multiaffine.
\end{itemize}
\end{lemma}

\proof
Suppose there exists $e \in E$ such that both $P_1$ and $P_2$ use
the variable $x_e$.  For $i=1,2$, let $d_i \ge 1$ be the degree of $P_i$
in the variable $x_e$, and let $Q_i \neq 0$ be the coefficient of
$x_e^{d_i}$ in $P_i$, considered as an element of the polynomial ring
$R[x_{\neq e}]$.  Then $Q_1 Q_2 \neq 0$ because $R[x_{\neq e}]$ is an
integral domain \cite[Theorem III.5.1 and Corollary III.5.7]{Hungerford_74}.
But this shows that the coefficient of $x_e^{d_1 + d_2}$ in $P_1 P_2$
is nonzero, contradicting the hypothesis that $P_1 P_2$ is multiaffine
(since $d_1 + d_2 \ge 2$).  This proves (a);  and (b) is an easy consequence.
\qed

\proofof{Proposition~\ref{prop.irreducible}}
If $|E|=1$, the result is trivial, so assume henceforth that $|E| \ge 2$.
Suppose that $P = P_1 P_2$ where $P_1$ and $P_2$ 
are non-constant polynomials over $R$.
By Lemma~\ref{lemma.irreducible},
there exist disjoint subsets $E_1,E_2 \subseteq E$
such that $P_i$ uses only the variables $\{x_e\}_{e \in E_i}$.
Moreover, each term $x^S$ in $P$ arises from a unique pair of terms
in $P_1$ and $P_2$ (namely, $x^{S \cap E_1}$ and $x^{S \cap E_2}$);
and conversely, each pair of nonzero terms in $P_1$ and $P_2$
gives rise to a nonzero term in $P$ (because $R$ is an integral domain).
Since every term in $P$ has degree $r(M)$ [the rank of $M$],
it follows that there exist integers $r_1,r_2 \ge 1$ such that
every term in $P_i$ has degree $r_i$ ($i=1,2$).

Let $B$ be a basis of $M$. 
For each element $i \in E\setminus B$, there is a unique circuit 
$C(i,B)$ of $M$ contained in $B \cup \{i\}$.
Consider the bipartite graph $G$ with vertex classes $B$ and $E\setminus B$
where a vertex $i \in E \setminus B$ is adjacent to a vertex $j \in B$
if and only if $j \in C(i,B) \setminus i$. 
By a theorem of Cunningham \cite{Cunningham_73} and Krogdahl \cite{Krogdahl_77}
(see \cite[Lemma 10.2.8]{Oxley_92}),
since $|E| \ge 2$ and $M$ is connected, $G$ is a 
connected graph. Thus, by interchanging $E_1$ and $E_2$ if necessary, we may 
assume that $E_1 \setminus B$ contains an element $i$ such that 
$C(i,B)$ meets $E_2 \cap B$.  Suppose $j \in C(i,B) \cap E_2 \cap B$. 
Then $(B \setminus \{j\}) \cup \{i\}$ is a basis of $M$.
But $|[(B \setminus \{j\}) \cup \{i\}]\cap E_1| > |B \cap E_1|$,
so $P_1$ has a term of degree exceeding $r_1$:  a contradiction.
\qed

{\bf Remark.}
This proof is a natural extension of the proof
given by Duffin and Morley \cite[Theorem 8]{Duffin_78}
of the irreducibility of the basis generating polynomial
of a connected {\em binary}\/ matroid.

\subsection{Parallel connection, series connection, and 2-sum}
  \label{sec_constructions.4} 

Let $P$ and $Q$ be polynomials such that $E(P)\cap E(Q)=\{e\}$.
The \emph{parallel connection of $P$ and $Q$} is
defined to be
\begin{equation}
   (P|\!|Q)(x)   \;=\;
      P^{\drop e}(x)Q^{/e}(x) \,+\, P^{/e}(x)Q^{\drop e}(x)
                              \,+\, x_{e}P^{/e}(x)Q^{/e}(x)  \;.
 \label{eq.def.parallel}
\end{equation}
Note that if both $P$ and $Q$ are affine in $e$, then so is $P|\!|Q$;
and if both $P$ and $Q$ are multiaffine, then so is $P|\!|Q$
(the latter assertion uses the fact that $E(P)$ and $E(Q)$ meet only in $e$).

\begin{PROP}[Parallel connection]
 \label{prop.parallel}
Let $P$ and $Q$ be polynomials such that $E(P)\cap E(Q)=\{e\}$,
and suppose further that both $P$ and $Q$ are affine in $e$.
If both $P$ and $Q$ have the half-plane property,
then $P|\!|Q$ also has the half-plane property.

If, in addition, both $P$ and $Q$ are homogeneous and multiaffine,
with $\supp(P)=\B(M)$ and $\supp(Q)=\B(N)$ for matroids $M$ and $N$,
and $e$ is not a loop in at least one of $P$ and $Q$,
then $\supp(P|\!|Q) =\B(M| \! |N)$,
the parallel connection of $M$ and $N$.\footnote{
   The standard notation for the parallel connection of $M$ and $N$
   is $P(M,N)$,
   but in order to avoid confusion with polynomials,
   we write $M|\!|N$ instead.
}

If $P$ and $Q$ are the basis generating polynomials for matroids
$M$ and $N$, and $e$ is not a loop in at least one of $P$ and $Q$,
then $P|\!|Q$ is the basis generating polynomial for $M|\!|N$.
\end{PROP}

\firstproof
If $e$ is a loop of $P$, then $(P|\!|Q)=P^{\drop e}Q^{/e}$;
while if $e$ is a coloop of $P$, then $P|\!|Q=P^{/e}Q$.
Propositions~\ref{prop3.1} and \ref{prop.product} imply that
$P^{\drop e}Q^{/e}$ and $P^{/e}Q$ have the half-plane property.
If $e$ is a loop or coloop of $Q$, we argue analogously.
Hence we may assume that $e$ is not a loop or coloop
in either $P$ or $Q$.

{}From Theorem~\ref{thm.local_hpp_single_e}(a)$\implies$(b),
both $(P,e)$ and $(Q,e)$ have the local half-plane property.
Let $x$ be such that $\re x_{f} >0$ for all $f\in (E(P)\cup E(Q)) \drop e$.
Then $P^{\drop e}(x)$, $P^{/e}(x)$, $Q^{\drop e}(x)$
and $Q^{/e}(x)$ are all nonzero, by the local half-plane
property. If $(P|\!|Q)(x)=0$, we can solve this equation
for $x_{e}$, yielding
\begin{equation}
   x_{e}  \;=\;  - \frac{P^{\drop e}(x)}{P^{/e}(x)}
                 - \frac{Q^{\drop e}(x)}{Q^{/e}(x)} \;,
\end{equation}
which has nonpositive real part as a consequence of the
local half-plane property for $(P,e)$ and $(Q,e)$.
Therefore, $P|\!|Q$ has the half-plane property.

According to \cite[Proposition 7.1.13(ii)]{Oxley_92},
the bases of $M|\!|N$ are exactly those sets that can be obtained from
some $B_1 \in \B(M)$ and some $B_2 \in \B(N)$
in one of the following ways:
\begin{itemize}
  \item[(i)]   $(B_1 \cup B_2) \drop e$, where $e \notin B_1$ and $e \in B_2$;
  \item[(ii)]  $(B_1 \cup B_2) \drop e$, where $e \in B_1$ and $e \notin B_2$;
  \item[(iii)] $B_1 \cup B_2$, where $e \in B_1$ and $e \in B_2$.
\end{itemize}
These three classes correspond precisely to the three terms of
(\ref{eq.def.parallel}).
Moreover, it is easy to see that
each basis $B$ of $M|\!|N$ is obtained in this result
from a unique pair $(B_1,B_2) \in \B(M)\times\B(N)$.
The claims about support and about basis generating polynomials
follow immediately.
\qed

\secondproof
Since $P$ and $Q$ are affine in $e$, we have
\begin{eqnarray}
   P(x)  & = &   P^{\drop e}(x_{\neq e}) + x_e P^{/e}(x_{\neq e})  \\[2mm]
   Q(x)  & = &   Q^{\drop e}(x_{\neq e}) + x_e Q^{/e}(x_{\neq e})
\end{eqnarray}
and hence
\be
  {\partial \over \partial x_e} (PQ)(x)  \;=\;
      P^{\drop e}(x)Q^{/e}(x) \,+\, P^{/e}(x)Q^{\drop e}(x)
                              \,+\, 2x_{e}P^{/e}(x)Q^{/e}(x)  \;.
 \label{eq.para2}
\ee
By Propositions~\ref{prop.product} and \ref{prop.derivs}(a),
the polynomial \reff{eq.para2} has the half-plane property.
Replacing $x_e$ by $x_e/2$ (which preserves the half-plane property),
we obtain $(P|\!|Q)(x)$.
The rest of the proposition is proven as before.
\qed

{\bf Remark.}
The assertion in Proposition~\ref{prop.parallel}
that $P|\!|Q$ has the half-plane property
does not require the hypothesis that $E(P)\cap E(Q)=\{e\}$.
But without this hypothesis, $P|\!|Q$ may not be multiaffine,
even if $P$ and $Q$ both are.

\begin{CORO}[2-sum]
 \label{cor.2sum}
Let $P$ and $Q$ be polynomials such that $E(P)\cap E(Q)=\{e\}$,
and suppose further that $P$ and $Q$ are both affine in $e$.
If both $P$ and $Q$ have the half-plane property,
then $(P|\!|Q)^{\drop e}$ also has the half-plane property.

If, in addition, both $P$ and $Q$ are homogeneous and multiaffine,
with $\supp(P)=\B(M)$ and $\supp(Q)=\B(N)$ for matroids $M$ and $N$,
and $e$ is not a loop or coloop in either $P$ or $Q$,
then $\supp((P|\!|Q)^{\drop e}) =\B(M \oplus_{2} N)$,
the 2-sum of $M$ and $N$.

If $P$ and $Q$ are the basis generating polynomials for matroids
$M$ and $N$, respectively, and $e$ is not a loop or coloop
in either $P$ or $Q$, then $(P|\!|Q)^{\drop e}$ is the basis
generating polynomial for $M \oplus_{2} N$.
\end{CORO}

\proof
This follows immediately from Propositions~\ref{prop3.1} and
\ref{prop.parallel} and the definition of 2-sum.
\qed

Let $P$ and $Q$ be polynomials such that $E(P)\cap E(Q)=\{e\}$.
The \emph{series connection of $P$ and $Q$} is
defined to be
\begin{equation}
   (P \& Q)(x)   \;=\;
      P^{\drop e}(x)Q^{\drop e}(x)  \,+\,
      x_e P^{\drop e}(x)Q^{/e}(x) \,+\, x_e P^{/e}(x)Q^{\drop e}(x)  \;.
\end{equation}
Note that if both $P$ and $Q$ are affine in $e$, then so is $P \& Q$,
and indeed we have $P \& Q = (P^{*e}|\!|Q^{*e})^{*e}$
where ${}^{*e}$ denotes duality with respect to $x_e$ only
[that is, $P^{*e}(x_e,x_{\neq e}) = x_e P(1/x_e,x_{\neq e})$;
 this operation preserves the half-plane property, by the same argument
 used in the proof of Proposition~\ref{prop.duality}].
If both $P$ and $Q$ are multiaffine, then so is $P \& Q$,
and we have $P \& Q = (P^{*}|\!|Q^{*})^{*}$.

\begin{PROP}[Series connection]
 \label{prop.series}
Let $P$ and $Q$ be polynomials such that $E(P)\cap E(Q)=\{e\}$,
and suppose further that $P$ and $Q$ are both affine in $e$.
If both $P$ and $Q$ have the half-plane property,
then $P \& Q$ also has the half-plane property.

If, in addition, both $P$ and $Q$ are homogeneous and multiaffine,
with $\supp(P)=\B(M)$ and $\supp(Q)=\B(N)$ for matroids $M$ and $N$,
and $e$ is not a coloop in at least one of $P$ and $Q$,
then $\supp(P \& Q) =\B(M\&N)$,
the series connection of $M$ and $N$.\footnote{
   The standard notation for the series connection of $M$ and $N$
   is $S(M,N)$,
   but in order to avoid confusion with polynomials,
   we write $M \& N$ instead.
}

If $P$ and $Q$ are the basis generating polynomials for matroids $M$ and $N$,
and $e$ is not a coloop in at least one of $P$ and $Q$,
then $P \& Q$ is the basis generating polynomial for $M\&N$.
\end{PROP}

\proof
The proof is analogous to that of Proposition~\ref{prop.parallel}.
\qed

\subsection{Principal truncation and principal extension} \label{sec.trunc}

Let $P$ be a polynomial with ground set $E$,
and let $\{ \lambda_e \} _{e \in E}$ be nonnegative real numbers.
Then the polynomial
\begin{equation}
   \tr_{\lambda}P(x) \;=\; \sum_{e\in E}\lambda_{e} P^{/e}(x)
\end{equation}
is called a \emph{(weighted) principal truncation} of $P$.
If $a$ is a new element (not occurring in $E$), then
\begin{equation}
   \ext_{\lambda}P(x_{\cup a}) \;=\; P(x)+x_{a}\tr_{\lambda} P(x)
\end{equation}
is called a \emph{(weighted) principal extension} of $P$;
here we have used the shorthand $x_{\cup a} = \{x_e\}_{e \in E \cup a}$.

\begin{PROP}[Principal truncation/extension]
  \label{prop.truncation}
Let $P$ be a polynomial with ground set $E$,
and let $\{ \lambda_e \} _{e \in E}$ be nonnegative real numbers.
If $P$ has the half-plane property,
then both $\tr_{\lambda}P$ and $\ext_{\lambda}P$ have the half-plane property.
\end{PROP}

\proof
Proposition~\ref{prop.derivs}(a) asserts that $\tr_{\lambda}P$
has the half-plane property.
Proposition~\ref{prop.derivs}(b) together with
Theorem~\ref{thm.local_hpp_single_e}(b)$\implies$(a)
imply that $\ext_{\lambda}P$ has the half-plane property.
\qed

%
%

The terms ``principal truncation'' and ``principal extension''
are taken from matroid theory;
they are justified by the following result:

\begin{PROP}
  \label{prop.truncation_support}
Let $P$ be a homogeneous multiaffine polynomial with ground set $E$,
and let $\{ \lambda_e \} _{e \in E}$ be nonnegative real numbers,
not all of which are zero.
Suppose that $\supp(P)=\B(M)$ for a matroid $M$ of nonzero rank,
and let $F$ be the closure in $M$ of the set
$S=\{e \colon\; \lambda_{e}>0\}$.
Then $\supp(\tr_{\lambda}P)=\B(\tr_{F}(M))$,
where $\tr_{F}(M)$ is the principal truncation of $M$ with respect to $F$.
Also, $\supp(\ext_{\lambda}P)=\B(M+_{F}a)$,
where $M+_{F}a$ is the principal extension of $M$ on $F$.
\end{PROP}

\proof
The bases of the truncation $\tr_{F}(M)$ are all the sets of the form
$B\drop f$ where $B\in \B(M)$ and $f\in B \cap F$.
The sets that appear in the support of $\tr_{\lambda}P$
are all the sets of the form
$B\drop g$ where $B\in \B(M)$ and $g\in B \cap S$.
Since $S\subseteq F$, we clearly have
$\supp(\tr_{\lambda}P)\subseteq\B(\tr_{F}(M))$.
Conversely, suppose $B\in \B(M)$ and $f\in B \cap F$,
and let $r$ be the rank of $M$.
Since $\rank(B\drop f)=r-1$ and $\rank((B\drop f)\cup S)=r$,
the independence augmentation axiom implies that
there is an element $g\in S$ such that $B'=(B\drop f)\cup\{g\}$
is a basis of $M$.
Since $B\drop f = B'\drop g$,
we conclude that $\B(\tr_{F}(M))\subseteq\supp(\tr_{\lambda}P)$,
as desired.

The bases of the extension $M+_{F}a$ are given by
\begin{equation}
   \B(M+_{F}a) \;=\;
  \B(M) \,\cup\,  \{B\cup\{a\} \colon\; B\in \B(\tr_{F}(M))\} \;.
\end{equation}
{}From this and the previous paragraph, the claim about the
support of $\ext_{\lambda}P$ follows readily.
\qed

In the situation of Propositions~\ref{prop.truncation}
and \ref{prop.truncation_support},
even if $P$ is the basis generating polynomial of the matroid $M$,
it need not be the case that $\tr_\lambda P$ is the
basis generating polynomial of $\tr_F(M)$,
since the coefficients need not all equal 1.
Rather, for this to be the case we must
choose the weights $\{\lambda_f\}_{f \in F}$ so as to have
\begin{equation}
   \sum_{f\in F \colon\, B\cup\{f\}\in\B(M)} \lambda_f  \;=\; 1
 \label{eq.nice_principal}
\end{equation}
for all $B \in \B(\tr_F(M))$.
Since $\B(\tr_F(M))$ is generally much larger than $F$,
this system of equations is usually overdetermined.
However, if this system does have a nonnegative solution
$\{\lambda_f\}_{f \in F}$, then for this choice of $\lambda$
it follows that $\tr_\lambda P$ is the basis generating polynomial
for $\tr_F(M)$,
and that $\ext_{\lambda} P$ is the basis generating polynomial for $M+_{F}a$.
When this happens we say that the principal truncation (or extension)
of $M$ by $F$ is \emph{nice}.

\begin{proposition}
 \label{prop.uniform_nice}
For every $1\leq r\leq n$, the uniform matroid $U_{r-1,n}$ is a nice
principal truncation of $U_{r,n}$.  Therefore, for all $0\leq r\leq n$,
the basis generating polynomial of $U_{r,n}$ has the half-plane property.
\end{proposition}

\proof
Fix $n$, and let $E = \{1,\ldots,n\}$ be the common ground set
of all the matroids $U_{r,n}$.
We have $U_{r-1,n}=\tr_E(U_{r,n})$.
Setting $\lambda_e = 1/(n-r+1)$ for all $e\in E$ gives a solution to the
the system \reff{eq.nice_principal}, so $U_{r-1,n}$ is a nice principal
truncation of $U_{r,n}$.  The basis generating polynomial of $U_{n,n}$
is $x^E$, which clearly has the half-plane property.  Reverse induction
on $r$ (from $n$ down to $0$) using Proposition~\ref{prop.truncation}
now shows that the basis generating polynomial of $U_{r,n}$
has the half-plane property for all $0\leq r\leq n$.
\qed

Here are a few more examples:

\bexam
  \label{exam.nice_principal_MK4plus}
The matroid $M(K_4)^+$ shown in Figure~\ref{fig.rank3.seven.4point}
(see Appendix~\ref{app_matroids})
is the principal extension of the graphic matroid $M(K_4)$
obtained by adding one point freely to one of the 3-point lines of $M(K_4)$,
say 123.  An easy computation shows that this principal extension is nice:
it suffices to take $\lambda_1 = \lambda_2 = \lambda_3 = 1/2$
and $\lambda_4 = \lambda_5 = \lambda_6 = 0$.
It follows (by Theorem~\ref{thm1.1} and Proposition~\ref{prop.truncation})
that $M(K_4)^+$ has the half-plane property.
\eexam

\bexam
  \label{exam.non-nice_principal_MK4pluse}
The matroid $M(K_4) + e$ (Figure~\ref{fig.rank3.seven.4point})
is the free extension of the graphic matroid $M(K_4)$,
i.e.\ the principal extension of $M(K_4)$
obtained by adding one point freely to the flat $F=E$.
An straightforward computation
shows that this principal extension is not nice:
the system \reff{eq.nice_principal} has no solution.

Let us remark that $M(K_4) + e$ fails the half-plane property
(Example~\ref{sec_counterexamples}.\ref{exam.counterexamples.4}).
This shows that non-nice principal extension does {\em not}\/ in general
preserve the half-plane property.
\eexam

\bexam
  \label{exam.nice_principal_F7m4pe}
The matroid $F_7^{-4} + e$ (Figure~\ref{fig_F7m4pe_etc})
is the free extension of the matroid $F_7^{-4}$ (Figure~\ref{fig_F7etc}),
i.e.\ the principal extension of $F_7^{-4}$
obtained by adding one point freely to the flat $F=E$.
An easy computation shows that this free extension is nice:
it suffices to take $\lambda_1 = 0$ and
$\lambda_2 = \ldots = \lambda_7 = 1/4$.
It follows that $F_7^{-4} + e$ has the half-plane property if $F_7^{-4}$ does.
Unfortunately, we have been unable to prove whether or not
$F_7^{-4}$ has the half-plane property;
but our numerical tests (Section~\ref{sec_numerical}) suggest that it does
(and that $F_7^{-4} + e$ does also).
\eexam

\bexam
  \label{exam.non-nice_principal_Q7}
The matroid $Q_7$ (Figure~\ref{fig.rank3.seven.4point})
is the free extension of the matroid $Q_7 \setminus 7$.
Curiously, the system \reff{eq.nice_principal} does have a solution,
but this (unique) solution fails to be nonnegative:
it is $\lambda_1 = -1/6$, $\lambda_2 = \lambda_3 = 1/2$,
$\lambda_4 = \lambda_5 = \lambda_6 = 1/3$.
So this free extension is not nice.
Nevertheless, it turns out that $Q_7$ has the half-plane property
(Corollary~\ref{cor.nice_transversal} and
 Example~\ref{sec_transversal}.\ref{exam.transversal.7}).
\eexam

\bigskip

We have written a {\sc Mathematica} program {\tt niceprincipal.m}
to test a principal extension for niceness
by solving the linear system \reff{eq.nice_principal};
it is available as part of the electronic version of this paper
at arXiv.org.


\subsection{Principal cotruncation and principal coextension} \label{sec.cotr}

Let $P$ be a {\em multiaffine}\/ polynomial with ground set $E$,
and let $\{ \lambda_e \} _{e \in E}$ be nonnegative real numbers.
Then the polynomial
\begin{equation}
   \cotr_{\lambda}P(x) \;=\;
   \sum_{e\in E} \lambda_{e} x_e P^{\drop e}(x_{\neq e})
\end{equation}
is called a \emph{(weighted) principal cotruncation} of $P$.
If $a$ is a new element (not occurring in $E$), then
\begin{equation}
   \coext_{\lambda}P(x_{\cup a}) \;=\; \cotr_{\lambda}P(x) \,+\, x_a P(x)
\end{equation}
is called a \emph{(weighted) principal coextension} of $P$.

Note that by Lemma~\ref{lemma.duality}, we have
\be
   (\cotr_{\lambda}P)^*  \;=\;
   \sum_{e\in E} \lambda_{e} (P^*)^{/e}(x_{\neq e})
   \;=\;  \tr_{\lambda} P^*
\ee
(to understand the first equality, observe that $P^{\drop e}$
 is considered as a polynomial on $E \drop e$,
 while the rest are considered as polynomials on $E$).
Likewise,
\be
   (\coext_{\lambda}P)^*  \;=\;
   x_a (\cotr_{\lambda} P)^* \,+\, P^*  \;=\;
   x_a \tr_{\lambda} P^* \,+\, P^* \;=\;
   \ext_{\lambda} P^*
\ee
(taking care again to notice that some of these are polynomials on $E \cup a$
 and others on $E$).
So principal cotruncation and coextension are indeed the duals of
principal truncation and extension.

\begin{PROP}[Principal cotruncation/coextension]
  \label{prop.cotruncation}
Let $P$ be a multiaffine polynomial with ground set $E$,
and let $\{ \lambda_e \} _{e \in E}$ be nonnegative real numbers.
If $P$ has the half-plane property,
then both $\cotr_{\lambda}P$ and $\coext_{\lambda}P$
have the half-plane property.
\end{PROP}

\proof
The result is trivial if $P \equiv 0$, so assume $P \not\equiv 0$.
We have
\begin{subeqnarray}
   {\cotr_{\lambda}P(x) \over P(x)}
   & = &
   \sum\limits_{e \in E} \lambda_e \, {x_e P^{\drop e}(x_{\neq e}) \over P(x)}
       \\[2mm]
   & = &
   \sum\limits_{e \in E} \lambda_e \,
      {x_e P^{\drop e}(x_{\neq e})  \over
       P^{\drop e}(x_{\neq e}) \,+\, x_e  P^{/e}(x_{\neq e})}
       \slabel{eq.cotr.b} \\[2mm]
   & = &
   \sum\limits_{e \in E} \lambda_e
       \left( {1 \over x_e} \,+\,
              {P^{/e}(x_{\neq e}) \over  P^{\drop e}(x_{\neq e})}
       \right)^{\! -1}
   \;,  \slabel{eq.cotr.c}
  \label{eq.cotr}
\end{subeqnarray}
which is real-part-positive by Theorem~\ref{thm.local_hpp_single_e}.
[Use \reff{eq.cotr.b} in place of \reff{eq.cotr.c} for those elements $e$
 having $P^{\drop e} \equiv 0$ or $P^{/e} \equiv 0$.]
So $\cotr_{\lambda}P$ has the half-plane property.
Moreover,
\be
   {\coext_{\lambda}P(x) \over P(x)}  \;=\;
   {\cotr_{\lambda}P(x) \over P(x)}  \,+\, x_a
\ee
is real-part-positive, so $\coext_{\lambda}P$ has the half-plane property.
\qed

The principal cotruncation of a matroid $M$ on a set $D$, which appears in the
next proposition, is also called the principal lift of $M$ on $D$
\cite[pp.~160--161]{Brylawski_86}.


\begin{PROP}
  \label{prop.cotruncation_support}
Let $P$ be a homogeneous multiaffine polynomial with ground set $E$,
and let $\{ \lambda_e \} _{e \in E}$ be nonnegative real numbers,
not all of which are zero.
Suppose that $\supp(P) = \B(M)$ for a matroid $M$ of rank less than $|E|$,
let $S = \{e\colon\;\lambda_e > 0\}$, and let $D$ be the subset of $E$ such
that $E \drop D$ is the union of all circuits contained in $E \drop S$.
Then $\supp(\cotr_{\lambda} P) = \B(\cotr_D(M))$,
where $\cotr_D(M)$ is the principal cotruncation of $M$ on $D$.
Also, $\supp(\coext_{\lambda} P) = \B(M \times_D a)$,
where $M \times_D a$ is the principal coextension of $M$ on $D$.
\end{PROP}

\proof
The bases of the cotruncation $\cotr_D(M)$ are all the sets of the form
$B \cup d$ where $B \in \B(M)$ and $d \in D \drop B$.  The sets
that appear in the support of $\cotr_{\lambda} P$ are all the sets
of the form $B \cup f$ where $B \in \B(M)$ and $f \in S \drop B$.
Since $S \subseteq D$, it follows that
$\supp(\cotr_{\lambda} P) \subseteq \B(\cotr_D(M))$.

Now suppose that $B \in \B(M)$ and $f \in D \drop B$. If $f \in S$, then
$B\cup f \in \B(\cotr_D(M))$. Now suppose that $f \notin S$.
Clearly $B \cup f$ contains a circuit $C$ of $M$ and $f \in C$. If $C$ avoids
$S$, then $C$ is a circuit contained in $E \drop S$,
so $C \subseteq E \drop D$, which is a contradiction.
Thus $C$ contains an element $g$ of $S$.
Then $(B \drop g) \cup f \in \B(M)$ and
$(B \drop g) \cup f \cup g = B \cup f \in \B(\cotr_D(M))$.
Hence
$\supp(\cotr_{\lambda} P) \supseteq \B(\cotr_D(M))$
and so equality holds.

The bases of the principal coextension $M\times_D a$ are given by
\be
   \B(M \times_D a)  \;=\;
      \{B \cup a \colon\; B \in \B(M)\} \,\cup\, \B(\cotr_D(M))
   \;.
\ee
{}From this and the previous paragraph, the claim about the support of
$\coext_{\lambda} P$ follows readily.
\qed

Note that even if $P$ is the basis generating polynomial of the matroid $M$,
it need not be the case that $\cotr_\lambda P$ is the
basis generating polynomial of $\cotr_D(M)$,
since the coefficients need not all equal 1.
Rather, for this to be the case we must
choose the weights $\{\lambda_d\}_{d \in D}$ so as to have
\begin{equation}
   \sum_{d\in D \colon\, B \drop d \in\B(M)} \lambda_d  \;=\; 1
 \label{eq.nice_principal_cotruncation}
\end{equation}
for all $B \in \B(\cotr_D(M))$.
When this system of equations has a nonnegative solution
$\{\lambda_d\}_{d \in D}$,
we say that the principal cotruncation (or coextension)
of $M$ by $D$ is \emph{nice}.

\begin{proposition}
 \label{prop.uniform_cotruncation_nice}
For every $1\leq r\leq n$, the uniform matroid $U_{r,n}$ is a nice
principal cotruncation of $U_{r-1,n}$.  Therefore, for all $0\leq r\leq n$,
the basis generating polynomial of $U_{r,n}$ has the half-plane property.
\end{proposition}

\proof
Fix $n$, and let $E = \{1,\ldots,n\}$ be the common ground set
of all the matroids $U_{r,n}$.
We have $U_{r,n}=\cotr_E(U_{r-1,n})$.
Setting $\lambda_e = 1/r$ for all $e\in E$ gives a solution to the
the system \reff{eq.nice_principal_cotruncation},
so $U_{r,n}$ is a nice principal cotruncation of $U_{r-1,n}$.
The basis generating polynomial of $U_{0,n}$ is $1$,
which clearly has the half-plane property.
Induction on $r$ using Proposition~\ref{prop.cotruncation}
shows that the basis generating polynomial of $U_{r,n}$
has the half-plane property for all $0\leq r\leq n$.
\qed

\medskip

{\bf Remark.}
The construction in this section was inspired by our analysis
of Fettweis' \cite{Fettweis_90} proof of the half-plane property
for the uniform matroids $U_{r,n}$ (see Section~\ref{sec_uniform}).

\subsection{Multiaffine part and full-rank matroid union}

Let $P(x) = \sum_{\bf m} a_{\bf m} x^{\bf m}$
be a polynomial with ground set $E$.
For any subset $A \subseteq E$, define
\begin{equation}
   P^{\flat A}(x)  \;=\;
   \sum_{{\bf m} \colon\, m_e \le 1 \,\forall e \in A}
   a_{\bf m} x^{\bf m}
   \;.
\end{equation}
When $A=E$ we shall write simply $P^\flat$,
and shall call it the {\em multiaffine part}\/ of $P$.
Clearly the map $P \mapsto P^{\flat A}$ is linear and idempotent.

\begin{PROP}
   \label{prop.flat}
If $P$ has the half-plane property,
then for any $A \subseteq E$,
$P^{\flat A}$ has the half-plane property.
\end{PROP}

\proof
Consider first the univariate case $|E|=1$ with $A=E$.
The case $P \equiv 0$ is trivial,
so let $P \not\equiv 0$ be a univariate polynomial of degree $k$
with the half-plane property:  we have
$P(x) = C \prod_{i=1}^k (x + \alpha_i)$ with $C \neq 0$
and $\re \alpha_i \ge 0$.
The cases $k=0,1$ are trivial, so assume $k \ge 2$.
Then $P^\flat(x) = a_0 + a_1 x$ where
\begin{subeqnarray}
   a_0  & = &  C \prod\limits_{i=1}^k \alpha_i  \\
   a_1  & = &  C \sum\limits_{j=1}^k \prod\limits_{i \neq j} \alpha_i
\end{subeqnarray}
If one or more of the $\alpha_i$ are zero, then $a_0 = 0$
and $P^\flat$ has the half-plane property.
Otherwise,
\begin{equation}
   {a_1 \over a_0}  \;=\;  \sum\limits_{j=1}^k  {1 \over \alpha_j}
\end{equation}
has nonnegative real part, so that $P^\flat$ again has the half-plane property.

Consider next the multivariate case with $|A| = 1$, i.e.\ $A = \{a\}$.
Let us write $P(x) = \sum_k P_k(x_{\neq a}) x_a^k$,
so that $P^{\flat a}(x) = P_0(x_{\neq a}) + P_1(x_{\neq a}) x_a$.
Applying the univariate case to the variable $x_a$,
we conclude that for each $x_{\neq a} \in H^{E \drop a}$,
the univariate polynomial $P^{\flat a}(\,\cdot\,, x_{\neq a})$
is either identically zero or else nonvanishing on $H$.
By Corollary~\ref{cor.Pcoeffs},
both $P_0$ and $P_1$ have the half-plane property,
hence each one is either identically zero
or else nonvanishing on $H^{E \drop a}$.
If both $P_0$ and $P_1$ are identically zero,
then $P^{\flat a}$ is identically zero and hence has the half-plane property.
If one or both of them is nonvanishing on $H^{E \drop a}$,
it follows that $P^{\flat a}(\,\cdot\,, x_{\neq a})$
cannot be identically zero for any $x_{\neq a} \in H^{E \drop a}$,
so $P^{\flat a}$ must be nonvanishing on $H^E$.

The general case with $|A| = \{a_1,\ldots,a_N\}$
is now obtained by successively applying the case $|A|=1$
to each $a_i$.
\qed

{\bf Remark.}  Note that $P^\flat$ can be identically zero
even if $P$ is not identically zero (and has the half-plane property),
e.g.\ $P(x) = x^2$.

\bigskip

For matroids $M$ and $N$, the \emph{union} $M\vee N$
is the matroid on ground set $E(M) \cup E(N)$
that has as its collection of independent sets
all sets of the form $I_{1}\cup I_{2}$, with $I_{1}$ independent in $M$
and $I_{2}$ independent in $N$.  The bases of $M\vee N$ are its
maximal independent sets, as always.
It is easy to see that $\rank(M\vee N) = \rank(M) + \rank(N)$
holds if and only if there exist $B_1 \in \B(M)$ and $B_2 \in \B(N)$
with $B_1 \cap B_2 = \emptyset$;
and in this case the bases of $M\vee N$ are precisely the sets
$B_1 \cup B_2$ with $B_1 \in \B(M)$, $B_2 \in \B(N)$
and $B_1 \cap B_2 = \emptyset$.
When $\rank(M\vee N) = \rank(M) + \rank(N)$ holds,
we call this {\em full-rank matroid union}\/.

\begin{PROP}[Matroid union]
  \label{prop.union}
Let $P$ and $Q$ be homogeneous multiaffine polynomials
with the same-phase property
such that $\supp(P)=\B(M)$ and $\supp(Q) =\B(N)$ for matroids $M$ and $N$.
\begin{itemize}
   \item[(a)] If $\rank(M\vee N) = \rank(M) + \rank(N)$,
        then $\supp((PQ)^\flat) = \B(M\vee N)$.
   \item[(b)] If $\rank(M\vee N) < \rank(M) + \rank(N)$,
        then $(PQ)^\flat \equiv 0$.
\end{itemize}
\end{PROP}

\proof
If $P(x)=\sum_{S}a_{S}x^{S}$ and $Q(x)=\sum_{S} b_{S}x^{S}$,
then $(PQ)^\flat(x) = \sum_{S} c_{S}x^{S}$ where
\begin{equation}
   c_S  \;=\;  \sum_{\begin{scarray}
                        T \cup U = S \\
                        T \cap U = \emptyset
                     \end{scarray}
                    }
                a_T b_U   \;.
 \label{eq.C_union}
\end{equation}
If $\rank(M\vee N) < \rank(M) + \rank(N)$,
then there are no pairs $T \in \supp(P)$ and $U \in \supp(Q)$
with $T \cap U = \emptyset$, so $(PQ)^\flat \equiv 0$.
If $\rank(M\vee N) = \rank(M) + \rank(N)$,
then the sets $S$ for which (\ref{eq.C_union})
contains at least one nonzero summand are precisely the bases of $M\vee N$;
and the same-phase property for $P$ and $Q$ ensures that all the
contributing terms $a_{T}b_{U}$ have the same phase, so that $c_{S}\neq 0$.
\qed

Note, however, that even if $P$ and $Q$ are the basis generating
polynomials of $M$ and $N$, it need not be the case that
$(PQ)^\flat$ is a multiple of the basis generating polynomial of $M\vee N$,
since the coefficients need not all be equal.
For example,
let $M$ and $N$ be matroids on the set $\{1,2,3,4\}$ with basis generating
polynomials $P(x) = x_1 + x_2 + x_3$ and $Q(x) = x_1 + x_2 + x_4$.
Then $(PQ)^{\flat}(x) = 2x_1x_2 + x_1x_3 + x_1x_4 +x_2x_3 +x_2x_4 + x_3x_4$.

\subsection{``Folding mod 2'' and convolution}  \label{sec.folding}

Let $P(x) = \sum_{\bf m} a_{\bf m} x^{\bf m}$
be a polynomial with ground set $E$.
For any subset $A \subseteq E$, define
\begin{equation}
   P^{\natural A}(x)  \;=\;
   \sum_{\bf m} a_{\bf m} \prod_{e \in A} x_e^{m_e \bmod 2}
                          \prod_{e \in E \drop A} x_e^{m_e}
\end{equation}
where $m_e \bmod 2 = 0$ or 1.
When $A=E$ we shall write simply $P^\natural$;
note that $P^\natural$ is multiaffine.
Clearly the map $P \mapsto P^{\natural A}$ is linear and satisfies
$(P_1 \cdots P_k)^{\natural A} =
 (P_1^{\natural A} \cdots P_k^{\natural A})^{\natural A}$.

\begin{PROP}
   \label{prop.natural}
If $P$ has the half-plane property and is not identically zero,
then for any $A \subseteq E$,
$P^{\natural A}$ has the half-plane property and is not identically zero.
\end{PROP}

\proof
Consider first the univariate case $|E|=1$.
Let $P \not\equiv 0$ be a univariate polynomial of degree $k$
with the half-plane property, so that
$P(x) = C \prod_{i=1}^k (x + \alpha_i)$ with $C \neq 0$
and $\re \alpha_i \ge 0$.
The cases $k=0,1$ are trivial, so assume $k \ge 2$.
By repeated use of the identity
$(P_1 P_2)^\natural = (P_1^\natural P_2^\natural)^\natural$
[where $P_1$ is a linear factor and $P_2$ is the product of the
 remaining factors],
it suffices to prove the case $k=2$.
In this case we have
$P^\natural(x) = C[(1+\alpha_1 \alpha_2) + (\alpha_1 + \alpha_2)x]$.
It is easy to see that $1+\alpha_1 \alpha_2$ and $\alpha_1 + \alpha_2$
cannot both be zero, so $P^\natural \not\equiv 0$.
If one of them is zero, then $P^\natural$ clearly has the
half-plane property.
If both of them are nonzero, then
\begin{equation}
   {1+\alpha_1 \alpha_2  \over  \alpha_1 + \alpha_2}
   \;=\;
   {1  \over  \alpha_1 + \alpha_2}
   \,+\,
   \left( {1 \over \alpha_1} + {1 \over \alpha_2} \right)^{-1}
\end{equation}
has nonnegative real part [this happens also if one of the $\alpha_i$
is zero], so $P^\natural$ again has the half-plane property.

Consider now the multivariate case, i.e.\ suppose that $P(x) \neq 0$
whenever $\re x_e > 0$ for all $e \in E$.
Let $A = \{a_1,\ldots,a_N\}$ and
define $A_j = \{a_1,\ldots, a_j\}$ for $1 \le j \le N$.
Applying the univariate case successively to each $x_{a_i}$
(with all other $x_e$ held fixed in the open right half-plane),
we produce a sequence of polynomials
$P^{\natural A_1}, P^{\natural A_2}, \ldots,
 P^{\natural A_N} \equiv P^{\natural A}$,
each of which has the half-plane property and is not identically zero.
\qed

{\bf Remarks.}
1.  There obviously exist polynomials $P \not\equiv 0$ for which
$P^\natural \equiv 0$, e.g.\ $P(x) = 1-x^2$.
Proposition~\ref{prop.natural} asserts that this cannot happen when
$P$ has the half-plane property.

2.  Suppose that in the definition of $P^{\natural A}$ we replace
``mod 2'' by ``mod $n$''.  Then the map $P \mapsto P^{\natural A}$ is
still linear and satisfies
$(P_1 \cdots P_k)^{\natural A} =
 (P_1^{\natural A} \cdots P_k^{\natural A})^{\natural A}$,
but for $n \ge 3$ it does not preserve the half-plane property.
A univariate counterexample is $P(x) = (1 + i\epsilon x)^n$
with $\epsilon$ real, for which
$P^\natural(x) = (1 + i\epsilon x)^n - (i\epsilon)^n (x^n -1)$.
For small $\epsilon \neq 0$, the roots of $P^\natural$ are
$x = \epsilon^{-1} i/(1-\omega) + O(\epsilon^{n-1})$
where $\omega \neq 1$ is an $n$th root of unity;
and for $n \ge 3$ at least one of these roots has strictly positive
real part.

3.  If $P(z_1,z_2) = a + bz_1 + cz_2 + dz_1 z_2$,
the {\em Asano contraction}\/ of $P$ is the univariate polynomial
$\widetilde{P}(z) = a + dz$.
The Asano contraction lemma \cite{Asano_70,Ruelle_71,Hinkkanen_97}
states that if $P$ is nonvanishing in the unit bidisc $|z_i| < 1$ ($i=1,2$),
then $\widetilde{P}$ is nonvanishing in the unit disc $|z| < 1$;
this lemma and its generalizations play an important role
in the derivation of Lee--Yang theorems in statistical mechanics
\cite{Ruelle_71,Ruelle_73,Lieb-Sokal_81,Hinkkanen_97}.
The Asano lemma can be obtained as a corollary of
Proposition~\ref{prop.natural},
using the M\"obius transformation
\begin{equation}
   z_i  \;=\;  {1-x_i \over 1+x_i}   \quad\Longleftrightarrow\quad
   x_i  \;=\;  {1-z_i \over 1+z_i}
 \label{eq.mobius}
\end{equation}
that maps the right half-plane $\re x_i > 0$
onto the unit polydisc $|z_i| < 1$.
It suffices to note that
\begin{eqnarray}
   Q(x_1,x_2)  & \equiv &
   (1+x_1)(1+x_2) P\Bigl( {1-x_1 \over 1+x_1}, {1-x_2 \over 1+x_2} \Bigr)
      \nonumber \\
   & = &
   (a+b+c+d) + (a-b+c-d) x_1 + (a+b-c-d) x_2 + (a-b-c+d) x_1 x_2
      \nonumber \\
   & \equiv &  \alpha + \beta x_1 + \gamma x_2 + \delta x_1 x_2
\end{eqnarray}
while
\be
   R(x) \;\equiv\;
   (1+x) \widetilde{P}\Bigl( {1-x \over 1+x} \Bigr)
   \;=\; (a+d) + (a-d) x
   \;=\; {\alpha + \delta \over 2} \,+\, {\beta + \gamma \over 2} x
   \;,
\ee
so that $2R$ can be obtained from $Q$ by setting $x_1 = x_2 = x$
and applying the $\natural$ operation.

\bigskip

Now let $P(x)=\sum_{S}a_{S}x^{S}$ and $Q(x)=\sum_{S} b_{S}x^{S}$
be multiaffine polynomials on the same ground set $E$.
The \emph{convolution} $P * Q$ of $P$ and $Q$
is defined to be
\begin{equation}
   (P * Q)(x)  \;=\;  \sum_{S,T}a_{S}b_{T}x^{S\triangle T} \;,
\end{equation}
where $\triangle$ denotes the symmetric difference of sets.\footnote{
   Strictly speaking, it is the {\em coefficients}\/
   $a = \{a_S\}$ and $b = \{b_S\}$ --- which are complex-valued
   functions on the group $2^E$ of subsets of $E$
   with respect to symmetric difference ---
   that are being convoluted here.
}
It is easily seen that convolution is a bilinear, commutative
and associative operation,
and that $P_1 * \cdots * P_k = (P_1 \cdots P_k)^\natural$.
It therefore follows immediately from
Propositions~\ref{prop.product} and \ref{prop.natural} that:

\begin{PROP}[Convolution]
 \label{prop.convolution}
Let $P$ and $Q$ be multiaffine polynomials with the half-plane property.
Then $P * Q$ also has the half-plane property.
Moreover, if $P$ and $Q$ are not identically zero,
then $P * Q$ is not identically zero.
\end{PROP}

{\bf Remarks.}
1.  The ring of multiaffine polynomials with respect to
convolution does have divisors of zero,
e.g.\ $P(x) = 1+x$ and $Q(x) = 1-x$ have $P*Q \equiv 0$.
But Proposition~\ref{prop.convolution} asserts that this cannot happen
when $P$ and $Q$ both have the half-plane property.

2. Proposition~\ref{prop.duality} (duality) is a special case
of Proposition~\ref{prop.convolution}, obtained by taking $Q(x) = x^E$.

3. Proposition~\ref{prop.convolution} is closely related to
Hinkkanen's composition theorem \cite{Hinkkanen_97},
which states that if $P(z)=\sum_{S}a_{S}z^{S}$ and $Q(z)=\sum_{S} b_{S}z^{S}$
are multiaffine polynomials (on the same ground set $E$)
that are nonvanishing in the unit polydisc,
then the \emph{Schur--Hadamard product}
\begin{equation}
   (P \circ Q)(z)  \;=\;  \sum_{S} a_{S} b_{S} z^{S}
\end{equation}
is either nonvanishing in the unit polydisc or else identically zero.
Indeed, Proposition~\ref{prop.convolution}
and Hinkkanen's composition theorem are interderivable,
using the M\"obius transformation \reff{eq.mobius}
that maps the right half-plane $\re x_i > 0$
onto the unit polydisc $|z_i| < 1$.
To see this, note first that, given a polynomial $P(x)=\sum_{S}a_{S}x^{S}$,
we can define a polynomial
\begin{equation}
   \widehat{P}(z)  \;=\;  (1+z)^E P\Bigl({1-z \over 1+z}\Bigr)
                   \;=\;  \sum_T \widehat{a}_T z^T
\end{equation}
where the coefficients $\{\widehat{a}_T\}$ are obtained from $\{a_S\}$
by Fourier transformation on the group $2^E$:
\begin{equation}
   \widehat{a}_T  \;=\;  \sum_{S \subseteq E} a_S \, (-1)^{|S \cap T|}
   \;.
\end{equation}
Clearly, $P$ has the half-plane property if and only if
$\widehat{P}$ has the ``unit polydisc property'';
and since Fourier transformation carries elementwise product
onto convolution (and vice versa), we have
$\widehat{P*Q} = \widehat{P} \circ \widehat{Q}$.
In this way, Proposition~\ref{prop.convolution}
can be given an alternate proof
as a corollary of Hinkkanen's composition theorem;
and, conversely, our proof of Proposition~\ref{prop.convolution}
yields an alternate proof of Hinkkanen's composition theorem.

\bigskip

In Proposition~\ref{prop.union} we showed that
matroid union can be obtained via the $\flat$ operation
in case $\rank(M\vee N) = \rank(M) + \rank(N)$.
Let us now show that, in this same full-rank situation,
matroid union can also be obtained via convolution
(followed by taking the leading part):

\begin{CORO}[Matroid union, revisited]
Let $P$ and $Q$ be homogeneous multiaffine polynomials
with the same-phase property
such that $\supp(P)=\B(M)$ and $\supp(Q) =\B(N)$ for matroids $M$ and $N$.
If $\rank(M\vee N) = \rank(M) + \rank(N)$,
then $(P * Q)^{\sharp} = (PQ)^\flat$
and hence $\supp((P * Q)^{\sharp})=\B(M\vee N)$.
\end{CORO}

\proof
As noted earlier,
the equality $\rank(M\vee N) = \rank(M) + \rank(N)$ implies that
there exist pairs $(S,T)$ with
$S \in \B(M)$, $T \in \B(N)$ and $S\cap T=\emptyset$;
moreover, the bases of $M\vee N$ are precisely the
corresponding sets $S \cup T$.
Now, for each $B \subseteq E$ of cardinality $\rank(M\vee N)$,
the coefficient $c_{B}$ of $B$ in $P * Q$
is the sum of $a_{S}b_{T}$ over all pairs $(S,T)\in\B(M)\times\B(N)$
for which $S\cup T=B$
(note that this implies $S \cap T =\emptyset$);
and the same-phase property for $P$ and $Q$
ensures that all the terms $a_{S}b_{T}$ have the same phase,
so that $c_{B}\neq 0$ whenever $B \in \B(M\vee N)$.
It follows that $P*Q$ indeed has degree $\rank(M\vee N)$
[rather than some lower degree],
so that $(P * Q)^{\sharp} = (PQ)^\flat$.
\qed

\section{Necessary and sufficient condition}  \label{sec_necessary}

\subsection{General case}  \label{sec_necessary.1}

Let $P$ be a polynomial in $n$ complex variables with complex coefficients.
Gregor \cite[Lemma 6]{Gregor_81} has given
a necessary and sufficient condition for $P$ to have the half-plane property:

\begin{proposition}[Gregor \protect\cite{Gregor_81}]
 \label{prop.Gregor}
Let $P$ be a polynomial in $n$ complex variables,
and fix $k \ge \deg P$.
For $x,y \in \R^n$, define the univariate
polynomial $\widehat{p}_{x,y}(\zeta) = \zeta^k P(\zeta x + \zeta^{-1} y)$.
Then the following are equivalent:
\begin{itemize}
   \item[(a)]  $P$ has the half-plane property and is not identically zero.
   \item[(b)]  For all $x,y \ge 0$ with $x+y > 0$,
       $\widehat{p}_{x,y}$ has the half-plane property
       and is not identically zero
       [that is,
       all the roots of $\widehat{p}_{x,y}$ lie in the closed left half-plane].
   \item[(c)]  For all $x,y > 0$,
       $\widehat{p}_{x,y}$ has the half-plane property
       and is not identically zero.
\end{itemize}
\end{proposition}

\proof
(a)$\implies$(b):
Let $x,y \ge 0$ with $x+y > 0$.
Then $\zeta \in H$ implies $\zeta x + \zeta^{-1} y \in H^n$,
so $\widehat{p}_{x,y}(\zeta) \neq 0$.

(b)$\implies$(c) is trivial.

(c)$\implies$(a):
Let $z = (z_1,\ldots,z_n)$ with with $\real z_i > 0$ for all $i$.
Choose $\zeta \in H$
in such a way that all the $z_i$ lie in the interior of the convex cone
generated by $\zeta$ and $\zeta^{-1}$
[i.e.\ choose $\zeta$ so that
 $\max |\arg z_i| < |\arg \zeta| < \pi/2$];
then there exist numbers $x_i, y_i > 0$
with $z_i = \zeta x_i + \zeta^{-1} y_i$.
It follows that $P(z) = \zeta^{-k} \widehat{p}_{x,y}(\zeta) \neq 0$.
\qed

If $P$ is homogeneous, this result can be simplified slightly:

\begin{proposition}
 \label{prop.generalrank}
Let $P$ be a homogeneous degree-$r$ polynomial in $n$ complex variables.
For $x,y \in \R^n$, define the univariate
polynomial $p_{x,y}(\zeta) = P(\zeta x + y)$.
Then the following are equivalent:
\begin{itemize}
   \item[(a)]  $P$ has the half-plane property and is not identically zero.
   \item[(b)]  For all $x,y \ge 0$ with $x+y > 0$,
       all the roots (real or complex) of $p_{x,y}$ lie in $(-\infty,0]$.
   \item[(c)]  For all $x,y > 0$,
       all the roots (real or complex) of $p_{x,y}$ lie in $(-\infty,0]$.
\end{itemize}
\end{proposition}

\proof
(a)$\implies$(b):
Let $x,y \ge 0$ with $x+y > 0$.
If $\zeta_0 \in \C \setminus (-\infty,0]$,
then we can find $\alpha,\beta \in H$ such that $\alpha/\beta = \zeta_0$.
But this means that
$p_{x,y}(\alpha/\beta) = \beta^{-r} P(\alpha x + \beta y) \neq 0$.

(b)$\implies$(c) is trivial.

(c)$\implies$(a):
Let $z = (z_1,\ldots,z_n)$ with $\real z_i > 0$ for all $i$.
Choose $\alpha,\beta \in H$
in such a way that all the $z_i$ lie in the interior of the convex cone
generated by $\alpha$ and $\beta$,
i.e.\ such that there exist numbers $x_i, y_i > 0$
with $z_i = \alpha x_i + \beta y_i$.
Then $P(z) = P(\alpha x + \beta y) = \beta^r p_{x,y}(\alpha/\beta) \neq 0$.
%
\qed

\bigskip
 
{\bf Remark.}  The necessary and sufficient conditions
given in Propositions~\ref{prop.Gregor} and \ref{prop.generalrank}
are algorithmically testable in time $c_1^{c_2^{n}}$
using cylindrical algebraic decomposition,
and in time $c^n$ using more recent algorithms
(see Section~\ref{sec.algorithms}).
But this computation seems thoroughly unfeasible in practice,
at least at present.

\subsection{Rank-2 case}

The drawback of the necessary-and-sufficient characterization
contained in Proposition~\ref{prop.generalrank} is, of course,
that the sufficient condition (b) or (c) is not easy to verify,
because of the universal quantification over $x,y \gtmaybeeq 0$.
However, in the rank-2 case
we can obtain an easily-checkable alternative condition.
Let $A = \{ a_{ij} \} _{i,j=1}^n$ be a symmetric matrix
with nonnegative real entries, i.e.\ $a_{ij} = a_{ji} \ge 0$.
(We shall see in Theorem~\ref{thm.same-phase}
 that the restriction to nonnegative real entries is no loss of generality.)
Define the homogeneous degree-2 polynomial in $n$ variables,
\be
 \label{def_PA}
   P_A(z)  \;=\;  \half \sum_{i,j=1}^n  a_{ij} z_i z_j
           \;=\;  \smhalf z^{\rm T} A z
   \;.
\ee
Let $\lambda_1(A) \ge \lambda_2(A) \ge \cdots \ge \lambda_n(A)$
be the eigenvalues of $A$.
We then have:

\begin{theorem}[due largely to \protect\cite{Gregor_81,Bapat_97}]
 \label{thm.rank2}
Let $A$ be a symmetric matrix with nonnegative real entries.
Then the following are equivalent:
\begin{itemize}
   \item[(a)] 
      Either $P_A(H^n) = \{0\}$ (i.e.\ $P_A \equiv 0$)
      or else $P_A(H^n) = \C \setminus (-\infty,0]$.
   \item[(b)] $P_A$ has the half-plane property.
   \item[(c)] $\lambda_2(A) \le 0$.
   \item[(d)] If $x,y \in \R^n$ with $y^{\rm T} A y \ge 0$,
       then $(x^{\rm T} Ay)^2 \ge (x^{\rm T} Ax) (y^{\rm T} Ay)$.
   \item[(e)] If $x,y \in \R^n$ with $x,y \ge 0$,
       then $(x^{\rm T} Ay)^2 \ge (x^{\rm T} Ax) (y^{\rm T} Ay)$.
\end{itemize}
\end{theorem}

\proof
If $A=0$, then (a)--(e) are all true;
so let us consider the nontrivial case $A \neq 0$
[hence $\lambda_1(A) > 0$ because $A$ is symmetric].
By the Perron--Frobenius theorem \cite[p.~66]{Gantmacher_59},
there exists an eigenvector $x^{(1)} \ge 0$ corresponding to
the eigenvalue $\lambda_1$.
Now extend this to an orthonormal basis of real eigenvectors
$x^{(1)}, \ldots, x^{(n)}$  satisfying $A x^{(k)} = \lambda_k x^{(k)}$.


(a)$\implies$(b) is trivial.

(b)$\implies$(c):
Suppose that $\lambda_2(A) > 0$.
Perturb $x^{(1)}, x^{(2)}$ slightly to get vectors
$\widetilde{x}^{(1)} > 0$ and $\widetilde{x}^{(2)}$
satisfying $\widetilde{x}^{(1) \rm T} A \widetilde{x}^{(2)} = 0$,
$\widetilde{x}^{(1) \rm T} A \widetilde{x}^{(1)} > 0$ and
$\widetilde{x}^{(2) \rm T} A \widetilde{x}^{(2)} > 0$.\footnote{
   For example, let $v$ be any vector with strictly positive components,
   and define $\widetilde{x}^{(1)} = x^{(1)} + \epsilon v$,
   $\widetilde{x}^{(2)} = x^{(2)} + \delta v$
   with $\delta = -\epsilon (\widetilde{x}^{(2) \rm T} A v) /
            (\widetilde{x}^{(1) \rm T} A v + \epsilon v^{\rm T} A v)$
   and $\epsilon > 0$ sufficiently small.
}
Consider the vector $z = \widetilde{x}^{(1)} + i\alpha \widetilde{x}^{(2)}$
where $\alpha \in \R$.  Then
\be
   P_A(z)  \;=\;  \smhalf (\widetilde{x}^{(1) \rm T} A \widetilde{x}^{(1)}
              \,-\, \alpha^2 \widetilde{x}^{(2) \rm T} A \widetilde{x}^{(2)})
   \;,
\ee
which vanishes when
$\alpha = \pm \sqrt{(\widetilde{x}^{(1) \rm T} A \widetilde{x}^{(1)}) /
                    (\widetilde{x}^{(2) \rm T} A \widetilde{x}^{(2)})}$.
So $P_A$ does not have the half-plane property.

(c)$\implies$(d):
If $x^{\rm T} Ax \le 0$ or $y^{\rm T} Ay = 0$, the assertion is trivial;
so we can assume that $x^{\rm T} Ax > 0$ and $y^{\rm T} Ay > 0$.
We must have $x^{(1)} \cdot x \neq 0$,
since otherwise [thanks to the negative-semidefiniteness of $A$
on the orthogonal complement of $x^{(1)}$]
we would have $x^{\rm T} Ax \le 0$;
and likewise we must have $x^{(1)} \cdot y \neq 0$.
Now define
\be
   g(\alpha)   \;=\;   (x+\alpha y)^{\rm T} A (x+\alpha y)   \;.
\ee
We have $g(0) = x^{\rm T} Ax > 0$.
On the other hand, for $\alpha = -(x^{(1)} \cdot x)/(x^{(1)} \cdot y)$
we have $x^{(1)} \cdot (x+\alpha y) = 0$ and hence $g(\alpha) \le 0$
[again by the negative-semidefiniteness of $A$ on $x^{(1)\perp}$].
So the quadratic equation $g(\alpha) = 0$ has a real solution,
which implies that its discriminant is nonnegative,
i.e.\ that $(x^{\rm T} Ay)^2 \ge (x^{\rm T} Ax) (y^{\rm T} Ay)$.

(d)$\implies$(e) is trivial.


(e)$\implies$(a):
Suppose that $P_A(z) = -C \le 0$ for some vector $z = (z_1,\ldots,z_n)$
with $\real z_i > 0$ for all $i$.
Choose $\theta$ so that $|\arg z_i| < \theta < \pi/2$ for all $i$.
Then all the $z_i$ lie in the interior of the convex cone
generated by $e^{i\theta}$ and $e^{-i\theta}$,
i.e.\ there exist vectors $x, y > 0$
such that $z = e^{i\theta} x + e^{-i\theta} y$.
Then
\be
   P_A(z) \;=\;  \smhalf [(x^{\rm T} A x) e^{2i\theta}  \,+\,
                          2 (x^{\rm T} A y) \,+\,
                          (y^{\rm T} A y) e^{-2i\theta}]
          \;=\;  -C  \;.
\ee
Defining $\zeta = e^{2i\theta}$, it follows that the quadratic equation
\be
 \label{dimpliesa.quadratic}
   (x^{\rm T} A x) \zeta^2 \,+\, 2 [(x^{\rm T} A y) + C] \zeta \,+\,
                                 (y^{\rm T} A y)   \;=\;  0
\ee
has a root $\zeta \in \C \setminus (-\infty,0]$.
Since all the coefficients of \reff{dimpliesa.quadratic} are positive,
this root $\zeta$ cannot be positive, so it must have a nonzero
imaginary part.  This means that the discriminant of \reff{dimpliesa.quadratic}
is negative, i.e.\ $[(x^{\rm T} Ay)+C]^2 < (x^{\rm T} Ax) (y^{\rm T} Ay)$.
Since $x^{\rm T} Ay \ge 0$ and $C \ge 0$,
it follows that $(x^{\rm T} Ay)^2 < (x^{\rm T} Ax) (y^{\rm T} Ay)$,
contradicting hypothesis (e).
Therefore $P_A(H^n) \subseteq \C \setminus (-\infty,0]$.

On the other hand, $P_A(1,\ldots,1) > 0$
and $P_A(\lambda,\ldots,\lambda) = \lambda^2 P_A(1,\ldots,1)$;
so as $\lambda$ runs over the open right half-plane,
$P_A(\lambda,\ldots,\lambda)$ runs over all of $\C \setminus (-\infty,0]$.
\qed

{\bf Remarks.}
1. The proof of (c)$\implies$(d) is taken from \cite[Theorem 4.4.6]{Bapat_97},
where it is also proven that (d)$\implies$(c).
After completing this proof, we learned that
the equivalence of (b), (c) and (e) had already been proven
by Fiedler and Gregor \cite[Theorems 5 and 8]{Gregor_81}.
Numerous equivalent characterizations are given in
\cite[Theorem 4.4.6]{Bapat_97}.

2. The reference in Theorem~\ref{thm.rank2}(c) to the eigenvalues
of the matrix $A$ (i.e.\ treating $A$ as a {\em linear operator}\/)
is in fact quite misleading.
What we really have here is a {\em quadratic form}\/
$Q(x) = \sum_{i,j=1}^n  a_{ij} x_i x_j$
on a real vector space $V$ ($= \R^n$)
equipped with a distinguished convex cone $C$
(i.e.\ the vectors with nonnegative components).
We then extend $Q$ bilinearly to the complexification $V+iV$,
and say that it has the ``half-plane property'' if it is either
identically zero or else is nonvanishing on the set
$C^\circ + iV$ (here $C^\circ$ denotes the interior of $C$).
Condition (c) is the statement that the inertia $(n_+,n_0,n_-)$
of the real quadratic form $Q$ satisfies $n_+ \le 1$.


\bigskip

Let us now specialize to the 0-1-valued case:
we are given a 2-uniform set system $\scrs$
on $[n]$ --- that is, a simple graph $G$ with vertex set $[n]$
and edge set $\scrs$ --- and we define
\begin{subeqnarray}
   a_{ij}  \;=\;  a_{ji}  & = &  \cases{1    & if $\{i,j\} \in \scrs$  \cr
                                        \noalign{\vskip 2pt}
                                        0    & if $\{i,j\} \notin \scrs$  \cr
                                       }
     \\[2mm]
   a_{ii}  & = & 0 \quad \hbox{for all $i$}
\end{subeqnarray}
Thus $A$ is simply the adjacency matrix of $G$,
and $P_A$ is the generating polynomial $P_\scrs$.
Theorem~\ref{thm.rank2} then becomes:

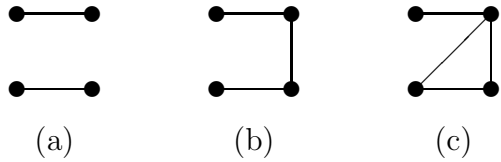
\begin{figure}
\setlength{\unitlength}{1cm}
\begin{center}
\begin{tabular}{c@{\qquad\qquad}c@{\qquad\qquad}c}
\begin{picture}(1,1)(0,0)
   \put(0,0){\circle*{0.2}}
   \put(0,1){\circle*{0.2}}
   \put(1,0){\circle*{0.2}}
   \put(1,1){\circle*{0.2}}
   \put(0,0){\line(1,0){1}}
   \put(0,1){\line(1,0){1}}
\end{picture}
&
\begin{picture}(1,1)(0,0)
   \put(0,0){\circle*{0.2}}
   \put(0,1){\circle*{0.2}}
   \put(1,0){\circle*{0.2}}
   \put(1,1){\circle*{0.2}}
   \put(0,0){\line(1,0){1}}
   \put(0,1){\line(1,0){1}}
   \put(1,0){\line(0,1){1}}
\end{picture}
&
\begin{picture}(1,1)(0,0)
   \put(0,0){\circle*{0.2}}
   \put(0,1){\circle*{0.2}}
   \put(1,0){\circle*{0.2}}
   \put(1,1){\circle*{0.2}}
   \put(0,0){\line(1,0){1}}
   \put(0,1){\line(1,0){1}}
   \put(1,0){\line(0,1){1}}
   \put(0,0){\line(1,1){1}}
\end{picture}
\\[3mm]
(a) & (b) & (c)
\end{tabular}
\end{center}
\vspace{-3mm}
\caption{
   Forbidden induced subgraphs in Corollary~\ref{cor.rank2}(d).
}
 \label{fig1}
\end{figure}

\begin{corollary}
 \label{cor.rank2}
Let $G$ be a simple graph on the vertex set $[n]$ with edge set $\scrs$,
and let $A$ be the corresponding adjacency matrix.
Then the following are equivalent:
\begin{itemize}
   \item[(a)]  $P_A$ has the half-plane property.
   \item[(b)]  $\lambda_2(A) \le 0$.
   \item[(c)]  $G$ is a complete multipartite graph (possibly empty)
      plus possible isolated vertices.
   \item[(d)]  $G$ contains none of the graphs in Figure~\ref{fig1}
      as induced subgraphs.
   \item[(e)]  $\scrs$ is the set of bases for a rank-2 matroid
      on the ground set $[n]$.
\end{itemize}
\end{corollary}

\proof
If $\scrs = \emptyset$ (i.e.\ $A=0$), then (a)--(e) are all true;
so let us consider the nontrivial case $\scrs \neq \emptyset$.
The equivalence of (a) and (b) follows from Theorem~\ref{thm.rank2}.
The equivalence of (b), (c) and (d)
is a result of Smith \cite[Theorem 1]{Smith_70}
(see also \cite[Theorem 6.7]{Cvetkovic_80} and \cite[Theorem 2.1]{Hoffman_72}).

Consider, finally, the basis exchange property for a matroid:
\begin{quote}
\begin{itemize}
   \item[(B2)]  If $B_1, B_2 \in \scrs$ and $i \in B_1 \setminus B_2$,
      then there exists $j \in B_2 \setminus B_1$ such that
      $(B_1 \setminus i) \cup j \in \scrs$.
\end{itemize}
\end{quote}
For a rank-2 matroid, this property holds trivially when
$|B_1 \cap B_2| = 1$ or 2.
When $B_1 \cap B_2 = \emptyset$ (so that $|B_1 \cup B_2| = 4$),
straightforward checking of cases shows that property (B2)
holds if and only if none of the configurations in Figure~\ref{fig1}
occurs on $B_1 \cup B_2$.
This proves the equivalence of (d) and (e).
\qed

\begin{corollary}
 \label{cor.rank2a}
All rank-2 matroids have the half-plane property,
as do all corank-2 matroids.
\end{corollary}

\proof
This is an immediate consequence of Corollary~\ref{cor.rank2}(e)$\implies$(a)
and Proposition~\ref{prop.duality}.
\qed


\section{Necessary conditions I: Same-phase property}
  \label{sec_same-phase}

In this section and the next, we derive some {\em necessary}\/ conditions
that a homogeneous polynomial with the half-plane property must satisfy.

Let us begin with a bit of motivation:
Many ``counting'' results in combinatorics find a more natural context
when they are generalized to allow nonnegative real weights,
rather than restricting all weights to be 0 or 1.
For example, the max-flow-min-cut theorem arises
by generalizing Menger's theorem in this way.
Now, in {\em some}\/ cases it is natural to go further and allow
complex weights, but in some cases it is not.
For example, in the Heilmann--Lieb theorem
(see Section~\ref{sec_Heilmann-Lieb} below),
the complex vertex weights $\{x_i\}_{i \in V}$ are the main point,
but it is essential that the edge weights $\{\lambda_e\}_{e \in E}$
be nonnegative.
Roughly speaking, it seems natural that {\em variables}\/
in multivariate polynomials be promoted to complex variables,
but it is often necessary that the {\em coefficients}\/
in the polynomial be nonnegative real numbers.\footnote{
   Of course, given a polynomial expression, it is not always obvious
   which objects should be considered to be ``variables''
   and which should be absorbed into the ``coefficients''.
   So this admittedly vague principle may in some cases be best understood
   backwards, i.e.\ as an injunction to treat as ``variables''
   those objects that can naturally be made complex,
   and as ``coefficients'' those that cannot.
}

With this background in mind, we shall prove that
if a {\em homogeneous}\/ polynomial
$P(x) = \sum_{\bf m} a_{\bf m} x^{\bf m}$ has the half-plane property,
then {\em necessarily}\/ all the coefficients $a_{\bf m}$ are nonnegative
modulo an overall multiplicative phase factor.

Let us say that a pair of nonzero complex numbers $a,b$
{\em have the same phase}\/ if $a/b \in (0,\infty)$.
Obviously, having the same phase is an equivalence relation
on $\C \setminus 0$;  and a collection $\{a_\alpha\}$
of nonzero complex numbers all have the same phase if and only if there exists
$\theta \in \R$ such that $e^{-i\theta} a_\alpha > 0$ for all $\alpha$.
We say that a polynomial $P(x) = \sum_{\bf m} a_{\bf m} x^{\bf m}$
has the {\em same-phase property}\/
if all the nonzero coefficients $a_{\bf m}$ have the same phase.

\begin{theorem}
   \label{thm.same-phase}
Let $P(x) = \sum_{\bf m} a_{\bf m} x^{\bf m}$
be a polynomial in $n$ complex variables that is homogeneous of degree $r$.
If $P$ has the half-plane property, then all the nonzero coefficients
$a_{\bf m}$ have the same phase.
\end{theorem}

\proof
The proof is by induction on $n$.
If $n=1$, we must have $P(x) = a_r x^r$, so the theorem is trivial.
So assume $n \ge 2$.
Let $M$ ($\le r$) be the degree of $x_n$ in $P(x)$,
and let us write
\be
   P(x)  \;=\;  \sum_{k=0}^M x_n^k P_k(x_1,\ldots,x_{n-1})
\ee
where of course
\be
   P_k(x_1,\ldots,x_{n-1})  \;=\;
   {1 \over k!} \, (\partial/\partial x_n)^k P(x) \Bigr| _{x_n=0}
   \;.
\ee
Clearly each coefficient of $P$ corresponds to exactly one coefficient
of exactly one $P_k$, and vice versa.
Now, each $P_k$ is a homogeneous polynomial (of degree $r-k$)
in $n-1$ variables;
and by Corollary~\ref{cor.Pcoeffs} it has the half-plane property.
So, by the inductive hypothesis, all its nonzero coefficients
have the same phase.

To complete the proof, we need only show that these phases
are the same for all $k$.
Applying Proposition~\ref{prop.generalrank} with
$x=(0,\ldots,0,1)$ and $y=(1,\ldots,1,0)$,
we get
\be
   p_{x,y}(\zeta)  \;=\;  P(1,\ldots,1,\zeta)  \;=\;
   \sum_{k=0}^M P_k(1,\ldots,1) \, \zeta^k
   \;.
\ee
Since $p_{x,y}$ has only real nonpositive zeros,
it can be written as $p_{x,y}(\zeta) = a \prod_{k=1}^M (\zeta + c_k)$
for some complex number $a$ and nonnegative real numbers $c_1,c_2,\ldots,c_M$;
so its coefficients $P_k(1,\ldots,1)$ must all have the same phase.
\qed

{\bf Remark.}
It is easy to see that if $P$ is homogeneous of degree 1,
then $P$ has the half-plane property
{\em if and only if}\/ all its nonzero coefficients have the same phase.

\bigskip

We can generalize Theorem~\ref{thm.same-phase} to a class of
not-necessarily-homogeneous polynomials.
Let us say that a polynomial
$P(x) = \sum_{\bf m} a_{\bf m} x^{\bf m}$
has {\em definite parity}\/
if all the nonzero monomials have total degree of the same parity
(i.e.\ $a_{\m}, a_{\m'} \neq 0$ implies that $|\m|\equiv|\m'|$ mod 2);
this is equivalent to saying that the polynomial $P$ is either even or odd.
We then have:

\begin{theorem}
   \label{thm.same-phase_parity}
Let $P(x) = \sum_{\bf m} a_{\bf m} x^{\bf m}$
be a polynomial in $n$ complex variables that has definite parity
and the half-plane property.
Then all the nonzero coefficients $a_{\bf m}$ have the same phase.
\end{theorem}

\noindent
We begin the proof of Theorem~\ref{thm.same-phase_parity} with a lemma:

\begin{lemma}
   \label{lemma.same-phase_parity}
Let $P$ and $Q$ be polynomials in $n$ complex variables,
and suppose that $P+cQ$ has the same-phase property for all real $c>0$.
Then either
\begin{itemize}
   \item[(a)] all the nonzero coefficients of $P$ and $Q$ have the same phase
\end{itemize}
or else
\begin{itemize}
   \item[(b)] $P$ and $Q$ each have the same-phase property,
      and $P = \alpha Q$ for some nonzero complex number $\alpha$.
\end{itemize}
\end{lemma}

\proof
Letting $c \to 0$ (resp.\ $c \to +\infty$),
we deduce that $P$ (resp.\ $Q$) has the same-phase property.
Let us write $P(x) = \sum_\m  a_{\m} x^\m$ and $Q(x) = \sum_\m  b_{\m} x^\m$.
Now there are two cases:
\begin{itemize}
   \item[1)]  There is a multi-index $\m$ such that both $a_{\m}$ and $b_{\m}$
      are nonzero.
   \item[2)]  There is no such multi-index $\m$.
\end{itemize}
In case 2, the fact that $P+Q$ has the same-phase property
implies that all the nonzero coefficients of $P$ and $Q$ have the same phase.

In case 1, choose a multi-index $\m$ such that both $a_{\m}$ and $b_{\m}$
are nonzero, and define $r_{\m} = a_{\m}/b_{\m}$.
If there is no multi-index $\n \neq \m$ with either $a_{\n}$ or $b_{\n}$
nonzero, then $P(x) = a_{\m} x^\m$ and $Q(x) = b_{\m} x^\m$,
which proves assertion (b) [with $\alpha = r_{\m}$].
So assume that there is at least one such multi-index $\n$,
and consider an arbitrary one of them.
Then, for all positive real numbers $c$,
the complex numbers $a_{\m} + cb_{\m}$ and $a_{\n} + cb_{\n}$
have the same phase whenever they are both nonzero.
If one of $a_{\n}$ and $b_{\n}$ is zero,
then by taking $c \to 0$ or $c \to +\infty$ we can conclude
(using the fact that $P$ and $Q$ have the same-phase property)
that $a_{\m}$, $b_{\m}$ and the nonzero member of $\{a_{\n}, b_{\n}\}$
all have the same phase,
which implies
(again using the fact that $P$ and $Q$ have the same-phase property)
that all the nonzero coefficients of $P$ and $Q$ have the same phase.
If, on the other hand, neither $a_{\n}$ nor $b_{\n}$ is zero,
let us define $r_{\n} = a_{\n}/b_{\n}$.
Then $(r_{\m}+c)b_{\m}$ and $(r_{\n}+c)b_{\n}$ have the same phase
for all positive real numbers $c \notin \{-r_{\m},-r_{\n}\}$.
Since $b_{\m}$ and $b_{\n}$ have the same phase,
so do $r_{\m}+c$ and $r_{\n}+c$.
Since this is valid for all positive real $c \notin \{-r_{\m},-r_{\n}\}$,
it follows that either $r_{\m}$ and $r_{\n}$ are both real and positive,
or else that $r_{\m} = r_{\n}$.
If the former case holds for at least one multi-index $\n$
for which $a_{\n}$ and $b_{\n}$ are both nonzero,
then all the nonzero coefficients of $P$ and $Q$ have the same phase.
If the latter case holds for all $\n$ for which
at least one of $a_{\n}$ and $b_{\n}$ is nonzero,
then $P = \alpha Q$ with $\alpha = r_{\m}$.
\qed

\proofof{Theorem~\ref{thm.same-phase_parity}}
We proceed by induction on $n$.
For $n=1$, since $P$ has definite parity,
its zeros are symmetric with respect to the origin.
Since $P$ also  has the half-plane property,
these zeros all lie symmetrically on the imaginary axis,
so that $P(x) = a x^m \prod_{i=1}^k (x^2 + c_i)$
for some complex number $a$, nonnegative integers $m$ and $k$,
and positive real numbers $c_1,c_2,\ldots,c_k$.
This implies the same-phase property for $P$.

For the induction step, suppose that $P(x_1,\ldots,x_{n+1})$
has definite parity and the half-plane property.
Let $P_j(x)$  be the coefficient of $x_{n+1}^j$  in $P(x)$.
Each $P_j$ has definite parity (this parity alternates with $j$)
and has the half-plane property (by Corollary~\ref{cor.Pcoeffs});
so by the inductive hypothesis each $P_j$ has the same-phase property.
We now want to show that the nonzero coefficients of the various $P_j$
all have the same phase.

Define $j_{\rm min} = \min\{j \colon\, P_j \not\equiv 0\}$ and
$j_{\rm max} = \max\{j \colon\, P_j \not\equiv 0\}$,
and suppose that $j_{\rm min} \le j < j_{\rm max}$.
If $P_j \equiv 0$ (which of course can happen only for $j>j_{\rm min}$),
then Theorem~\ref{thm.fettweis} states that
$P_{j+1}/P_{j-1}$ is a strictly positive constant,
so that the nonzero coefficients of $P_{j-1}$ and $P_{j+1}$
have the same phase.
If $P_j \not\equiv 0$ and $P_{j+1} \not\equiv 0$,
then Theorem~\ref{thm.fettweis} states that
$P_{j+1}/P_j$ is real-part-positive on $H^n$.
Thus, by applying Theorem~\ref{thm.local_hpp_single_e}
to the polynomial $P_j(x) + y P_{j+1}(x)$,
where $y$ is a new variable,
we deduce that this polynomial has the half-plane property
(in the variables $x_1,\ldots,x_n,y$).
It follows that, for any nonnegative real numbers $c_1,\ldots,c_n$,
the polynomial $P_j(x) + (c_1 x_1 + \ldots + c_n x_n) P_{j+1}(x)$
also has the half-plane property
(in the variables $x_1,\ldots,x_n$).
Since this latter polynomial has definite parity,
by the inductive hypothesis it has the same-phase property.
Now apply Lemma~\ref{lemma.same-phase_parity}
with $P=P_j$ and $Q(x) = x_1 P_{j+1}(x)$.
In case (a), we conclude that
all the nonzero coefficients of $P_j$ and $P_{j+1}$ have the same phase.
In case (b), we conclude that $P_j/P_{j+1} = \alpha x_1$
for some nonzero complex number $\alpha$;
but since $P_j/P_{j+1}$ is real-part-positive on $H^n$,
$\alpha$ must be a positive real number,
so we again conclude that
all the nonzero coefficients of $P_j$ and $P_{j+1}$ have the same phase.

The foregoing observations show
that the nonzero coefficients of all the $P_j$ have the same phase,
and hence that $P$ has the same-phase property.
This completes the induction step, and hence the proof.
\qed

\bigskip

Finally, let us take this opportunity to discuss the image of $H^n$
under a homogeneous polynomial $P$.

\begin{proposition}
   \label{prop.imageofHn}
Let $P \not\equiv 0$ be a homogeneous polynomial of degree $r \ge 1$
in $n$ complex variables.
If $P$ does not have the half-plane property, then $P(H^n) = \C$.
If $P$ does have the half-plane property, then there exists $\theta \in \R$
such that
\be
        P(H^n)  \;=\;  e^{i\theta} H^r  \;=\;
        \cases{ e^{i\theta} H                         & if $r=1$ \cr
                \noalign{\vskip 2mm}
                \C \setminus e^{i\theta} (-\infty,0]  & if $r=2$ \cr
                \noalign{\vskip 2mm}
                \C \setminus 0                        & if $r \ge 3$ \cr
              }
 \label{eq.prop.imageofHn}
\ee
\end{proposition}        

\noindent
[We apologize for the inconsistent notations
 $H^n = H \times \ldots \times H$ and
 $H^r = \{ z^r \colon\, z \in H \}$.
 We trust that this will not cause any confusion.]

\proof
Clearly $P$ is nonconstant, so by the open mapping theorem, $P(H^n)$ is open.
By homogeneity, $P(H^n)$ is a cone (i.e.\ invariant under multiplication
by any $\lambda > 0$).
If $P$ does not have the half-plane property, then $P(H^n) \ni 0$;
and an open cone containing 0 must be all of $\C$.
Now suppose that $P$ does have the half-plane property,
so that $P(H^n) \not\ni 0$.
Then $P(1,\ldots,1) \neq 0$, so let $\theta = \arg P(1,\ldots,1)$.
Since $P(\zeta,\ldots,\zeta) = \zeta^r P(1,\ldots,1)$,
it follows that $P(\zeta,\ldots,\zeta)$ covers the set $e^{i\theta} H^r$
as $\zeta$ ranges over $H$.
Therefore, $P(H^n) \supseteq e^{i\theta} H^r$.

Let us now prove the reverse containment.
Suppose first that $r=1$ and $P(x) = \sum_{i=1}^n a_i x_i$.
Then, as remarked above, $P$ has the half-plane property
if and only if all the nonzero $a_i$ have the same phase,
i.e.\ there exists $\theta \in \R$ such that
$\real (e^{-i\theta} a_i) \ge 0$ for all $i$.
But in this case $P(H^n) = e^{i\theta} H$.

Suppose next that $r=2$ and $P(x) = \sum_{i,j=1}^n a_{i,j} x_i x_j$.
By Theorem~\ref{thm.same-phase},
there exists $\theta \in \R$ such that
$\real (e^{-i\theta} a_{ij}) \ge 0$ for all $i,j$.
But then Theorem~\ref{thm.rank2}(b)$\implies$(a)
tells us that $P(H^n) = \C \setminus e^{i\theta} (-\infty,0]$.

If $r \ge 3$, the containment $P(H^n) \subseteq \C \setminus 0$ is trivial.
\qed

\begin{corollary}
  \label{cor.disjoint}
Fix $N \ge 2$.
Let $\{E_i\}_{i=1}^N$ be {\em disjoint}\/ finite sets,
and for each $i$ let $P_i$ be a polynomial of degree $r_i$
on the ground set $E_i$.
[Note that $P_i$ need not be homogeneous.]
Suppose that $P = \sum_{i=1}^N P_i$, considered as a polynomial
on the ground set $E = \bigcup_{i=1}^N E_i$, has the half-plane property.
Then there is at most one index $i$ with $r_i \ge 2$ and $P_i \not\equiv 0$.
\end{corollary}

\proof
Suppose that there are two such indices $i,j$.
By setting $x_e = 0$ for all $e \in E \setminus (E_i \cup E_j)$,
we can assume without loss of generality that $N=2$ and $E = E_1 \cup E_2$.
Let us write $x^{(1)} = \{x_e\}_{e \in E_1}$ and
$x^{(2)} = \{x_e\}_{e \in E_2}$, and let us define
\be
   Q_\lambda(x)  \;=\; \lambda^{-1} P(\lambda^{1/r_1} x^{(1)},
                                      \lambda^{1/r_2} x^{(2)})
   \;.
\ee
Then $Q_\lambda$ has the half-plane property for all $\lambda > 0$;
and by Hurwitz's theorem, so does
\be
   \lim_{\lambda\to +\infty} Q_\lambda(x)
   \;=\;
   P_1^\sharp(x^{(1)}) \,+\, P_2^\sharp(x^{(2)})
\ee
(recall from Section~\ref{sec_leadingpart} that ${}^\sharp$ denotes
 leading part).
So we can assume without loss of generalize that $P_1$ and $P_2$
are homogeneous.
Applying Proposition~\ref{prop.imageofHn},
we conclude that there exist $\theta_1, \theta_2 \in \R$
such that $P_i(H^{E_i}) \supseteq \C \setminus e^{i\theta_i} (-\infty,0]$
for $i=1,2$.
But then we can choose $x^{(1)} \in H^{E_1}$ and $x^{(2)} \in H^{E_2}$
so that $P_1(x^{(1)}) + P_2(x^{(2)}) = 0$,
showing that $P$ does not have the half-plane property.
\qed

{\bf Remark.}
We can also allow the sets $E_i$ to be non-disjoint,
provided that $\deg P_1$, $\deg P_2 \ge |E_1 \cap E_2| + 2$:
it suffices to fix $x_e \in H$ for $e \in E_1 \cap E_2$
and then apply Corollary~\ref{cor.disjoint} in the remaining variables.

\bigskip

Let us conclude by mentioning an alternate proof of
Theorem~\ref{thm.same-phase},
which is more lengthy than the one given here but also more elementary
(because it does not require Proposition~\ref{prop.generalrank}).
First one proves Theorem~\ref{thm.same-phase} in the multiaffine case
by induction on $r$, deducing the same-phase property for $P$
from that of its contractions $P^{/i}$;
the key tool here (ensuring the necessary ``connectedness'')
is Corollary~\ref{cor.disjoint}.
Then one uses a Grace--Walsh--Szeg\"o construction
to reduce the case of a non-multiaffine polynomial
to that of a multiaffine polynomial in a larger number of variables.

\section{Necessary conditions II: Matroidal support}  \label{sec_matroidal}

\subsection{The exchangeability theorem}

In this section we shall prove that
if $P \not\equiv 0$ is a homogeneous multiaffine polynomial
with the half-plane property,
then $\supp(P)$ is necessarily the collection of bases of a matroid.

\begin{theorem}
   \label{thm.matroidal}
Let $P(x) = \sum_{S \subseteq E, |S|=r} a_S x^S$ be a multiaffine
homogeneous degree-$r$ polynomial (on the ground set $E$)
that has the half-plane property and is not identically zero.
Then there exists a rank-$r$ matroid $M$ on the ground set $E$
such that $\supp(P) = \scrb(M)$.
\end{theorem}

\noindent
When $\supp(P) = \scrb(M)$, we call the (obviously unique) matroid $M$
the {\em support matroid}\/ of $P$.

Indeed, with only a little more work we can prove a generalization
to homogeneous polynomials that are not necessarily multiaffine.
First, a bit of notation:
For multi-indices $\m$ and $\p$, let $\m\wedge\p$ be
the multi-index defined by $(\m\wedge\p)(e)=\min[\m(e), \p(e)]$
for all $e\in E$.
For $e \in E$, let us define the multi-index $\d_e$ by
\be
   \d_e(f)  \;=\;  \cases{1   & if $f=e$       \cr
                          0   & if $f \neq e$  \cr
                         }
\ee
If $\mathcal{S}$ is a set of multi-indices and $\m,\p\in\mathcal{S}$,
we say that \emph{$\m$ is exchangeable towards $\p$ in $\mathcal{S}$}
(and write ``$\m\rightarrow\p$ in $\mathcal{S}$'')
to denote the following condition:
for every $e\in E$ such that $\m(e)>\p(e)$, there
exists an $f\in E$ such that $\m(f)<\p(f)$ and
$\m-\d_{e}+\d_{f}\in\mathcal{S}$.
This is an obvious generalization of the basis-exchange property for matroids,
to which it reduces if all the multi-indices are 0-1-valued.

\begin{THM}
   \label{thm.matroidal_non-multiaffine}
Let $P(x) = \sum_{|{\bf m}|=r} a_{\bf m} x^{\bf m}$ be a
homogeneous degree-$r$ polynomial (on the ground set $E$)
that has the half-plane property and is not identically zero.
Then, for every $\m,\m' \in \supp(P)$,
$\m$ is exchangeable towards $\m'$ in $\supp(P)$.
\end{THM}

\noindent
Clearly, by specializing Theorem~\ref{thm.matroidal_non-multiaffine}
to the multiaffine case, we obtain Theorem~\ref{thm.matroidal}.

\proofof{Theorem~\ref{thm.matroidal_non-multiaffine}}
Since $P$ is homogeneous, by Theorem~\ref{thm.same-phase}
we may assume without loss of generality
that $a_{\m} > 0$ for all $\m\in\supp(P)$.

We proceed by induction on $r$.  The cases $r=0$ and $r=1$ are
trivial.  The induction step fails when passing from $r=1$ to $r=2$,
for reasons that will be explained below.
Therefore, we treat $r=2$ as the base case of the induction,
and deduce it from Theorem~\ref{thm.rank2} as follows:
Let $A=(a_{ij})_{i,j\in E}$
be the unique symmetric real matrix such that $P=P_{A}$.
If the conclusion of the theorem fails,
then there exist $\m,\m'\in\supp(P)$ and $e\in E$ with $\m(e)>\m'(e)$
such that for any $f\in E$ with $\m(f)<\m'(f)$ we have
$\m-\d_{e}+\d_{f}\not\in\supp(P)$.
Since $r=2$ we may write $\m=\d_{e}+\d_{g}$ and
$\m'=\d_{b}+\d_{c}$ (some of $b,c,e,g$ may coincide).
Since neither $\d_{g}+\d_{b}$ nor $\d_{g}+\d_{c}$ is in $\supp(P)$,
we have $a_{bg}=a_{gb}=a_{cg}=a_{gc}=0$.
Let $x=\d_{e}+\rho\d_{g}$ for some $\rho>0$, and let $y=\d_{b}+\d_{c}$.
We see that
\begin{subeqnarray}
   x^{\rm T}Ay   & = &  (a_{eb}+a_{ec})/2   \\[2mm]
   y^{\rm T}Ay   & = &  a_{bc}+(a_{bb}+a_{cc})/2 \;>\; 0  \\[2mm]
   x^{\rm T}Ax   & = &   \rho a_{eg}+(a_{ee}+\rho^{2}a_{gg})/2 \;>\; 0
\end{subeqnarray}
where $a_{eg}, a_{bc} > 0$ since $\m,\m' \in \supp(P)$].
As $\rho\rightarrow\infty$,
inequality (e) of Theorem~\ref{thm.rank2}
eventually fails, so $P$ does not have the half-plane property.

For the induction step, assume that $r\geq 3$ and that the
theorem holds for homogeneous polynomials of degree $r-1$.
We shall employ an inner induction on $n=|E|$.
The base case $n=1$ is trivial.
So assume the result for homogeneous polynomials of degree $r$
and ground sets of size strictly less than $n$.
Let $P$ be a homogeneous polynomial of degree $r$
on an $n$-element ground set $E$,
and let $\m,\m' \in \supp(P)$.
We shall show that $\m\rightarrow\m'$ in $\supp(P)$,
i.e.\ that for every $e\in E$ such that $\m(e)>\m'(e)$, there
exists an $f\in E$ such that $\m(f)<\m'(f)$ and
$\m-\d_{e}+\d_{f}\in\supp(P)$.
We divide the argument into three cases:
\begin{enumerate}
   \item There exists $b \in E$ such that $\m(b)=\m'(b)=0$.
   \item There exists $b \in E$ such that $\m(b)>0$ and $\m'(b)>0$.
   \item For every $b\in E$, exactly one of $\m(b)>0$ or $\m'(b)>0$ holds.
\end{enumerate}

\textit{Case 1.}\
Suppose there exists an element $b$ such that $\m(b)=\m'(b)=0$.
Then consider the polynomial $P^{\drop b}$,
which by Proposition~\ref{prop3.1}
has the half-plane property.
The ground set $E\drop b$ of $P^{\drop b}$
is strictly smaller than the ground set $E$ of $P$,
so the hypothesis of the inner induction applies.
Both $\m$ and $\m'$ are in the support of $P^{\drop b}$,
so we have $\m\rightarrow\m'$ in $\supp(P^{\drop b})$.
Since $\supp(P^{\drop b})\subseteq \supp(P)$,
it follows that $\m\rightarrow\m'$ in $\supp(P)$.

\medskip

\textit{Case 2.}\
Suppose there exists an element $b$ such that $\m(b)>0$ and $\m'(b)>0$.
Then consider the polynomial $P^{/b}$,
which by Proposition~\ref{prop3.1}
has the half-plane property.
This polynomial is homogeneous of degree $r-1$,
so the hypothesis of the outer induction applies.
Both $\m-\d_{b}$ and $\m'-\d_{b}$ are in the support of $P^{/b}$,
so we have $(\m-\d_{b}) \rightarrow(\m'-\d_{b})$ in $\supp(P^{/b})$.
Now let $e\in E$ be such that $\m(e)>\m'(e)$.
Since $(\m-\d_{b})(e) >(\m'-\d_{b})(e)$,
there is an $f\in E$ with $\m(f)<\m'(f)$ and
$\m-\d_{b}-\d_{e}+\d_{f}\in\supp(P^{/b})$.
This implies that $\m-\d_{e}+\d_{f}\in\supp(P)$.

This argument shows that whenever $\p,\q\in\supp(P)$
are such that $\p\wedge\q>0$, we have $\p\rightarrow\q$ in $\supp(P)$.
This fact will be used repeatedly in the remainder of the proof.

\medskip

\textit{Case 3.}\
Suppose that for every $b\in E$, exactly one of $\m(b)>0$
or $\m'(b)>0$ holds.
We proceed in a series of steps:

(a)
We claim that there is a multi-index $\k\in\supp(P)$
such that $\m\wedge\k\neq 0$ and $\m'\wedge\k\neq 0$.
Suppose not, and specialize
$x_{b}=y$ if $\m(b)>0$ and $x_{b}=z$ if $\m'(b)>0$.
The result is a polynomial of the form $Ay^{r}+Bz^{r}$
for some numbers $A,B>0$, which has the half-plane property.
However, upon substituting the values $y=B^{1/r}e^{i\pi/2r}$ and
$z=A^{1/r}e^{-i\pi/2r}$ this polynomial vanishes, a contradiction
(since $r \ge 2$).
Hence there must exist a multi-index $\k\in\supp(P)$ as claimed.

(b)
Next, we claim that there is a multi-index $\k' \in\supp(P)$
such that $|\m\wedge\k'|\geq 2$ and $|\m'\wedge\k'|\geq 1$.
(This is what fails in the case $r=2$.)
Perhaps the $\k$ found in step (a) already satisfies this condition.
If not, then $|\m\wedge\k|=1$ and $|\m'\wedge\k|\geq 2$.
Let $b\in E$ be such that $(\m'\wedge\k)(b)>0$.
Since $\k\rightarrow\m$ in $\supp(P)$ [as shown in Case 2],
there is a $g\in E$ with $\k(g)<\m(g)$ such that
$\k'=\k-\d_{b}+\d_{g}\in\supp(P)$.
This $\k'$ meets the required conditions.

(c)
Next, we claim that there is a multi-index $\p\in\supp(P)$
such that $|\m\wedge\p|\geq 2$ and $|\m'\wedge\p|\geq 1$
and, moreover, for all $b\in E$ such that $\m(b)>0$ we have $\p(b)\leq\m(b)$.
To show this, we define
\begin{equation}
   \nu(\p,\m)  \;=\;  \sum_{b\in E:\ \m(b)>0}  \max[0,\p(b)-\m(b)]
\end{equation}
and choose $\k' \in\supp(P)$ as in step (b).
If $\nu(\k',\m)=0$ then we are done,
so assume that $\nu(\k',\m)>0$.  Let $b\in E$ be such that
$\m(b)>0$ and $\k'(b)>\m(b)$.  Since $\k'\rightarrow\m$ in $\supp(P)$,
there is a $g\in E$ with $\k'(g)<\m(g)$ such that
$\p'=\k'-\d_{b}+\d_{g} \in \supp(P)$.
We have $\nu(\p',\m)<\nu(\k',\m)$ and $|\m\wedge\p'|=|\m\wedge\k'| \geq 2$
and $|\m'\wedge\p'|=|\m'\wedge\k'| \geq 1$.
Repeating this argument inductively as required
produces a multi-index $\p\in\supp(P)$ satisfying the claim.

(d)
Let $e\in E$ be such that $\m(e)>\m'(e)$,
and consider two subcases:
\begin{quote}
\begin{itemize}
   \item[(i)] $\m(e)>\p(e)$.
   \item[(ii)] $\m(e)=\p(e)$.
\end{itemize}
\end{quote}

Subcase (i):
If $\m(e)>\p(e)$, then since $\m\rightarrow\p$ in $\supp(P)$,
there is an $f\in E$ with $\m(f)<\p(f)$ such that
$\m-\d_{e}+\d_{f} \in \supp(P)$.
Since $\m(f)<\p(f)$ we must have $\m(f)=0$, so that $\m(f)<\m'(f)$.

Subcase (ii):
Here $\m(e)=\p(e)>0=\m'(e)$.  Therefore, since
$\p\rightarrow\m'$ in $\supp(P)$, there is a $g\in E$ with $\p(g)<\m'(g)$
such that $\q=\p-\d_{e}+\d_{g}\in\supp(P)$.
Now we have $\m(e)>\q(e)$ and $\m\wedge\q>0$ and, moreover,
for all $b\in E$ such that $\m(b)>0$ we have $\q(b)\leq\m(b)$.
We may therefore repeat the argument of subcase (i)
with $\q$ in place of $\p$.
\qed

{\bf Remark.}  In step (a) of Case 3 we needed to invoke the
same-phase property (Theorem~\ref{thm.same-phase})
in order to ensure that $A,B \neq 0$.
In the multiaffine case this is unnecessary,
because only $\m$ and $\m'$ contribute to $A$ and $B$,
so $A,B \neq 0$ is guaranteed;
and this fact is enough to imply that $Ay^r + Bz^r$
cannot have the half-plane property (by Corollary~\ref{cor.disjoint}).
In most other respects, however,
the proof of Theorem~\ref{thm.matroidal_non-multiaffine}
is not much more complicated than the proof for the multiaffine case.

\bigskip

Theorem~\ref{thm.matroidal_non-multiaffine}
shows that the support of a homogeneous polynomial
with the half-plane property satisfies a multi-analogue of
the matroid basis-exchange axiom.  Such structures were introduced
by Bouchet and Cunningham \cite{Bouchet_95} and are called ``jump systems''
(see also \cite{Lovasz_97,Geelen_nd}).
For $\p\in\ZZ^{E}$, let $\|\p\|=\sum_{e\in E}|\p(e)|$.
A nonempty subset $\scrj\subseteq\ZZ^{E}$ is called a
\emph{jump system} (with ground set $E$) in case it satisfies
the following condition:
\begin{itemize}
   \item[(J)]  For any $\m,\m'\in\scrj$ and $\k\in\ZZ^{E}$ such that
       $\|\m-\k\|=1$ and $\|\k-\m'\|=\|\m-\m'\|-1$, either $\k\in\scrj$
       or there exists $\k'\in\scrj$ with $\|\k-\k'\|=1$ and
       $\|\k'-\m'\|=\|\m-\m'\|-2.$
\end{itemize}
A jump system $\scrj$ is said to have \emph{constant sum} when
$|\m|=|\m'|$ for all $\m,\m'\in\scrj$.  In this case, the condition
defining a jump system is equivalent to requiring that
for all $\m,\m'\in\scrj$, both $\m\rightarrow\m'$ in $\scrj$ and
$(-\m)\rightarrow(-\m')$ in $-\scrj \equiv \{-\k\colon\, \k\in\scrj\}$.
(The definition of exchangeability extends in the natural way from
$\NN^{E}$ to $\ZZ^{E}$.)

\begin{CORO}
 \label{cor.jumpsystem}
Let $P \not\equiv 0$ be a homogeneous polynomial with the half-plane
property.  Then $\supp(P)$ is a jump system with constant sum.
\end{CORO}

\proof
Since $P$ is homogeneous, $\supp(P)$ has constant sum.
For any $\m,\m'\in\supp(P)$, we have $\m\rightarrow\m'$ in $\supp(P)$
by Theorem \ref{thm.matroidal_non-multiaffine}.
For each $e\in E$, let $\p(e)=\deg_{e} P$,
and consider $Q(x)=x^{\p}P(1/x)$.
Since the open right half-plane
is invariant under the transformation $z\mapsto 1/z$,
$Q$ is also a homogeneous polynomial with the half-plane property.
Thus, by Theorem \ref{thm.matroidal_non-multiaffine} again,
$(\p-\m)\rightarrow(\p-\m')$ in $\supp(Q)$.
But this means that $(-\m)\rightarrow(-\m')$ in $-\supp(P)$.
\qed

A nonempty intersection of a jump system in $\ZZ^{E}$ with the set
$\{0,1\}^{E}$ is known as a {\em delta-matroid}\/, and is usually
regarded as a set system (by identifying a set with its
characteristic vector).  When a delta-matroid is uniform
(i.e.\ all sets have the same size), then it is the collection of bases
of a matroid (as is easily seen),
so that Theorem~\ref{thm.matroidal} is a special case of
Corollary~\ref{cor.jumpsystem}.

\bigskip

It is natural to wonder what can be said about the support of
a polynomial $P$ with the half-plane property,
if one omits (or weakens) the assumption that $P$ is homogeneous.
The fact that the same-phase property need not hold in this case
causes great complications;  one might start by considering
the subclass of polynomials $P$ that have the half-plane property
{\em and}\/ the same-phase property.

\begin{question}
 \label{quest.matroidal.1}
If $P$ is multiaffine and has the half-plane property,
is $\supp(P)$ a delta-matroid?  What if $P$ also has the same-phase property?
\end{question}

\begin{question}
 \label{quest.matroidal.2}
Assume that $P$ has the half-plane property and
has definite parity
(i.e.\ $\m,\m'\in\supp(P)$ implies that $|\m|\equiv|\m'|$ mod 2).
Is $P$ then a jump system?
(Recall from Theorem~\ref{thm.same-phase_parity}
 that all such polynomials have the same-phase property.)
\end{question}


\subsection{The weak half-plane property}   \label{sec_matroidal.weakHPP}

Recall that a set system $\scrs$ (resp.\ a matroid $M$)
is said to have the half-plane property
if its generating polynomial $P_\scrs$
(resp.\ its basis generating polynomial $P_{\scrb(M)}$) does.
Let us now say that a set system $\scrs$ (resp.\ a matroid $M$)
has the {\em weak half-plane property}\/
if there exists a polynomial $P$ with the half-plane property
for which $\supp(P) = \scrs$ [resp.\ $\supp(P) = \scrb(M)$].
This is obviously a weaker condition, since the nonzero coefficients
of $P$ need not all be equal.
(By Theorem~\ref{thm.same-phase} they must, however,
 all have the same phase in the matroid case or,
 more generally, when the set system is $r$-uniform.)
Theorem~\ref{thm.matroidal} can then be rephrased as saying that
if an $r$-uniform set system has the weak half-plane property,
it is necessarily the collection of bases of a matroid.
A central open question is the converse:

\begin{question}
  \label{question.weak_HPP}
Does every matroid $M$ have the weak half-plane property?
And if not, which ones do?
\end{question}

\noindent
We shall soon show (Corollary~\ref{cor.determinant})
that all matroids representable over $\C$ have the
weak half-plane property.
But we are totally in the dark about matroids not representable over $\C$,
such as the Fano matroid $F_7$ and its coextensions $AG(3,2)$ and $S_8$
(which are representable only over fields of characteristic 2),
the matroids $T_8$, $R_9$, $S(5,6,12)$ and $PG(2,3)$
(which are representable only over fields of characteristic 3),
%
or the V\'amos and non-Pappus matroids
(which are not representable over any field).

\bigskip

{\bf Remark.}  The weak half-plane property for an $n$-element
matroid $M$ is algorithmically testable in time $c_1^{c_2^{2n+|\scrb(M)|}}$
using cylindrical algebraic decomposition,
and in time $c^{n|\scrb(M)|}$
using more recent algorithms (see Section~\ref{sec.algorithms}).
But this is thoroughly unfeasible in practice!

\section{The determinant condition and $\sqrt[6]{1}$-matroids}
   \label{sec_determinant}

\subsection{The determinant condition}

Let $A$ be an $r \times n$ matrix with complex entries, and define
\be
   Q_A(x)   \;=\;   \det(A X A^*)
 \label{def_QA}
\ee
where $X = \diag(x_1,\ldots,x_n)$
and ${}^*$ denotes Hermitian conjugate.
Clearly $Q_A$ is multiaffine and is homogeneous of degree $r$.
We have the following straightforward generalization of Theorem~\ref{thm1.1}:

\begin{theorem}
 \label{thm.QA}
Let $A$ be an arbitrary complex $r \times n$ matrix, and define
$Q_A$ by \reff{def_QA}.  Then:
\begin{itemize}
   \item[(a)]  $Q_A$ has the half-plane property.
   \item[(b)]  $Q_A(x) =  \sum\limits_{\begin{scarray}
                                          S \subseteq [n] \\
                                          |S| = r
                                       \end{scarray}}
                          |\det(A \restrict S)|^2 \, x^S   \,$,
      where $A \restrict S$ denotes the (square) submatrix of $A$
      using the columns indexed by the set $S$.
\end{itemize}
\end{theorem}

\proof
Note first that if $\rank A < r$,
then $Q_A \equiv 0$ and so trivially has the half-plane property;
so let us assume that $\rank A = r$.
Then $\ker A^* = 0$,
i.e.\ for every nonzero $\psi \in \C^r$,
we have $A^* \psi \neq 0$ in $\C^n$.
It follows that, for each $\psi \neq 0$, the quantity
\be
   \psi^* A X A^* \psi  \;=\;
   \sum_{i=1}^n |(A^* \psi)_i|^2 \, x_i
\ee
has strictly positive real part whenever $\real x_i > 0$ for all $i$;
so in particular $A X A^* \psi \neq 0$.
Therefore the matrix $A X A^*$ is nonsingular,
and so has a nonzero determinant.

The identity (b) is an immediate consequence of the
Cauchy--Binet theorem.\footnote{
   The Cauchy--Binet theorem states that if
   $A$ is an $m \times n$ matrix and $B$ is an $n \times m$ matrix,
   where $m \le n$, then
   $$ \det(AB)  \;=\;  \sum\limits_{\begin{scarray}
                                       S \subseteq [n] \\
                                       |S| = m
                                    \end{scarray}}
          \left( \det A\bigl[ [m] | S \bigr] \right)
          \left( \det B\bigl[ S | [m] \bigr] \right)
      \;.
   $$
   Here the sum runs over all $m$-element subsets
   $S \subseteq [n] \equiv \{1,2,\ldots,n\}$,
   $A\bigl[ [m] | S \bigr]$ denotes the submatrix of $A$
   consisting of the columns from $S$ (taken in order),
   and $B\bigl[ S | [m] \bigr]$ denotes the submatrix of $B$
   consisting of the rows from $S$ (taken in order).
   See e.g.\ \cite[pp.~128--129]{Marcus_65} for a proof.
}
\qed

{\bf Remark.}
This proof is a direct generalization of the proof of Theorem~\ref{thm1.1},
to which it reduces if we take $A$ to be the directed vertex-edge incidence
matrix for any orientation of $G$ with the $i_0$th row suppressed
(here $i_0$ is an arbitrary vertex of $G$).

\begin{corollary}
 \label{cor.determinant}
\begin{itemize}
   \item[(a)]  Every matroid representable over $\C$ has the
      weak half-plane property.
   \item[(b)]  Let $M$ be a rank-$r$ matroid on $n$ elements
that can be represented over $\C$ by an $r \times n$ matrix $A$
for which every $r \times r$ subdeterminant is either zero
or else of modulus 1.  Then $M$ has the half-plane property.
\end{itemize}
\end{corollary}

\proof
(a) This is an immediate consequence of Theorem~\ref{def_QA},
since $\det(A \restrict S) \neq 0$ if and only if $S \in \scrb(M)$.

(b) Clearly $|\det(A \restrict S)|^2 = 1$ if $S \in \scrb(M)$
and 0 if $S \notin \scrb(M)$.
So, by Theorem~\ref{def_QA}(b) we have $Q_A = P_{\scrb(M)}$,
and by Theorem~\ref{def_QA}(a) $Q_A$ has the half-plane property.
\qed

\subsection{$(F,G)$-representability and $\sqrt[6]{1}$-matroids} \label{sec_FG}

Which matroids have the property of Corollary~\ref{cor.determinant}(b)?
Certainly regular matroids do,
as they can be represented over $\R$ (hence also over $\C$)
by a totally unimodular $r \times n$ matrix $A$
(i.e.\ one for which every square subdeterminant is either 0, $+1$ or $-1$).
But, in fact, Corollary~\ref{cor.determinant}(b) applies to
a larger class of matroids.
Let us pose these questions  more generally in the following context.

Let $F$ be a field, let $F^*$ be the multiplicative group $F \setminus 0$,
and let $G$ be a subgroup of $F^*$.
If $A$ is a matrix over $F$,
let us call $A$ an {\em $(F,G)$-matrix}\/
if every nonzero subdeterminant of $A$ lies in $G$
(note, in particular, that every nonzero entry of $A$ must lie in $G$).
And let us call $A$ a {\em weak $(F,G)$-matrix}\/
if every nonzero $r \times r$ subdeterminant of $A$ lies in $G$,
where $A$ has rank $r$.
Finally, let us call a matroid {\em $(F,G)$-representable}\/
if it is representable over $F$ by an $(F,G)$-matrix,
and {\em weakly $(F,G)$-representable}\/
if it is representable over $F$ by a weak $(F,G)$-matrix.

The concept of an $(F,G)$-representable matroid
was introduced by Whittle \cite{Whittle_97}
and studied further, under the hypothesis that $-1 \in G$,
by Semple and Whittle \cite{Semple_96}
within the more general framework of partial fields.\footnote{
   These authors call it a {\em $(G,F)$-matroid}\/.
}
(We shall not need this more general notion here,
and we refer the reader to Semple and Whittle's paper for a discussion of it.)
As these authors note, many important classes of
matroids are special cases of $(F,G)$-representable matroids:
\begin{itemize}
   \item $F$-representable matroids:  $G=F^*$.
   \item Regular matroids (also known as unimodular matroids  \cite{White_87b}):
      Let $F$ be any field of characteristic zero
      (e.g.\ $F=\Q,\R$ or $\C$), and let $G = \{-1,1\}$.
   \item $k$-regular matroids \cite{Whittle_95,Whittle_97,Semple_98,Oxley_00}:
      Let $F$ be the field $\Q(\alpha_1,\ldots,\alpha_k)$
      obtained by extending the rationals
      by $k$ algebraically independent transcendental elements
      $\alpha_1,\ldots,\alpha_k$, and let $G$ be the set of all
      products of integer powers of differences of distinct members
      of $\{0,1,\alpha_1,\ldots,\alpha_k\}$.
      Thus, 0-regular matroids are just the regular matroids;
       1-regular matroids are also called near-regular.
       A matroid is  $\omega$-regular if it is $k$-regular for some $k$.
   \item Dyadic matroids \cite{Whittle_95,Whittle_97}:
      Let $F$ be any field of characteristic zero,
      and let $G = \{ \pm 2^k \colon\;  k \in \Z \}$.
   \item Sixth-root-of-unity matroids ($\sqrt[6]{1}$-matroid for short)
      \cite{Whittle_97}:
      Let $F=\C$ and $G$ be the multiplicative group $\Z_6$
      of complex sixth roots of unity.
   \item Complex unimodular matroids:  Let $F=\C$ and
      $G$ be the multiplicative group $U(1)$ of complex numbers of modulus 1.
\end{itemize}

Most ``naturally arising'' examples of $(F,G)$-representability ---
including those relevant to this paper --- have $-1 \in G$,
and (as we shall see) the theory takes a simpler form under this hypothesis.
Nevertheless, there is also some interesting theory that can be developed
for the case $-1 \notin G$:
for example, in Appendix~\ref{app_F1rep},
we characterize $(F,\{1\})$-representable matroids.
We shall therefore comment briefly on which results seem to need
$-1 \in G$ and which do not.

We begin with some easy lemmas:

\begin{lemma}
 \label{lemma_matrix_operations}
Let $A$ be an $(F,G)$-matrix, and let $B$ be an $F$-matrix.
\begin{itemize}
   \item[(a)]  If $B$ is obtained from $A$ by a sequence of row or column
      deletions, then $B$ is an $(F,G)$-matrix.
   \item[(b)]  If $B$ is obtained from $A$ by a sequence of row or column
      scalings by factors lying in $G$, then $B$ is an $(F,G)$-matrix.
   \item[(c)]  If $-1 \in G$, and $B$ is obtained from $A$ by a sequence
      of row or column interchanges, then $B$ is an $(F,G)$-matrix.
   \item[(d)]
      If $-1 \in G$,
      and $B$ is obtained from $A$ by a sequence of pivots,
      then $B$ is an $(F,G)$-matrix.
\end{itemize}
\end{lemma}

\proof
(a)--(c) are trivial.
The proof of (d) is essentially identical to that of
\cite[Proposition 3.3]{Semple_96}.
\qed

\begin{lemma}
 \label{lemma_FG_operations}
The class of $(F,G)$-representable matroids is closed under deletion
and direct sum.  Moreover, if $-1 \in G$, then the class is closed under
contraction, duality, 2-sum, and series and parallel connection.
\end{lemma}

\proof
The proof of the first sentence is straightforward.
When $-1 \in G$, closure under contraction follows from
Lemma~\ref{lemma_matrix_operations}(d)
by mimicking the argument of \cite[Proposition 3.2.6]{Oxley_92}.
As noted in \cite[Proposition 4.2]{Semple_96},
closure under duality, 2-sum, and series and parallel connection
follows by mimicking the standard arguments for fields that appear in
\cite[Theorem 2.2.8 and Proposition 7.1.21]{Oxley_92}.
\qed

It follows from the characterization of $(F,\{1\})$-representable matroids
in Appendix~\ref{app_F1rep} that, when $-1 \not \in G$, the class of
$(F,G)$-representable matroids need not be closed under duality, 2-sum,
series connection, or parallel connection.
One question we have been unable to answer is the following:

\begin{question}
   \label{quest.contraction}
When $-1 \notin G$,
is the class of $(F,G)$-representable matroids always
closed under contraction?
\end{question}

\noindent
It follows from Theorem~\ref{F1} that the answer is affirmative
when $G=\{1\}$, but we do not know whether this is so in general.

Our first nontrivial result is that weak $(F,G)$-representability
is no more general than $(F,G)$-representability, at least when $-1 \in G$.
We do not know whether this result holds without the assumption that
$-1 \in G$.

\begin{proposition}
 \label{prop.weakFG}
Let $F$ be a field, and let $G$ be a subgroup of the multiplicative group
$F^*$  that contains $-1$.
Let $A$ be a rank-$r$ $m \times n$ matrix over $F$
that is a weak $(F,G)$-matrix.
Then there is a rank-$r$ $r \times n$ matrix $B$ over $F$
that is an $(F,G)$-matrix and for which $M[A] = M[B]$.
\end{proposition}

\proof
Let $M = M[A]$. Since $A$ has rank $r$, we may delete
$m-r$ rows from $A$ leaving an $r \times n$ matrix $A_1$  whose rows are
linearly independent. Moreover, $A_1$ is a weak $(F,G)$-matrix
representing $M$. Now the following row and column operations leave a
determinant unchanged
except for possibly multiplying it by $-1$:
\begin{enumerate}
\item[(i)] Interchanging two rows or two columns.
\item[(ii)] Adding a multiple of one row to another.
\end{enumerate}
By operations (i) and (ii), $A_1$ can be transformed
into a matrix $A_2$ of the form $[D | Z]$ where $D$ is a diagonal matrix
all of whose diagonal
entries are nonzero.

The matrix $[D | Z]$ can be further transformed as follows.
Let the top-left $2\times 2$ submatrix of $D$ be
$\left[\!
\begin{array}{cc}
d_1 & 0 \\
0 & d_2
\end{array}
\!\right]$. In $[D | Z]$, perform the following operations:
\begin{enumerate}
\item[(i)] replace row 2 by row 2 plus $d_1^{-1}$ times  row 1;
\item[(ii)] replace row 1 by row 1 minus $d_1$ times  row 2;
\item[(iii)] replace row 2 by row 2 plus $d_1^{-1}$ times  row 1; and
\item[(iv)] interchange rows 1 and 2.
\end{enumerate}
We now have, in the top-left corner, the submatrix
$\left[\!\begin{array}{cc}
1 & 0 \\
0 & -d_1d_2
\end{array}
\!\right]$.
We can repeat this process until all the diagonal entries in $D$ except
possibly the last, $d_r$, are ones.
The resulting matrix $[D' | Z']$ is a weak $(F,G)$-matrix
representing $M$.  Since $\det D' = d_r$, it follows that $d_r \in G$.
Thus, we can multiply the last row of $[D' | Z']$ by $d_r^{-1}$ to
obtain a weak $(F,G)$-matrix $A_3$ representing $M$ of the form $[I_r|Y]$.

We shall show next that $A_3$ is an  $(F,G)$-matrix. Let $A'_3$ be a
$k \times k$ submatrix of $A_3$ for some $k$ with $1 \le k \le r-1$.
The identity submatrix $I_r$ of $A_3$ has an $r \times (r-k)$ submatrix
$B_1$ each of whose columns contains all zeros in the rows of $A_3$
that meet  $A'_3$. Let $A''_3$ be the submatrix of $A_3$ consisting of
those columns that meet some column of $A'_3$.  As $A_3$ is a weak
$(F,G)$ matrix, the determinant of the matrix $[B_1|A''_3]$ is in $G$.
But the value of this determinant is $\pm \det A_3''$.  Since $-1 \in G$,
we deduce that $\det A_3'' \in G$.  We conclude that $A_3$ is
indeed an $(F,G)$-matrix representing $M$.
\qed

\begin{proposition}
 \label{prop.FGstandardform}
Let $F$ be a field, and let $G$ be a subgroup of $F^*$ that contains $-1$.
Let $A$ be an $(F,G)$-matrix of rank $r$.
Then there is an $(F,G)$-matrix $B = [I_r|D]$ such that
$M[A] \simeq M[B]$ and $B$ is obtainable from $A$ by a sequence of
row or column interchanges, pivots, and deletions of zero rows.
\end{proposition}

\proof
This follows from Lemma~\ref{lemma_matrix_operations}
by well-known arguments:
see e.g.\ \cite[Propositions 3.5 and 4.1]{Semple_96}.
\qed

If $D$ is a matrix over $F$,
let us denote by $\scrg(D^\sharp)$ the simple bipartite graph
whose vertices correspond to the rows and columns of $D$
and whose edges correspond to the nonzero entries of $D$
\cite[pp.~190, 194]{Oxley_92}.
If $H$ is a subgroup of $F^*$,
let us call $D$ {\em $H$-normalized}\/ if there exists
a basis $B$ for the cycle matroid of $\scrg(D^\sharp)$
[i.e.\ a maximal spanning forest in $\scrg(D^\sharp)$]
such that each entry in $D$ corresponding to an edge in $B$ lies in $H$.

We aim next to show that,
given an $(F,G)$-matrix $[I_r|D]$,
we can perform a sequence of row and column scalings by elements of $G$
(thereby maintaining an $(F,G)$-matrix)
to yield a matrix $[I_r|D']$ for which $D'$ is $\{1\}$-normalized.
In fact, a much stronger result is true:

\begin{proposition}
\label{scaling}
Let $F$ be a field, and let $G$ be a subgroup of $F^*$.
Consider a matrix $[I_r|D]$ over $F$, in which all the nonzero elements
of $D$ lie in $G$.
Let $B=\{b_1,...,b_k\}$ be a basis for the cycle matroid
of $\scrg(D^\sharp)$,
and let $\{\theta_1,...,\theta_k\}$ be elements of $G$.
Then, by a sequence of row and column scalings with scale factors in $G$,
one can obtain from $[I_r|D]$ a matrix $[I_r|D']$
in which the entry of $D'$ corresponding to $b_i$ is $\theta_i$.
In particular, every subdeterminant of $[I_r|D']$
differs from the corresponding subdeterminant of $[I_r|D]$
by a factor lying in $G$.
\end{proposition}

\proof
The proof of \cite[Theorem 6.4.7]{Oxley_92} carries through with very minor
modifications to prove this proposition.
\qed

By definition, every $\sqrt[6]{1}$-matroid is a complex unimodular matroid.
In fact, the two classes coincide:

\begin{theorem}
 \label{thm.sixthroot}
The classes of complex unimodular matroids and $\sqrt[6]{1}$-matroids
are equal.
\end{theorem}

\proof
Evidently we need only
prove that every complex unimodular matroid is  a $\sqrt[6]{1}$-matroid.
By Propositions~\ref{prop.FGstandardform} and \ref{scaling},
we can represent any complex unimodular matroid $M$
by a matrix $[I_r|D]$ over $\C$ in which all nonzero subdeterminants
have modulus one and in which $D$ is $\{1\}$-normalized.
The desired result then follows from the following two lemmas:

\begin{lemma}
 \label{lemma.sixth.1}
Let $D$ be a $\Z_6$-normalized matrix over $\C$ in which all nonzero
subdeterminants have modulus one.
Then all the nonzero entries of $D$ are sixth roots of unity.
\end{lemma}

\begin{lemma}
 \label{lemma.sixth.2}
Let $A$ be a square matrix over $\C$, whose determinant has modulus 1,
such that all nonzero entries are sixth roots of unity.
Then $\det(A)$ is in fact a sixth root of unity.
\end{lemma}

Let us prove these two lemmas in reverse order:

\proofof{Lemma~\ref{lemma.sixth.2}}
Since $\det(A)$ is a linear combination, with coefficients $\pm 1$,
of products of entries in $A$, it follows that
$\det(A)$ belongs to the additive subgroup of $\C$
generated by the sixth roots of unity.
In the complex plane, the elements of this subgroup are
the vertices of a triangular lattice,
which intersects the unit circle precisely at the sixth roots of unity.
\qed

\proofof{Lemma~\ref{lemma.sixth.1}}
Let $B$ be a basis for the cycle matroid of $\scrg(D^\sharp)$
such that each entry in $D$ corresponding to an edge in $B$
is a sixth root of unity.
%
%
For each nonzero entry $d$ of $D$,
let $e_d$ be the corresponding edge of $\scrg(D^\sharp)$.
If $e_d \not \in B$, then adding $e_d$ to the subgraph of $\scrg(D^\sharp)$
induced by $B$ creates a unique cycle $C_d$,
whose length is some even integer $2k \ge 4$.
Assume that there exists a nonzero entry $d$
that is not a sixth root of unity,
and choose one for which $k$ is minimal.

Let $D_d$ be the $k \times k$ submatrix of $D$ induced by the
vertices of $C_d$.
By construction, all the entries of $D_d$ corresponding to edges in
$C_d \setminus e_d$ are sixth roots of unity.
Moreover, if $D_d$ contains any nonzero entry $d'$ besides those of $C_d$,
then the corresponding edge $e_{d'}$ of $\scrg(D^\sharp)$
is a diagonal of $C_d$
(note that this cannot occur in the minimal case $2k=4$).
It follows that $|C_{d'}| < |C_d|$,
so the choice of $d$ implies that $d'$ is a sixth root of unity.
Hence every nonzero entry of $D_d$ except possibly $d$
is a sixth root of unity.

Consider now the subgraph of $\scrg(D^\sharp)$
induced by the vertices of $C_d$;
among the cycles of this graph containing $e_d$,
choose a cycle $C'$ of shortest length $2j \le 2k$.
Let $D'_d$ be the $j \times j$ submatrix of $D_d$
induced by the vertices of $C'$.
Then each row and column of $D'_d$ has exactly two nonzero entries
corresponding to edges of $C'$,
and no other nonzero entries (by the minimality of $C'$).
Moreover, all the nonzero entries of $D'_d$, except possibly $d$,
are sixth roots of unity.

Because each row and column of $D'_d$ has exactly two nonzero entries,
the expansion of $\det(D'_d)$ in permutations
has exactly two nonzero terms, exactly one of which contains $d$.
We conclude that $\det(D'_d) = a-bd$ where $a$ and $b$ are
sixth roots of unity.
If $\det(D'_d)=0$, then $d$ is a sixth root of unity.
If $\det(D'_d)$ has modulus 1, then so does $c \equiv a/b - d$.
Hence $d$ is a complex number of modulus 1,
which differs from a sixth root of unity ($a/b)$
by another complex number  ($c$) of modulus 1.
It follows by simple geometry that $d$ is a sixth root of unity.
\qed

%
%
%

It is well known that the regular matroids are characterized by the
following list of equivalent properties \cite{White_87b}:

\begin{theorem}
The following are equivalent for a matroid $M$:
\begin{itemize}
   \item[(a)]  $M$ is $(F,\{\pm 1\})$-representable for some field $F$
       of characteristic 0.
   \item[(a)]  $M$ is weakly $(F,\{\pm 1\})$-representable for some field $F$
       of characteristic 0.
   \item[(c)]  $M$ is $(\Q,\{\pm 1\})$-representable.
   \item[(d)]  $M$ is weakly $(\Q,\{\pm 1\})$-representable.
   \item[(e)]  $M$ is representable over $GF(2)$ and over at least one
       field of characteristic not equal to $2$.
   \item[(f)]  $M$ is representable over all fields.
\end{itemize}
\end{theorem}
%
%
\noindent
More recently, Whittle \cite{Whittle_97} has given an analogous
characterization of the $\sqrt[6]{1}$-matroids.
Combining his (deep) results with our (comparatively elementary)
Proposition~\ref{prop.weakFG} and Theorem~\ref{thm.sixthroot},
we obtain:

\begin{theorem}[Whittle \protect\cite{Whittle_97}]
The following are equivalent for a matroid $M$:
\begin{itemize}
   \item[(a)]  $M$ is representable over $GF(3)$ and at least one field
      of characteristic $2$.
   \item[(b)]  $M$ is representable over $GF(3)$ and $GF(4)$.
   \item[(c)]  $M$ is representable over $GF(3)$ and all fields $GF(2^{2k})$
      for $k$ integer.
   \item[(d)]  $M$ is representable over all fields $F$ that contain a root
      of the polynomial $x^2 - x + 1$.
      [In particular, $M$ is representable over all fields $GF(q)$ for which
       $q$ is a power of 3, $q$ is a square,
       or $q \equiv 1 \!\!\!\! \pmod{3}$.]
    \item[(e)] $M$ is representable over all fields $GF(q)$
       for which $q \not\equiv 2 \!\!\!\! \pmod{3}$.
    \item[(f)]  $M$ is $(\C,\Z_6)$-representable.
   \item[(g)]  $M$ is $(\C,U(1))$-representable.
   \item[(h)]  $M$ is weakly $(\C,U(1))$-representable.
\end{itemize}
\end{theorem}

The excluded-minor characterization of regular matroids
is well known \cite[Theorems 13.1.1 and 13.1.2]{Oxley_92}:

\begin{theorem}[Tutte]
\quad

\vspace*{-2mm}
\begin{itemize}
\item[(a)] A binary matroid $M$ is regular if and only if
it has no minor isomorphic to $F_7$ or $F_7^*$.
\item[(b)] A matroid $M$ is regular if and only if
it has no minor isomorphic to $U_{2,4}$, $F_7$ or $F_7^*$.
\end{itemize}
\end{theorem}

\noindent
More recently, Geelen, Gerards and Kapoor \cite{Geelen_00a},
as a corollary of their important work in determining the
excluded minors for $GF(4)$-representable matroids,
have given an analogous excluded-minor characterization
of $\sqrt[6]{1}$-matroids:

\begin{theorem}[Geelen--Gerards--Kapoor \protect\cite{Geelen_00a}]
\label{GGK}
\quad

\vspace*{-2mm}
\begin{itemize}
\item[(a)] A ternary matroid $M$ is a $\sqrt[6]{1}$-matroid if and only if
$M$ has no minor isomorphic to $F_7^-$, $(F_7^-)^*$ or $P_8$.
\item[(b)] A matroid $M$ is a $\sqrt[6]{1}$-matroid if and only if $M$ has no
minor
isomorphic to $U_{2,5}$, $U_{3,5}$, $F_7$, $F_7^*$,
$F_7^-$, $(F_7^-)^*$ or $P_8$.
\end{itemize}
\end{theorem}

\noindent
It will be shown in Section~\ref{sec_counterexamples}
that none of $F_7$, $F_7^*$, $F_7^-$, $(F_7^-)^*$ and $P_8$
has the half-plane property.
The next two corollaries come from combining this fact
with the last two theorems and Corollary~\ref{cor.determinant}:

\begin{corollary}  A binary matroid has the half-plane property
if and only if it is regular.
\end{corollary}

\begin{corollary} A ternary matroid has the half-plane property
if and only if it is a $\sqrt[6]{1}$-matroid.
\end{corollary}

To conclude this section, we remark that a much less elementary proof of
Theorem~\ref{thm.sixthroot} than the one given above comes from using
Theorem~\ref{GGK}(b) together with the fact, which is not difficult to verify,
that none of $U_{2,5}$, $U_{3,5}$, $F_7$, $F_7^*$,
$F_7^-$, $(F_7^-)^*$ and $P_8$ is a complex unimodular matroid.

\section{Uniform matroids}  \label{sec_uniform}

\subsection{Half-plane property}

The basis generating polynomial of the uniform matroid $U_{r,n}$
is the elementary symmetric polynomial
\be
   E_{r,n}(x_1,\ldots,x_n)   \;=\;
   \sum_{1 \le i_1 < i_2 < \ldots < i_r \le n}  x_{i_1} x_{i_2} \cdots x_{i_r}
\ee
(we set $E_{0,n} \equiv 1$).
We have the following fundamental result:

\begin{theorem}
 \label{thm.uniform}
Let $0 \le r \le n$.  Then:
\begin{itemize}
   \item[(a)]  The elementary symmetric polynomial $E_{r,n}$
       has the half-plane property.
   \item[(b)]  For $r \ge 1$, the rational function
       $F_{r,n} \equiv E_{r,n}/E_{r-1,n}$
       is strictly real-part-positive on $H^n$.
\end{itemize}
\end{theorem}

Note first that (a) implies (b),
by Proposition~\ref{prop.derivs}(b) with
$\lambda_i = 1/(n-r+1)$ for all $i$;
and (b) implies (a), by Lemmas~\ref{lemma.RPP_implies_HPP_2}
and \ref{lemma.RPP_implies_HPP}
and the fact (which is a special case of Proposition~\ref{prop.irreducible})
that the elementary symmetric polynomials
$E_{r,n}$ with $r < n$ are irreducible over any field.
[Alternatively, the truth of (b) for {\em all}\/ $r \ge 1$
allows one to prove (a) by induction on $r$,
starting from $E_{0,n} \equiv 1$
and using Lemma~\ref{lemma.RPP_implies_HPP}.]
So we can try to prove whichever half seems more convenient.

We have already given two proofs of Theorem~\ref{thm.uniform}(a):
one based on nice principal truncation (i.e.\ differentiation)
starting from the easy case $r=n$ (Proposition~\ref{prop.uniform_nice}),
and another based on nice principal cotruncation
starting from the easy case $r=0$
(Proposition~\ref{prop.uniform_cotruncation_nice}).
We shall give a third proof in Section~\ref{sec_transversal},
based on showing that $U_{r,n}$ is a nice transversal matroid
(Corollary~\ref{cor.nice_transversal} and
 Example~\ref{sec_transversal}.\ref{exam.transversal.1} following it).
Let us here give two more proofs:
one based on the Grace--Walsh--Szeg\"o coincidence theorem,
and one based on ``series-parallel'' identities for the
rational functions $F_{r,n}$.

\fourthproof
We have $E_{r,n}(\xi,\ldots,\xi) = {n \choose r} \xi^r$,
which is manifestly novanishing for $\xi$ in the open right half-plane $H$.
Since $E_{r,n}$ is symmetric and multiaffine,
the Grace--Walsh--Szeg\"o coincidence theorem (Theorem~\ref{thm.grace})
implies that $E_{r,n}$ is nonvanishing in $H^n$.
\qed

\par\medskip\noindent{\sc Fifth proof} \cite{Fettweis_90}.
We shall prove (b) by induction on $r$.
The case $r=1$ is trivial for all $n$.
Now the following identities are easy to prove (for $1 \le r \le n$):
\begin{eqnarray}
   E_{r,n}(x)  & = &  E_{r,n-1}(x_{\neq i}) \,+\, x_i E_{r-1,n-1}(x_{\neq i})
     \qquad\hbox{for any index } i  \\[2mm]
   E_{r,n}(x)  & = &  {1 \over r} \sum_{i=1}^n x_i E_{r-1,n-1}(x_{\neq i})
\end{eqnarray}
We therefore have
\begin{subeqnarray}
   F_{r,n}(x)  & = &  {1 \over r} \sum_{i=1}^n
                      {x_i E_{r-1,n-1}(x_{\neq i})  \over  E_{r-1,n}(x)}
        \\[2mm]
   & = &  {1 \over r} \sum_{i=1}^n
          {x_i E_{r-1,n-1}(x_{\neq i})
           \over
           E_{r-1,n-1}(x_{\neq i}) \,+\, x_i E_{r-2,n-1}(x_{\neq i})
          }
        \\[2mm]
   & = &  {1 \over r} \sum_{i=1}^n
          {x_i F_{r-1,n-1}(x_{\neq i})
           \over
           F_{r-1,n-1}(x_{\neq i}) \,+\, x_i
          }
        \\[2mm]
   & = &  {1 \over r} \sum_{i=1}^n
          \left( {1 \over x_i} \,+\, {1 \over F_{r-1,n-1}(x_{\neq i})}
          \right) ^{\! -1}
          \;,
 \slabel{eq.Frn_synthesis.d}
 \label{eq.Frn_synthesis}
\end{subeqnarray}
so that the strict real-part-positivity of $F_{r,n}$ follows from that
of $F_{r-1,n-1}$.
\qed

{\bf Remarks.}
1.  The computation \reff{eq.Frn_synthesis} is really just \reff{eq.cotr}
specialized to the case of uniform matroids.
So this is just an explicit version of the ``nice cotruncation'' proof
of Proposition~\ref{prop.uniform_cotruncation_nice}.
Indeed, it was our analysis of the proof \reff{eq.Frn_synthesis}
that led us to abstract the idea of weighted principal cotruncation
(Section~\ref{sec.cotr}).

2. The identity \reff{eq.Frn_synthesis.d} shows how $F_{r,n}(x)$
can by synthesized as the admittance
of a 2-terminal series-parallel network
whose elementary branch admittances
are positive multiples of the $x_i$.
To begin with, $F_{1,n}(x) = x_1 + \cdots + x_n$
is the admittance of branches $x_1,\ldots,x_n$ placed in parallel.
Then $F_{r,n}(x)$ is the admittance of $n$ branches in parallel,
the $i$th of which consists of admittances $x_i/r$
and $F_{r-1,n-1}(x_{\neq i})/r$ in series.
[Of course, the admittance $F_{r-1,n-1}(x_{\neq i})/r$
 is obtained from the network giving $F_{r-1,n-1}(x_{\neq i})$
 by dividing each branch admittance by $r$.]

Now let $\scrg_{r,n}$ be the graph obtained by this construction:
it contains approximately $n^r$ edges,
and each edge $e$ is assigned an admittance $\alpha_e x_{i(e)}$
for some positive constant $\alpha_e$ and some index $i(e)$.
Then the spanning-tree polynomial $T_{\scrg_{r,n}}(\{x_e\})$,
with its arguments specialized as $x_e \to \alpha_e x_{i(e)}$,
contains the polynomial $E_{r,n}(x_1,\ldots,x_n)$ as a factor.
In this way, Theorem~\ref{thm.uniform}(a) can be obtained
as a corollary of Theorem~\ref{thm1.1}.

More generally, one might try to prove the half-plane property
for a polynomial $P$ by starting from the spanning-tree polynomial $T_G$
of a large graph $G$, specializing to fewer variables,
and then extracting a factor.
We do not know whether this method can be applied fruitfully
in other cases.

\subsection{Brown--Colbourn property}

Let us now prove that for $r < n$,
the uniform matroids $U_{r,n}$ not only have the half-plane property
but have the Brown--Colbourn property
(which is stronger, by virtue of Corollary~\ref{cor.shiftedHPP}).

\begin{theorem}
 \label{thm.uniform.BC}
Let $0 \le r \le n-1$.  Then the uniform matroid $U_{r,n}$ has
the ``Brown--Colbourn property'',
i.e.\ $\real x_i < -1/2$ for all $e$ implies
$P_{\scri(U_{r,n})}(x) \neq 0$,
where $\scri(U_{r,n}) = \{ S \subseteq \{1,\ldots,n\} \colon\, |S| \le r \}$
is the collection of independent sets of $U_{r,n}$.
\end{theorem}

\proof
Since
$P_{\scri(U_{r,n})}(x_1,\ldots,x_n) = \sum_{k=0}^r E_{k,n}(x_1,\ldots,x_n)$
is manifestly symmetric and multiaffine,
it follows from
the Grace--Walsh--Szeg\"o coincidence theorem (Theorem~\ref{thm.grace})
that it suffices to prove the Brown--Colbourn property
for the univariate polynomial
\be
   I_{r,n}(z)  \;=\;  P_{\scri(U_{r,n})}(z,\ldots,z)
               \;=\; \sum_{k=0}^r {n \choose k} \, z^k
   \;.
\ee
A slightly sharper result than this
is given in Proposition~\ref{prop.uniform.BC} below.
\qed

As preparation for the statement and proof
of Proposition~\ref{prop.uniform.BC},
let us note that the univariate reliability polynomial [cf. \reff{def.relpoly}]
for the uniform matroid $U_{r,n}$ is given by
\be
   {\rm Rel}_{r,n}(q)  \;=\; (1-q)^n I_{r,n}\Bigl( {q \over 1-q} \Bigr)
                       \;=\; \sum_{k=0}^r {n \choose k} \, q^k (1-q)^{n-k}
\ee
and hence can be written as ${\rm Rel}_{r,n}(q) = (1-q)^{n-r} H_{r,n}(q)$ with
\begin{subeqnarray}
  H_{r,n}(q)  & = &  \sum\limits_{k=0}^r {n \choose k} \, q^k (1-q)^{r-k}    \\
      & = &  \sum\limits_{k=0}^r \sum\limits_{j=0}^{r-k}
                 (-1)^j {n \choose k} {r-k \choose j} \, q^{j+k}             \\
      & = &  \sum\limits_{\ell=0}^r q^\ell \sum\limits_{j=0}^{\ell}
                 (-1)^j {n \choose \ell-j} {r-\ell+j \choose j}           \\
      & = &  \sum\limits_{\ell=0}^r {n-r-1+\ell \choose \ell} \, q^\ell
\end{subeqnarray}
where we define ${n-r-1+\ell \choose \ell} = \delta_{\ell 0}$
in case $r=n$.\footnote{
   The binomial-coefficient identity used in the last step
   can be deduced from \cite[eqns.~(1.7) and (1.9)]{Greene_90}.
   Alternatively, it is a specialization of
   \cite[Exercise 4.15(a)]{Wilf_94} followed by \cite[eqn.~(1.7)]{Greene_90}.
}

\begin{proposition}[Wagner \protect\cite{Wagner_00}]
   \label{prop.uniform.BC}
Let $0 \le r \le n-1$.  Then all the zeros of $H_{r,n}(q)$ lie in the annulus
\be
   {1 \over n-r}  \;\le\;  |q|  \;\le\;  {r \over n-1}  \;.
\ee
In particular, all the zeros of $H_{r,n}(q)$ lie in $|q| \le 1$,
so that all the zeros of $I_{r,n}(z)$ lie in $\re z \ge -1/2$.
\end{proposition}

\proof
The ratios of successive coefficients of $H_{r,n}$ are
$\lambda_\ell = {n-r-1+\ell \choose \ell} \big/ {n-r+\ell \choose \ell+1}
 = (\ell+1)/(n-r+\ell)$,
which is nondecreasing as $\ell$ runs from 0 to $r-1$.
Thus, by the Enestr\"om--Kakeya theorem
\cite[Theorem 30.3 and Exercise 2]{Marden_66}
(see also \cite{Anderson_79})
it follows that all the zeros of $H_{r,n}$
lie in the annulus $\lambda_0 \le |q| \le \lambda_{r-1}$.
\qed

\section{The permanent condition: Transversal and cotransversal matroids}
   \label{sec_transversal}

\subsection{Heilmann--Lieb theorem}  \label{sec_Heilmann-Lieb}

We begin by recalling the Heilmann--Lieb \cite{Heilmann_72} theorem
on the zeros of matching polynomials.
This theorem is most often quoted in its univariate version
(see e.g.\ \cite[Section 8.5]{Lovasz_86}),
but it is the multivariate result
\cite[Theorem 4.6 and Lemma 4.7]{Heilmann_72} that is truly fundamental.

Let $G=(V,E)$ be a loopless graph, and let us define the matching polynomial
with edge weights $\{\lambda_e\}_{e \in E}$ and
vertex weights $\{x_i\}_{i \in V}$:
\be
   M_G(x;\lambda)  \;=\; \sum_{{\rm matchings}\, M} \;
                         \prod_{e = ij \in M}  \lambda_e x_i x_j
   \;.
\ee
(Here $e=ij$ means that the endpoints of $e$ are $i,j$;
 if $G$ is not simple, this is a slight abuse of notation
 but unlikely to cause any confusion.
 Note also that we could, if we wanted, restrict attention to simple graphs
 by replacing each set of parallel edges $e_1,\ldots,e_k$
 with a single edge $e'$ of weight $\lambda_{e'} = \sum_{i=1}^k \lambda_{e_i}$.)
Define also the complementary matching polynomial
\be
   \widetilde{M}_G(x;\lambda)  \;=\;   x^V \, M_G(1/x;\lambda)
\ee
where $x^V = \prod_{i \in V} x_i$.
If $e = ij \in E$, we have the fundamental recursion relation
\be
   M_G(x;\lambda)  \;=\;  M_{G \setminus e}(x;\lambda)  \,+\,
                       \lambda_e x_i x_j M_{G-i-j}(x;\lambda)
\ee
or equivalently
\be
   \widetilde{M}_G(x;\lambda)  \;=\;
             \widetilde{M}_{G \setminus e}(x;\lambda)  \,+\,
             \lambda_e \widetilde{M}_{G-i-j}(x;\lambda)
   \;.
 \label{eq.transversal.star1}
\ee
Repeated application of \reff{eq.transversal.star1} to all edges incident on
a vertex $i$, followed by deletion of $i$, yields the key identity
\be
   \widetilde{M}_G(x;\lambda)  \;=\;
             x_i \widetilde{M}_{G-i}(x;\lambda)  \,+\,
             \sum\limits_{\begin{scarray}
                             e \sim i \\
                             e=ij
                          \end{scarray}}
               \lambda_e \widetilde{M}_{G-i-j}(x;\lambda)
 \label{eq.transversal.star2}
\ee
where $e \sim i$ denotes that $e$ is incident on $i$.
The Heilmann--Lieb theorem asserts that
if the edge weights are nonnegative,
then the polynomials $M_G$ and $\widetilde{M}_G$ have the half-plane property:

\begin{theorem}[Heilmann and Lieb \protect\cite{Heilmann_72}]
 \label{thm.heilmann-lieb}
Let $G=(V,E)$ be a loopless graph,
and let $\{\lambda_e\}_{e \in E}$ be nonnegative edge weights.
If $\real x_i > 0$ for all $i \in V$, then
\begin{itemize}
   \item[(a)]  $\widetilde{M}_G(x;\lambda) \neq 0$.
   \item[(b)]  For every $i \in V$, $\widetilde{M}_{G-i}(x;\lambda) \neq 0$.
   \item[(c)]  For every $i \in V$,
      $\real {\displaystyle \widetilde{M}_G(x;\lambda) \over
              \displaystyle \widetilde{M}_{G-i}(x;\lambda)} > 0$.
   \item[(d)]  $M_G(x;\lambda) \neq 0$.
\end{itemize}
\end{theorem}


\proof
We shall prove (a)--(c) by induction on $|V|$.
The theorem is trivial if $|V|=0$ or 1,
since $\widetilde{M}_G = 1$ for $|V|=0$
and $\widetilde{M}_G(x_1) = x_1$ for $|V|=1$.
So assume that (a)--(c) hold for all graphs with fewer than $n$ vertices,
and let $|V|=n$.
Then (b) holds for $G$ because (a) holds for $G-i$.
It then follows from \reff{eq.transversal.star2} that
\be
   {\widetilde{M}_G(x;\lambda) \over \widetilde{M}_{G-i}(x;\lambda)}
   \;=\;
   x_i \,+\,
   \sum\limits_{\begin{scarray}
                             e \sim i \\
                             e=ij
                          \end{scarray}}
               \lambda_e
   {\widetilde{M}_{G-i-j}(x;\lambda) \over \widetilde{M}_{G-i}(x;\lambda)}
   \;.
\ee
The right-hand side has positive real part because $\real x_i > 0$,
$\lambda_e \ge 0$ and
\linebreak
$\real[\widetilde{M}_{G-i-j}(x;\lambda) / \widetilde{M}_{G-i}(x;\lambda)] > 0$
by hypothesis (c) applied to $G-i$
(using the fact that the function $z \mapsto 1/z$
 maps the right half-plane to itself).
We conclude that (c) holds for $G$.  Hence so does (a).

Finally, (d) follows from (a) by Proposition~\ref{prop.duality}.
\qed

{\bf Remarks.}
1.  The same conclusion obviously holds [with a reversal of sign in (c)]
if $\real x_i < 0$ for all $i$.
But the proof does not apply to other rotated half-planes,
as we need invariance under the function $z \mapsto 1/z$.
Indeed, the result for other half-planes is false already for $G=K_2$,
which has $\widetilde{M}_G(x_1,x_2;\lambda) = \lambda_{12} + x_1 x_2$.

2.  This proof is basically the same as that of Heilmann and Lieb
\cite[Lemma 4.7]{Heilmann_72}, but is slightly simpler and stronger:
they prove (c) only for pairs $(G',G'-i)$ belonging to their set $\scre_G$.

3.  How could one have guessed the key inductive hypothesis (c)?
It suffices to note, from \reff{eq.transversal.star2},
that $\partial \widetilde{M}_G/\partial x_i = \widetilde{M}_{G-i}$;
and of course $\widetilde{M}_{G-i} \not\equiv 0$
because it has constant term 1.
Therefore, by either Proposition~\ref{prop.derivs}(b)
or Theorem~\ref{thm.local_hpp_single_e}(c),
if $\widetilde{M}_G$ has the half-plane property,
then necessarily hypothesis (c) must hold.

4. Here is a simple alternate proof of Theorem~\ref{thm.heilmann-lieb}(d),
based on the concept of multiaffine part (Proposition~\ref{prop.flat}):
The polynomial $P_G(x) = \prod_{e=ij \in E} (1 + \lambda_e x_i x_j)$
manifestly has the half-plane property whenever $\lambda_e \ge 0$ for all $e$.
But $P_G^\flat$ (the multiaffine part of $P_G$) is precisely $M_G$.

\subsection{Application to transversal matroids}

Now let us specialize the Heilmann--Lieb theorem to the case of
a bipartite graph $G=(V,E)$ with bipartition $V = A \cup B$,
setting $x_j = 1$ for $j \in B$
and considering $M_G(x;\lambda)$ as a polynomial in $\{x_i\}_{i \in A}$.
Then the restricted matching polynomial
\be
   \bar{M}_G(x;\lambda)  \;=\; \sum_{{\rm matchings}\, M} \;
                         \prod_{\begin{scarray}
                                  e = ij \in M \\
                                  i \in A  \\
                                  j \in B
                                \end{scarray}}
                                 \lambda_e x_i
\ee
also has the half-plane property,
provided that the edge weights are nonnegative.

Consider now the transversal matroid $M[G,A]$ with ground set $A$
defined by the bipartite graph $G$, in which a subset $S \subseteq A$
is declared independent if it can be matched into $B$.
Defining the weighted sum of such matchings,
\be
   c(S;\lambda)  \;=\;  \sum_{\begin{scarray}
                                 {\rm matchings}\, M \\
                                 V(M) \cap A = S
                              \end{scarray}} \;
                        \prod_{e \in M} \lambda_e
  \;,
\ee
we find immediately that
\be
   \bar{M}_G(x;\lambda) \;=\;  \sum_{S \in \scri(M[G,A])}  c(S;\lambda) \, x^S
   \;.
\ee
So the restricted matching polynomial $\bar{M}_G(x;\lambda)$,
which has the half-plane property,
is {\em almost}\/ the independent-set generating polynomial for $M[G,A]$:
the only trouble comes from the weights $c(S;\lambda)$,
which need not be equal.
Likewise, the restricted {\em maximum-}\/matching polynomial
\be
   \bar{M}_G^\sharp(x;\lambda) \;=\;
      \sum_{\begin{scarray}
               S \in \scri(M[G,A])  \\
               |S| = {\rm rank} \, M[G,A]
            \end{scarray}}
      c(S;\lambda) \, x^S
\ee
also has the half-plane property (by Proposition~\ref{prop.shiftedHPP}),
and is almost the basis generating polynomial for $M[G,A]$,
again modulo the problem of possibly unequal weights.

Note also that we may, if we wish, let $G$ be the {\em complete}\/
bipartite graph on the vertex set $A \cup B$,
since any undesired edges can always be given weight $\lambda_e = 0$.
In particular, if we take $A=[n]$, $B=[r]$ and $G=K_{n,r}$,
we have the following permanental analogue of Theorem~\ref{thm.QA}
(see \cite{Minc_78} for the definition and properties of permanents):

\begin{theorem}
 \label{thm.permanent}
Let $\Lambda$ be an arbitrary nonnegative $r \times n$ matrix
($r \le n$), and define $P_\Lambda$ by
\be
   P_\Lambda(x)  \;=\;  \per(\Lambda X)
\ee
where $X = \diag(x_1,\ldots,x_n)$.  Then:
\begin{itemize}
   \item[(a)]  $P_\Lambda$ has the half-plane property.
   \item[(b)]  $P_\Lambda(x) =  \sum\limits_{\begin{scarray}
                                          S \subseteq [n] \\
                                          |S| = r
                                       \end{scarray}}
                          \per(\Lambda \restrict S) \, x^S   \,$,
      where $\Lambda \restrict S$ denotes the (square) submatrix of $\Lambda$
      using the columns indexed by the set $S$.
\end{itemize}
\end{theorem}

\proof
Part (b) is an immediate consequence of the definition
of the permanent of a non-square matrix.
Part (a) follows immediately from the Heilmann--Lieb theorem and the fact that
\be
   P_\Lambda(x)  \;=\;
   \bar{M}_G^\sharp(x;\lambda)  \;=\;
   \lim_{\alpha\to\infty}  \alpha^{-r} \bar{M}_G(\alpha x;\lambda)
\ee
with $G=K_{n,r}$.
\qed

Let us now call the pair $(G,A)$ {\em nice}\/
if there exists a collection $\{\lambda_e\}_{e\in E}$
of nonnegative edge weights
so that $c(S;\lambda)$ has the same nonzero value
for all bases $S$ of $M[G,A]$.
And let us call the transversal matroid $M$ nice
if there exists a nice pair $(G,A)$ such that $M \simeq M[G,A]$.
The foregoing results can be rephrased as follows:

\begin{corollary}
 \label{cor.nice_transversal}
\quad

\vspace*{-2mm}
\begin{itemize}
   \item[(a)] Every transversal matroid has the
      weak half-plane property.
   \item[(b)] Every nice transversal matroid has the half-plane property.
\end{itemize}
\end{corollary}

\noindent
Unfortunately, Corollary~\ref{cor.nice_transversal}(a)
adds nothing to Corollary~\ref{cor.determinant}(a),
since every transversal matroid is representable
over all fields of sufficiently large cardinality,
and hence in particular over $\C$ \cite[Corollary 12.2.17]{Oxley_92}.
But Corollary~\ref{cor.nice_transversal}(b) is powerful,
as we shall soon see.

Our next task should be to characterize nice pairs $(G,A)$
and nice transversal matroids $M$.
Unfortunately, we are unable to do this in general,
so we shall content ourselves with giving
some examples of nice and non-nice transversal matroids.
As a preliminary, let us observe that the restriction of a nice pair $(G,A)$
to any subset $A' \subseteq A$ (and the same set $B$)
yields a nice pair $(G',A')$.
In particular, any restriction of a nice transversal matroid
is a nice transversal matroid.

Let us also note that, to prove that a matroid $M$ has the half-plane property,
it suffices by Proposition~\ref{prop.duality}
to show that either $M$ {\em or its dual}\/
is a nice transversal matroid.
If $M^*$ is transversal, we say that $M$ is {\em cotransversal}\/
(also called a {\em strict gammoid}\/);
if $M^*$ is nice transversal, we say that $M$ is {\em co-nice cotransversal}\/.

\subsection{Examples}  \label{sec_transversal_examples}

In this subsection we give
some examples of nice and non-nice transversal matroids.
Sometimes, instead of specifying the pair $(G,A)$,
we shall find it more convenient to specify the family
$\scra = \{A_j\}_{j \in B}$ of subsets of $A$
defined by $i \in A_j$ if and only if $ij \in E$.
The family $\scra$ is called a {\em presentation}\/ for $M[G,A]$.

\bexam
  \label{exam.transversal.1}
Every uniform matroid $U_{r,n}$ is nice.
Indeed, it suffices to take the obvious presentation,
namely the one induced by the complete bipartite graph $G=K_{n,r}$
with $|A|=n$ and $|B|=r$, and to set $\lambda_e = 1$ for all $e$.
\eexam

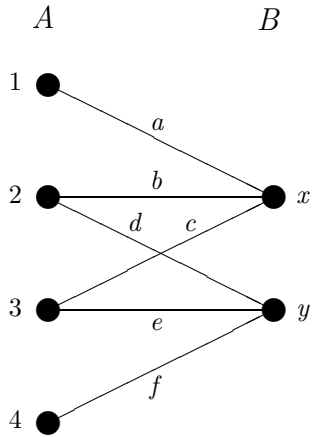
\begin{figure}[t]
\setlength{\unitlength}{1.5cm}
\begin{center}
\begin{picture}(2,4)(0,0)
   \put(0,0){\circle*{0.2}}  
   \put(0,1){\circle*{0.2}}  
   \put(0,2){\circle*{0.2}}  
   \put(0,3){\circle*{0.2}}  
   \put(2,1){\circle*{0.2}}  
   \put(2,2){\circle*{0.2}}  
   \put(0,0){\line(2,1){2}}  
   \put(0,1){\line(2,0){2}}  
   \put(0,1){\line(2,1){2}}  
   \put(0,2){\line(2,-1){2}} 
   \put(0,2){\line(2,0){2}}  
   \put(0,3){\line(2,-1){2}} 
   \put(-0.35,-0.05){\footnotesize 4}
   \put(-0.35,0.95){\footnotesize 3}
   \put(-0.35,1.95){\footnotesize 2}
   \put(-0.35,2.95){\footnotesize 1}
   \put(2.2,0.95){\footnotesize\it y}
   \put(2.2,1.95){\footnotesize\it x}
   \put(-0.15,3.5){\normalsize\it A}
   \put(1.85,3.48){\normalsize\it B}
   \put(0.9,2.6){\footnotesize\it a}
   \put(0.9,2.07){\footnotesize\it b}
   \put(1.2,1.7){\footnotesize\it c}
   \put(0.7,1.7){\footnotesize\it d}
   \put(0.9,0.82){\footnotesize\it e}
   \put(0.9,0.25){\footnotesize\it f}
\end{picture}
\end{center}
\caption{
   A bipartite graph $G=(V,E)$ with bipartition $V = A \cup B$.
}
 \label{fig_bipartite}
\end{figure}

\bexam
  \label{exam.transversal.2}
Not all pairs $(G,A)$ are nice.
To see this, consider the bipartite graph $G$ shown in
Figure~\ref{fig_bipartite},
with $A=\{1,2,3,4\}$ and $B=\{x,y\}$.
Then every 2-element subset of $A$ is matchable into $B$,
and the weights are
\begin{subeqnarray}
   c(12) & = & ad \\
   c(13) & = & ae \\
   c(14) & = & af \\
   c(23) & = & be + cd \\
   c(24) & = & bf \\
   c(34) & = & cf
\end{subeqnarray}
where for conciseness we have used the name of the edge $e$
 in place of $\lambda_e$.
Suppose that all these subset weights are equal and nonzero.
Then all the edge weights $a,\ldots,f$ must be nonzero,
and moreover we must have $d=e=f$ and $a=b=c$.
But then $c(23) = 2 c(24)$, a contradiction.
So this pair $(G,A)$ is not nice.
Note, nevertheless, that $M[G,A] \simeq U_{2,4}$,
which has an alternative (nonisomorphic) presentation
that {\em is}\/ nice (Example~\ref{sec_transversal}.\ref{exam.transversal.1}).
\eexam

\bexam
  \label{exam.transversal.3}
Not all transversal matroids are nice.
To see this, consider the rank-3 whirl $\scrw^3$,
which has ground set $\{1,2,3,4,5,6\}$
and 3-point lines 123, 345 and 561.
It is transversal with
maximal presentation
$\scra = \{456,\, 126,\, 234\}$
(i.e.\ the sets of the presentation are the complements of the
three 3-point lines).
By a general result (\cite{Mason_70} or \cite[Theorem 5.2.6]{Brualdi_87}),
this is the unique maximal presentation.
Now suppose that we could weight the edges in the bipartite graph
to make $\scrw^3$ nice.
The bases 413, 513 and 613 all occur exactly once as a transversal;
therefore, all the edges associated with the first set
in the presentation must receive the same weight, call it $\lambda(1)$.
Likewise, 513, 523 and 563 all occur once,
so all the edges associated with the second set
in the presentation must receive the same weight, call it $\lambda(2)$.
Finally, 512, 513 and 514 all occur once,
so all the edges associated with the third set
in the presentation must receive the same weight, call it $\lambda(3)$.
Thus the basis 135 gets weight $\lambda(1)\lambda(2)\lambda(3)$.
But the basis 246 gets weight $2\lambda(1)\lambda(2)\lambda(3)$.
Therefore $\scrw^3$ is not nice.

It follows that any transversal matroid having $\scrw^3$ as a restriction
is also non-nice:  this includes $\scrw^{3+}$ and $\scrw^3 + e$
(see Figures~\ref{fig.rank3.seven.4point} and \ref{fig_F7etc}
in Appendix~\ref{app_matroids})
and their free extensions $\scrw^{3+} + e$ and $\scrw^3 + e + f$
(see Figure~\ref{fig_F7m4pe_etc}).
\eexam

\bexam
  \label{exam.transversal.4}
For $n_1,n_2,n_3 \ge 3$, let $L_{n_1,n_2,n_3}$ be the rank-3 matroid
consisting of three nonintersecting lines containing
$n_1$, $n_2$ and $n_3$ points, respectively.
This matroid has a presentation consisting of
the complements of the three lines,
in which each basis occurs exactly twice as a transversal.
Therefore this matroid is nice (taking all weights $\lambda_e$ equal).

It follows that the matroids $P_6$, $S_7$, $F_7^{-6}$ and $P'''_7$
(see Figures~\ref{fig.rank3.le6}--\ref{fig_P7etc}
 in Appendix~\ref{app_matroids})
are nice, since they can be obtained by deleting elements
from $L_{3,3,3}$ or $L_{4,3,3}$.
\eexam

\bexam
  \label{exam.transversal.5}
More generally, let $L_{n_1,n_2,n_3;n'}$
(with $n_1,n_2,n_3 \ge 3$ and $n' \ge 0$)
be the rank-3 matroid consisting of three nonintersecting lines
plus $n'$ freely added points.
Let us show that the smallest example with freely added points,
namely $M = L_{3,3,3;1}$, is {\em not}\/ nice:

Let the 3-point lines be
$\hbox{\it abc}$, $\hbox{\it def}$ and $\hbox{\it ghi}$,
and call the freely added point 0.
Then $M$ has a presentation
consisting of the complements of the three 3-point lines, i.e.\
$\scra = \{0\hbox{\it abcdef}, 0\hbox{\it abcghi}, 0\hbox{\it defghi}\}$.
Indeed, this is the unique presentation.\footnote{
   {\sc Proof.}  The presentation $\scra$ is a maximal presentation,
   since adding an element to any set will result in one of
   $\hbox{\it abc}$\/, $\hbox{\it def}$\/ or $\hbox{\it ghi}$\/
   occurring as a transversal.
   By a general result (\cite{Mason_70} or \cite[Theorem 5.2.6]{Brualdi_87}),
   this is the unique maximal presentation.
   The complements of the three sets in the presentation are the three flats
   $\hbox{\it ghi}$\/, $\hbox{\it def}$\/ and $\hbox{\it abc}$\/.
   Any other presentation must come from the maximal one
   by deleting elements from the sets in the presentation.
   This will increase the complements of these sets.
   But the complements of the sets in any presentation
   must be flats in $M$, and there are no flats other than $E(M)$
   that contain any one of
   $\hbox{\it abc}$\/, $\hbox{\it def}$\/ or $\hbox{\it ghi}$\/.
   This proves the claim.
}
Now suppose that we could weight the edges in the bipartite graph
to make $M$ nice.
For each element $x$ of $0\hbox{\it defghi}$,
the basis $\hbox{\it abx}$ occurs exactly twice as a transversal:
as $\hbox{\it abx}$ and as $\hbox{\it bax}$.
This means that all the edges associated with the third set
in the presentation must receive the same weight, call it $\lambda(3)$.
By the same argument, all the edges associated with the
second set in the presentation must receive the same weight $\lambda(2)$,
and all the edges associated with the first set in the presentation
must receive the same weight $\lambda(1)$.
Thus the basis $\hbox{\it abd}$
gets weight $2\lambda(1)\lambda(2)\lambda(3)$.
But the basis $0\hbox{\it ae}$
gets weight $3\lambda(1)\lambda(2)\lambda(3)$.
Therefore $M$ is not nice.
\eexam

\bexam
  \label{exam.transversal.6}
Let $L_{n_1;n'}$ (with $n_1 \ge 3$ and $n' \ge 1$)
be the rank-3 matroid consisting of one $n_1$-point line
plus $n'$ freely added points.  
By deletion from $L_{n_1,3,3}$
(Example~\ref{sec_transversal}.\ref{exam.transversal.4}),
we can conclude that $L_{n_1;n'}$ is nice whenever $n' \le 4$.
Let us show that the smallest remaining case, $L_{3;5}$, is {\em not}\/ nice:

$L_{3;5}$ is transversal with unique maximal presentation 
$\scra = \{12345678, 12345678, 45678\}$.
Let us denote the weights associated with the first (resp.\ second, third)
set in the presentation as $\lambda_i$ (resp.\ $\lambda'_i, \lambda''_i$).
We shall use the combinations
$\alpha_{ij} = \lambda_i \lambda'_j + \lambda_j \lambda'_i$
and
$\beta_{ij} = \lambda_i \lambda'_j - \lambda_j \lambda'_i$.

For each pair of distinct elements $i,j \in \{1,2,3\}$
and each $x \in \{4,5,6,7,8\}$,
the basis {\it ijx}\/ occurs exactly twice as a transversal,
always with the same two ways of selecting $i$ and $j$ from the first two sets.
Fixing $i$ and $j$ and letting $x$ vary,
we conclude that all the edges associated with
the third set in the presentation must receive the same weight,
which without loss of generality we may take equal to 1. 
It then follows that $\alpha_{ij} = 1$
for all distinct $i,j \in \{1,2,3\}$.

For each $i \in \{1,2,3\}$ and
each triplet of distinct elements $x,y,z \in \{4,5,6,7,8\}$,
the bases {\it ixy}\/, {\it ixz}\/ and {\it iyz}\/
each occur exactly four times as transversals.
Thus $\alpha_{ix} + \alpha_{iy} = 1$,
$\alpha_{ix} + \alpha_{iz} = 1$
and $\alpha_{iy} + \alpha_{iz} = 1$.
By comparing the first equation with the difference of the last two,
we deduce that $\alpha_{ix} = 1/2$
for all  $i \in \{1,2,3\}$ and $x \in \{4,5,6,7,8\}$.


For each triplet of distinct elements $x,y,z \in \{4,5,6,7,8\}$,
the basis {\it xyz}\/ occurs exactly six times as a transversal,
so we have $\alpha_{xy} + \alpha_{xz} + \alpha_{yz} = 1$.
Therefore, for each quadruplet of distinct elements
$x,y,z,w \in \{4,5,6,7,8\}$, we have
$\alpha_{xy} + \alpha_{xz} + \alpha_{yz}  = 1$,
$\alpha_{xy} + \alpha_{xw} + \alpha_{yw}  = 1$,
$\alpha_{xz} + \alpha_{xw} + \alpha_{zw}  = 1$ and
$\alpha_{yz} + \alpha_{yw} + \alpha_{zw}  = 1$.
Subtracting the sum of the last two equations from the sum of the first two,
we conclude that $\alpha_{xy} = \alpha_{zw}$.
Since there are {\em five}\/ elements in $\{4,5,6,7,8\}$,
we can conclude that
$\alpha_{xy} = 1/3$ for all distinct $x,y \in \{4,5,6,7,8\}$.

Finally, for any distinct $i,j \in \{1,2,3\}$
and distinct $x,y \in \{4,5,6,7,8\}$, we have
$\alpha_{ij}\alpha_{xy} - \alpha_{iy}\alpha_{jx} = 1/3 - 1/4 = 1/12$.
Simplifying this, we find
$\beta_{ix} \beta_{jy} = -1/12$ for all such $i,j,x,y$.
But this implies that all the numbers $\beta_{ix}$
($1 \le i \le 3$, $4 \le x \le 8$)
are equal and take the value $\pm \sqrt{-1/12}$.
In particular, there are no {\em real}\/ solutions.
[The niceness equations do have {\em complex}\/ solutions of the form
$\lambda_1 = \lambda_2 = \lambda_3 = e^{\pm i \pi/6} /\sqrt{2}$,
$\lambda'_1 = \lambda'_2 = \lambda'_3 = e^{\mp i \pi/6} /\sqrt{2}$,
$\lambda_4 = \ldots = \lambda_8 = 1/\sqrt{6}$,
$\lambda'_4 = \ldots = \lambda'_8 = 1/\sqrt{6}$
and
$\lambda''_4 = \ldots = \lambda''_8 = 1$.]
\eexam

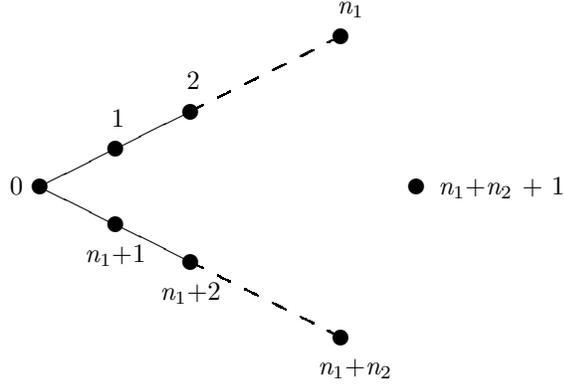
\begin{figure}[t]
\setlength{\unitlength}{1cm}
\begin{center}
\begin{picture}(5,5)(0,-2.5)
   \put(0,0){\circle*{0.2}}    
   \put(1,0.5){\circle*{0.2}}  
   \put(2,1){\circle*{0.2}}    
   \put(4,2){\circle*{0.2}}    
   \put(1,-0.5){\circle*{0.2}} 
   \put(2,-1){\circle*{0.2}}   
   \put(4,-2){\circle*{0.2}}   
   \put(5,0){\circle*{0.2}}    
   \drawline(0,0)(2,1)
   \dashline{0.2}(2,1)(4,2)
   \drawline(0,0)(2,-1)
   \dashline{0.2}(2,-1)(4,-2)
   \put(-0.4,-0.1){\footnotesize 0}
   \put(0.95,0.8){\footnotesize 1}
   \put(1.95,1.3){\footnotesize 2}
   \put(3.95,2.3){\footnotesize\it n$_1$}
   \put(0.6,-1.0){\footnotesize\it n$_1$\rm +1}
   \put(1.6,-1.5){\footnotesize\it n$_1$\rm +2}
   \put(3.7,-2.5){\footnotesize\it n$_1$\rm +\footnotesize\it n$_2$}
   \put(5.3,-0.1){\footnotesize\it n$_1$\rm +\footnotesize\it n$_2$ \rm+ 1}
\end{picture}
\end{center}
\caption{
   The matroid $M_{n_1,n_2}$.
}
 \label{fig_Mn}
\end{figure}

\bexam
  \label{exam.transversal.7}
For $n_1,n_2 \ge 2$, let $M_{n_1,n_2}$ be the rank-3 matroid on
the $(n_1+n_2+2)$-element ground set $\{0,1,\ldots,n_1+n_2+1\}$
such that there are lines $\{0,1,2,\ldots,n_1\}$
and $\{0,n_1+1,n_1+2,\ldots,n_1+n_2\}$
but no other 3-element circuits (see Figure~\ref{fig_Mn}).
This matroid is transversal with $A=\{0,1,\ldots,n_1+n_2+1\}$
and $B=\{x,y,z\}$, and has as a presentation
\begin{subeqnarray}
   A_x & = &  \{0,1,2,\ldots,n_1+n_2\}  \\
   A_y & = &  \{1,2,\ldots,n_1,n_1+n_2+1\}  \\
   A_z & = &  \{n_1+1,n_1+2,\ldots,n_1+n_2,n_1+n_2+1\}
\end{subeqnarray}
Now put weights
\begin{subeqnarray}
   \lambda_{0x} & = &  1  \\
   \lambda_{ix} & = &  1/2 \quad\hbox{for } 1 \le i \le n_1+n_2  \\
   \lambda_{iy} & = &  1 \qquad\hbox{for all } i \in A_y \\
   \lambda_{iz} & = &  1 \qquad\hbox{for all } i \in A_z
\end{subeqnarray}
It is easily checked that $c(S;\lambda) = 1$
for every basis $S$ of $M_{n_1,n_2}$.
So this is a nice presentation,
and $M_{n_1,n_2}$ is a nice transversal matroid.
In particular, $Q_6 \simeq M_{2,2}$, $P_6 \simeq M_{2,3} \setminus 0$,
$Q_7 \simeq M_{2,3}$, $S_7 \simeq M_{3,4} \setminus \{0,1\}$ and
$P'''_7 \simeq M_{3,3} \setminus 0$ are nice transversal matroids.
\eexam

\begin{proposition}
 \label{prop.nle6}
All matroids on a ground set of at most 6 elements have the
half-plane property.
\end{proposition}

\proof
The half-plane property is preserved by
direct sums (Proposition~\ref{prop.product})
and 2-sums (Corollary~\ref{cor.2sum}).
Therefore, we can restrict attention to 3-connected matroids.
Since all rank-1 and rank-2 matroids have the half-plane property
(Corollary~\ref{cor.rank2a}), we can restrict attention to
matroids of rank $\ge 3$.
Finally, since the half-plane property is invariant under duality
(Proposition~\ref{prop.duality}), we can restrict attention to
matroids of rank $\le \lfloor n/2 \rfloor$,
where $n$ is the number of elements.
So it suffices to consider 3-connected rank-3 matroids on 6 elements.
All these matroids are shown in Figure~\ref{fig.rank3.le6}
(see Appendix~\ref{app_matroids}).
The uniform matroid $U_{3,6}$ has the half-plane property,
as we have shown many times (Section~\ref{sec_uniform}).
The graphic matroid $M(K_4)$ and the whirl $\scrw^3$
are $\sqrt[6]{1}$-matroids, hence have the half-plane property
by Corollary~\ref{cor.determinant}(b).
Finally, we have just shown
(Examples~\ref{sec_transversal}.\ref{exam.transversal.4}
 and \ref{sec_transversal}.\ref{exam.transversal.7})
that $Q_6$ and $P_6$ are nice transversal matroids,
hence also have the half-plane property.
\qed

\bexam
  \label{exam.transversal.8}
The matroid $F_7^{-5}$
(see Figure~\ref{fig_F7etc} in Appendix~\ref{app_matroids})
is transversal with unique maximal presentation
$\scra = \{2347, 4567, 1234567\}$.
The bases 251, 351, 451 and 751 all occur exactly once as a transversal;
therefore, all the edges associated with the first set
in the presentation must receive the same weight, call it $\lambda(1)$.
By symmetry,
all the edges associated with the second set
in the presentation must receive the same weight, call it $\lambda(2)$.
But then 471 occurs twice as a tranversal,
so it gets twice the weight that 251 gets.
It follows that $F_7^{-5}$ is not nice.
\eexam

\bexam
  \label{exam.transversal.9}
The matroid $(F_7^{-3})^*$
(see Figure~\ref{fig_F7etc})
is transversal with unique maximal presentation
$\scra = \{123, 147, 156, 345\}$.
The bases 1465, 2465 and 3465 all occur exactly once as a transversal,
so all the edges associated with the first set
in the presentation must receive the same weight, call it $\lambda(1)$.
By symmetry, all the edges associated with the second set
in the presentation must receive the same weight, call it $\lambda(2)$,
and all the edges associated with the third set
in the presentation must receive the same weight, call it $\lambda(3)$.
Finally, the bases 1763, 1764 and 1765 all occur exactly once as a transversal,
so all the edges associated with the fourth set
in the presentation must receive the same weight, call it $\lambda(4)$.
But then 1345 occurs three times as a tranversal,
so it gets three times the weight that 1456 gets.
So $(F_7^{-3})^*$ is not nice.
\eexam

\bexam
  \label{exam.transversal.10}
The matroid $(F_7^{-4})^*$
(see Figure~\ref{fig_F7etc})
is transversal with unique maximal presentation
$\scra = \{123, 147, 156, 1234567\}$.
The bases 1456, 2456 and 3456 all occur exactly twice as transversals,
all with the same two ways of selecting 456 from the last three sets
in the presentation;
so all the edges associated with the first set
in the presentation must receive the same weight, call it $\lambda(1)$.
By symmetry, all the edges associated with the second set
in the presentation must receive the same weight, call it $\lambda(2)$,
and all the edges associated with the third set
in the presentation must receive the same weight, call it $\lambda(3)$.
Without loss of generality we can set
$\lambda(1) = \lambda(2) = \lambda(3) = 1$.
Let us denote the weights associated with the fourth set in the
presentation as $\lambda_i$ ($1 \le i \le 7$).
Bases then get weights as follows:
\begin{center}
\begin{tabular}{ll}
   1245: & $\lambda_1 + \lambda_2 + \lambda_4 + \lambda_5$ \\
   1345: & $\lambda_1 + \lambda_3 + \lambda_4 + \lambda_5$ \\
   1246: & $\lambda_1 + \lambda_2 + \lambda_4 + \lambda_6$ \\
   1346: & $\lambda_1 + \lambda_3 + \lambda_4 + \lambda_6$ \\
   1276: & $\lambda_1 + \lambda_2 + \lambda_7 + \lambda_6$ \\
   1376: & $\lambda_1 + \lambda_3 + \lambda_7 + \lambda_6$ \\
   1234: & $\lambda_2 + \lambda_3$ \\
   1475: & $\lambda_4 + \lambda_7$ \\
   1562: & $\lambda_5 + \lambda_6$
\end{tabular}
\end{center}
It follows from the first six bases that
$\lambda_2 = \lambda_3$, $\lambda_5 = \lambda_6$ and $\lambda_4 = \lambda_7$.
It then follows from the last three bases that
$\lambda_2 = \lambda_3 = \ldots = \lambda_7$.
But then we must have $\lambda_1 = -\lambda_2$.
Therefore, the equations have no {\em nonnegative}\/ solution,
so $(F_7^{-4})^*$ is not nice.

Now consider the matroid $(F_7^{-4} + e)^*$
(see Figure~\ref{fig_F7m4pe_etc})
which is transversal with unique maximal presentation
$\scra = \{123, 147, 156, 12345678, 12345678\}$.
The first part of the argument just given for $(F_7^{-4})^*$
carries over with slight modification to $(F_7^{-4} + e)^*$:
the bases 14568, 24568 and 34568 all occur exactly {\em four}\/ times
as transversals, all with the same four ways of selecting 4568
from the last four sets in the presentation,
so all the edges associated with the first set in the
presentation must receive the same weight $\lambda(1)$;
and by symmetry, all the edges associated with the second (resp.\ third) set
in the presentation must receive the same weight
$\lambda(2)$ [resp.\ $\lambda(3)$].
Without loss of generality we can set
$\lambda(1) = \lambda(2) = \lambda(3) = 1$.
Let us denote the weights associated with the fourth (resp.\ fifth)
set in the presentation as $\lambda_i$ (resp.\ $\lambda'_i$)
for $1 \le i \le 8$.
Bases then get weights as follows:
\begin{center}
\begin{tabular}{ll}
   14568: & $(\lambda_5 + \lambda_6) \lambda'_8 \,+\,
             (\lambda'_5 + \lambda'_6) \lambda_8$ \\
   14578: & $(\lambda_4 + \lambda_7) \lambda'_8 \,+\,
             (\lambda'_4 + \lambda'_7) \lambda_8$ \\
   12458: & $(\lambda_1 + \lambda_2 + \lambda_4 + \lambda_5) \lambda'_8 \,+\,
             (\lambda'_1 + \lambda'_2 + \lambda'_4 + \lambda'_5) \lambda_8$ \\
   12678: & $(\lambda_1 + \lambda_2 + \lambda_6 + \lambda_7) \lambda'_8 \,+\,
             (\lambda'_1 + \lambda'_2 + \lambda'_6 + \lambda'_7) \lambda_8$
\end{tabular}
\end{center}
All bases must get equal (and nonzero) weight.
Thus, by comparing the sum of the weights of the first two bases
with the sum of the weights of the last two bases,
we conclude that $\lambda_1 \lambda'_8 = \lambda_2 \lambda'_8 =
                  \lambda'_1 \lambda_8 = \lambda'_2 \lambda_8 = 0$
(since the weights are nonnegative).
But, by symmetry, the element 2 can be replaced by 3, 4, 5, 6 or 7.
Hence $\lambda_i \lambda'_8 = \lambda'_i \lambda_8 = 0$
for $1 \le i \le 7$, a contradiction.
Therefore $(F_7^{-4} + e)^*$ is not nice.
\eexam

\bexam
  \label{exam.transversal.11}
The matroid $(\scrw^3 + e)^*$
(see Figure~\ref{fig_F7etc})
is transversal with unique maximal presentation
$\scra = \{123, 345, 156, 1234567\}$.
The bases 1467, 2467 and 3467 all occur exactly once as a transversal;
so all the edges associated with the first set
in the presentation must receive the same weight, call it $\lambda(1)$.
By symmetry, all the edges associated with the second set
in the presentation must receive the same weight, call it $\lambda(2)$,
and all the edges associated with the third set
in the presentation must receive the same weight, call it $\lambda(3)$.
But then 1357 occurs twice as a transversal,
so it gets twice the weight that 1467 does.
So $(\scrw^3 + e)^*$ is not nice.

Now consider the matroid $(\scrw^3 + e + f)^*$
(see Figure~\ref{fig_F7m4pe_etc})
which is transversal with unique maximal presentation
$\scra = \{123, 345, 156, 12345678, 12345678\}$.
The argument just given for $(\scrw^3 + e)^*$
carries over verbatim to $(\scrw^3 + e + f)^*$
if we replace each basis $\hbox{\it abc}7$ by $\hbox{\it abc}78$
and we change ``all occur exactly once as a transversal'' to
``all occur exactly twice as a transversal, all with the same two ways
of selecting 78 from the last two sets in the presentation''.
It follows that $(\scrw^3 + e + f)^*$ is not nice.
\eexam   

\bexam
  \label{exam.transversal.12}
The matroid $(F_7^{-5})^*$
(see Figure~\ref{fig_F7etc})
is transversal with (non-maximal) presentation
$\scra = \{123, 156, 2467, 3457\}$,
in which each basis occurs exactly twice as a transversal.
Therefore this matroid is nice (taking all weights $\lambda_e$ equal).
\eexam


\bexam
  \label{exam.transversal.13}
The matroid $(\scrw^{3+})^*$
(see Figure~\ref{fig.rank3.seven.4point})
is transversal with unique maximal presentation
$\scra = \{156, 345, 1237, 1237\}$.
The bases 1427, 5427 and 6427 all occur exactly twice as transversals,
all with the same two ways of selecting 27 from the last two sets
in the presentation;
so all the edges associated with the first set
in the presentation must receive the same weight, call it $\lambda(1)$.  
By symmetry, all the edges associated with the second set
in the presentation must receive the same weight, call it $\lambda(2)$.
Let us denote the weights associated with the third (resp.\ fourth)
set in the presentation as $\lambda_i$ (resp.\ $\lambda'_i$) for $i=1,2,3,7$.
The bases 6417, 6427 and 6437 all occur
exactly twice as transversals,
from which we can conclude that
$\lambda_i \lambda'_7 + \lambda'_i \lambda_7 = \alpha > 0$ for $i=1,2,3$.
But then 1357 gets weight $2\lambda(1) \lambda(2) \alpha$,
while 1427 gets weight $\lambda(1) \lambda(2) \alpha$.
So $(\scrw^{3+})^*$ is not nice.

Now consider the matroid $(\scrw^{3+} + e)^*$
(see Figure~\ref{fig_F7m4pe_etc})
which is transversal with unique maximal presentation
$\scra = \{156, 345, 1237, 1237, 12345678\}$.
The argument just given for $(\scrw^{3+})^*$
carries over verbatim to $(\scrw^{3+} + e)^*$
if we replace each basis $\hbox{\it abcd}$ by $\hbox{\it abcd\/}8$.
So $(\scrw^{3+} + e)^*$  is not nice.
\eexam

\bexam
  \label{exam.transversal.14}
The matroid $(Q_7)^*$
(see Figure~\ref{fig.rank3.seven.4point})
is transversal with unique maximal presentation
$\scra = \{123, 1456, 1456, 1234567\}$.
The bases 1457, 2457 and 3457 all occur exactly twice as transversals,
all with the same two ways of selecting 457 from the last three sets
in the presentation;
so all the edges associated with the first set
in the presentation must receive the same weight,
which without loss of generality we take to be 1.
Let us denote the weights associated with the second (resp.\ third, fourth)
set in the presentation as $\lambda_i$ (resp.\ $\lambda'_i, \lambda''_i$).
For each $x \in \{4,5,6\}$,
the basis $21x7$ occurs twice as a transversal (as $21x7$ and $2x17$);
so $\lambda_1 \lambda'_x + \lambda_x \lambda'_1$
takes the same value $\alpha \neq 0$ for $x=4,5,6$.
For each pair of distinct $x,y \in \{4,5,6\}$,
the basis $1xy7$ occurs twice as a transversal (as $1xy7$ and $1yx7$);
so $\lambda_x \lambda'_y + \lambda_y \lambda'_x = \alpha$ as well.
Bases then get weights as follows:
\begin{center}
\begin{tabular}{ll}
   1452: & $\alpha (\lambda''_1 + \lambda''_2 + \lambda''_4 + \lambda''_5)$ \\ 
   1453: & $\alpha (\lambda''_1 + \lambda''_3 + \lambda''_4 + \lambda''_5)$ \\ 
   1463: & $\alpha (\lambda''_1 + \lambda''_3 + \lambda''_4 + \lambda''_6)$ \\ 
   1563: & $\alpha (\lambda''_1 + \lambda''_3 + \lambda''_5 + \lambda''_6)$ \\ 
   2143: & $\alpha (\lambda''_2 + \lambda''_3)$ \\
   2456: & $\alpha (\lambda''_4 + \lambda''_5 + \lambda''_6)$ \\
\end{tabular}
\end{center}
Comparing the first four bases, we conclude that
$\lambda''_2 = \lambda''_3 \equiv \beta$
and $\lambda''_4 = \lambda''_5 = \lambda''_6 \equiv \gamma$.
But then $\lambda''_1 + \beta + 2\gamma = 2\beta = 3\gamma$,
which implies that $\lambda''_1 = -\smhalf \gamma$.
Therefore, the equations have no {\em nonnegative}\/ solution,
so $(Q_7)^*$ is not nice.
\eexam

\bexam
  \label{exam.transversal.15}
The matroid $P''_7$
(see Figure~\ref{fig_P7etc})
is transversal with unique maximal presentation
$\scra = \{1267, 1345, 4567\}$.
The bases 134, 234, 634 and 734 all occur exactly once as a transversal;
therefore, all the edges associated with the first set
in the presentation must receive the same weight.
Likewise, 216, 236, 246 and 256 (resp.\ 134, 135, 136 and 137) all occur once,
so all the edges associated with the second (resp.\ third) set
in the presentation must receive the same weight.
But then 147 gets twice the weight of 134.
So $P''_7$ is not nice.
\eexam

\bexam
  \label{exam.transversal.16}
The matroid $(P'_7)^*$
(see Figure~\ref{fig_P7etc})
is transversal with unique maximal presentation
$\scra = \{123, 156, 267, 345\}$.
The bases 1574, 2574 and 3574 all occur exactly once as a transversal;
therefore, all the edges associated with the first set
in the presentation must receive the same weight.
Analogous arguments using 2174, 2574, 2674
(resp.\ 3524, 3564, 3574 or 2673, 2674, 2675)
show that all the edges associated with the second (resp.\ third or fourth)
set in the presentation must receive the same weight.
But then 1263 gets twice the weight of 1574.
So $(P'_7)^*$ is not nice.

Now consider the matroid $(P'_7 + e)^*$
(see Figure~\ref{fig_F7m4pe_etc})
which is transversal with unique maximal presentation
$\scra = \{123, 156, 267, 345, 12345678\}$.
The argument just given for $(P'_7)^*$
carries over verbatim to $(P'_7 + e)^*$
if we replace each basis $\hbox{\it abcd}$ by $\hbox{\it abcd\/}8$.
So $(P'_7 + e)^*$  is not nice.   
\eexam

\bexam
  \label{exam.transversal.17}
The matroid $(P''_7)^*$
(see Figure~\ref{fig_P7etc})
is transversal with (non-maximal) presentation
$\scra = \{123, 267, 345, 14567\}$,
in which each basis occurs exactly twice as a transversal.
Therefore $(P''_7)^*$ is nice (taking all weights $\lambda_e$ equal).
\eexam

\bexam
  \label{exam.transversal.18}
More generally,
for each $n \ge 1$, let $N_n$ be the rank-$(n+1)$ transversal matroid on
the ground set $\{1,\ldots,2n+1\}$ given by the presentation
\begin{subeqnarray}
   A_{\hat{0}} & = &  \{1,2,4,6,8,\ldots,2n-4,2n-2,2n,2n+1\}  \\
   A_{\hat{j}} & = &  \{2j-1,2j,2j+1\} \quad\hbox{for } 1 \le j \le n
\end{subeqnarray}
Geometrically, the dual of $N_n$
consists of 3-point lines 123, 345, 567, \ldots,
$(2n-1)(2n)(2n+1)$ joined together in general position in rank $n$.
Let $S$ be a matchable subset of $\{1,\ldots,2n+1\}$.
By using induction on $n$ and considering separately the cases when
$S$ does and does not contain two consecutive members of $\{1,\ldots,2n+1\}$,
we can show (after some non-trivial work)
that $S$ arises from exactly two matchings.
So, putting equal weights on all edges,
we conclude that this is a nice presentation,
and $N_n$ is a nice transversal matroid.
We have $N_1 \simeq U_{2,3}$, $N_2 \simeq U_{3,5}$ and $N_3 \simeq (P''_7)^*$.
\eexam

\bexam
  \label{exam.transversal.19}
The matroid (non-Pappus$\setminus 1)^*$ is transversal with maximal presentation 
$\scra = \{247, 269, 348, 359, 456\}$. 
For each $x \in \{4,5,6\}$, the basis $7283x$ occurs exactly once as a 
transversal, so all edges associated with the last set in the presentation 
receive the same nonzero weight.
For each $y \in \{3,5,9\}$, the basis $728y6$ occurs exactly once as a 
transversal, so all edges associated with the fourth set in the presentation 
receive the same nonzero weight.
For each $z \in \{2,6,9\}$, the basis $7z834$ occurs exactly once as a 
transversal, so all edges associated with the second set in the presentation 
receive the same nonzero weight.
Now consider the bases 29356 and 26354.
The first occurs exactly once as a transversal,
while the second occurs exactly twice (as 26354 and as 42356).
Since the transversals 29356 and 26354 get the same nonzero weight,
the transversal 42356 must get weight 0.
Therefore the edge joining the element 4 to the first set in the 
presentation has weight 0.  So we may delete 4 from the first set and still 
have a presentation.  This implies that $\{2,7\}$ contains a cocircuit of 
(non-Pappus$\setminus 1)^*$ and hence contains a circuit of
non-Pappus$\setminus 1$, which is false.
So (non-Pappus$\setminus 1)^*$ is not nice.
\eexam


\bigskip

We have written a {\sc Mathematica} program {\tt nicetransversal.m}
to test a presentation $\scra$ for niceness,
using the {\sc Mathematica} function {\tt Solve}
to solve the niceness equations $c(S;\lambda) = 1$ for $S \in \scrb(M[\scra])$.
This program is available as part of the electronic version of this paper
at arXiv.org.
However, because the niceness equations are a polynomial system
of degree $r$ in a large number of variables,
this program sometimes crashes for lack of memory
or fails to complete even in several days of CPU time.

\subsection{Rank-3 transversal matroids}  \label{sec_rank3_transversal}


Let us now make a systematic study of the niceness
of rank-3 transversal matroids.
We begin by recalling Brylawski's \cite{Brylawski_75}
algorithm for constructing transversal matroids.
A {\em principal transversal matroid}\/
(also called a {\em fundamental transversal matroid}\/) 
of rank $r$ is obtained by beginning with a distinguished basis
$B = \{1,2,\ldots,r\}$,
which we view geometrically as a simplex with vertices $1,2,\ldots,r$,
and adding elements as follows:
First one can add elements freely to the flat of rank 0 (i.e.\ add loops).
Then one can add elements freely to the flats of rank 1 (i.e.\ add elements
parallel to the basis elements).
Then one can add elements freely to the lines spanned by two basis elements:
such elements pick up the dependencies that are forced by virtue of their
lying on that line, but have no other dependencies.
One continues this process by adding elements freely
on the flats spanned by three basis elements,
then by four basis elements, \ldots\ and finally by $r$ basis elements
(the last are elements that are free in the matroid).

\begin{theorem}[Brylawski \protect\cite{Brylawski_75}]
 \label{thm.brylawski.1}
Every transversal matroid of rank $r$ is a restriction of a
principal transversal matroid of rank $r$,
obtained by deleting some elements of the original basis $B$.
\end{theorem}

\begin{theorem}[Brylawski \protect\cite{Brylawski_75}]
 \label{thm.brylawski.2}
A principal transversal matroid has a unique presentation,
which consists of the $r$ cocircuits
that are the complements of the flats spanned by the sets $B \setminus i$
($1 \le i \le r$).
\end{theorem}

Applying this construction in the case $r=3$, we conclude that
the most general {\em simple}\/ rank-3 principal transversal matroid
is $C_{n_1,n_2,n_3;n'}$ with $n_1,n_2,n_3 \ge 2$
and $n' \ge 0$, defined as the rank-3 matroid on $n_1 + n_2 + n_3 + n' - 3$
elements consisting of consisting of three distinguished points
(namely, the simplex vertices), three lines that join pairs of these vertices
and contain $n_1, n_2$ and $n_3$ points, respectively,
together with $n'$ points freely added in the plane spanned by the vertices.
All simple rank-3 transversal matroids can then be obtained
by deleting zero or more simplex vertices of $C_{n_1,n_2,n_3;n'}$.
By considering all the possibilities
for simplex vertices to be deleted or not,
and for the resulting lines to be nontrivial or trivial,
we obtain the following classes of matroids:
\begin{itemize}
   \item[1)]  
      $C_{n_1,n_2,n_3;n'}$ with $n_1,n_2,n_3 \ge 3$ and $n' \ge 0$.
   \item[2)]  
      $D_{n_1,n_2;n'}$ with $n_1,n_2 \ge 3$ and $n' \ge 0$,
      consisting of two intersecting lines containing $n_1$ and $n_2$ points,
      respectively, together with $n'$ freely added points.
   \item[3)]  
      $E_{n_1,n_2,n_3;n'}$ with $n_1,n_2,n_3 \ge 3$ and $n' \ge 0$,
      consisting of an $n_2$-point line that is met by nonintersecting lines
      with $n_1$ and $n_3$ points, respectively,
      together with $n'$ freely added points.
   \item[4)]  
      $F_{n_1,n_2,n_3;n'}$ with $n_1,n_2,n_3 \ge 3$ and $n' \ge 0$,
      consisting of an intersecting pair of lines containing $n_1$ and $n_2$
      points, respectively, and a third line containing $n_3$ points that
      does not meet the first two, together with $n'$ freely added points.
   \item[5)]  
      $L_{n_1,n_2,n_3;n'}$ with $n_1,n_2,n_3 \ge 3$ and $n' \ge 0$,
      consisting of three nonintersecting lines
      having $n_1, n_2$ and $n_3$ points, respectively,
      together with $n'$ freely added points.
   \item[6)]  
      $L_{n_1,n_2;n'}$ with $n_1,n_2 \ge 3$ and $n' \ge 0$,
      consisting of two nonintersecting lines having $n_1$ and $n_2$ points,
      respectively, together with $n'$ freely added points.
   \item[7)]  
      $L_{n_1;n'}$ with $n_1 \ge 3$ and $n' \ge 1$,
      consisting of one $n_1$-point line together with $n'$ freely added points.
   \item[8)]  The uniform matroid $U_{3,n'}$ with $n' \ge 3$,
      consisting of $n'$ points and no non-trivial lines.
\end{itemize}
We can now determine which of these matroids are nice:

{\em Class 1.}\/
The smallest case $C_{3,3,3;0} \simeq \scrw^3$ is not nice
(Example~\ref{sec_transversal}.\ref{exam.transversal.3} above).
And since, starting from any larger matroid $C_{n_1,n_2,n_3;n'}$
one can obtain $C_{3,3,3;0}$ by deleting elements,
it follows that all such matroids are non-nice.

{\em Class 2.}\/
The matroid $D_{n_1,n_2;n'}$ is nice when $n'=0$ or 1
(Example~\ref{sec_transversal}.\ref{exam.transversal.7}).
If $n' \ge 2$, the smallest case is $D_{3,3;2} \simeq F_7^{-5}$,
which is not nice
(Example~\ref{sec_transversal}.\ref{exam.transversal.8}).
Therefore, all cases with $n' \ge 2$ are non-nice.

{\em Class 3.}\/
The smallest case is $E_{3,3,3;0} \simeq P''_7$,
which is not nice
(Example~\ref{sec_transversal}.\ref{exam.transversal.15}).

{\em Class 4.}\/  The smallest case is $F_{3,3,3;0}$,
which is not nice;
the argument is similar to that of
Example~\ref{sec_transversal}.\ref{exam.transversal.5}
but is longer.

{\em Class 5.}\/  The matroid $L_{n_1,n_2,n_3;0}$
is nice (Example~\ref{sec_transversal}.\ref{exam.transversal.4}).
All matroids $L_{n_1,n_2,n_3;n'}$ with $n' \ge 1$
are not nice (Example~\ref{sec_transversal}.\ref{exam.transversal.5}).

{\em Class 6.}\/  By deletion from $L_{n_1,n_2,n_3;0}$,
we can conclude that $L_{n_1,n_2;n'}$ is nice when $n' \le 2$.
The smallest remaining case is $L_{3,3;3}$, which is not nice;
the argument is similar to that of
Example~\ref{sec_transversal}.\ref{exam.transversal.5}
but is longer.

{\em Class 7.}\/  By deletion from $L_{n_1,n_2,n_3;0}$,
we can conclude that $L_{n_1;n'}$ is nice when $n' \le 4$.
The smallest remaining case is $L_{3;5}$, which is not nice
(Example~\ref{sec_transversal}.\ref{exam.transversal.6}).

{\em Class 8.}\/  All uniform matroids are nice
(Example~\ref{sec_transversal}.\ref{exam.transversal.1}).

\subsection{A wild speculation}  \label{sec_wild_transversal}

Not all transversal matroids are nice;
but this means only that our {\em method for proving}\/
the half-plane property fails in these cases,
not that the half-plane property itself fails.
Indeed, we do not know a single example of a transversal matroid
that fails the half-plane property.
Moreover, we have conducted extensive numerical experiments
(see Section~\ref{sec_numerical})
on the rank-3 transversal matroids discussed in the previous subsection,
and it seems plausible that all rank-3 transversal matroids
have the half-plane property.
Might {\em all}\/ transversal matroids (and hence all gammoids)
have the half-plane property?
Unfortunately, we have no idea how to prove this.

\bigskip

{\bf Remark.}
The strongly base-orderable matroids \cite{Ingleton_77,Brualdi_87}
form a minor-closed class that contains all the transversal matroids;
so one might entertain the even stronger conjecture that
``every strongly base-orderable matroid has the half-plane property''.
However, this conjecture is false,
since the matroid $F_7^{-3}$ is strongly base-orderable
but does not have the half-plane property
(Example~\ref{sec_counterexamples}.\ref{exam.counterexamples.7} below).
To see that $F_7^{-3}$ is strongly base-orderable,
it suffices to note that a rank-3 matroid is strongly base-orderable
$\;\Longleftrightarrow\;$ it is base-orderable
$\;\Longleftrightarrow\;$ it has no restriction isomorphic to $M(K_4)$
\cite{Ingleton_77}.

\section{Counterexamples}   \label{sec_counterexamples}

In this section we use Proposition~\ref{prop.generalrank}(a)$\implies$(b)
to show that certain polynomials do {\em not}\/ have the half-plane property.

\bexam
  \label{exam.counterexamples.1}
Consider the Fano matroid $F_7$ with the ground set numbered
as shown in Figure~\ref{fig_F7etc} (see Appendix~\ref{app_matroids}).
Its basis set $\scrb(F_7)$ consists of all 3-element subsets of [7]
except $\{1,2,3\}$, $\{3,4,5\}$, $\{1,5,6\}$, $\{1,4,7\}$, $\{2,5,7\}$,
$\{3,6,7\}$ and $\{2,4,6\}$.
Let $P$ be the basis generating polynomial $P_{\scrb(F_7)}$.
If we take $x = \chi_{\{1,2,4,5\}}$ and $y = \chi_{\{3,6,7\}}$
(where $\chi_A$ denotes characteristic function of the set $A$),
we obtain $p_{x,y}(\zeta) = 4\zeta^3 + 12\zeta^2 + 12 \zeta$,
whose roots are $\zeta = 0, (-3 \pm \sqrt{3}\,i)/2$.
So the Fano matroid does not possess the half-plane property.

In fact, the Fano matroid is minor-minimal for failing to have
the half-plane property.
Indeed, it follows from Proposition~\ref{prop.nle6}
that {\em every}\/ 7-element matroid lacking the half-plane property
is minor-minimal.
%
\eexam

\bexam
  \label{exam.counterexamples.2}
Consider next the non-Fano matroid $F_7^-$ (Figure~\ref{fig_F7etc}),
which is obtained from $F_7$ by relaxing the circuit-hyperplane $\{2,4,6\}$,
so that $\scrb(F_7^-) = \scrb(F_7) \cup \{ \{2,4,6\} \}$.
With the same choices of $x$ and $y$,
we obtain $p_{x,y}(\zeta) = 4\zeta^3 + 13\zeta^2 + 12 \zeta$,
whose roots are $\zeta = 0, (-13 \pm \sqrt{23}\,i)/8$.
So the non-Fano matroid does not possess the half-plane property either;
and, as noted above, it is minor-minimal.
%
\eexam

\bexam
  \label{exam.counterexamples.3}
Consider next the matroid $F_7^{--}$ obtained from $F_7^-$
by relaxing $\{1,4,7\}$ (Figure~\ref{fig_F7etc}).
Choosing $x=\chi_{\{1,4,7\}}$ and $y = \chi_{\{2,3,5,6\}}$,
we obtain $p_{x,y}(\zeta) = \zeta^3 + 12\zeta^2 + 13\zeta + 4$,
whose roots are $\zeta \approx -10.834170$ and
$\zeta \approx -0.582915 \pm 0.171501\,i$.
So $F_7^{--}$ does not possess the half-plane property either;
and, as noted above, it is minor-minimal.
%
\eexam

\bexam
  \label{exam.counterexamples.4}
A similar approach shows that the matroid $M(K_4) + e$
obtained from $F_7^{--}$
by relaxing $\{3,4,5\}$ (Figure~\ref{fig_F7etc})
also fails to possess the half-plane property.
But we can no longer choose $x$ and $y$ to be characteristic functions.
Instead, let us take a unified approach to $F_7$, $F_7^-$, $F_7^{--}$
and $M(K_4) + e$ that gives additional insight into {\em why}\/
the half-plane property fails.

Let us start with the graphic matroid $M(K_4)$,
which of course {\em does}\/ have the half-plane property
(by Theorem~\ref{thm1.1}).
Let $K_4$ have vertex set $\{1,2,3,4\}$.
If we take $x = \chi_{\{12,13,14\}}$ and $y = \chi_{\{23,24,34\}}$,
we find $p_{x,y}(\zeta) \equiv P_{\scrb(M(K_4))}(\zeta x+y) = \zeta(\zeta+3)^2$,
which has a double root at $\zeta = -3$.
So $M(K_4)$ satisfies the condition of Proposition~\ref{prop.generalrank}(b),
but ``just barely'':  by suitable perturbations we may be able to split
the double root into a pair of complex-conjugate roots,
thereby proving the failure of the half-plane property
for the perturbed matroid.

In the matroids $F_7$ {\em et al}\/.,
let us consider 4 to be the ``new element'',
i.e.\ $F_7 \setminus 4 \simeq F_7^- \setminus 4 \simeq \ldots \simeq M(K_4)$.
Then the above choice of $x$ and $y$ corresponds to taking
$y$ to be the characteristic function of a 3-point line in $M(K_4)$
and $x$ to be the characteristic function of the complementary set,
e.g.\ $x = \chi_{\{1,2,5\}}$ and $y = \chi_{\{3,6,7\}}$.
So let us perturb this slightly, by taking
$x = \chi_{\{1,2,5\}} + \epsilon \chi_{\{4\}}$
and $y = \chi_{\{3,6,7\}} + a\epsilon \chi_{\{4\}}$
with $\epsilon, a \ge 0$.
(Our previous choices for $F_7$ and $F_7^-$
 correspond to $\epsilon = 1$ and $a=0$.)
We then have:
\begin{quote}
\begin{itemize}
   \item[$F_7$:] $p_{x,y}(\zeta) = (1+3\epsilon)\zeta^3 +
                     [6+(6+3a)\epsilon]\zeta^2 +
                     [9+(3+6a)\epsilon]\zeta + 3a\epsilon$ \\[1mm]
     Roots $\zeta = -(a/3)\epsilon \,+\, O(\epsilon^{3/2})$ and
           $\zeta = -3 \pm 2 \sqrt{a-3} \, \epsilon^{1/2} \,+\, O(\epsilon)$.
   \item[$F_7^-$:] $p_{x,y}(\zeta) = (1+3\epsilon)\zeta^3 +
                     [6+(7+3a)\epsilon]\zeta^2 +
                     [9+(3+7a)\epsilon]\zeta + 3a\epsilon$ \\[1mm]
     Roots $\zeta = -(a/3)\epsilon \,+\, O(\epsilon^{3/2})$ and
           $\zeta = -3 \pm \sqrt{3(a-3)} \, \epsilon^{1/2} \,+\, O(\epsilon)$.
   \item[$F_7^{--}$:] $p_{x,y}(\zeta) = (1+3\epsilon)\zeta^3 +
                     [6+(8+3a)\epsilon]\zeta^2 +
                     [9+(3+8a)\epsilon]\zeta + 3a\epsilon$ \\[1mm]
     Roots $\zeta = -(a/3)\epsilon \,+\, O(\epsilon^{3/2})$ and
           $\zeta = -3 \pm \sqrt{2(a-3)} \, \epsilon^{1/2} \,+\, O(\epsilon)$.
   \item[$M(K_4) + e$:] $p_{x,y}(\zeta) = (1+3\epsilon)\zeta^3 +
                     [6+(9+3a)\epsilon]\zeta^2 +
                     [9+(3+9a)\epsilon]\zeta + 3a\epsilon$ \\[1mm]
     Roots $\zeta = -(a/3)\epsilon \,+\, O(\epsilon^{3/2})$ and
           $\zeta = -3 \pm \sqrt{a-3} \, \epsilon^{1/2} \,+\, O(\epsilon)$.
\end{itemize}
\end{quote}
So any choice of $a \in [0,3)$ will yield nonreal roots for
small $\epsilon > 0$ for all four matroids.

Let us conclude this example by remarking that
a second (nonisomorphic) choice of $x,y$
in $M(K_4)$ also yields a double root:
for $x = \chi_{\{12,34\}}$ and $y = \chi_{\{13,14,23,24\}}$,
we find $p_{x,y}(\zeta) = 4(\zeta+1)^2$,
which has a double root at $\zeta = -1$.
An analogous perturbation of this choice also yields nonreal roots
for $F_7$ {\em et al}\/., provided that $a \in [0,1)$.
%
\eexam

\bexam
  \label{exam.counterexamples.5}
Consider the polynomial
\be
   P_\mu(z_1,\ldots,z_7)  \;=\;  P_{\scrb(F_7)}(z)  \,+\,  \mu z_2 z_4 z_6
   \;.
\ee
For $\mu=0$ (resp.\ $\mu=1$) this is the basis generating polynomial
of $F_7$ (resp.\ $F_7^-$).
For $\mu=4$ it is the polynomial $Q_A(z) = \det(A Z A^{\rm T})$
obtained from $Z = \diag(z_1,\ldots,z_7)$ and the matrix
\be
   A  \;=\;  \left[ \begin{array}{ccccccc}
                        1 & 1 & 0 & 0 & 0 & 1 & 1 \\
                        0 & 1 & 1 & 1 & 0 & 0 & 1 \\
                        0 & 0 & 0 & 1 & 1 & 1 & 1
                    \end{array}
             \right]
   \;,
\ee
which represents $F_7^-$ over any field of characteristic $\neq 2$
(in particular, over $\C$);
by Theorem~\ref{thm.QA}, $Q_A$ has the half-plane property.
As in Example~\ref{sec_counterexamples}.\ref{exam.counterexamples.4},
let us choose $x = \chi_{\{1,2,5\}} + \epsilon \chi_{\{4\}}$
and $y = \chi_{\{3,6,7\}} + a\epsilon \chi_{\{4\}}$
with $\epsilon, a \ge 0$.
We then obtain
$p_{x,y}(\zeta) = (1+3\epsilon)\zeta^3 + [6+(6+\mu+3a)\epsilon]\zeta^2 +
                     [9+(3+6a+\mu a)\epsilon]\zeta + 3a\epsilon$,
whose roots are $\zeta = -(a/3)\epsilon \,+\, O(\epsilon^{3/2})$ and
$\zeta = -3 \pm \sqrt{(4-\mu)(a-3)} \, \epsilon^{1/2} \,+\, O(\epsilon)$.
So, if $-\infty < \mu < 4$ (resp.\ $\mu > 4$),
then any choice $0 \le a < 3$ (resp.\ $a > 3$)
will yield nonreal roots for small $\epsilon > 0$.
So $P_\mu$ has the half-plane property {\em only}\/ for $\mu = 4$!
\eexam

\bexam
  \label{exam.counterexamples.6}
Consider, more generally, the polynomial
\be
   P_{\mu,\nu,\rho}(z_1,\ldots,z_7)  \;=\;
   P_{\scrb(F_7)}(z)  \,+\,  \mu z_2 z_4 z_6  \,+\,  \nu z_1 z_4 z_7  \,+\,
                             \rho z_3 z_4 z_5
   \;.
\ee
For $(\mu,\nu,\rho) = (0,0,0)$ this is $P_{\scrb(F_7)}$;
for $(\mu,\nu,\rho) = (1,0,0)$ it is $P_{\scrb(F_7^-)}$;
for $(\mu,\nu,\rho) = (1,1,0)$ it is $P_{\scrb(F_7^{--})}$;
for $(\mu,\nu,\rho) = (1,1,1)$ it is $P_{\scrb(M(K_4) + e)}$.

Some cases of the polynomial $P_{\mu,\nu,\rho}$
are of the form $Q_A(z) = \det(A Z A^*)$
[where $Z = \diag(z_1,\ldots,z_7)$] for a suitable complex matrix $A$,
and thus have the half-plane property by Theorem~\ref{thm.QA}.
Consider, for instance, the matrix
\be
   A  \;=\;  \left[ \begin{array}{ccccccc}
                        1 & 1 & 0 & a & 0 & 1 & 1 \\
                        0 & 1 & 1 & 1 & 0 & 0 & 1 \\
                        0 & 0 & 0 & 1 & 1 & b & b
                    \end{array}
             \right]
\ee
For $a=0$ and $b = e^{i\theta}$,
we have $Q_A = P_{\mu,\nu,\rho}$ with
$(\mu,\nu,\rho) = (2 + 2\cos\theta, 2 - 2\cos\theta, 0)$
[that is, $\mu,\nu \ge 0$ with $\mu+\nu=4$ and $\rho=0$].\footnote{
   We remark that the matrix $A$ represents $F_7^{--}$ over any field of
   at least four elements, provided that $a=0$ and $b \notin \{0,1,-1\}$.
}
For $a = e^{\pm \pi i/3}$ and $b=1$,
we have $Q_A = P_{\mu,\nu,\rho}$ with $(\mu,\nu,\rho) = (3,0,1)$;
and for $a = e^{\pm \pi i/3}$ and $b=e^{\mp 2\pi i/3}$,
we have $Q_A = P_{\mu,\nu,\rho}$ with $(\mu,\nu,\rho) = (0,3,1)$.
Moreover, by appropriately permuting columns of the matrix $A$
({\em without}\/ permuting the column labels),
the triplet $(\mu,\nu,\rho)$ can be permuted at will.\footnote{
   This is because there are automorphisms of $F_7$
   that realize an arbitrary permutation of the lines 246, 147 and 345.
}
We do not know whether any additional cases of $P_{\mu,\nu,\rho}$
can be obtained from a determinant.

For $0 \le \mu,\nu,\rho \le 2$ with $\mu+\nu+\rho = 4$,
the polynomial $P_{\mu,\nu,\rho}$ arises from the principal extension of
$F_7 \setminus 4 \simeq M(K_4)$ by the new element 4
with weights
\begin{subeqnarray}
   \lambda_2 \;=\; \lambda_6 & = & 1 - \mu/2  \\
   \lambda_1 \;=\; \lambda_7 & = & 1 - \nu/2  \\
   \lambda_3 \;=\; \lambda_5 & = & 1 - \rho/2
\end{subeqnarray}
and thus has the half-plane property by Proposition~\ref{prop.truncation}.

Finally, as in Example~\ref{sec_counterexamples}.\ref{exam.counterexamples.5},
let us choose $x = \chi_{\{1,2,5\}} + \epsilon \chi_{\{4\}}$
and $y = \chi_{\{3,6,7\}} + a\epsilon \chi_{\{4\}}$
with $\epsilon, a \ge 0$.
We then obtain the same $p_{x,y}$
as in Example~\ref{sec_counterexamples}.\ref{exam.counterexamples.5},
except that $\mu$ is replaced by $\mu+\nu+\rho$.
So we obtain nonreal roots for small $\epsilon > 0$
whenever $\mu + \nu + \rho \neq 4$.
It follows (using also Theorem~\ref{thm.same-phase})
that $P_{\mu,\nu,\rho}$ has the half-plane property {\em only}\/
when $\mu,\nu,\rho \ge 0$ with $\mu+\nu+\rho=4$.
But we do not know whether these necessary conditions are sufficient.
\eexam

\bexam
  \label{exam.counterexamples.7}
Consider the matroid $F_7^{-3}$
obtained from $F_7^-$ by relaxing the circuit-hyperplanes
$\{2,5,7\}$ and $\{3,6,7\}$ (Figure~\ref{fig_F7etc}).
Let us take $x = \chi_{\{1,4,5\}} + \epsilon\chi_{\{2\}}$
and $y = \chi_{\{3,6,7\}}$.
Then $p_{x,y}(\zeta) = (1+3\epsilon)\zeta^3 + (6+8\epsilon)\zeta^2 +
  (9+3\epsilon)\zeta + 1$,
which has nonreal roots whenever
$0.090685 \ltapprox \epsilon \ltapprox 0.494485$.
So $F_7^{-3}$ does not have the half-plane property;
and, as noted above, it is minor-minimal.
%

Note, finally, that $F_7^{-3}$ is a relaxation of $P_7$,
which is a $\sqrt[6]{1}$-matroid and hence has the half-plane property
(see Section~\ref{sec_determinant}).
So relaxation does {\em not}\/ in general preserve the half-plane property.
\eexam

\bexam
  \label{exam.counterexamples.8}
Consider the rank-4 matroid $P_8$ represented over any field
of characteristic $\neq 2$ by the matrix
\be
\left[
   \begin{array}{rrrr|rrrr}
      1 & 0 & 0 & 0 &    0 & 1 & 1 &  2  \\
      0 & 1 & 0 & 0 &    1 & 0 & 1 &  1  \\
      0 & 0 & 1 & 0 &    1 & 1 & 0 &  1  \\
      0 & 0 & 0 & 1 &    2 & 1 & 1 &  0
   \end{array}
\right]
\ee
Its basis set $\scrb(P_8)$ consists of all 4-element subsets of [8] except
$\{1,2,3,8\}$, $\{1,2,4,7\}$, $\{1,3,4,6\}$, $\{2,3,4,5\}$, $\{1,4,5,8\}$,
$\{2,3,6,7\}$, $\{1,5,6,7\}$, $\{2,5,6,8\}$, $\{3,5,7,8\}$ and $\{4,6,7,8\}$.
If we take $x = \chi_{\{1,4,5,8\}}$ and $y = \chi_{\{2,3,6,7\}}$,
we obtain $p_{x,y}(\zeta) = 16\zeta^3 + 28\zeta^2 + 16\zeta$,
whose roots are $\zeta = 0, (-7 \pm \sqrt{15}\,i)/8$.
So $P_8$ does not possess the half-plane property.
Moreover, it is minor-minimal,
since all the single-element contractions (resp.\ deletions) of $P_8$
are isomorphic to $P_7$ (resp.\ $P_7^*$),
which are $\sqrt[6]{1}$-matroids and hence have the half-plane property.

Observe now that $\{1,4,5,8\}$ and $\{2,3,6,7\}$
form the unique pair of disjoint circuit-hyperplanes in $P_8$.
Let us denote by $P'_8$ the matroid obtained from $P_8$ by relaxing
one of these circuit-hyperplanes (say, $\{1,4,5,8\}$),
and by $P''_8$ the matroid obtained by relaxing both of them.
Let us make the same choices of $x$ and $y$ as for $P_8$.
Then, for $P'_8$ one has
$p_{x,y}(\zeta) = \zeta^4 + 16\zeta^3 + 28\zeta^2 + 16\zeta$,
whose roots are $\zeta = 0$ and
$\zeta \approx -14.093869, -0.953065 \pm 0.476353\,i$.
And for $P''_8$ one has
$p_{x,y}(\zeta) = \zeta^4 + 16\zeta^3 + 28\zeta^2 + 16\zeta + 1$,
whose roots are $\zeta \approx -14.093459, -0.070955,
  -0.917793 \pm 0.397059\,i$.
So neither $P'_8$ nor $P''_8$ possesses the half-plane property.
We suspect (but are unable to prove) they are minor-minimal:
all the single-element contractions (resp.\ deletions) of $P'_8$
are isomorphic to $P_7$ or $P'_7$ (resp.\ $P_7^*$ or $(P'_7)^*$),
and all the single-element contractions (resp.\ deletions) of $P''_8$
are isomorphic to $P'_7$ (resp.\ $(P'_7)^*$);
and our numerical results (Section~\ref{sec_numerical})
suggest that $P'_7$ (and hence also its dual) has the half-plane property.
%
%
%
\eexam


\bexam
  \label{exam.counterexamples.9}
Consider the Pappus and non-Pappus matroids
(Figure~\ref{fig_pappus} in Appendix~\ref{app_matroids}).
Take $x = 2\chi_{\{1\}} + \chi_{\{3,4,6,7\}} + \epsilon \chi_{\{9\}}$
and $y = \chi_{\{2,5,8\}}$.
We then have:
\begin{quote}
\begin{quote}
\begin{itemize}
   \item[Pappus:]
      $p_{x,y}(\zeta) = (16+14\epsilon)\zeta^3 + (33+15\epsilon)\zeta^2 +
                             (18+3\epsilon)\zeta + 1$  \\[1mm]
      Roots $\zeta = -{1 \over 16} + O(\epsilon)$
         and $\zeta = -1 \pm (i\sqrt{2/15}) \epsilon^{1/2} + O(\epsilon)$.
   \item[non-Pappus:]
      $p_{x,y}(\zeta) = (16+14\epsilon)\zeta^3 + (33+16\epsilon)\zeta^2 +
                             (18+3\epsilon)\zeta + 1$  \\[1mm]
      Roots $\zeta = -{1 \over 16} + O(\epsilon)$
         and $\zeta = -1 \pm (i/\sqrt{15}) \epsilon^{1/2} + O(\epsilon)$.
\end{itemize}
\end{quote}
\end{quote}
So the Pappus and non-Pappus matroids do not have the half-plane property.
We do not know whether they are minor-minimal;
but we suspect that they are, since our numerical experiments
(Section~\ref{sec_numerical})
suggest that their two non-isomorphic deletions,
non-Pappus$\setminus 1$ and Pappus$\setminus e \simeq$ non-Pappus$\setminus 9$,
do have the half-plane property.
(Of course, all their contractions are rank-2 and hence have the
 half-plane property.)
\eexam


\bigskip

It is worth noting that non-Pappus$\setminus 9$
is ``borderline'' for the half-plane property
in the same sense that $M(K_4)$ is
(Example~\ref{sec_counterexamples}.\ref{exam.counterexamples.4}),
i.e.\ $p_{x,y}(\zeta)$ has a double root at $\zeta=-1$
for a suitable choice of $x,y$.
This fact suggests that other extensions of non-Pappus$\setminus 9$
might also fail to have the half-plane property.
This is indeed the case:

\bexam
  \label{exam.counterexamples.10}
Consider the matroid (non-Pappus$\setminus 9)+e$,
obtained from non-Pappus$\setminus 9$ by adding a new element freely
(let us call this new element 9).
Then its bases are those of non-Pappus plus $\{2,6,9\}$ and $\{3,5,9\}$.
Take $x = 2\chi_{\{1\}} + \chi_{\{3,4,6,7\}}$
and $y = \chi_{\{2,5,8\}} + \epsilon \chi_{\{9\}}$.
We then have:
\begin{quote}
\begin{itemize}
   \itemindent 2.2cm
   \item[\qquad(non-Pappus$\setminus 9)+e$:]
      $p_{x,y}(\zeta) = 16\zeta^3 + (33+14\epsilon)\zeta^2 +
                             (18+18\epsilon)\zeta + (1+3\epsilon)$  \\[1mm]
      \hspace*{2.05cm}
      Roots $\zeta = -{1 \over 16} + O(\epsilon)$
         and $\zeta = -1 \pm (i/\sqrt{15}) \epsilon^{1/2} + O(\epsilon)$.
\end{itemize}
\end{quote}
So (non-Pappus$\setminus 9)+e$
does not have the half-plane property.
\eexam

\section{Numerical experiments}
   \label{sec_numerical}

Given a homogeneous multiaffine polynomial $P(x_1,\ldots,x_n)$,
we have searched numerically for counterexamples to the half-plane property
using two methods:

\medskip

{\em Elementary method.}\/
Choose $x_1,\ldots,x_{n-1}$ uniformly at random in the rectangle
$(0,1) + (-1,1)\,i$;
then solve $P(x_1,\ldots,x_n) = 0$ for $x_n$.
If $\Re x_n > 0$, we have found a counterexample to the half-plane property.
(Thanks to homogeneity, there is no loss of generality in choosing
 $x_1,\ldots,x_{n-1}$ in the specified rectangular subset of $\C$.)

\medskip

{\em Method using Proposition~\ref{prop.generalrank}.}\/
Choose vectors $a,b$ uniformly at random in $[0,1]^n$,
and compute the roots of the univariate polynomial
$p_{a,b}(\zeta) = P(\zeta a + b)$.
If at least one of these roots has a nonzero imaginary part,
we have found a counterexample to the half-plane property.
(Thanks to homogeneity, there is no loss of generality in choosing
 $a,b$ in $[0,1]^n \subset \R^n$.)

\medskip

\begin{table}[tp]
\centering
\begin{tabular}{|l|c|r|c|r|}
\cline{2-5}
 \multicolumn{1}{c}{\mbox{}}
        & \multicolumn{2}{|c|}{Elementary} &
          \multicolumn{2}{|c|}{Method using} \\
 \multicolumn{1}{c}{\mbox{}}
        & \multicolumn{2}{|c|}{Method} &
          \multicolumn{2}{|c|}{Proposition~\ref{prop.generalrank}} \\
\hline
\multicolumn{1}{|c|}{Matroid} &
     \# tries & \# counterexamples & \# tries & \# counterexamples \\
\hline\hline
$F_7$         &  $10^8$  & 4152 &  $10^6$  &   81566 \\ 
$F_7^-$       &  $10^8$  &  687 &  $10^6$  &   42620 \\ 
$F_7^{--}$    &  $10^8$  &   34 &  $10^6$  &   12794 \\ 
$M(K_4)+e$     &  $10^8$  &   0 &  $10^6$  &    1060 \\ 
$F_7^{-3}$     &  $10^8$  &   0 &  $10^6$  &     695 \\ 
$F_7^{-4}$     &  $10^8$  &   0 &  $10^8$  &       0 \\ 
$\scrw^3 + e$  &  $10^8$  &   0 &  $10^8$  &       0 \\ 
\hline
$\scrw^{3+}$   &  $10^8$  &   0 &  $10^8$  &       0 \\ 
\hline
$P'_7$         &  $10^8$  &   0 &  $10^8$  &   0 \\ 
\hline
$F_7^{-4} + e$ &   ---    & --- &  $10^7$  &       0 \\ 
$\scrw^3 +e+f$ &   ---    & --- &  $10^7$  &       0 \\ 
$\scrw^{3+}+e$ &   ---    & --- &  $10^7$  &       0 \\ 
$P'_7 + e$     &   ---    & --- &  $10^7$  &       0 \\ 
\hline
$P_8$          &  $10^8$  & 278 &  $10^6$  &   10930 \\ 
$P'_8$         &  $10^8$  & 114 &  $10^6$  &    5590 \\ 
$P''_8$        &  $10^8$  &  40 &  $10^6$  &    2723 \\ 
\hline
$V_8$ (V\'amos) & $10^8$  &   0 &  $10^7$  &       0 \\ 
\hline
Pappus         &  $10^8$  &   0 & $10^6$  &   544 \\    
non-Pappus     &  $10^8$  &   0 & $10^6$   &   17 \\    
(non-Pappus $\drop\,9)+e$
               &  $10^8$  &   0 & $10^7$   &    6 \\    
non-Pappus $\drop$ 1 &  $10^8$  &   0 & $10^7$  &  0 \\ 
non-Pappus $\drop$ 9 &  $10^8$  &   0 & $10^8$  &  0 \\
\hline
$C_{3,3,3;3}$  &  $10^7$  &   0 &  $10^7$  &       0 \\ 
$C_{4,4,4;4}$  &  $10^7$  &   0 &  $10^6$  &       0 \\ 
$C_{5,5,5;5}$  &  $10^7$  &   0 &  $10^6$  &       0 \\ 
$C_{6,6,6;6}$  &  $10^6$  &   0 &  $10^5$  &       0 \\ 
$C_{7,7,7;7}$  &  $10^6$  &   0 &  $10^5$  &       0 \\ 
\hline
\end{tabular}
\vspace{2mm}
\caption{
   Results of numerical experiments to test the half-plane property.
}
 \label{table_numerical}
\end{table}

We first applied these methods to some matroids for which
the half-plane property is known to fail ($F_7$, $F_7^-$, \ldots),
in order to get a rough feeling for the rate of finding counterexamples.
We then applied them to some matroids for which
the half-plane property is an open question
($F_7^{-4}$, $\scrw^3 + e$, $\scrw^{3+}$, $P'_7$, \ldots).
The results of our numerical experiments are shown in
Table~\ref{table_numerical}.
As can be seen, the method using Proposition~\ref{prop.generalrank}
is vastly more powerful than the elementary method:
in those cases where the half-plane property fails,
the proportion of counterexamples is larger by a factor of
$\approx 2000$--20000 or even more.
Roughly speaking, one pair $a,b$ in Proposition~\ref{prop.generalrank}
corresponds to a large set of counterexamples in the elementary method.
Indeed, if we had used only the elementary method
we would have failed to detect the failure of the half-plane property
for $M(K_4)+e$, $F_7^{-3}$, Pappus, non-Pappus and (non-Pappus $\drop\,9)+e$,
even with $10^8$ tries.

Our numerical results suggest, first of all, that the matroids
$F_7^{-4}$, $\scrw^3 + e$, $\scrw^{3+}$, $P'_7$
and their free extensions,
as well as $V_8$, non-Pappus $\drop$ 1 and non-Pappus $\drop$ 9,
probably {\em do}\/ all have the half-plane property.
One urgent problem, therefore, is to try to prove (or disprove)
these alleged facts;
the proofs may well require new tools.
We find it particularly surprising that
the V\'amos matroid $V_8$ has (or appears to have) the half-plane property:
since this matroid is not representable over any field (much less over $\C$),
neither the determinant condition nor the permanent condition applies to it
--- nor can the constructions in Section~\ref{sec_constructions} lead to it,
starting from a polynomial whose support is the collection of bases
of a representable matroid ---
so we have no idea why it should have
even the {\em weak}\/ half-plane property.

Finally, we made a systematic study of the half-plane property for
rank-3 transversal matroids,
all of which are restrictions of $C_{n,n,n;n}$ for some $n \ge 3$
(see Section~\ref{sec_rank3_transversal}).
The results in Table~\ref{table_numerical} strongly suggest that
{\em all}\/ rank-3 transversal matroids have the half-plane property.
And they suggest the bold conjecture that perhaps all transversal matroids
of {\em any}\/ rank have the half-plane property.

\section{Open questions}
   \label{sec_open}

We conclude by presenting some open questions raised by,
or related to, the results of this paper.

\subsection{The strong half-plane property}

Let us say that a polynomial $P$ in $n$ complex variables has
the {\em strong half-plane property}\/
if 
\be
   \left. \partial^{\bf m}_x |P(x+iy)|^2 \right|_{x=0}  \;\ge\; 0
 \label{eq.strongHPP}
\ee
for all multi-indices ${\bf m}$ and all $y \in \R^n$.
It is easy to show that \reff{eq.strongHPP} implies the half-plane property:
Fix $y \in \R^n$;  then either
$\left. \partial^{\bf m}_x |P(x+iy)|^2 \right|_{x=0} > 0$
for at least one multi-index ${\bf m}$,
which implies $P(x+iy) \neq 0$ for all vectors $x > 0$,
or else
$\left. \partial^{\bf m}_x |P(x+iy)|^2 \right|_{x=0} = 0$ for all ${\bf m}$,
which implies $P(x+iy) = 0$ for all vectors $x \ge 0$
and hence $P \equiv 0$ by analytic continuation from a real environment.

\begin{question}
Does the half-plane property imply the strong half-plane property?
\end{question}

\noindent
It is easy to see that the answer is affirmative when $n=1$
(just factor $P$);  hence it is also true for general $n$
if one considers only multi-indices ${\bf m}$ with a single nonzero component.
But we do not know whether it is true in general.
This question was inspired by Newman's ``strong Lee--Yang theorem''
\cite[Section 3]{Newman_74}.

\subsection{Same-phase theorem}

\begin{question}
Can the same-phase theorem
(Theorems~\ref{thm.same-phase} and \ref{thm.same-phase_parity})
be extended to a larger class of polynomials?
\end{question} 

\noindent
The Fettweis--Basu lemma (Lemma~\ref{lemma.fettweis}),
applied in many variables successively, may be a useful tool.

\subsection{Support of polynomials with the half-plane property}

In Theorems~\ref{thm.matroidal} and \ref{thm.matroidal_non-multiaffine}
we provided some {\em necessary}\/ conditions for a subset $S \subseteq \N^E$
to be the support of a homogeneous polynomial with the half-plane property.
It is natural to seek a generalization to the non-homogeneous case:

\begin{problem}
Find necessary conditions for a subset $S \subseteq \N^E$
to be the support of a polynomial with the half-plane property.
\end{problem}

\noindent
Two special cases of this problem were posed earlier:

\begin{question}[= Question~\protect\ref{quest.matroidal.1}]
If $P$ is multiaffine and has the half-plane property,
is $\supp(P)$ a delta-matroid?  What if $P$ also has the same-phase property?
\end{question}
 
\begin{question}[= Question~\protect\ref{quest.matroidal.2}]
Assume that $P$ has the half-plane property and
has definite parity
(i.e.\ $\m,\m'\in\supp(P)$ implies that $|\m|\equiv|\m'|$ mod 2).
Is $\supp(P)$ then a jump system?
(Recall from Theorem~\ref{thm.same-phase_parity}
 that all such polynomials have the same-phase property.)
\end{question}

More generally, one can pose the probably-much-more-difficult problem:

\begin{problem}
 \label{prob.support.2}
Find necessary and sufficient conditions for a subset $S \subseteq \N^E$
to be the support of a polynomial with the half-plane property.
\end{problem}

\noindent
Already in the homogeneous multiaffine case
this question looks difficult,
and constitutes the converse of Theorem~\ref{thm.matroidal}:

\begin{question}[= Question~\protect\ref{question.weak_HPP}]
Does every matroid $M$ have the weak half-plane property?
And if not, which ones do?
(For instance, does the Fano matroid $F_7$ have the weak half-plane property?)
\end{question}   

\noindent
By Corollary~\ref{cor.determinant},
every matroid representable over $\C$ has the weak half-plane property.
But we know {\em nothing}\/ beyond this.
Indeed, both of our {\em ab initio}\/ methods for constructing
polynomials with the half-plane property
--- the determinant method (Theorem~\ref{thm.QA})
and the permanent method (Theorem~\ref{thm.permanent}) ---
lead always to polynomials $P$ whose support matroid is $\C$-representable;
and all the constructions of Section~\ref{sec_constructions},
applied to a polynomial $P$ whose support matroid is $\C$-representable,
lead to another polynomial with the same property.\footnote{
   The operations of
   Sections~\ref{sec_constructions.1}--\ref{sec_constructions.4},
   applied to a matroid that is representable over a field $F$,
   always produce another $F$-representable matroid:
   see \cite[Proposition 3.2.4]{Oxley_92} for deletion and contraction,
   \cite[Corollary 2.2.9]{Oxley_92} for duality,
   \cite[Proposition 4.2.15]{Oxley_92} for direct sum,
   \cite[Proposition 7.1.21]{Oxley_92} for
      parallel connection and series connection,
   and \cite[Proposition 7.1.23]{Oxley_92} for 2-sum.
   The other operations of Section~\ref{sec_constructions}
   (principal extension, truncation, coextension and cotruncation,
    and full-rank matroid union)
   can lead out of the class of $F$-representable matroids,
   but always produce a matroid that is representable over
   some finite extension field of $F$:
   see \cite[Propositions 7.3.5 and 7.4.17]{Brylawski_86}
   and \cite{Piff_70}.
   In particular, applied to a $\C$-representable matroid,
   they always yield another $\C$-representable matroid.
}
So even the following is an open problem:

\begin{problem}
Construct a polynomial $P$ with the half-plane property
whose support matroid is not representable over $\C$.
(Or prove that it is impossible.)
\end{problem}

\noindent
Note that our numerical results (Section~\ref{sec_numerical})
suggest that the V\'amos matroid $V_8$,
which is not representable over any field, does have the half-plane property.

\bigskip

Theorems~\ref{thm.matroidal} and \ref{thm.matroidal_non-multiaffine}
have the general form:  if $P(x) = \sum_{{\bf m}} a_{\bf m} x^{\bf m}$ 
has the half-plane property and certain coefficients $a_{\bf m}$, $a_{\bf m'}$
are nonzero, then so are (one or more of) certain other coefficients
$a_{\bf m''}$.
It is natural to ask whether this qualitative result can be extended
to a {\em quantitative}\/ lower bound on some combination of
the corresponding coefficients $a_{\bf m''}$.

\begin{question}
Can Theorems~\ref{thm.matroidal} and \ref{thm.matroidal_non-multiaffine}
be extended to quantitative inequalities on the coefficients?
\end{question} 

\noindent
This question is vaguely reminiscent of Kung's discussion of the
heuristic analogy between determinantal identities
and basis-exchange properties in matroids \cite{Kung_86,Kung_01}.
One approach might be to use
Proposition~\ref{prop.generalrank}(a)$\implies$(b);
if successful, this would be a nice generalization
of Theorem~\ref{thm.rank2} to the higher-rank case.

\subsection{$(F,G)$-representability of matroids}

\begin{question}[= Question~\protect\ref{quest.contraction}]
When $-1 \notin G$,
is the class of $(F,G)$-representable matroids always
closed under contraction?
\end{question}

\begin{question}
If $-1 \notin G$,
does weak $(F,G)$-representability necessarily imply $(F,G)$-representability?
\end{question}

\begin{question}
Is Theorem~\ref{thm.sixthroot} the prototype of a more general theorem
asserting that, in certain cases,
every $(F,G)$-representable matroid is in fact $(F,G')$-representable
for some specified subgroup $G' \subsetneq G$?
\end{question}

\subsection{Determinant class, permanent class, and the half-plane property}

In this paper we have given two distinct
methods for constructing polynomials with the half-plane property:
the determinant construction (Theorem~\ref{thm.QA})
and the permanent construction (Theorem~\ref{thm.permanent}).
It is natural to ask what is the relation (if any)
between these two classes of polynomials,
and between them and the class of all polynomials
with the half-plane property.
Let us pose this question as follows:

Fix integers $r$ and $n$ ($0 \le r \le n$),
and let $\scrp_{r,n}$ be the set of degree-$r$ homogeneous multiaffine
polynomials in $n$ variables with complex coefficients,
$P(x) = \sum_{S \subseteq [n], |S|=r} a_S x^S$.
Let $\scrp_{r,n}^+$ be the subset of $\scrp_{r,n}$
consisting of polynomials with nonnegative coefficients.
By identifying a polynomial with its coefficients,
the spaces $\scrp_{r,n}$ and $\scrp_{r,n}^+$ can be thought of as
$\C^{n \choose r}$ and $[0,\infty)^{n \choose r}$, respectively.
Let us now define three subsets of $\scrp_{r,n}^+$
(thanks to Theorem~\ref{thm.same-phase},
 there is no loss of generality in restricting attention to $\scrp_{r,n}^+$):
\begin{itemize}
   \item $\scrh_{r,n}$:  the polynomials with the half-plane property.
   \item $\scrd_{r,n}$:  the polynomials of ``determinant class'',
      i.e.\ $P(x) = \det(A X A^*)$
      [or equivalently $a_S = |\det(A \restrict S)|^2$]
      for some $r \times n$ complex matrix $A$.
   \item $\scrm_{r,n}$:  the polynomials of ``permanent (or matching) class'',
      i.e.\ $P(x) = \per(\Lambda X)$
      [or equivalently $a_S = \per(\Lambda \restrict S)$]
      for some $r \times n$ nonnegative matrix $\Lambda$.
\end{itemize}
We have shown that
$\scrp_{r,n}^+ \supseteq \scrh_{r,n} \supseteq \scrd_{r,n} \cup \scrm_{r,n}$,
and these containments are in general strict
(see Example~\ref{sec_open}.\ref{exam_open.1} below).

\begin{question}
What is the relation between the spaces
   $\scrp_{r,n}^+$, $\scrh_{r,n}$, $\scrd_{r,n}$ and $\scrm_{r,n}$?
\end{question}

Most ambitiously, we can ask for a complete characterization of the
determinant class $\scrd_{r,n}$ and the permanent class $\scrp_{r,n}$:

\begin{problem}
 \label{prob.determinant}
Find necessary and sufficient conditions for a set of
nonnegative real numbers $\{ a_S \}_{S \subseteq [n], |S|=r}$
to be representable in the form $a_S = |\det(A \restrict S)|^2$
for some $r \times n$ complex matrix $A$.
\end{problem}

\begin{problem}
 \label{prob.permanent}
Find necessary and sufficient conditions for a set of
nonnegative real numbers $\{ a_S \}_{S \subseteq [n], |S|=r}$
to be representable in the form $a_S = \per(\Lambda \restrict S)$
for some $r \times n$ nonnegative matrix $\Lambda$.
\end{problem}

\noindent
For the analogous problem with $a_S = \det(A \restrict S)$,
the necessary and sufficient condition is given by the
Grassmann--Pl\"ucker syzygies
\cite[Proposition 1.6.1]{White_87a} \cite[Section 1.2]{Bokowski_89}.
But the presence of the modulus-square seems to make
Problem~\ref{prob.determinant} quite difficult.
Indeed, it seems nontrivial even in the rank-2 case.
Here is an even more special case that illustrates some of the complexity:

\bexam
 \label{exam_open.1}
Consider the homogeneous degree-2 polynomial in $n \ge 3$ variables
$P(z) = \half \sum_{i,j=1}^n  a_{ij} z_i z_j$ where
\be
   a_{ij}  \;=\;
   \cases{\mu  & if $(i,j) = (1,2)$ or $(2,1)$ \cr
          \noalign{\vskip 1mm}
          1    & in all other cases of $i \neq j$ \cr
          \noalign{\vskip 1mm}
          0    & if $i=j$
         }
 \label{eq.exam_open.1}
\ee
One then finds, after some calculation, that:
\begin{itemize}
   \item  $P$ has the half-plane property if and only if
     $0 \le \mu \le (2n-4)/(n-3)$.\footnote{
        {\sc Sketch of proof:}
        Let $A$ be the matrix defined in \reff{eq.exam_open.1}.
        Then it can be shown by induction on $n$ that
        $$\det(\lambda I - A) \;=\;
          (\lambda+1)^{n-3} (\lambda+\mu)
          [\lambda^2 - (\mu+n-3)\lambda + (n-3)\mu - (2n-4)]  \;.$$
        So $A$ has are $n-3$ eigenvalues $-1$,
        one eigenvalue $-\mu$, and a pair of eigenvalues
        $$ {\mu+n-3 \pm \sqrt{\mu^2 - 2(n-3)\mu + (n^2+2n-7)}  \over 2}  \;.$$
        The claim then follows from Theorem~\ref{thm.rank2}.
}
   \item $P$ belongs to the determinant class if and only if:
      \begin{itemize}
         \item For $n=3$, $0 \le \mu < \infty$.
         \item For $n=4$, $0 \le \mu \le 4$.
         \item For $n=5$, $\mu=0$ or $\mu=3$.
         \item For $n \ge 6$, never.
      \end{itemize}
   \item $P$ belongs to the permanent class
         if and only if $0 \le \mu \le 2$.\footnote{
         {\sc Sketch of proof:}
         Suppose that there exist nonnegative numbers
         $\{\lambda^{(1)}_i\}_{i=1}^n$ and $\{\lambda^{(2)}_i\}_{i=1}^n$
         such that
         $$a_{ij} \;=\; \lambda^{(1)}_i \lambda^{(2)}_j +
                        \lambda^{(1)}_j \lambda^{(2)}_i
          \qquad\hbox{for all $i \neq j$}   \eqno(*)$$
         By considering separately the cases $\lambda^{(1)}_1 =0$
         and $\lambda^{(1)}_1 \neq 0$ (and likewise for $\lambda^{(2)}_1$),
         one shows in each case that
         $\lambda^{(1)}_3 = \ldots = \lambda^{(1)}_n \equiv \alpha$ and
         $\lambda^{(2)}_3 = \ldots = \lambda^{(2)}_n \equiv \beta$.
         Using this fact together with (*),
         one shows (after some algebra) that $0 \le \mu \le 2$,
         and that every $\mu$ in this interval is attainable.
}
\end{itemize}
\eexam

\bigskip

One can ask also about general properties of the spaces
$\scrh_{r,n}$, $\scrd_{r,n}$ and $\scrm_{r,n}$.
For example, are they convex?
The answer is no in general;
indeed, a convex combination of polynomials in $\scrd_{r,n} \cap \scrm_{r,n}$
need not even lie in $\scrh_{r,n}$.
To see this, just let $P$ and $Q$ be polynomials
depending on disjoint subsets of variables,
e.g.\ $P(x_1,\ldots,x_4) = x_1 x_2$ and $Q(x_1,\ldots,x_4) = x_3 x_4$,
both of which lie in $\scrd_{2,4} \cap \scrm_{2,4}$;
but by Corollary~\ref{cor.disjoint},
$(P+Q)/2$ does not even have the half-plane property.

On the other hand, by Theorem~\ref{thm.matroidal} we can
``stratify'' $\scrh_{r,n}$ as
$\scrh_{r,n} = \bigcup\limits_M \scrh_{r,n}(M)$ where
\be
   \scrh_{r,n}(M)  \;=\;  \{P \in \scrh_{r,n} \colon\; \supp(P) = \B(M) \}
\ee
and the union runs over all rank-$r$ matroids $M$ on the ground set $[n]$;
and likewise for $\scrd_{r,n}$ and $\scrm_{r,n}$.
%
%
%
One might then ask whether the sets
$\scrh_{r,n}(M)$, $\scrd_{r,n}(M)$ and $\scrm_{r,n}(M)$ are convex
for fixed $M$.
But this too is false [at least for $\scrh_{r,n}(M)$ and $\scrd_{r,n}(M)$]:
in fact, a convex combination of polynomials in $\scrd_{r,n}(M)$
need not even lie in $\scrh_{r,n}(M)$.
To see this, let $P$ be the polynomial of
Example~\ref{sec_open}.\ref{exam_open.1},
and let $Q$ be the analogous polynomial with the indices $\{1,2\}$
replaced by $\{3,4\}$.
Then, for $n=4$, $P$ and $Q$ belong to $\scrd_{2,4}(U_{2,4})$
for $0 < \mu \le 4$;
but an easy calculation using Theorem~\ref{thm.rank2}
shows that $(P+Q)/2$ has the half-plane property only for $0 \le \mu \le 3$.

Finally, analogous examples show that the sets
$\scrh_{r,n}(M)$ and $\scrd_{r,n}(M)$
need not be log-convex or harmonic-mean-convex.
For example, let $P$ be the polynomial of
Example~\ref{sec_open}.\ref{exam_open.1},
and let $Q$ be the analogous polynomial with the indices $\{1,2\}$
replaced by $\{1,3\}$.
Then, for $n=4$, $P$ and $Q$ belong to $\scrd_{2,4}(U_{2,4})$
for $0 < \mu \le 4$;
but a straightforward (though messy) calculation using Theorem~\ref{thm.rank2}
shows that the polynomial $R$
obtained from $P$ and $Q$ by coefficientwise geometric mean
(resp.\ coefficientwise harmonic mean\footnote{
   The harmonic mean of $a$ and $b$ is $2ab/(a+b)$.
})
has the half-plane property only for $\mu \ge 1/16$ (resp.\ $\mu \ge 1/7$).

\subsection{Half-plane property for transversal matroids}

As noted in Section~\ref{sec_wild_transversal},
we have proven the half-plane property only for a subclass of
transversal matroids (the ``nice'' ones),
but we do not know a single example of a transversal matroid
that fails the half-plane property.
Our numerical experiments on rank-3 transversal matroids
(Section~\ref{sec_numerical})
suggest the following conjecture:

\begin{conjecture}
All rank-3 transversal matroids have the half-plane property.
\end{conjecture}

\noindent
If this conjecture is indeed true,
its proof will very likely require new techniques,
which could potentially shed light also on other related problems.
More ambitiously, we can raise the following question:

\begin{question}
Might {\em all} transversal matroids (and hence all gammoids)
have the half-plane property?
\end{question}

\noindent
Either a proof or a counterexample would be of considerable interest.

\subsection{Half-plane property for 7-element rank-3 matroids}

\begin{table}[tp]
\hspace*{-9mm}
\begin{tabular}{|l|l|l|}
\hline
 \multicolumn{1}{|c}{\bf Have HPP} & 
 \multicolumn{1}{|c}{\bf Do not have HPP} & 
 \multicolumn{1}{|c|}{\bf HPP unknown}  \\
\hline\hline
 $Q_7$ (Example~\ref{sec_transversal}.\ref{exam.transversal.7}) &
   $F_7$ (Example~\ref{sec_counterexamples}.\ref{exam.counterexamples.1}) &
   $F_7^{-4}$ \\
 $S_7$ (Examples~\ref{sec_transversal}.\ref{exam.transversal.4}
        and \ref{sec_transversal}.\ref{exam.transversal.7})  &
   $F_7^-$ (Example~\ref{sec_counterexamples}.\ref{exam.counterexamples.2}) &
   $\scrw^3 + e$ \\
 $M(K_4)^+$ (Example~\ref{sec_constructions}.\ref{exam.nice_principal_MK4plus}
       and Corollary~\ref{cor.determinant})
   & $F_7^{--}$ (Example~\ref{sec_counterexamples}.\ref{exam.counterexamples.3})
   & $\scrw^{3+}$   \\
 $F_7^{-5}$ (Example~\ref{sec_transversal}.\ref{exam.transversal.12}) &
   $M(K_4)+e$ (Example~\ref{sec_counterexamples}.\ref{exam.counterexamples.4})
   & $P'_7$  \\
 $F_7^{-6}$ (Example~\ref{sec_transversal}.\ref{exam.transversal.4}) &
   $F_7^{-3}$ (Example~\ref{sec_counterexamples}.\ref{exam.counterexamples.7})
   &     \\
 $P_7$ (Corollary~\ref{cor.determinant})  &    &      \\
 $P''_7$ (Examples~\ref{sec_transversal}.\ref{exam.transversal.17}
          and \ref{sec_transversal}.\ref{exam.transversal.18})  &   &   \\
 $P'''_7$ (Examples~\ref{sec_transversal}.\ref{exam.transversal.4}
        and \ref{sec_transversal}.\ref{exam.transversal.7})  &    &    \\
 $U_{3,7}$ (Theorem~\ref{thm.uniform})    &    &     \\
\hline
\end{tabular}
\vspace{2mm}
\caption{
   The 7-element rank-3 3-connected matroids,
   divided according to whether or not they have the
   half-plane property (HPP).
}
 \label{table_HPP_rank3}
\end{table}

We have shown that all matroids of rank or corank at most 2 have the
half-plane property (Corollary~\ref{cor.rank2a}), as do all matroids
on a ground set of at most 6 elements (Proposition~\ref{prop.nle6}).
So the first nontrivial case arises with 7-element rank-3 matroids;
we would like to know which ones have, and which ones do not have,
the half-plane property.
In Table~\ref{table_HPP_rank3} we divide the
7-element rank-3 3-connected matroids (see Appendix~\ref{subsec_rank3_7el})
into three categories:
those we have proven to have the half-plane property,
those we have proven not to have the half-plane property,
and those for which we have no proof either way.
There are exactly four matroids in the latter category:
$F_7^{-4}$, $\scrw^3 + e$, $\scrw^{3+}$ and $P'_7$.
Our numerical experiments (Section~\ref{sec_numerical})
suggest that these latter four matroids probably {\em do}\/ have
the half-plane property;
but proving it may well require new techniques.

\subsection{Algorithms}

As noted in Sections~\ref{sec.algorithms}, \ref{sec_necessary.1} and
\ref{sec_matroidal.weakHPP},
both the half-plane property and the weak half-plane property
are algorithmically testable, using quantifier-elimination methods.
But, at least with existing algorithms and currently available
computer hardware, these computations do not seem to be feasible in practice
for any interesting matroids (e.g.\ 7-element rank-3 matroids).

\begin{problem}
Find algorithms for testing the half-plane property
(or, more ambitiously, the weak half-plane property)
that are feasible in practice.
\end{problem}

\noindent
More modestly, one can ask:

\begin{problem}
Find heuristic numerical methods for testing the half-plane property
that are ``more powerful'' than the method based on
Proposition~\ref{prop.generalrank}
(see Section~\ref{sec_numerical}).
\end{problem}

\section*{Acknowledgments}

We wish to thank Sankar Basu, Jason Brown, Jim Geelen, Russ Lyons,
Richard Stanley, Dirk Vertigan and Neil White
for helpful comments and suggestions.
We also wish to thank James Renegar and Adam Strzebonski
for correspondence concerning quantifier-elimination algorithms;
J\"urgen Bokowski, Andreas Dress and Bernd Sturmfels
for correspondence concerning determinants and Grassmann--Pl\"ucker identities;
and Sankar Basu, Charles Desoer and Paul Penfield
for giving us pointers to the engineering literature.
Finally, we wish to thank Marc Noy and Dominic Welsh for organizing the
Workshop on Tutte Polynomials and Related Topics (Barcelona, September 2001),
at which this work was first presented.

This research was supported in part by
NSA grants MDA904--99--1--0030 and MDA904--01--1--0026 (J.G.O.),
NSF grants PHY--9900769 and PHY--0099393 (A.D.S.),
and an operating grant from NSERC (D.G.W.).

\appendix
\section{Matroids considered in this paper}  \label{app_matroids}

The half-plane property is preserved by
direct sums (Proposition~\ref{prop.product})
and 2-sums (Corollary~\ref{cor.2sum});
moreover, it is trivially preserved by adjoining loops or parallel elements.
Therefore, we can restrict attention to simple 3-connected matroids.
(In fact, every 3-connected matroid with at least 4 elements
 is automatically simple \cite[Proposition 8.1.6]{Oxley_92}.)
Since all rank-1 and rank-2 matroids have the half-plane property
(Corollary~\ref{cor.rank2a}), we can restrict attention to
matroids of rank $\ge 3$.
Finally, since the half-plane property is invariant under duality
(Proposition~\ref{prop.duality}), we can restrict attention to
matroids of rank $\le \lfloor n/2 \rfloor$,
where $n$ is the number of elements.

In this appendix we list all 3-connected rank-3 matroids
on 6 or 7 elements, as well as a few larger matroids
that will play a role in this paper.
We follow where possible the notation of \cite[Appendix]{Oxley_92};
further information on many of these matroids can be found there.

\subsection{Rank-3 matroids on 6 elements}  \label{subsec_rank3_6el}

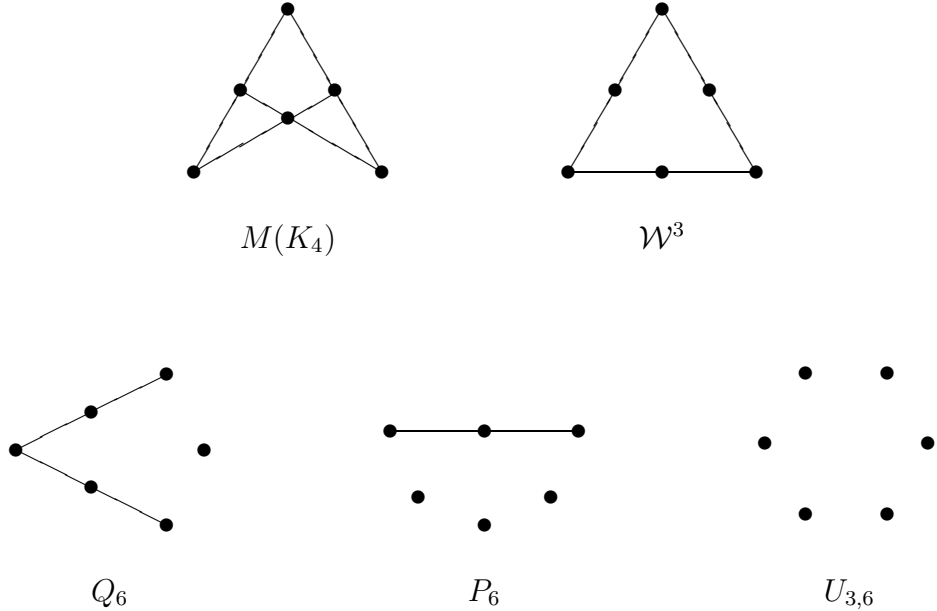
\begin{figure}
\setlength{\unitlength}{2.5cm}
\begin{center}
\begin{tabular}{c@{\qquad\qquad\qquad}c}
\begin{picture}(1,1)(0,0)
   \put(0,0){\circle*{0.07}}
   \put(1,0){\circle*{0.07}}
   \put(0.5,0.866){\circle*{0.07}}
   \put(0.5,0.289){\circle*{0.07}}
   \put(0.25,0.433){\circle*{0.07}}
   \put(0.75,0.433){\circle*{0.07}}
   \begin{drawjoin}
      \jput(0,0){\ }
      \jput(0.5,0.866){\ }
   \end{drawjoin}
   \begin{drawjoin}
      \jput(0.5,0.866){\ }
      \jput(1,0){\ }
   \end{drawjoin}
   \begin{drawjoin}
      \jput(0,0){\ }
      \jput(0.75,0.417){\ }
   \end{drawjoin}
   \begin{drawjoin}
      \jput(0.25,0.433){\ }
      \jput(1,0){\ }
   \end{drawjoin}
\end{picture}
&
\begin{picture}(1,1)(0,0)
   \put(0,0){\circle*{0.07}}
   \put(0.5,0){\circle*{0.07}}
   \put(1,0){\circle*{0.07}}
   \put(0.5,0.866){\circle*{0.07}}
   \put(0.25,0.433){\circle*{0.07}}
   \put(0.75,0.433){\circle*{0.07}}
   \put(0,0){\line(1,0){1}}
   \begin{drawjoin}
      \jput(0,0){\ }
      \jput(0.5,0.866){\ }
   \end{drawjoin}
   \begin{drawjoin}
      \jput(0.5,0.866){\ }
      \jput(1,0){\ }
   \end{drawjoin}
\end{picture}
\\[5mm]
$M(K_4)$ & $\scrw^3$ \\[10mm]
\end{tabular}
\begin{tabular}{c@{\qquad\qquad\qquad}c@{\qquad\qquad\qquad}c}
\begin{picture}(1,1)(0,0)
   \put(0,0.4){\circle*{0.07}}
   \put(0.4,0.2){\circle*{0.07}}
   \put(0.4,0.6){\circle*{0.07}}
   \put(0.8,0){\circle*{0.07}}
   \put(0.8,0.8){\circle*{0.07}}
   \put(1,0.4){\circle*{0.07}}
   \begin{drawjoin}
      \jput(0,0.4){\ }
      \jput(0.8,0.8){\ }
   \end{drawjoin}
   \begin{drawjoin}
      \jput(0,0.4){\ }
      \jput(0.8,0){\ }
   \end{drawjoin}
\end{picture}
&
\begin{picture}(1,1)(0,0)
   \put(0,0.5){\circle*{0.07}}
   \put(0.5,0.5){\circle*{0.07}}
   \put(1,0.5){\circle*{0.07}}
   \put(0,0.5){\line(1,0){1}}
   \put(0.5,0){\circle*{0.07}}
   \put(0.147,0.147){\circle*{0.07}}
   \put(0.853,0.147){\circle*{0.07}}
\end{picture}
&
%
\begin{picture}(0.9,1)(0,0)
   \put(0,0.433){\circle*{0.07}}
   \put(0.216,0.808){\circle*{0.07}}
   \put(0.65,0.808){\circle*{0.07}}
   \put(0.866,0.433){\circle*{0.07}}
   \put(0.65,0.058){\circle*{0.07}}
   \put(0.216,0.058){\circle*{0.07}}
\end{picture}
\\[5mm]
$Q_6$ & $P_6$ & $U_{3,6}$ \\[10mm]
\end{tabular}
\end{center}
\vspace{-3mm}
\caption{
   The 3-connected rank-3 matroids on 6 elements.
   All of these are obtained from $M(K_4)$ by a sequence of relaxations.
}
 \label{fig.rank3.le6}
\end{figure}

The 3-connected rank-3 matroids on 6 elements
are shown in Figure~\ref{fig.rank3.le6}.

\begin{quote}
\begin{itemize}
   \item[$M(K_4)$:]  Regular (in fact graphic and cographic).
       Self-dual. 
       Not transversal or cotransversal.
       Has half-plane property (Theorem~\ref{thm1.1}).
   \item[$\scrw^3$:]
       $F$-representable if and only if $|F| \ge 3$.
       Sixth-root-of-unity but not regular.
       Self-dual. 
       Transversal but not nice;
       cotransversal but not co-nice
       (Example~\ref{sec_transversal}.\ref{exam.transversal.3}).
       Has half-plane property (Corollary~\ref{cor.determinant}).
   \item[$Q_6$:]
       $F$-representable if and only if $|F| \ge 4$.
       Not sixth-root-of-unity.
       Self-dual. 
       Nice transversal and co-nice cotransversal
       (Example~\ref{sec_transversal}.\ref{exam.transversal.7}).
       Has half-plane property (Corollary~\ref{cor.nice_transversal}).
   \item[$P_6$:]
       $F$-representable if and only if $|F| \ge 5$.
       Not sixth-root-of-unity.
       Self-dual. 
       Nice transversal and co-nice cotransversal
       (Examples~\ref{sec_transversal}.\ref{exam.transversal.4}
        and \ref{sec_transversal}.\ref{exam.transversal.7}).
       Has half-plane property (Corollary~\ref{cor.nice_transversal}).
   \item[$U_{3,6}$:]
       $F$-representable if and only if $|F| \ge 4$.
       Not sixth-root-of-unity.
       Self-dual. 
       Nice transversal and co-nice cotransversal
       (Example~\ref{sec_transversal}.\ref{exam.transversal.1}).
       Has half-plane property (Theorem~\ref{thm.uniform}).
\end{itemize}
\end{quote}

\subsection{Rank-3 matroids on 7 elements}  \label{subsec_rank3_7el}

\begin{figure}
\setlength{\unitlength}{2.5cm}
\begin{center}
\begin{tabular}{c@{\qquad\qquad\qquad}c}
\begin{picture}(1,1)(0,0)
   \put(0,0.4){\circle*{0.07}}
   \put(0.266,0.533){\circle*{0.07}}
   \put(0.533,0.666){\circle*{0.07}}
   \put(0.266,0.266){\circle*{0.07}}
   \put(0.533,0.133){\circle*{0.07}}
   \put(0.8,0){\circle*{0.07}}
   \put(1,0.4){\circle*{0.07}}
   \begin{drawjoin}
      \jput(0,0.4){\ }
      \jput(0.533,0.666){\ }
   \end{drawjoin}
   \begin{drawjoin}
      \jput(0,0.4){\ }
      \jput(0.8,0){\ }
   \end{drawjoin}
   \put(-0.15,0.35){\footnotesize 1}
   \put(0.235,0.60){\footnotesize 2}
   \put(0.50,0.73){\footnotesize 3}
   \put(0.235,0.1){\footnotesize 4}
   \put(0.50,-0.03){\footnotesize 5}
   \put(0.77,-0.16){\footnotesize 6}
   \put(1.08,0.35){\footnotesize 7}
\end{picture}
&
\begin{picture}(1,1)(0,0)
   \put(0,0.5){\circle*{0.07}}
   \put(0.333,0.5){\circle*{0.07}}
   \put(0.666,0.5){\circle*{0.07}}
   \put(1,0.5){\circle*{0.07}}
   \put(0,0.5){\line(1,0){1}}
   \put(0.5,0){\circle*{0.07}}
   \put(0.147,0.147){\circle*{0.07}}
   \put(0.853,0.147){\circle*{0.07}}
   \put(-0.03,0.57){\footnotesize 1}
   \put(0.30,0.57){\footnotesize 2}
   \put(0.635,0.57){\footnotesize 3}
   \put(0.97,0.57){\footnotesize 4}
   \put(0.47,-0.16){\footnotesize 6}
   \put(0.12,-0.02){\footnotesize 5}
   \put(0.82,-0.02){\footnotesize 7}
\end{picture}
\\[5mm]
$Q_7$ & $S_7$ \\[10mm]
\end{tabular}
\begin{tabular}{c@{\qquad\qquad\qquad}c}
\begin{picture}(1,1)(0,0)
   \put(0,0){\circle*{0.07}}
   \put(1,0){\circle*{0.07}}
   \put(0.5,0.866){\circle*{0.07}}
   \put(0.5,0.289){\circle*{0.07}}
   \put(0.125,0.217){\circle*{0.07}}
   \put(0.25,0.433){\circle*{0.07}}
   \put(0.75,0.433){\circle*{0.07}}
   \begin{drawjoin}
      \jput(0,0){\ }
      \jput(0.5,0.866){\ }
   \end{drawjoin}
   \begin{drawjoin}
      \jput(0.5,0.866){\ }
      \jput(1,0){\ }
   \end{drawjoin}
   \begin{drawjoin}
      \jput(0,0){\ }
      \jput(0.75,0.417){\ }
   \end{drawjoin}
   \begin{drawjoin}
      \jput(0.25,0.433){\ }
      \jput(1,0){\ }
   \end{drawjoin}
   \put(-0.15,-0.03){\footnotesize 3}
   \put(-0.025,0.21){\footnotesize 7}
   \put(1.07,-0.03){\footnotesize 5}
   \put(0.47,0.95){\footnotesize 1}
   \put(0.465,0.12){\footnotesize 4}
   \put(0.1,0.45){\footnotesize 2}
   \put(0.8,0.45){\footnotesize 6}
\end{picture}
&
\begin{picture}(1,1)(0,0)
   \put(0,0){\circle*{0.07}}
   \put(0.5,0){\circle*{0.07}}
   \put(1,0){\circle*{0.07}}
   \put(0.5,0.866){\circle*{0.07}}
   \put(0.166,0.289){\circle*{0.07}}
   \put(0.333,0.577){\circle*{0.07}}
   \put(0.75,0.433){\circle*{0.07}}
   \put(0,0){\line(1,0){1}}
   \begin{drawjoin}
      \jput(0,0){\ }
      \jput(0.5,0.866){\ }
   \end{drawjoin}
   \begin{drawjoin}
      \jput(0.5,0.866){\ }
      \jput(1,0){\ }
   \end{drawjoin}
   \put(-0.15,-0.03){\footnotesize 3}
   \put(0.02,0.32){\footnotesize 7}
   \put(1.07,-0.03){\footnotesize 5}
   \put(0.47,0.95){\footnotesize 1}
   \put(0.465,-0.15){\footnotesize 4}
   \put(0.18,0.60){\footnotesize 2}
   \put(0.8,0.45){\footnotesize 6}
\end{picture}
\\[5mm]
$M(K_4)^+$ & $\scrw^{3+}$
\end{tabular}
\end{center}
\vspace{-3mm}
\caption{
   The 3-connected rank-3 matroids on 7 elements with a 4-point line.
}
 \label{fig.rank3.seven.4point}
\end{figure}
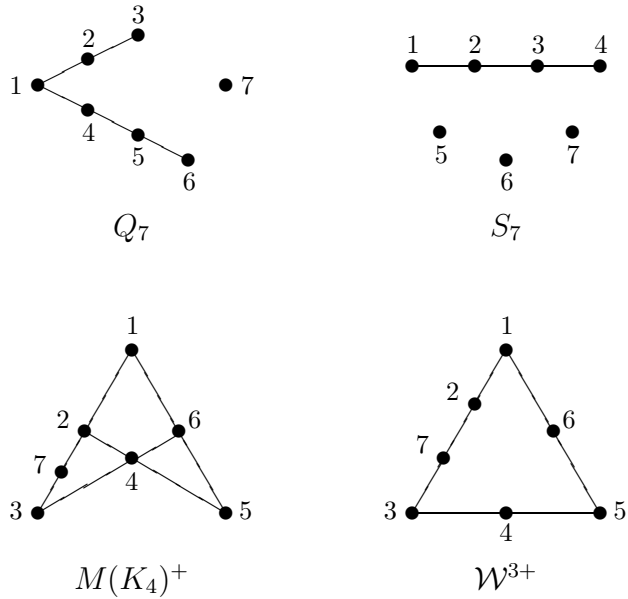

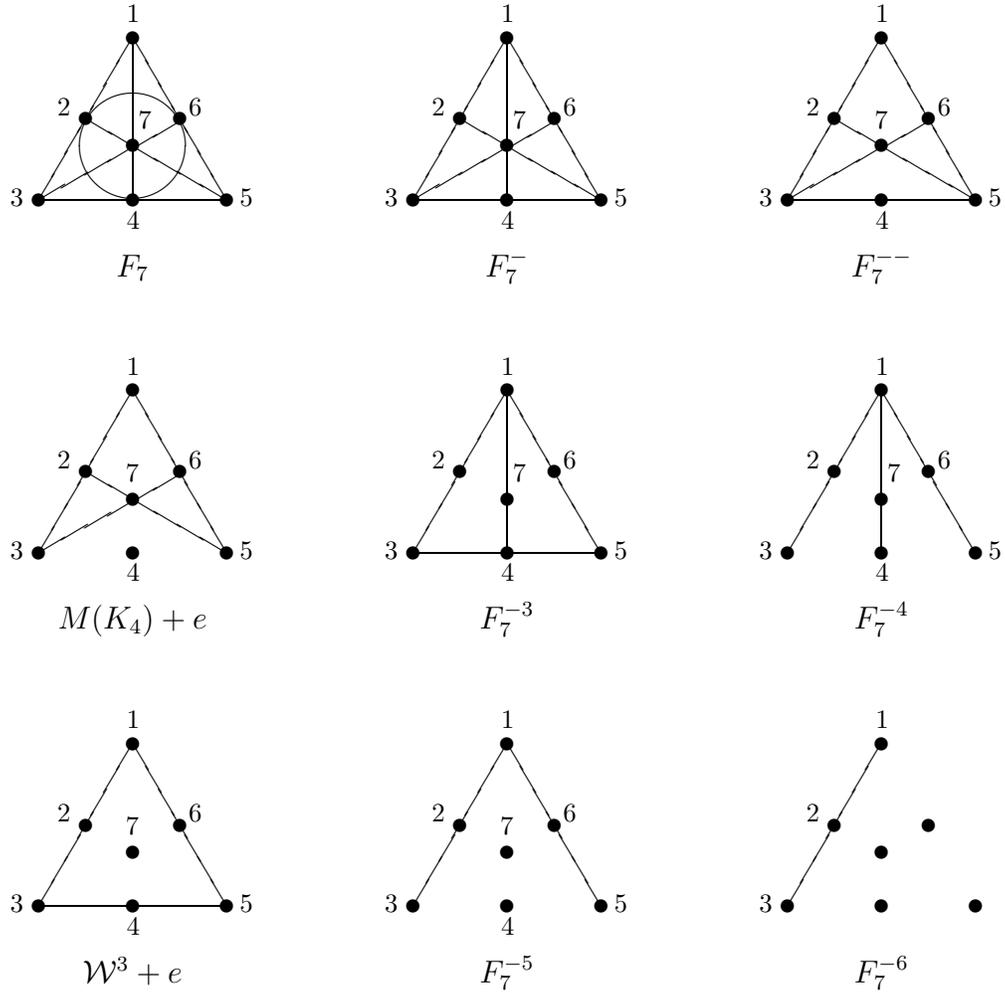
\begin{figure}
\setlength{\unitlength}{2.5cm}
\begin{center}
\begin{tabular}{c@{\qquad\qquad\qquad}c@{\qquad\qquad\qquad}c}
\begin{picture}(1,1)(0,0)
   \put(0,0){\circle*{0.07}}
   \put(-0.15,-0.03){\footnotesize 3}
   \put(0.5,0){\circle*{0.07}}
   \put(0.47,-0.15){\footnotesize 4}
   \put(1,0){\circle*{0.07}}
   \put(1.07,-0.03){\footnotesize 5}
   \put(0.5,0.866){\circle*{0.07}}
   \put(0.47,0.95){\footnotesize 1}
   \put(0.5,0.289){\circle*{0.07}}
   \put(0.53,0.38){\footnotesize 7}
   \put(0.25,0.433){\circle*{0.07}}
   \put(0.1,0.45){\footnotesize 2}
   \put(0.75,0.433){\circle*{0.07}}
   \put(0.8,0.45){\footnotesize 6}
   \put(0,0){\line(1,0){1}}
   \put(0.5,0){\line(0,1){0.866}}
   \begin{drawjoin}
      \jput(0,0){\ }
      \jput(0.5,0.866){\ }
   \end{drawjoin}
   \begin{drawjoin}
      \jput(0.5,0.866){\ }
      \jput(1,0){\ }
   \end{drawjoin}
   \begin{drawjoin}
      \jput(0,0){\ }
      \jput(0.75,0.417){\ }
   \end{drawjoin}
   \begin{drawjoin}
      \jput(0.25,0.433){\ }
      \jput(1,0){\ }
   \end{drawjoin}
   \put(0.5,0.289){\circle{0.578}}
\end{picture}
&
\begin{picture}(1,1)(0,0)
   \put(0,0){\circle*{0.07}}
   \put(-0.15,-0.03){\footnotesize 3}
   \put(0.5,0){\circle*{0.07}}
   \put(0.47,-0.15){\footnotesize 4}
   \put(1,0){\circle*{0.07}}
   \put(1.07,-0.03){\footnotesize 5}
   \put(0.5,0.866){\circle*{0.07}}
   \put(0.47,0.95){\footnotesize 1}
   \put(0.5,0.289){\circle*{0.07}}
   \put(0.53,0.38){\footnotesize 7}
   \put(0.25,0.433){\circle*{0.07}}
   \put(0.1,0.45){\footnotesize 2}
   \put(0.75,0.433){\circle*{0.07}}
   \put(0.8,0.45){\footnotesize 6}
   \put(0,0){\line(1,0){1}}
   \put(0.5,0){\line(0,1){0.866}}
   \begin{drawjoin}
      \jput(0,0){\ }
      \jput(0.5,0.866){\ }
   \end{drawjoin}
   \begin{drawjoin}
      \jput(0.5,0.866){\ }
      \jput(1,0){\ }
   \end{drawjoin}
   \begin{drawjoin}
      \jput(0,0){\ }
      \jput(0.75,0.417){\ }
   \end{drawjoin}
   \begin{drawjoin}
      \jput(0.25,0.433){\ }
      \jput(1,0){\ }
   \end{drawjoin}
\end{picture}
&
\begin{picture}(1,1)(0,0)
   \put(0,0){\circle*{0.07}}
   \put(-0.15,-0.03){\footnotesize 3}
   \put(0.5,0){\circle*{0.07}}
   \put(0.47,-0.15){\footnotesize 4}
   \put(1,0){\circle*{0.07}}
   \put(1.07,-0.03){\footnotesize 5}
   \put(0.5,0.866){\circle*{0.07}}
   \put(0.47,0.95){\footnotesize 1}
   \put(0.5,0.289){\circle*{0.07}}
   \put(0.465,0.38){\footnotesize 7}
   \put(0.25,0.433){\circle*{0.07}}
   \put(0.1,0.45){\footnotesize 2}
   \put(0.75,0.433){\circle*{0.07}}
   \put(0.8,0.45){\footnotesize 6}
   \put(0,0){\line(1,0){1}}
   \begin{drawjoin}
      \jput(0,0){\ }
      \jput(0.5,0.866){\ }
   \end{drawjoin}
   \begin{drawjoin}
      \jput(0.5,0.866){\ }
      \jput(1,0){\ }
   \end{drawjoin}
   \begin{drawjoin}
      \jput(0,0){\ }
      \jput(0.75,0.417){\ }
   \end{drawjoin}
   \begin{drawjoin}
      \jput(0.25,0.433){\ }
      \jput(1,0){\ }
   \end{drawjoin}
\end{picture}
\\[5mm]
$F_7$ & $F_7^-$ & $F_7^{--}$ \\[10mm]
\end{tabular}
\begin{tabular}{c@{\qquad\qquad\qquad}c@{\qquad\qquad\qquad}c}
\begin{picture}(1,1)(0,0)
   \put(0,0){\circle*{0.07}}
   \put(-0.15,-0.03){\footnotesize 3}
   \put(0.5,0){\circle*{0.07}}
   \put(0.47,-0.15){\footnotesize 4}
   \put(1,0){\circle*{0.07}}
   \put(1.07,-0.03){\footnotesize 5}
   \put(0.5,0.866){\circle*{0.07}}
   \put(0.47,0.95){\footnotesize 1}
   \put(0.5,0.289){\circle*{0.07}}
   \put(0.465,0.38){\footnotesize 7}
   \put(0.25,0.433){\circle*{0.07}}
   \put(0.1,0.45){\footnotesize 2}
   \put(0.75,0.433){\circle*{0.07}}
   \put(0.8,0.45){\footnotesize 6}
   \begin{drawjoin}
      \jput(0,0){\ }
      \jput(0.5,0.866){\ }
   \end{drawjoin}
   \begin{drawjoin}
      \jput(0.5,0.866){\ }
      \jput(1,0){\ }
   \end{drawjoin}
   \begin{drawjoin}
      \jput(0,0){\ }
      \jput(0.75,0.417){\ }
   \end{drawjoin}
   \begin{drawjoin}
      \jput(0.25,0.433){\ }
      \jput(1,0){\ }
   \end{drawjoin}
\end{picture}
&
\begin{picture}(1,1)(0,0)
   \put(0,0){\circle*{0.07}}
   \put(-0.15,-0.03){\footnotesize 3}
   \put(0.5,0){\circle*{0.07}}
   \put(0.47,-0.15){\footnotesize 4}
   \put(1,0){\circle*{0.07}}
   \put(1.07,-0.03){\footnotesize 5}
   \put(0.5,0.866){\circle*{0.07}}
   \put(0.47,0.95){\footnotesize 1}
   \put(0.5,0.289){\circle*{0.07}}
   \put(0.53,0.38){\footnotesize 7}
   \put(0.25,0.433){\circle*{0.07}}
   \put(0.1,0.45){\footnotesize 2}
   \put(0.75,0.433){\circle*{0.07}}
   \put(0.8,0.45){\footnotesize 6}
   \put(0,0){\line(1,0){1}}
   \put(0.5,0){\line(0,1){0.866}}
   \begin{drawjoin}
      \jput(0,0){\ }
      \jput(0.5,0.866){\ }
   \end{drawjoin}
   \begin{drawjoin}
      \jput(0.5,0.866){\ }
      \jput(1,0){\ }
   \end{drawjoin}
\end{picture}
&
\begin{picture}(1,1)(0,0)
   \put(0,0){\circle*{0.07}}
   \put(-0.15,-0.03){\footnotesize 3}
   \put(0.5,0){\circle*{0.07}}
   \put(0.47,-0.15){\footnotesize 4}
   \put(1,0){\circle*{0.07}}
   \put(1.07,-0.03){\footnotesize 5}
   \put(0.5,0.866){\circle*{0.07}}
   \put(0.47,0.95){\footnotesize 1}
   \put(0.5,0.289){\circle*{0.07}}
   \put(0.53,0.38){\footnotesize 7}
   \put(0.25,0.433){\circle*{0.07}}
   \put(0.1,0.45){\footnotesize 2}
   \put(0.75,0.433){\circle*{0.07}}
   \put(0.8,0.45){\footnotesize 6}
   \put(0.5,0){\line(0,1){0.866}}
   \begin{drawjoin}
      \jput(0,0){\ }
      \jput(0.5,0.866){\ }
   \end{drawjoin}
   \begin{drawjoin}
      \jput(0.5,0.866){\ }
      \jput(1,0){\ }
   \end{drawjoin}
\end{picture}
\\[5mm]
$M(K_4)+e$ & $F_7^{-3}$ & $F_7^{-4}$ \\[10mm]
\end{tabular}
\begin{tabular}{c@{\qquad\qquad\qquad}c@{\qquad\qquad\qquad}c}
\begin{picture}(1,1)(0,0)
   \put(0,0){\circle*{0.07}}
   \put(-0.15,-0.03){\footnotesize 3}
   \put(0.5,0){\circle*{0.07}}
   \put(0.47,-0.15){\footnotesize 4}
   \put(1,0){\circle*{0.07}}
   \put(1.07,-0.03){\footnotesize 5}
   \put(0.5,0.866){\circle*{0.07}}
   \put(0.47,0.95){\footnotesize 1}
   \put(0.5,0.289){\circle*{0.07}}
   \put(0.465,0.38){\footnotesize 7}
   \put(0.25,0.433){\circle*{0.07}}
   \put(0.1,0.45){\footnotesize 2}
   \put(0.75,0.433){\circle*{0.07}}
   \put(0.8,0.45){\footnotesize 6}
   \put(0,0){\line(1,0){1}}
   \begin{drawjoin}
      \jput(0,0){\ }
      \jput(0.5,0.866){\ }
   \end{drawjoin}
   \begin{drawjoin}
      \jput(0.5,0.866){\ }
      \jput(1,0){\ }
   \end{drawjoin}
\end{picture}
&
\begin{picture}(1,1)(0,0)
   \put(0,0){\circle*{0.07}}
   \put(-0.15,-0.03){\footnotesize 3}
   \put(0.5,0){\circle*{0.07}}
   \put(0.47,-0.15){\footnotesize 4}
   \put(1,0){\circle*{0.07}}
   \put(1.07,-0.03){\footnotesize 5}
   \put(0.5,0.866){\circle*{0.07}}
   \put(0.47,0.95){\footnotesize 1}
   \put(0.5,0.289){\circle*{0.07}}
   \put(0.465,0.38){\footnotesize 7}
   \put(0.25,0.433){\circle*{0.07}}
   \put(0.1,0.45){\footnotesize 2}
   \put(0.75,0.433){\circle*{0.07}}
   \put(0.8,0.45){\footnotesize 6}
   \begin{drawjoin}
      \jput(0,0){\ }
      \jput(0.5,0.866){\ }
   \end{drawjoin}
   \begin{drawjoin}
      \jput(0.5,0.866){\ }
      \jput(1,0){\ }
   \end{drawjoin}
\end{picture}
&
\begin{picture}(1,1)(0,0)
   \put(0,0){\circle*{0.07}}
   \put(-0.15,-0.03){\footnotesize 3}
   \put(0.5,0){\circle*{0.07}}
   \put(1,0){\circle*{0.07}}
   \put(0.5,0.866){\circle*{0.07}}
   \put(0.47,0.95){\footnotesize 1}
   \put(0.5,0.289){\circle*{0.07}}
   \put(0.25,0.433){\circle*{0.07}}
   \put(0.1,0.45){\footnotesize 2}
   \put(0.75,0.433){\circle*{0.07}}
   \begin{drawjoin}
      \jput(0,0){\ }
      \jput(0.5,0.866){\ }
   \end{drawjoin}
\end{picture}
\\[5mm]
$\scrw^3 + e$ & $F_7^{-5}$ & $F_7^{-6}$
\end{tabular}
\end{center}
\vspace{-3mm}
\caption{
   The 3-connected rank-3 non-uniform matroids on $7$ elements
   that have no 4-point line and are not $\omega$-regular.
   All of these are obtained from the Fano matroid $F_7$
   by a sequence of relaxations.
}
 \label{fig_F7etc}
\end{figure}

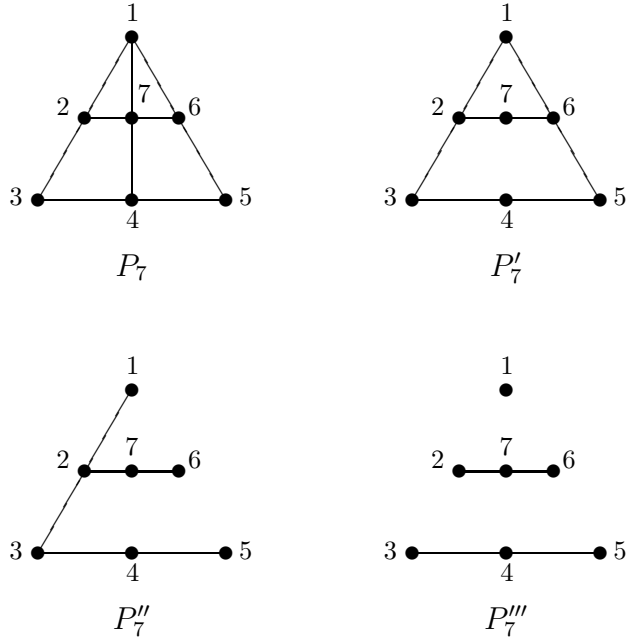
\begin{figure}
\setlength{\unitlength}{2.5cm}
\begin{center}
\begin{tabular}{c@{\qquad\qquad\qquad}c}
\begin{picture}(1,1)(0,0)
   \put(0,0){\circle*{0.07}}
   \put(-0.15,-0.03){\footnotesize 3}
   \put(0.5,0){\circle*{0.07}}
   \put(0.47,-0.15){\footnotesize 4}
   \put(1,0){\circle*{0.07}}
   \put(1.07,-0.03){\footnotesize 5}
   \put(0.5,0.866){\circle*{0.07}}
   \put(0.47,0.95){\footnotesize 1}
   \put(0.5,0.433){\circle*{0.07}}
   \put(0.53,0.52){\footnotesize 7}
   \put(0.25,0.433){\circle*{0.07}}
   \put(0.1,0.45){\footnotesize 2}
   \put(0.75,0.433){\circle*{0.07}}
   \put(0.8,0.45){\footnotesize 6}
   \put(0,0){\line(1,0){1}}
   \put(0.5,0){\line(0,1){0.866}}
   \put(0.25,0.433){\line(1,0){0.5}}
   \begin{drawjoin}
      \jput(0,0){\ }
      \jput(0.5,0.866){\ }
   \end{drawjoin}
   \begin{drawjoin}
      \jput(0.5,0.866){\ }
      \jput(1,0){\ }
   \end{drawjoin}
\end{picture}
&
\begin{picture}(1,1)(0,0)
   \put(0,0){\circle*{0.07}}
   \put(-0.15,-0.03){\footnotesize 3}
   \put(0.5,0){\circle*{0.07}}
   \put(0.47,-0.15){\footnotesize 4}
   \put(1,0){\circle*{0.07}}
   \put(1.07,-0.03){\footnotesize 5}
   \put(0.5,0.866){\circle*{0.07}}
   \put(0.47,0.95){\footnotesize 1}
   \put(0.5,0.433){\circle*{0.07}}
   \put(0.465,0.52){\footnotesize 7}
   \put(0.25,0.433){\circle*{0.07}}
   \put(0.1,0.45){\footnotesize 2}
   \put(0.75,0.433){\circle*{0.07}}
   \put(0.8,0.45){\footnotesize 6}
   \put(0,0){\line(1,0){1}}
   \put(0.25,0.433){\line(1,0){0.5}}
   \begin{drawjoin}
      \jput(0,0){\ }
      \jput(0.5,0.866){\ }
   \end{drawjoin}
   \begin{drawjoin}
      \jput(0.5,0.866){\ }
      \jput(1,0){\ }
   \end{drawjoin}
\end{picture}
\\[5mm]
$P_7$ & $P'_7$ \\[10mm]
\end{tabular}
\begin{tabular}{c@{\qquad\qquad\qquad}c}
\begin{picture}(1,1)(0,0)
   \put(0,0){\circle*{0.07}}
   \put(-0.15,-0.03){\footnotesize 3}
   \put(0.5,0){\circle*{0.07}}
   \put(0.47,-0.15){\footnotesize 4}
   \put(1,0){\circle*{0.07}}
   \put(1.07,-0.03){\footnotesize 5}
   \put(0.5,0.866){\circle*{0.07}}
   \put(0.47,0.95){\footnotesize 1}
   \put(0.5,0.433){\circle*{0.07}}
   \put(0.465,0.52){\footnotesize 7}
   \put(0.25,0.433){\circle*{0.07}}
   \put(0.1,0.45){\footnotesize 2}
   \put(0.75,0.433){\circle*{0.07}}
   \put(0.8,0.45){\footnotesize 6}
   \put(0,0){\line(1,0){1}}
   \put(0.25,0.433){\line(1,0){0.5}}
   \begin{drawjoin}
      \jput(0,0){\ }
      \jput(0.5,0.866){\ }
   \end{drawjoin}
\end{picture}
&
\begin{picture}(1,1)(0,0)
   \put(0,0){\circle*{0.07}}
   \put(-0.15,-0.03){\footnotesize 3}
   \put(0.5,0){\circle*{0.07}}
   \put(0.47,-0.15){\footnotesize 4}
   \put(1,0){\circle*{0.07}}
   \put(1.07,-0.03){\footnotesize 5}
   \put(0.5,0.866){\circle*{0.07}}
   \put(0.47,0.95){\footnotesize 1}
   \put(0.5,0.433){\circle*{0.07}}
   \put(0.465,0.52){\footnotesize 7}
   \put(0.25,0.433){\circle*{0.07}}
   \put(0.1,0.45){\footnotesize 2}
   \put(0.75,0.433){\circle*{0.07}}
   \put(0.8,0.45){\footnotesize 6}
   \put(0,0){\line(1,0){1}}
   \put(0.25,0.433){\line(1,0){0.5}}
\end{picture}
\\[5mm]
$P''_7$ & $P'''_7$
\end{tabular}
\end{center}
\vspace{-3mm}
\caption{
   The 3-connected rank-3 non-uniform matroids on $7$ elements
   that have no 4-point line and are $\omega$-regular.
   All of these are obtained from $P_7$ by a sequence of relaxations.
}
 \label{fig_P7etc}
\end{figure}

Besides the uniform matroid $U_{3,7}$,
we classify the 3-connected rank-3 matroids on 7 elements as follows
\cite{Semple_98}:
\begin{itemize}
   \item[(a)] Those with a 4-point line (Figure~\ref{fig.rank3.seven.4point}).
   \item[(b)] Those with no 4-point line:
   \begin{itemize}
       \item[(b1)] Those that are not $\omega$-regular
            (Figure~\ref{fig_F7etc}).
            All of these are obtained from the Fano matroid $F_7$
            by a sequence of relaxations.
       \item[(b2)] Those that are $\omega$-regular
            (Figure~\ref{fig_P7etc}).
            All of these are obtained from $P_7$ by a sequence of relaxations.
   \end{itemize}
\end{itemize}

\subsubsection{Matroids with a 4-point line}

\begin{quote}
\begin{itemize}
   \item[$Q_7$:]
       $F$-representable if and only if $|F| \ge 5$.
       Not sixth-root-of-unity.
       Nice transversal
       (Example~\ref{sec_transversal}.\ref{exam.transversal.7});
       co-transversal but not co-nice
       (Example~\ref{sec_transversal}.\ref{exam.transversal.14}).
       Non-nice principal extension of $Q_6$;
       non-nice free extension of $Q_7 \setminus 7$
       (Example~\ref{sec_constructions}.\ref{exam.non-nice_principal_Q7}).
       Has half-plane property (Corollary~\ref{cor.nice_transversal}).
   \item[$S_7$:]
       $F$-representable if and only if $|F| \ge 7$.
       Not sixth-root-of-unity.
       Nice transversal
       (Example~\ref{sec_transversal}.\ref{exam.transversal.7});
       co-transversal
       (we don't know whether or not it is co-nice).
       Non-nice principal extension of $P_6$;
       non-nice free extension of $S_7 \setminus 7$.
       Has half-plane property (Corollary~\ref{cor.nice_transversal}).
   \item[$M(K_4)^+$:]
       $F$-representable if and only if $|F| \ge 3$.
       Sixth-root-of-unity but not regular.
       Not transversal or cotransversal.
       Nice principal extension of $M(K_4)$ with
       $\lambda_1 = \lambda_2 = \lambda_3 = 1/2$
       (Example~\ref{sec_constructions}.\ref{exam.nice_principal_MK4plus}).
       This matroid is called $O_7$ in \cite{Geelen_00a}.
       Has half-plane property
       (Example~\ref{sec_constructions}.\ref{exam.nice_principal_MK4plus}
       or Corollary~\ref{cor.determinant}).
   \item[$\scrw^{3+}$:]
       $F$-representable if and only if $|F| \ge 4$.
       Not sixth-root-of-unity.
       Transversal but not nice
       (Example~\ref{sec_transversal}.\ref{exam.transversal.3});
       co-transversal but not co-nice
       (Example~\ref{sec_transversal}.\ref{exam.transversal.13}).
       Non-nice principal extension of $\scrw^3$;
       non-nice principal extension of $Q_7 \setminus 7$.
       Not known whether it has half-plane property.
\end{itemize}
\end{quote}

\subsubsection{The Fano matroid $F_7$ and its relaxations}

\begin{quote}
\begin{itemize}
   \item[$F_7$:]
       $F$-representable if and only if $F$ has characteristic 2.
       Not sixth-root-of-unity.
       Not transversal or cotransversal.
       Doubly transitive automorphism group.
       Does not have half-plane property
       (Example~\ref{sec_counterexamples}.\ref{exam.counterexamples.1}).
   \item[$F_7^-$:]
       $F$-representable if and only if $F$ has characteristic $\neq 2$.
       Not sixth-root-of-unity.
       Not transversal or cotransversal.
       Does not have half-plane property
       (Example~\ref{sec_counterexamples}.\ref{exam.counterexamples.2}).
   \item[$F_7^{--}$:]
       $F$-representable if and only if $|F| \ge 4$
       (Example~\ref{sec_counterexamples}.\ref{exam.counterexamples.6}).
       Not sixth-root-of-unity.
       Not transversal or cotransversal.
       Non-nice principal extension of $M(K_4)$.
       Does not have half-plane property
       (Example~\ref{sec_counterexamples}.\ref{exam.counterexamples.3}).
   \item[$M(K_4)+e$:]
       $F$-representable if and only if $|F| \ge 5$.
       Not sixth-root-of-unity.
       Not transversal or cotransversal.
       Non-nice free extension of $M(K_4)$
       [Example~\ref{sec_constructions}.\ref{exam.non-nice_principal_MK4pluse}].
       Does not have half-plane property
       (Example~\ref{sec_counterexamples}.\ref{exam.counterexamples.4}).
   \item[$F_7^{-3}$:]
       $F$-representable if and only if $|F| \ge 5$.
       Not sixth-root-of-unity.
       Not transversal;
       cotransversal but not co-nice
       (Example~\ref{sec_transversal}.\ref{exam.transversal.9}).
       Non-nice principal extension of $\scrw^3$.
       Does not have half-plane property
       (Example~\ref{sec_counterexamples}.\ref{exam.counterexamples.7}).
   \item[$F_7^{-4}$:]
       $F$-representable if and only if $|F| \ge 4$.
       Not sixth-root-of-unity.
       Not transversal;  cotransversal but not co-nice
       (Example~\ref{sec_transversal}.\ref{exam.transversal.10}).
       Non-nice principal extension of $Q_6$.
       Not known whether it has half-plane property.
   \item[$\scrw^3 + e$:]
       $F$-representable if and only if $|F| \ge 5$.
       Not sixth-root-of-unity.
       Transversal but not nice
       (Example~\ref{sec_transversal}.\ref{exam.transversal.3});
       cotransversal but not co-nice
       (Example~\ref{sec_transversal}.\ref{exam.transversal.11}).
       Non-nice free extension of $\scrw^3$;
       non-nice principal extension of $Q_6$.
       Not known whether it has half-plane property.
   \item[$F_7^{-5}$:]
       $F$-representable if and only if $|F| \ge 5$.
       Not sixth-root-of-unity.
       Transversal but not nice
       (Example~\ref{sec_transversal}.\ref{exam.transversal.8});
       co-nice cotransversal
       (Example~\ref{sec_transversal}.\ref{exam.transversal.12}).
       Non-nice free extension of $Q_6$;
       non-nice principal extension of $P_6$.
       Has half-plane property (Corollary~\ref{cor.nice_transversal}).
   \item[$F_7^{-6}$:]
       $F$-representable if and only if $|F| \ge 7$.
       Not sixth-root-of-unity.
       Nice transversal
       (Example~\ref{sec_transversal}.\ref{exam.transversal.4});
       co-transversal
       (we don't know whether or not it is co-nice).
       Non-nice free extension of $P_6$;
       non-nice principal extension of $U_{3,6}$.
       Has half-plane property (Corollary~\ref{cor.nice_transversal}).
\end{itemize}
\end{quote}

\subsubsection{$P_7$ and its relaxations}

\begin{quote}
\begin{itemize}
   \item[$P_7$:]
       $F$-representable if and only if $|F| \ge 3$.
       Sixth-root-of-unity but not regular.
       $k$-regular for all $k \ge 1$ (see \cite{Semple_99}).
       Not transversal or cotransversal.
       Has half-plane property (Corollary~\ref{cor.determinant}).
   \item[$P'_7$:]
       $F$-representable if and only if $|F| \ge 4$.
       $k$-regular for all $k \ge 2$ (see \cite{Semple_99}).
       Not transversal; cotransversal but not co-nice
       (Example~\ref{sec_transversal}.\ref{exam.transversal.16}).
       Non-nice principal extension of $\scrw^3$.
       Not known whether it has half-plane property.
   \item[$P''_7$:]
       $F$-representable if and only if $|F| \ge 5$.
       $k$-regular for all $k \ge 3$ (see \cite{Semple_99}).
       Transversal but not nice
       (Example~\ref{sec_transversal}.\ref{exam.transversal.15});
       co-nice cotransversal
       (Examples~\ref{sec_transversal}.\ref{exam.transversal.17}
        and \ref{sec_transversal}.\ref{exam.transversal.18}).
       Non-nice principal extension of $Q_6$;
       non-nice principal extension of $R_6$.
       Has half-plane property (Corollary~\ref{cor.nice_transversal}).
   \item[$P'''_7$:]
       $F$-representable if and only if $|F| \ge 7$.
       $k$-regular for all $k \ge 4$ (see \cite{Semple_99}).
       Nice transversal
       (Examples~\ref{sec_transversal}.\ref{exam.transversal.4}
       and \ref{sec_transversal}.\ref{exam.transversal.7});
       co-transversal
       (we don't know whether or not it is co-nice).
       Non-nice free extension of $R_6$;
       non-nice principal extension of $P_6$.
       Has half-plane property (Corollary~\ref{cor.nice_transversal}).
\end{itemize}
\end{quote}

\subsection{Some rank-3 matroids on 8 or 9 elements}

The matroids $F_7^{-4}$, $\scrw^3 + e$, $\scrw^{3+}$ and $P'_7$
play a special role in this paper,
as they are the only rank-3 7-element matroids
for which we are unable to prove whether or not they have
the half-plane property.
Since our numerical experiments (Section~\ref{sec_numerical})
suggest that they probably {\em do}\/ have the half-plane property,
we have also investigated some single-element extensions
of these matroids in an (unsuccessful) effort to find one
that fails the half-plane property.
For brevity we discuss here only the four matroids
that are obtained as {\em free}\/ extensions
(see Figure~\ref{fig_F7m4pe_etc}):

\begin{figure}
\setlength{\unitlength}{2.5cm}
\begin{center}
\begin{tabular}{c@{\qquad\qquad\qquad}c}
\begin{picture}(1,1)(0,0)
   \put(0,0){\circle*{0.07}}
   \put(-0.15,-0.03){\footnotesize 3}
   \put(0.5,0){\circle*{0.07}}
   \put(0.47,-0.15){\footnotesize 4}
   \put(1,0){\circle*{0.07}}
   \put(1.07,-0.03){\footnotesize 5}
   \put(0.5,0.866){\circle*{0.07}}
   \put(0.47,0.95){\footnotesize 1}
   \put(0.5,0.289){\circle*{0.07}}
   \put(0.53,0.38){\footnotesize 7}
   \put(0.25,0.433){\circle*{0.07}}
   \put(0.1,0.45){\footnotesize 2}
   \put(0.75,0.433){\circle*{0.07}}
   \put(0.8,0.45){\footnotesize 6}
   \put(0.5,0){\line(0,1){0.866}}
   \put(1.05,0.433){\circle*{0.07}}
   \put(1.02,0.51){\footnotesize 8}
   \begin{drawjoin}
      \jput(0,0){\ }
      \jput(0.5,0.866){\ }
   \end{drawjoin}
   \begin{drawjoin}
      \jput(0.5,0.866){\ }
      \jput(1,0){\ }
   \end{drawjoin}
\end{picture}
&
\begin{picture}(1,1)(0,0)
   \put(0,0){\circle*{0.07}}
   \put(-0.15,-0.03){\footnotesize 3}
   \put(0.5,0){\circle*{0.07}}
   \put(0.47,-0.15){\footnotesize 4}
   \put(1,0){\circle*{0.07}}
   \put(1.07,-0.03){\footnotesize 5}
   \put(0.5,0.866){\circle*{0.07}}
   \put(0.47,0.95){\footnotesize 1}
   \put(0.5,0.50){\circle*{0.07}}
   \put(0.47,0.58){\footnotesize 7}
   \put(0.5,0.30){\circle*{0.07}}
   \put(0.465,0.15){\footnotesize 8}
   \put(0.25,0.433){\circle*{0.07}}
   \put(0.1,0.45){\footnotesize 2}
   \put(0.75,0.433){\circle*{0.07}}
   \put(0.8,0.45){\footnotesize 6}
   \put(0,0){\line(1,0){1}}
   \begin{drawjoin}
      \jput(0,0){\ }
      \jput(0.5,0.866){\ }
   \end{drawjoin}
   \begin{drawjoin}
      \jput(0.5,0.866){\ }
      \jput(1,0){\ }
   \end{drawjoin}
\end{picture}
\\[5mm]
$F_7^{-4} + e$ & $\scrw^3 + e + f$ \\[10mm]
\end{tabular}
\begin{tabular}{c@{\qquad\qquad\qquad}c}
\begin{picture}(1,1)(0,0)
   \put(0,0){\circle*{0.07}}
   \put(0.5,0){\circle*{0.07}}
   \put(1,0){\circle*{0.07}}
   \put(0.5,0.866){\circle*{0.07}}
   \put(0.166,0.289){\circle*{0.07}}
   \put(0.333,0.577){\circle*{0.07}}
   \put(0.75,0.433){\circle*{0.07}}
   \put(0,0){\line(1,0){1}}
   \begin{drawjoin}
      \jput(0,0){\ }
      \jput(0.5,0.866){\ }
   \end{drawjoin}
   \begin{drawjoin}
      \jput(0.5,0.866){\ }
      \jput(1,0){\ }
   \end{drawjoin}
   \put(-0.15,-0.03){\footnotesize 3}
   \put(0.02,0.32){\footnotesize 7}
   \put(1.07,-0.03){\footnotesize 5}
   \put(0.47,0.95){\footnotesize 1}
   \put(0.465,-0.15){\footnotesize 4}
   \put(0.18,0.60){\footnotesize 2}
   \put(0.8,0.45){\footnotesize 6}
   \put(0.5,0.289){\circle*{0.07}}
   \put(0.465,0.38){\footnotesize 8}
\end{picture}         
&
\begin{picture}(1,1)(0,0)
   \put(0,0){\circle*{0.07}}
   \put(-0.15,-0.03){\footnotesize 3}
   \put(0.5,0){\circle*{0.07}}
   \put(0.47,-0.15){\footnotesize 4}
   \put(1,0){\circle*{0.07}}
   \put(1.07,-0.03){\footnotesize 5}
   \put(0.5,0.866){\circle*{0.07}}
   \put(0.47,0.95){\footnotesize 1}
   \put(0.5,0.433){\circle*{0.07}}
   \put(0.465,0.52){\footnotesize 7}
   \put(0.25,0.433){\circle*{0.07}}
   \put(0.1,0.45){\footnotesize 2}
   \put(0.75,0.433){\circle*{0.07}}
   \put(0.8,0.45){\footnotesize 6}
   \put(0,0){\line(1,0){1}}
   \put(0.25,0.433){\line(1,0){0.5}}
   \begin{drawjoin}
      \jput(0,0){\ }
      \jput(0.5,0.866){\ }
   \end{drawjoin}
   \begin{drawjoin}
      \jput(0.5,0.866){\ }
      \jput(1,0){\ }
   \end{drawjoin}
   \put(1.05,0.433){\circle*{0.07}}
   \put(1.02,0.51){\footnotesize 8}
\end{picture}
\\[5mm]
$\scrw^{3+} + e$ & $P'_7 + e$
\end{tabular}
\end{center}
\vspace{-3mm}
\caption{
   Some rank-3 8-element matroids that are free extensions of
   rank-3 7-element matroids.
}
 \label{fig_F7m4pe_etc}
\end{figure}
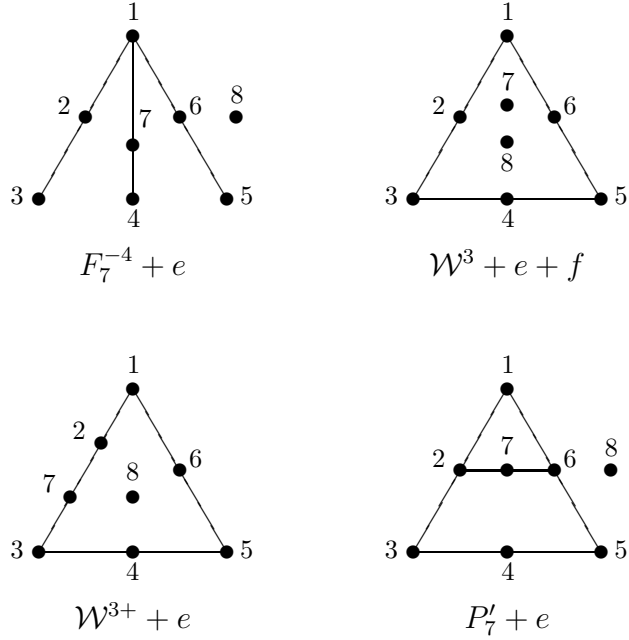

\begin{quote}
\begin{itemize}
   \item[$F_7^{-4} + e$:]
       $F$-representable if and only if $|F| \ge 7$.
       Not sixth-root-of-unity.
       Not transversal;
       cotransversal but not co-nice
       (Example~\ref{sec_transversal}.\ref{exam.transversal.10}).
       Nice principal extension of $F_7^{-4}$ with
       $\lambda_1 = 0$, $\lambda_2 = \ldots = \lambda_7 = 1/4$
       (Example~\ref{sec_constructions}.\ref{exam.nice_principal_F7m4pe}).
       Every single-element deletion is isomorphic to
       $F_7^{-4}$, $F_7^{-5}$ or $U_{3,7}$.
       Not known whether it has half-plane property.
   \item[$\scrw^3 + e + f$:]
       $F$-representable if and only if $|F| \ge 7$.
       Not sixth-root-of-unity.
       Transversal but not nice
       (Example~\ref{sec_transversal}.\ref{exam.transversal.3});
       cotransversal but not co-nice
       (Example~\ref{sec_transversal}.\ref{exam.transversal.11}).
       Non-nice free extension of $\scrw^3 + e$.
       Every single-element deletion is isomorphic to
       $\scrw^3 + e$, $F_7^{-5}$ or $F_7^{-6}$.
       Not known whether it has half-plane property.
   \item[$\scrw^{3+} + e$:]
       $F$-representable if and only if $|F| \ge 7$.
       Not sixth-root-of-unity.
       Transversal but not nice
       (Example~\ref{sec_transversal}.\ref{exam.transversal.3});
       cotransversal but not co-nice
       (Example~\ref{sec_transversal}.\ref{exam.transversal.13}).
       Non-nice free extension of $\scrw^{3+}$.
       Every single-element deletion is isomorphic to
       $Q_7$, $S_7$, $\scrw^{3+}$, $\scrw^3 + e$ or $F_7^{-5}$.
       Not known whether it has half-plane property.
   \item[$P'_7 + e$:]
       $F$-representable if and only if $|F| \ge 7$.
       Not sixth-root-of-unity.
       Not transversal;
       cotransversal but not co-nice
       (Example~\ref{sec_transversal}.\ref{exam.transversal.16}).
       Non-nice free extension of $P'_7$.
       Every single-element deletion is isomorphic to
       $\scrw^3 + e$, $F_7^{-5}$, $P'_7$ or $P'''_7$.
       Not known whether it has half-plane property.
\end{itemize}
\end{quote}

The Pappus and non-Pappus matroids have rank 3 and 9 elements,
and are shown in Figure~\ref{fig_pappus}.
We also consider some related matroids:
the deletions non-Pappus$\setminus 1$ and non-Pappus$\setminus 9$,
and the free extension (non-Pappus $\drop\,9)+e$.
These matroids have the following properties:

\begin{figure}[t]
\setlength{\unitlength}{0.7cm}
\begin{center}
\begin{tabular}{c@{\qquad\qquad\qquad}c}
\begin{picture}(6,4)(0,0)
   \put(0,0){\circle*{0.2}}  
   \put(3,0){\circle*{0.2}}  
   \put(6,0){\circle*{0.2}}  
   \put(0,4){\circle*{0.2}}  
   \put(3,4){\circle*{0.2}}  
   \put(6,4){\circle*{0.2}}  
   \put(1.5,2){\circle*{0.2}}  
   \put(3,2){\circle*{0.2}}    
   \put(4.5,2){\circle*{0.2}}  
   \put(0,0){\line(1,0){6}}  
   \put(0,4){\line(1,0){6}}  
   \put(0,0){\line(3,4){3}}  
   \put(3,0){\line(3,4){3}}  
   \put(0,0){\line(3,2){6}}  
   \put(0,4){\line(3,-2){6}} 
   \put(0,4){\line(3,-4){3}} 
   \put(3,4){\line(3,-4){3}} 
   \put(1.5,2){\line(1,0){3}}  
   \put(-0.1,-0.5){\footnotesize 4}
   \put(2.9,-0.5){\footnotesize 5}
   \put(5.9,-0.5){\footnotesize 6}
   \put(-0.1,4.3){\footnotesize 1}
   \put(2.9,4.3){\footnotesize 2}
   \put(5.9,4.3){\footnotesize 3}
   \put(1.0,1.9){\footnotesize 7}
   \put(2.9,2.3){\footnotesize 8}
   \put(4.8,1.9){\footnotesize 9}
\end{picture}
&
\begin{picture}(6,4)(0,0)
   \put(0,0){\circle*{0.2}}  
   \put(3,0){\circle*{0.2}}  
   \put(6,0){\circle*{0.2}}  
   \put(0,4){\circle*{0.2}}  
   \put(3,4){\circle*{0.2}}  
   \put(6,4){\circle*{0.2}}  
   \put(1.5,2){\circle*{0.2}}  
   \put(3,2){\circle*{0.2}}    
   \put(4.5,2){\circle*{0.2}}  
   \put(0,0){\line(1,0){6}}  
   \put(0,4){\line(1,0){6}}  
   \put(0,0){\line(3,4){3}}  
   \put(3,0){\line(3,4){3}}  
   \put(0,0){\line(3,2){6}}  
   \put(0,4){\line(3,-2){6}} 
   \put(0,4){\line(3,-4){3}} 
   \put(3,4){\line(3,-4){3}} 
   \put(-0.1,-0.5){\footnotesize 4}
   \put(2.9,-0.5){\footnotesize 5}
   \put(5.9,-0.5){\footnotesize 6}
   \put(-0.1,4.3){\footnotesize 1}
   \put(2.9,4.3){\footnotesize 2}
   \put(5.9,4.3){\footnotesize 3}
   \put(1.0,1.9){\footnotesize 7}
   \put(2.9,2.3){\footnotesize 8}
   \put(4.8,1.9){\footnotesize 9}
\end{picture}
\\[5mm]
Pappus & non-Pappus \\[10mm]
\end{tabular}
\begin{tabular}{c}
\begin{picture}(6,4)(0,0)
   \put(0,0){\circle*{0.2}}  
   \put(3,0){\circle*{0.2}}  
   \put(6,0){\circle*{0.2}}  
   \put(0,4){\circle*{0.2}}  
   \put(3,4){\circle*{0.2}}  
   \put(6,4){\circle*{0.2}}  
   \put(1.5,2){\circle*{0.2}}  
   \put(3,2){\circle*{0.2}}    
   \put(4.5,2){\circle*{0.2}}  
   \put(0,0){\line(1,0){6}}  
   \put(0,4){\line(1,0){6}}  
   \put(0,0){\line(3,4){3}}  
   \put(0,0){\line(3,2){6}}  
   \put(0,4){\line(3,-2){6}} 
   \put(0,4){\line(3,-4){3}} 
   \put(-0.1,-0.5){\footnotesize 4}
   \put(2.9,-0.5){\footnotesize 5}
   \put(5.9,-0.5){\footnotesize 6}
   \put(-0.1,4.3){\footnotesize 1}
   \put(2.9,4.3){\footnotesize 2}
   \put(5.9,4.3){\footnotesize 3}
   \put(1.0,1.9){\footnotesize 7}
   \put(2.9,2.3){\footnotesize 8}
   \put(4.8,1.9){\footnotesize 9}
\end{picture}
\\[5mm]
(non-Pappus $\drop\,9)+e$
\end{tabular}
\end{center}
\caption{
   The Pappus and non-Pappus matroids,
   along with the free extension (non-Pappus $\drop\,9)+e$.
}
 \label{fig_pappus}
\end{figure}
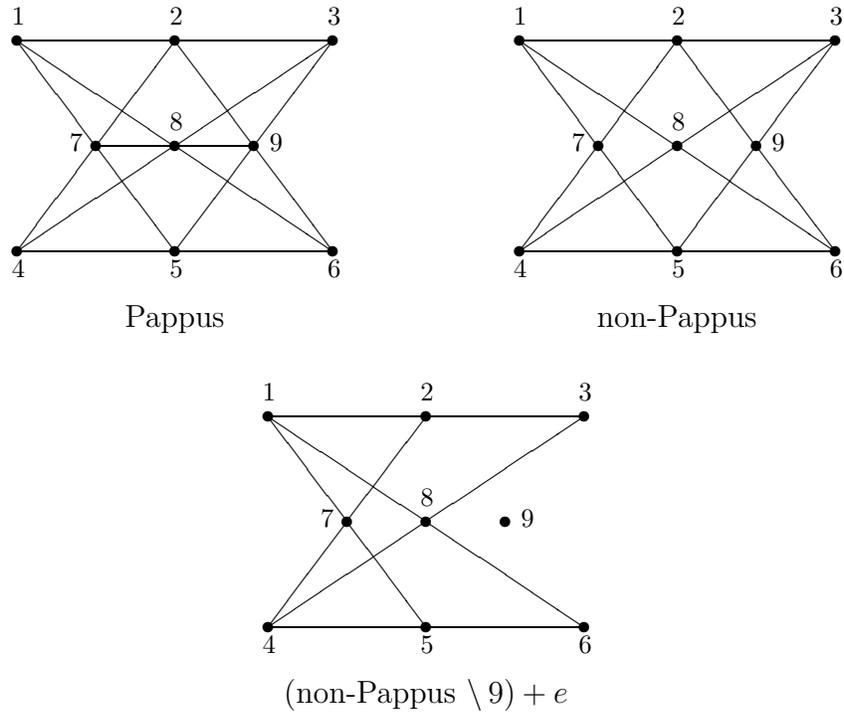

\begin{quote}
\begin{itemize}
   \item[Pappus:]
       $F$-representable if and only if $|F| = 4$ or $|F| \ge 7$.\footnote{
   The assertion on \cite[p.~516]{Oxley_92} is in error.
}
       Not sixth-root-of-unity.
       Not transversal or cotransversal.
       Transitive automorphism group.
       Every single-element deletion is isomorphic to non-Pappus$\setminus 9$.
       Does not have half-plane property
       (Example~\ref{sec_counterexamples}.\ref{exam.counterexamples.9}).
   \item[non-Pappus:]
       Not representable over any field.
       Not sixth-root-of-unity.
       Not transversal or cotransversal.
       Every single-element deletion is isomorphic to non-Pappus$\setminus 1$
          or non-Pappus$\setminus 9$.
       Does not have half-plane property
       (Example~\ref{sec_counterexamples}.\ref{exam.counterexamples.9}).
   \item[non-Pappus$\setminus 1$:]
       $F$-representable if and only if $|F| \ge 5$.
       Not sixth-root-of-unity.
       Not transversal;  cotransversal but not co-nice
       (Example~\ref{sec_transversal}.\ref{exam.transversal.19}).
       Every single-element deletion is isomorphic to 
       $F_7^{-4}$, $\scrw^3 + e$, $F_7^{-5}$, $P'_7$ or $P''_7$.
       Not known whether it has half-plane property.
   \item[non-Pappus$\setminus 9$:]
       $F$-representable if and only if $|F| \ge 4$.
       Not sixth-root-of-unity.
       Not transversal or cotransversal.
       Every single-element deletion is isomorphic to $F_7^{-4}$ or $P'_7$.
       Not known whether it has half-plane property.
   \item[(non-Pappus $\drop\,9)+e$:]
       $F$-representable if $|F| \ge 7$.
       Not sixth-root-of-unity.
       Not transversal or cotransversal.
       Every single-element deletion is isomorphic to
          $F_7^{-4} + e$, $P'_7 + e$ or non-Pappus$\setminus 9$.
       Does not have half-plane property
       (Example~\ref{sec_counterexamples}.\ref{exam.counterexamples.10}).
\end{itemize}
\end{quote}

\subsection{Some rank-4 matroids on 8 elements}

The rank-4 matroid $P_8$ is represented over any field
of characteristic $\neq 2$ by the matrix\footnote{
   The $GF(3)$-representation of $P_8$ on \cite[p.~512]{Oxley_92}
   has a misprint:  the bottom right element should be 0, not 1.
}
\be
\left[
   \begin{array}{rrrr|rrrr}
      1 & 0 & 0 & 0 &    0 & 1 & 1 &  2  \\
      0 & 1 & 0 & 0 &    1 & 0 & 1 &  1  \\
      0 & 0 & 1 & 0 &    1 & 1 & 0 &  1  \\
      0 & 0 & 0 & 1 &    2 & 1 & 1 &  0
   \end{array}
\right]
\ee
A geometric representation of $P_8$ over the reals can be obtained
by starting from a 3-dimensional cube and then rotating one face
of the cube by $45^\circ$ in its plane \cite[Figure 14]{Geelen_00a}.
There is a unique pair of disjoint circuit-hyperplanes in $P_8$,
namely $\{1,4,5,8\}$ and $\{2,3,6,7\}$.
We denote by $P'_8$ (resp.\ $P''_8$)
the matroid obtained from $P_8$ by relaxing
one (resp.\ both) of these circuit-hyperplanes.

The V\'amos matroid $V_8$ also has rank 4;
it is pictured in \cite[p.~76, Figure 2.4]{Oxley_92}.

\begin{quote}
\begin{itemize}
   \item[$P_8$:]
       $F$-representable if and only if $F$ has characteristic $\neq 2$.
       Not sixth-root-of-unity.
       Self-dual.
       Not transversal or cotransversal.
       Transitive automorphism group.
       Every single-element contraction is isomorphic to $P_7$;
       every single-element deletion is isomorphic to $(P_7)^*$.
       See \cite{Geelen_00a,Geelen_00b} for more information.
       Does not have half-plane property
       (Example~\ref{sec_counterexamples}.\ref{exam.counterexamples.8}).
   \item[$P'_8$:]
       $F$-representable if and only if $|F| \ge 4$.
       Not sixth-root-of-unity.
       Self-dual.
       Not transversal or cotransversal.
       Every single-element contraction is isomorphic to $P_7$ or $P'_7$;
       every single-element deletion is isomorphic to $(P_7)^*$ or $(P'_7)^*$.
       Does not have half-plane property
       (Example~\ref{sec_counterexamples}.\ref{exam.counterexamples.8}).
   \item[$P''_8$:]
       $F$-representable if and only if $|F| \ge 5$.
       Not sixth-root-of-unity.
       Self-dual.
       Not transversal or cotransversal.
       Transitive automorphism group.
       Every single-element contraction is isomorphic to $P'_7$;
       every single-element deletion is isomorphic to $(P'_7)^*$.
       See \cite{Geelen_00a,Geelen_00b} for more information.
       Does not have half-plane property
       (Example~\ref{sec_counterexamples}.\ref{exam.counterexamples.8}).
   \item[$V_8$:]
       Not representable over any field.
       Not sixth-root-of-unity.
       Self-dual.
       Not transversal or cotransversal.
       Every single-element contraction is isomorphic to
           $F_7^{-4}$ or $F_7^{-5}$;
       every single-element deletion is isomorphic to
           $(F_7^{-4})^*$ or $(F_7^{-5})^*$.
       Not known whether it has half-plane property.
\end{itemize}
\end{quote}

\subsection{Some transversal matroids}

Some interesting transversal matroids are introduced in
Sections~\ref{sec_transversal_examples} and \ref{sec_rank3_transversal}.
In particular, all rank-3 transversal matroids are classified
in Section~\ref{sec_rank3_transversal};
and the rank-($n+1$) transversal matroid $N_n$ is defined in
Example~\ref{sec_transversal}.\ref{exam.transversal.18}.

\section{$(F,\{1\})$-representability of matroids}  \label{app_F1rep}


\newcommand{\de}{\backslash}
\newcommand{\ag}{AG(2,3)}
\newcommand{\ba}{\backslash}
\newcommand{\sixrt}{\sqrt[6]{1}}
\newcommand{\co}{{\rm co}}
\newcommand{\si}{{\rm si}}
\newcommand{\cl}{{\rm cl}}
\newcommand{\utf}{U_{2,4}}
\newcommand{\mkf}{M(K_4)}

We defined $(F,G)$-representations of matroids in Section~\ref{sec_FG}
and noted that, when $-1 \in G$, such representations are covered under
the theory
of partial fields developed by Semple and Whittle \cite{Semple_96}.
When $-1 \notin G$, the class of $(F,G)$-representable matroids
is less well-behaved.
For example, neither reordering the columns, reordering the rows,
nor pivoting is guaranteed to produce another $(F,G)$-representation.
In this appendix, we consider $(F,G)$-representable matroids in
the case that $G$ is the trivial group $\{1\}$,
and we characterize such matroids.
If $F$ has characteristic 2, then it is clear that a matroid $M$ is
$(F,\{1\})$-representable if and only if $M$ is binary.
When the characteristic of $F$ is not 2,
the determination of all $(F,\{1\})$-representable matroids,
which is contained in Theorem~\ref{F1}, is more difficult.

We now inductively define when a graph is a {\it chain of cycles}\/.
\begin{itemize}
\item[(i)] A single cycle $C$ of length at least three is a chain
of cycles, which we write as $(C)$.
\item[(ii)] Suppose that $(C_1,C_2,\ldots,C_n)$ is a chain of
cycles and that each of $C_1,C_2,\ldots,C_n$ is a cycle.
Let $e$ be an edge of $C_n$ that is in none of
$C_1,C_2,\ldots,C_{n-1}$ and let $C_{n+1}$ be a cycle of length
at least three that contains $e$ but no other edge of
$(C_1,C_2,\ldots,C_n)$.  Let $(C_1,C_2,\ldots,C_{n+1})$
be the parallel connection of $(C_1,C_2,\ldots,C_n)$
and $C_{n+1}$ across the edge $e$. Then
$(C_1,C_2,\ldots,C_{n+1})$ is a chain of cycles.
\end{itemize}
A graph $G'$ is an {\em augmented chain of cycles}\/
if $G'$ can be obtained from a chain of cycles $G$ by adding
some (possibly empty) set of edges so that each added edge is
parallel to some edge of $G$.

The next theorem is the main result of this appendix. The graph
$G_6$  has 6 vertices and 9 edges
and is obtained by attaching a
triangle via parallel connection to each edge of a triangle.

\begin{theorem}
\label{F1}
The following statements are equivalent for a matroid $M$ and a field $F$
of characteristic other than 2.
\begin{enumerate}
\item[(i)] $M$ is $(F,\{1\})$-representable.
\item[(ii)] Each  component of $M$ is a loop, a rank-1 uniform matroid,
or the cycle matroid of an augmented chain of cycles.
\item[(iii)] $M$ has no minor isomorphic to $U_{2,4}, M(K_4), M(K_{2,3})$,
or $M(G_6)$.
\end{enumerate}
\end{theorem}

A consequence of Theorem~\ref{F1} is that, when the characteristic of $F$
is not 2, the class of $(F,\{1\})$-representable matroids is not closed
under duality.
For example, if $G$ is the graph that is obtained from $G_6$
by contracting one edge incident with a degree-2 vertex, then
$G$ is an augmented chain of cycles, so $M(G)$ is $(F,\{1\})$-representable.
However, $M^*(G)$ has $M(K_{2,3})$ as a minor and so is not
$(F,\{1\})$-representable.
We can also use $M(G)$ to show that, when the characteristic of $F$ is not 2,
the class of $(F,\{1\})$-representable matroids is not closed under 2-sum,
series connection, or parallel connection.  To see this, we observe that
$M(G_6)$ is the 2-sum of $M(G)$ and the graph obtained by adding an edge
parallel to one edge of a triangle; $M(G_6)$ is also the series connection of
$M(G)$ and $U_{1,2}$; and, finally, $M(G_6)$ is the parallel connection of the
simplification of $M(G)$ and a triangle.

\bigskip

The proof of Theorem~\ref{F1} will use the following definition.
Fix an integer $r \ge 1$;
and, for $1 \le i \le r$, let ${\bf e}_i$ be the column vector of length $r$
with a 1 in the $i$th row and zeros elsewhere.
An $r \times n$ matrix $A$ will be called {\it special}\/ if:
\begin{itemize}
          \item[(a)] the first column is ${\bf e}_1$,
        and the last column is ${\bf e}_r$;
   \item[(b)] each column is either ${\bf e}_i$ or
        ${\bf e}_i + {\bf e}_{i+1}$;
   \item[(c)] a column ${\bf e}_i$ is always followed by either ${\bf e}_i$ or
       ${\bf e}_i + {\bf e}_{i+1}$; and
   \item[(d)] a column ${\bf e}_i + {\bf e}_{i+1}$ is always followed by
       either ${\bf e}_i + {\bf e}_{i+1}$, ${\bf e}_{i+1}$,
       or ${\bf e}_{i+1} + {\bf e}_{i+2}$.
\end{itemize}
Note that a special matrix may have repeated columns;
but if two columns are equal to $c$, say, then all
intermediate columns must also be equal to $c$.

The next result is the main step in the proof of Theorem~\ref{F1}.

\begin{proposition}
\label{P}
Let $A$ be an $(F,\{1\})$-representation of a connected loopless matroid $M$.
If the number of rows of $A$ equals $r(M)$, then  $A$ is special.
\end{proposition}

\proof
Let $A$ be the
$r \times n$ matrix $[a_{ij}]$.
Assume that the columns of $A$ are indexed, in order, by $1,2,\ldots, n$ where
$E(M) = \{1,2,\ldots,n\}$.
The proposition certainly holds if $A$ has one row, so
assume that $r \ge 2$.
Note that $A$ has no zero columns, because $M$ is loopless;
and $A$ has no zero rows, because $r = r(M)$.
Since every subdeterminant of $A$ is in $\{0,1\}$, it follows that
$A$ has no $2 \times 2$
submatrices of the form
$\left[\!
\begin{array}{cc}
* & 1 \\
1 & 0
\end{array}
\!\right]$
or
$\left[\!
\begin{array}{cc}
0 & 1 \\
1 & *
\end{array}
\!\right]$.
Therefore:
\begin{noname}
\label{first}
If $a_{ij} = 0$, then
\begin{enumerate}
\item[(i)] $a_{ij'} = 0$ for all $j'$ with $1 \le j' \le j$, or
$a_{i'j} = 0$ for all $i'$ with $1 \le i' \le i$; and
\item[(ii)] $a_{ij'} = 0$ for all $j'$ with $j \le j' \le n$, or
$a_{i'j} = 0$ for all $i'$ with $i \le i' \le r$.
\end{enumerate}
\end{noname}
We rewrite these conditions as:
\begin{noname}
\label{nwse} If $a_{ij} = 0$, then
\begin{enumerate}
\item[(i)] either all the entries north of $a_{ij}$ are zero, or
all the entries west of $a_{ij}$ are zero; and
\item[(ii)] either all the entries south of $a_{ij}$ are zero, or
all the entries east of $a_{ij}$ are zero.
\end{enumerate}
\end{noname}

\begin{lemma}
\label{nesw}
Let $a_{ij} = 0$. Then either
\begin{enumerate}
\item[(i)] $a_{i'j'} = 0$ for all $i'$ and $j'$ with $1 \le i' \le i$ and
$j \le j' \le n$; or
\item[(ii)] $a_{i'j'} = 0$ for all $i'$ and $j'$ with $i \le i' \le r$ and
$1 \le j' \le j$.
\end{enumerate}
\end{lemma}

\proof
Since $a_{ij} = 0$, by (\ref{nwse})(i), either
\begin{enumerate}
\item[(a)] all entries north of $a_{ij}$ are 0; or
\item[(b)] all entries west of $a_{ij}$ are 0.
\end{enumerate}
Suppose that (a) holds.  Since column $j$ is non-zero, if $1 \le i' \le i$,
then, by
(\ref{nwse})(ii), all entries east of $a_{i'j}$ are 0. Thus (i) holds.
Now suppose that (b) holds.  Then, since row $i$ is non-zero,
if $1 \le j' \le j$, then, by (\ref{nwse})(ii), all entries south of $a_{ij'}$
are zero, and (ii) holds.
\qed

Again it will be convenient to rewrite the last lemma as:
\begin{noname}
\label{nesw2} If $a_{ij} = 0$, then  either
\begin{enumerate}
\item[(i)] all the entries north-east of $a_{ij}$ are zero; or
\item[(ii)]  all the entries south-west of $a_{ij}$ are zero.
\end{enumerate}
\end{noname}

By combining (\ref{first}) with the fact that $A$ has no zero rows and
no zero columns, we deduce that no column of $A$ has $[1~ 0~ 1]^{\rm T}$
as a submatrix and no row of $A$ has $[1~ 0~ 1]$
as a submatrix. Thus every column of $A$ and every row of $A$ consists of a
(possibly empty) sequence of zeros, followed by a non-empty sequence of ones,
followed by a (possibly empty) sequence of zeros.

Let $D_3$ be the matrix
$\displaystyle\left[\!\begin{array}{ccc}
1 & 1 & 0 \\
1 & 1 & 1 \\
0 & 1 & 1
\end{array}\!\right]$.
Then $D_3$ has determinant $-1$ and so the next lemma is immediate.

\begin{lemma}
\label{D3}
The matrix $A$ does not have $D_3$ as a submatrix.
\end{lemma}

\begin{lemma}
\label{I}
If $a_{ij} = 1$ and $a_{(i+1)j} = 0$, then
$j < n$ and  $a_{i(j+1)} = 1$.
\end{lemma}

\proof
If $j =n$, then, since $a_{ij} = 1$, it follows by (\ref{nwse}) that row $i+1$
of $A$ is zero; a contradiction.  Thus $j < n$.  Suppose that
$a_{i(j+1)} = 0$. Then, by (\ref{nesw2}), all entries north-east of
$a_{i(j+1)}$ are zero and all entries south-west of $a_{(i+1)j}$ are zero.
Thus $A$ has the block form
$\left[\!
\begin{array}{cc}
A_1 & 0 \\
0 & A_2
\end{array}
\!\right]$
and so $M$ is disconnected; a contradiction.
\qed

By a similar argument to that just given, we get:
\begin{lemma}
\label{II}
If $a_{ij} = 0$ and $a_{(i+1)j} = 1$, then
$j \ge 2$ and  $a_{(i+1)(j-1)} = 1$.
\end{lemma}

\begin{lemma}
\label{nz}
$a_{r1} = 0 = a_{1n}$.
\end{lemma}

\proof
Suppose that $a_{r1} = 1$.
If $a_{i1} = 0$ for some $i < r$, then, by (\ref{nesw}), row $i$ is zero;
a contradiction. Thus column 1 of $A$ consists of all ones and, similarly,
row $r$ consists of all ones. Since $r \ge 2$, there is a column of $A$
that does not consist of all ones.  Let $j$ be the first such column. Then
$a_{j1} = 0$ and $\{1,2,\ldots,j-1\}$ is the ground set of a rank-1
component of $M$; a contradiction. We conclude that $a_{r1} = 0$, and
a similar argument shows that $a_{1n} = 0$.
\qed

Assume now that Proposition~\ref{P} is false,
and let $A$ be a counterexample having the minimum number of columns.
In particular,  no two consecutive  columns of $A$ are equal.
The following lemmas derive some properties of this minimal counterexample,
culminating in a contradiction.

\begin{lemma}
Every column of $A$ has at most two ones.
\end{lemma}

\proof
Suppose that, for some $j$ with $1 < j < n$, column $j$ has a unique one, which
occurs in row $i$ say. Then $1 < i < r$, and every entry north-east of $a_{(i-1)j}$ is
zero, as is every entry south-west of $a_{(i+1)j}$. Thus
$M/j$ is disconnected, so $M$ is a parallel
connection with basepoint $j$.  Indeed, $M$ is the parallel connection
of $M[A_1]$ and $M[A_2]$ where $A_1$ is the submatrix of $A$ that consists of
$a_{ij}$ and all entries to its north-west, and $A_2$ is the submatrix of $A$
that consists of
$a_{ij}$ and all entries to its south-east. Since $M$ is connected, so are both
$M[A_1]$ and $M[A_2]$.  It follows that both $A_1$ and $A_2$ are special
matrices and, therefore, so is $A$; a contradiction. Thus every column of $A$
except possibly the first or the last has at least two ones.

Now suppose that $A$ has a single one in its first column. Then   $a_{11} = 1$.
 Suppose that $M/1$ is connected.  Let $A_1$ be the matrix obtained from $A$ by
deleting the first row and column.  Then $A_1$  is special.
By Lemma~\ref{nz}, $a_{1n} = 0$.  Let the first zero entry in row 1 of $A$
be $a_{1m}$.  Then $m \ge 3$ otherwise 1 is a coloop of $M$. Now
$a_{11} = a_{12} =\ldots = a_{1(m-1)} = 1$.
Moreover,
$a_{22} = a_{23} =\ldots = a_{2(m-1)} = 1$ otherwise $M[A_1]$ has a loop. Since
$A$ has no consecutive equal columns and $A_1$ is special,
$m \le 4$. If
$m = 3$, then  columns 1
and 2 of $A$ are
$[1~0~0~\ldots~0]^{\rm T}$ and  $[1~1~0~\ldots~0]^{\rm T}$,
and it follows that $A$
is special; a contradiction. We may now suppose that $m=4$.
Then $a_{24} = 0$ otherwise  the submatrix of $A$ determined by the first three rows
and columns 2, 3, and 4 is $D_3$; a contradiction. Thus $a_{2j} = 0$ for all
$j \ge 4$. Now the matrix $A$ represents $M$ over $F$. If we subtract
row 2 from row 1 in $A$, we get another $F$-representation for $M$, which shows
that 1 is a coloop of $M$; a contradiction.

We may now assume that $M/1$ is disconnected. Then $M\backslash 1$ is
connected. It follows that the matrix $A \backslash 1$, which is obtained from
$A$ by deleting the first column, is special. Thus column 2 of $A$ equals
$e_1$.  Therefore column 1 is parallel to column 2; a contradiction.

We conclude that the first column of $A$ has at least two ones, and, by
symmetry,  the last column of $A$ has at least two ones.
\qed

A submatrix $C$ of a matrix $D$ will be called {\it consecutive} if the rows of
$C$ are consecutive rows of $D$ and the columns of $C$ are consecutive columns
of $D$.

\begin{lemma}
$A$ has
$\left[\!
\begin{array}{cc}
1 & 1 \\
0 & 1
\end{array}
\!\right]$
as a consecutive submatrix.
\end{lemma}

\proof
Proceed north from $a_{r1}$, which we know from Lemma~\ref{nz} is 0, until the
first 1 is reached.
Then we have
$\left[\!
\begin{array}{c}
1  \\
0
\end{array}
\!\right]$
as a submatrix meeting consecutive rows.  Then proceed east from the 0 in this
submatrix until the first 1 is met. By Lemma~\ref{I}, we deduce that $A$ has
$\left[\!
\begin{array}{ccccc}
1 & 1 & \ldots & 1 & 1 \\
0 & 0 & \ldots & 0 & 1
\end{array}
\!\right]$
as a consecutive submatrix, and the lemma follows.
\qed

Let $J_k'$ be the $k \times k$ matrix all of whose entries are one except for
that in the south-west corner, which is zero.

\begin{lemma}
\label{23}
Every consecutive $J_2'$-submatrix of $A$
is contained in a consecutive $J_3'$-submatrix of $A$.
\end{lemma}

\proof
Suppose that
$\displaystyle\left[\!\begin{array}{ccc}
c_{11} & c_{12} & c_{13}\\
1 & 1 & c_{23}\\
0 & 1 & c_{33}
\end{array}\!\right]$
is a consecutive submatrix $C$ of $A$ containing
$\left[\!
\begin{array}{cc}
1 & 1 \\
0 & 1
\end{array}
\!\right]$
as a submatrix.  The additional row must exist and $c_{11} = 1$
otherwise $A$ has a column with a single one. The additional column must exist
and $c_{33} = 1$
otherwise $A$ has a row with a single one and so has a coloop.

Suppose $c_{12} = 0$.  Then adjoin an additional column $x$ to $A$ that is zero
everywhere except in the row corresponding to the second row of $C$.  This
column should be added between the first and second columns of $C$ with the
resulting matrix being $A'$. The matroid $M[A']/x$ is disconnected. Thus $M[A]$
is the 2-sum of two matroids with basepoint $x$. Let $A_1'$ be the submatrix of
$A'$ consisting of the non-zero entry in column $x$ and all entries north-west
of it, and $A_2'$ be the submatrix of $A'$ consisting of the non-zero entry in
column $x$ and all entries south-east of it. Then $M[A']$ is the parallel
connection of $M[A_1']$ and $M[A_2']$ and so each of these is connected.  Since
$A$ has a submatrix equal to $A_1'$, the last matrix is an
$(F,\{1\})$-representation for a connected matroid. Since $A_2'$ is obtained by
adjoining the column $[1~0~0~\ldots~0]^{\rm T}$
to the beginning of a submatrix of
$A$, it too is an $(F,\{1\})$-representation for a connected matroid. Thus each
of $A_1'$ and $A_2'$ is special. It follows that $A'$ is special and so too is
$A$; a contradiction.  We conclude that $c_{12} = 1$.

Suppose $c_{23} = 0$. Let $y$ be the column of $A$ corresponding to the second
column of $C$.  Then $A\backslash y$ has the form
$\left[\!
\begin{array}{cc}
A_1 & 0 \\
0 & A_2
\end{array}
\!\right]$
and so $M\backslash y$ is disconnected and the rank of $A$ is the
sum of the ranks, $r_1$ and $r_2$, say, of $A_1$ and $A_2$. We may apply row
operations in $A$ to transform $A_2$ into a matrix that has a permutation of
$I_{r_2}$ as a submatrix.  These row operations will not affect the first $r_1$
rows of $A$.  Then, by contracting the elements of $M$ corresponding to the
columns of $I_{r_2}$, we see that the submatrix $A_1''$ consisting of $A_1$ and
the first $r_1$ rows of column $y$ represents a connected matroid.  This is
because $M$ is the series connection of this matroid and another with respect
to the basepoint $y$.  Thus $A_1''$ is special.   But the last column of
$A_1''$ has more than one one; a contradiction.  We conclude that $c_{23} = 1$.

Finally, since $D_3$ is not a submatrix of $A$, we deduce that $c_{13} = 1$.
\qed

\begin{lemma}
For all $k \ge 2$, the matrix $A$ has
$J_k'$ as a consecutive submatrix.
\end{lemma}

\proof We shall argue by induction.
By the last two lemmas, $A$ has $J_2'$ and $J_3'$ as consecutive submatrices.
Assume that $A$ has $J_m'$ as a consecutive submatrix, for some $m \ge 3$.
The rows $i$ and $i+1$ of $A$ that correspond to rows 1 and 2 of $J_m'$ must
differ.  Therefore, in $A$, the submatrix induced by rows $i$ and $i+1$ has
either
\begin{enumerate}
\item[(i)] in the consecutive columns immediately after $J_m'$, a $2 \times t$
matrix $H_1$ of the form
$\left[\!
\begin{array}{ccccc}
1 & 1 & \ldots & 1 & 0 \\
1 & 1 & \ldots & 1 & 1
\end{array}
\!\right]$
for some $t \ge 1$;
\item[(ii)] in the consecutive columns immediately before $J_m'$, a $2 \times
s$ matrix $H_2$ of the form
$\left[\!
\begin{array}{ccccc}
1 & 1 & \ldots & 1 & 1 \\
0 & 1 & \ldots & 1 & 1
\end{array}
\!\right]$
for some $s \ge 1$.
\end{enumerate}

In case (i), the column of $A$ corresponding to the last column of $H_1$ must
have ones in the $m - 1$ rows immediately below the unique zero otherwise
(\ref{nesw2}) is contradicted.  It follows that $A$ has $D_3$ as a submatrix; a
contradiction.

In case (ii), the unique zero in $H_2$ is in a consecutive submatrix of $A$
equal to $J_2'$.  Thus, by Lemma~\ref{23}, the $(i-1)$st row of $A$ must have
ones in all the columns corresponding to those of $J_m'$ otherwise $A$ has
$D_3$ as a submatrix. We deduce that $A$ has, as a consecutive submatrix, the
matrix $J_m''$ that is obtained by adjoining a row of ones to the beginning of
$J_m'$.

The columns $j$ and $j+1$ of $A$ that correspond to the last two columns of
$J_m''$ must differ. If they differ in some row before row $i-1$ of $A$, let
row $u$ be the highest indexed such row where they differ. This row has
$[1~ 0]$ in columns $j$ and $j+1$ of $A$.  Then, by Lemma~\ref{nesw},
$a_{uj} = a_{u(j-1)} = \ldots = a_{u(j-m+2)} = 1$ and so $A$ has $D_3$ as a
submatrix; a contradiction.  We deduce that columns $j$ and $j+1$ of $A$ differ
in some row after row $i + m-1$.  Let row $v$ be the smallest indexed such row
 where these columns differ. Then
$$
   \left[\!\begin{array}{cc}
a_{(v-1)j} & a_{(v-1)(j+1)} \\
a_{vj} & a_{v(j+1)}
\end{array}\!\right]
 \;=\; J_2'  \;.
$$
By Lemma~\ref{23}, this $J_2'$ is contained in a consecutive submatrix of $A$
equal to $J_3'$. To avoid having $D_3$ as a submatrix of $A$, it follows that
$a_{(i-1)(j+2)} = a_{i(j+2)} = \ldots = a_{(i+m-1)(j+2)} = 1$. Thus $A$ has
$J_{m+1}'$ as a consecutive submatrix and the lemma follows by induction.
\qed

Proposition~\ref{P} follows immediately from
the last lemma since it implies that
the matrix $A$ is infinite.
\qed

\begin{lemma}
\label{A}
If $A$ is a special matrix, then every square submatrix of $A$ has determinant
in $\{0,1\}$.
\end{lemma}

\proof
We argue by induction on the number $n$ of columns of $A$.  If $n=1$, then the
result is immediate. Assume the result holds for $n = k$ and let $n = k+1$. Now
let $A'$ be a square submatrix of $A$. We may assume that the columns of $A$
are distinct since if $A'$ has repeated columns, its determinant is 0. If $A'$
meets the first column of $A$, then $\det A' = 0$ unless $A'$ meets the first
row of $A$. In the exceptional case, $\det A'$ equals the determinant of the
submatrix $A'_1$ of $A'$ that is obtained by deleting the first row and column.
 Now $A_1'$ is a submatrix of $A_1$, the matrix obtained by deleting the first
row and column of $A$. Moreover, it is not difficult to see that $A_1$ is
special. Thus, by induction, $\det A' \in \{0,1\}$. We may now assume that $A'$
does not meet the first column of $A$. If $A'$ does not meet the first row of
$A$, then $A'$ is a submatrix of the special matrix $A_1$ and the result
follows by the induction assumption. Thus we may assume that $A'$ contains the
entry $a_{12}$ of $A$ otherwise $\det A' = 0$. Since $a_{12}$ is the only
non-zero entry of the first row of $A'$, we have $\det A' = \det A_1'$. But
$A_1'$ is a submatrix of the special matrix $A_1$ and again the result follows
by the induction assumption.
\qed

We are now ready to prove Theorem~\ref{F1}.

\proofof{Theorem~\ref{F1}}
We shall show first that (ii) implies (iii). Suppose that $M$ satisfies (ii).
Then every minor of $M$ also satisfies (ii).  The construction of $M$
guarantees that $M$ is the cycle matroid of a series-parallel network $G$. Thus
$M$ has no minor isomorphic to $U_{2,4}$ or $M(K_4)$.  Moreover, $G$ is an
outerplanar graph, so $M$ has no minor isomorphic to $M(K_{2,3})$. Finally, we
observe that $M$ does not have $M(G_6)$ as a minor since $G_6$ cannot be
written as a chain of cycles. Thus (iii) holds.

Next we show that (iii) implies (ii). Suppose that $M$ is a simple connected
matroid that satisfies (iii) and has rank at least two, and assume that $M$ is
a minor-minimal matroid that does not satisfy (ii). The fact that $M$ has none
of $U_{2,4}, M(K_4),$ and $M(K_{2,3})$  as a minor ensures that $M \cong M(G)$
where $G$  is an outerplanar graph with at least three vertices.
Take an outerplanar embedding of $G$.  The boundary of the infinite face is a
Hamilton cycle $C$ of $G$ and all edges not in $C$ are chords of $C$. Certainly
no chords of $C$ cross. Now take a chord $e$ of $C$ with endpoints $x$ and $y$,
say, such that there is an $xy$-path $P_1$ in $C$ such that every vertex of
$V(P_1) - \{x,y\}$ has degree 2 in $G$. Then $M\backslash P_1$ satisfies (iii)
and hence also satisfies (ii).  Thus $M\backslash P_1$ is the cycle matroid of
a chain of cycles $(C_1,C_2,\ldots,C_n)$, and $M$ is obtained by taking the
parallel connection of $M\backslash P_1$ and the circuit with ground set $P_1
\cup e$. If $e$ is an edge of $C_1$ that is in no other $C_i$, then $M$ is the
cycle matroid of a chain of cycles. Therefore, by symmetry, we may now assume that either $e$
is an edge that is common to $C_i$ and $C_{i+1}$ for some $i$, or $e$ is in
exactly one $C_j$ for some $j$ in $\{2,3,\ldots,n-1\}$. The first case  cannot
arise because $G$ is outerplanar.  The second case implies without difficulty
that $M(G)$ has a minor isomorphic to $M(G_6)$.  This contradiction implies
that $M$ satisfies (ii).

Now suppose that $M$ satisfies (ii). We shall show that $M$ satisfies (i). If
the components $M_1,M_2,\ldots,M_m$ of $M$ have $(F,\{1\})$-representations
$A_1,A_2,\ldots,A_m$, respectively, then the matrix whose block form has
$A_1,A_2,\ldots,A_m$ on the main diagonal and zeros elsewhere is an
$(F,\{1\})$-representation for $M$.  Thus it suffices to prove that each
component of $M$ satisfies (i). Clearly a loop and a coloop are
$(F,\{1\})$-representable. Moreover, if the simple matroid associated with $M$ is
$(F,\{1\})$-representable, then so too is $M$. Thus we may assume that $M$
is simple. Let $C$ be a $k$-cycle for some $k \ge 3$.
Then the $(k-1) \times k$ matrix  with columns
${\bf e}_1,{\bf e}_1+{\bf e}_2,{\bf e}_2+{\bf e}_3,\ldots,
{\bf e}_{k-2} + {\bf e}_{k-1},{\bf e}_{k-1}$ is an
$F$-representation for $M[C]$.  The construction of a representation for a
parallel connection of two matroids is straightforward (see, for example,
 \cite[Proposition 7.1.21]{Oxley_92}). Using this, it is not difficult to show
that  the cycle matroid of a chain of cycles can be represented by a special
matrix.  It follows by Lemma~\ref{A} that $M$ satisfies (i).

Finally, we show that (i) implies (ii).  Assume that $M$ is
$(F,\{1\})$-representable and consider a  component $M'$
of $M$ of rank at
least two. Let $A$ be an $(F,\{1\})$-representation for $M'$.
We may assume that $A$ has $r(M')$ rows otherwise we can delete some row from
$A$ and still retain an $(F,\{1\})$-representation for $M'$.
Then,
by Proposition~\ref{P}, $A$ is special. From the definition of a special
matrix, it is not difficult to see that $M'$ is the cycle matroid of
an augmented chain of cycles.
\qed

\end{document}